\documentclass[11 pt, reqno]{article}

\usepackage{amsmath,amsfonts,amssymb,amsthm}
\usepackage[colorlinks=true,hyperindex=true]{hyperref}
\usepackage{cancel,bbm}
\usepackage{mathabx} 
\usepackage{comment}
\usepackage{enumitem}
\usepackage{cite}

\usepackage{chngcntr}
\counterwithin*{equation}{section}

\usepackage[left=2.5cm, right=2.5cm, top=3 cm, bottom= 3 cm]{geometry}
\usepackage{color}
\usepackage{fancyhdr}
\usepackage{latexsym}

\usepackage[T1]{fontenc}

\usepackage{bm}

\newtheorem{ccounter}{ccounter}[section]
\newtheorem{thm}[ccounter]{Theorem}
\newtheorem{lem}[ccounter]{Lemma}
\newtheorem{cor}[ccounter]{Corollary}
\newtheorem{defn}[ccounter]{Definition}
\newtheorem{prop}[ccounter]{Proposition}
\newtheorem{ass}[ccounter]{Assumption}
\newtheorem{ex}[ccounter]{Example}

\theoremstyle{definition}
\newtheorem*{remark}{Remark}

\def\bet{\begin{thm}}
\def\eet{\end{thm}}
\def\bel{\begin{lem}}
\def\eel{\end{lem}}
\def\bas{\begin{ass}}
\def\eas{\end{ass}}
\def\bec{\begin{cor}}
\def\eec{\end{cor}}
\def\bed{\begin{defn}}
\def\eed{\end{defn}}
\def\bep{\begin{prop}}
\def\eep{\end{prop}}
\def\beq{\begin{equation}}
\def\eeq{\end{equation}}
\def\proof{\noindent {\bf Proof.}\ \ }
\def\bea{\begin{equation*}}
\def\eea{\end{equation*}}
\def\tr{\mathrm{Tr}}
\def\bex{\begin{ex}}
\def\eex{\end{ex}}

\def\rr{\mathbb{R}}

\def\cc{\mathbb{C}}
\def\ff{\mathbb{F}}
\def\1{\boldsymbol{1}}

\let\Im\relax
\DeclareMathOperator{\Im}{Im}

\let\Re\relax
\DeclareMathOperator{\Re}{Re}

\DeclareMathOperator{\sgn}{sgn}

\let\le\relax
\def\le{\leq}

\let\ge\relax
\def\ge{\geq}

\def\e{\mathrm{e}}
\def\i{\mathrm{i}}
\def\del{\partial}
\def\d{\mathrm{d}}
\def\eps{\varepsilon}
\renewcommand\leq\varleq
\renewcommand\geq\vargeq
\def\ee{\mathrm{E}}

\def\F{\mathcal{F}}
\def\O{\mathcal{O}}

\def\ee{\mathbb{E}}

\def\pp{\mathbb{P}}

\def\rhosc{\rho_{\mathrm{sc}}}

\def\mfa{\mathfrak{m}}
\def\A{\mathcal{A}}

\def\Var{\mathrm{Var}}
\def\bx{\mathbf{x}}
\def\by{\mathbf{y}}
\def\S{\mathcal{S}}

\def\mfa{\mathfrak{a}}

\def\mfb{\mathfrak{b}}

\def\tilZ{\tilde{Z}}

\def\nn{\mathbb{N}}
\def\mfc{\mathfrak{c}}

\def\dd{\mathbb{D}}
\def\mfx{\varkappa}
\def\GMC{\mathrm{GMC}}

\def\nc{\normalcolor}
\def\hatgam{\hat{\gamma}}
\def\tilmu{\tilde{\mu}}
\def\tilW{\tilde{W}}

\def\tilL{\tilde{L}}
\def\calL{\mathcal{L}}
\def\tilG{\tilde{G}}
\def\tilZ{\tilde{Z}}
\def\blam{\boldsymbol{\lambda}}
\def\bz{\boldsymbol{z}}
\def\eeps{\epsilon}
\def\hatZ{\hat{Z}}
\def\hatL{\hat{L}}
\def\hatcalL{\hat{\mathcal{L}}}

\def\Pone{P^{(1)}}
\def\Ptwo{P^{(2)}} 
\def\tilF{\tilde{\mathcal{F}}}
\def\GB{\mathcal{G}}
\def\yY{\mathcal{Y}}
\def\kB{\mathcal{K}}
\def\xX{\mathcal{X}}
\def\bmE{\bm{E}}
\def\hC{\hat{\mathcal{C}}}
\def\cC{\mathcal{C}}
\def\vV{\mathcal{V}}
\def\bw{\boldsymbol{w}}
\def\bnu{\boldsymbol{\nu}}
\def\dz{\mathrm{d}^2 z}
\def\dw{\mathrm{d}^2 w}
\def\Wla{W^{\lambda}_A}

\def\ddb{\mathbb{D}^\beta}
\def\kapb{\kappa_4}
\def\Kb{K^\beta} 

\def\GMC{\mathrm{GMC}}
\def\Eone{\mathcal{E}_1}
\def\Etwo{\mathcal{E}_2}

\def\ddtwo{\mathbb{D}^{(2)}}
\def\EE{\mathcal{E}}
\def\hnu{\hat{\nu}}
\def\ddhb{\hat{\mathbb{D}}^{\beta}}
\def\cD{\mathcal{D}}
\def\hxX{\hat{\mathcal{X}}}
\def\rcC{\mathcal{C}_\rr} 
\def\Wlac{W^{\lambda}_{A, \cc}}
\def\Wlar{W^{\lambda}_{A, \rr}}
\def\Wlaf{W^{\lambda}_{A, \ff}}
\def\Qlac{Q^{\lambda}_{A, \cc}}
\def\Qlar{Q^{\lambda}_{A, \rr}}
\def\Qlaf{Q^{\lambda}_{A, \ff}}
\def\Nfb{N^{\mfb}}
\def\mfe{\mathfrak{e}}
\def\mfh{\mathfrak{h}}

\newcommand\blfootnote[1]{%
  \begingroup
  \renewcommand\thefootnote{}\footnote{#1}%
  \addtocounter{footnote}{-1}%
  \endgroup
}

\title{?????????????????/}

\begin{document}

\begin{table}
\centering

\begin{tabular}{c}

\multicolumn{1}{c}{\parbox{12cm}{\begin{center}\Large{\bf Gaussian Multiplicative Chaos for i.i.d. matrices}\end{center}}}\\
\\
\end{tabular}
\begin{tabular}{ c c c  }
Giorgio Cipolloni $^*$
& \phantom{blah} & 
Benjamin Landon $^\dagger$
 \\
  & & \\
\end{tabular}
\\
\begin{tabular}{c}
\multicolumn{1}{c}{\today}\\
\\
\end{tabular}

\begin{tabular}{p{15 cm}}
\small{{\bf Abstract:} We consider $N\times N$ matrices $X$ with i.i.d. entries, and prove that the random measure $\mu_N(z):=\nc\frac{ | \det (X-z)|^\gamma}{\ee[ | \det (X-z)|^\gamma]} \d z\d \overline{z}\nc$ converges to the Gaussian Multiplicative Chaos (GMC) on the unit disc in the full subcritical regime $\gamma \in (0, 2 \sqrt{2})$ as $N \to \infty$. Our result holds for both symmetry classes and in particular is new even for real Ginibre matrices. This result is the first of its kind for any non-invariant ensemble of random matrices. In the case of real matrices, we further prove that the restriction of $\mu_N(z)$ to the unit interval converges to a one-dimensional GMC. We establish the asymptotics for the $K$-point function of $| \det (X-z)|$ at any collection of mesoscopically separated points $z_i$. Our methods are analytic and probabilistic in nature, relying in part on the dynamical approach based on Dyson Brownian motion. }
\end{tabular}
\end{table} 
\blfootnote{$^*$ Email: cipolloni@axp.mat.uniroma2.it,   
Mathematics Department, University of Rome Tor Vergata, Rome, Italy. 
 ORCID: 0000-0002-4901-7992.    }
 \blfootnote{$^\dagger$ Corresponding author.  
Email: blandon@math.toronto.edu, 
Department of Mathematics, University of Toronto, Toronto, Canada.
ORCID: 0000-0003-2090-4886}
\blfootnote{MSC code (2020): 60B20}
\tableofcontents
\section{Introduction}

A fundamental problem of random matrix theory is determining the  behavior of large random matrix spectra. One natural choice of a random matrix $X$ is one with independent, identically distributed (i.i.d.) centered entries. The \emph{circular law} of  Girko \cite{girko1985circular} and Bai \cite{bai1997circular} answers the most basic asymptotic question about the eigenvalues  of $X$, by showing that in the limit of large dimension $N \to \infty$, the empirical eigenvalue measure converges  to the uniform measure on the unit disc $\mathbb{D}$.
\beq
\lim_{N\to+\infty}\frac{1}{N}\sum_i f(\sigma_i)= \frac{1}{\pi}\int_\mathbb{D}f(z)\,\d z \d \bar{z}.
\eeq 
 Here, we normalize $X$ so that
 $\ee[ |X_{ij}|^2]= \frac{1}{N}$ and denote its eigenvalues by $\{ \sigma_i\}_{i=1}^N$.

The next essential probabilistic question is to determine the fluctuations  around the limiting circular law. The formula $\sum_i f(\sigma_i)=\frac{1}{2\pi}\int_\cc \Delta f(z)\log |\det (X-z)|\, \d z \d \bar{z}$ naturally leads one to consider the fluctuations of the (asymptotically centered) log-characteristic polynomial,
\beq \label{eqn:phiN-int-def}
\Phi_N (z) := \log |\det (X-z) | + N\frac{(1- |z|^2 )_+}{2},
\eeq
which is a fundamental observable given its relationship with the spectrum of $X$. 

This question of fluctuations was first answered in the seminal work of Rider and Vir\'ag \cite{rider2007noise} that proved that $\Phi_N(z)$ converges to the Gaussian Free Field (GFF)  on the unit disc. The GFF $Y(z)$ is a  centered Gaussian process with covariance kernel\footnote{In the case of real matrices there is an additional term on the RHS which we omit in the introduction. In fact, \cite{rider2007noise} only proves the result for the complex Ginibre ensemble, where the entries are Gaussian. This  was later extended to normal matrices \cite{ameur2011random, ameur2015random}, general matrices with i.i.d. entries \cite{cipolloni2021fluctuation, cipolloni2023central, cipolloni2023mesoscopic}, and to a $(2+1)$-dimensional field  in \cite{bourgade2024fluctuations}.}
\beq \label{eqn:int-GFF-cov}
\mathrm{Cov} (Y(z) , Y(w) ) = - \frac{1}{2} \log |z-w| .
\eeq
It is clear from \eqref{eqn:int-GFF-cov} that $Y(z)$ has to be suitably interpreted as a generalized function and so 
the convergence of \cite{rider2007noise} is only in a weak sense. In particular, this form of convergence is insufficient to determine the more delicate properties of $\Phi_N (z)$. For example, the work \cite{rider2007noise} leaves aside any possible \emph{pointwise} convergence of $\Phi_N(z)$, or any understanding of the structure of its extreme values. 

Investigation of the extreme values of $\Phi_N (z)$ is intimately related to the fact that \cite{rider2007noise} connects  i.i.d. matrices to \emph{log-correlated fields.}  Log-correlated  fields such as the GFF appear  in diverse areas of mathematics, including statistical mechanics, SPDEs, conformal field theory, and analytic number theory. They also arise in the context of disordered systems and phase transitions \cite{carpentier2001glass,fyodorov2008freezing,fyodorov2009statistical} in the physics literature, as well as in mathematical finance \cite{duchon2012forecasting,bacry2008log}. A notable area of intense activity is the study of the extreme values or \emph{high points} of log-correlated fields. This is motivated in part by desires to extend the classical extreme value theory to random variables with strong correlation \cite{arguin2016extrema}.

Of particular relevance here are the  works of Fyodorov, Hiary and Keating \cite{fyodorov2012freezing,fyodorov2014freezing}, which have precipitated a series of developments in analytic number theory \cite{arguin2019maximum, arguin2020fyodorov, arguin2023fyodorov, harper2019riemann} and random matrix theory \cite{arguin2017maximum, chhaibi2019circle, paquette2018maximum, paquette2022extremal}. Emphasizing the log-correlated structure they proposed the characteristic polynomial of a random matrix as a model for the statistical properties of the Riemann $\zeta$-function, leading to a series of predictions about its high points. 
Among the first works to put the random matrix calculations on rigorous footing is a remarkable paper of Webb, proving that fractional powers of the characteristic polynomial of a Haar-distributed unitary matrix converge to an object called the Gaussian Multiplicative Chaos (GMC) \cite{webb2015characteristic}.\footnote{At least in a certain range of parameters, discussed in detail later.}

The foundations of the theory of the GMC were built by Kahane in the 80s \cite{kahane1985chaos}. Kahane was motivated in part by the modelling of intermittency in fluid mechanics \cite{kolmogorov1941local, kolmogorov1962refinement, mandelbrot1974multiplications, mandelbrot1974intermittent, obukhov1962some}. Given a  Gaussian field $Y(x)$ on $\rr^d$ with covariance similar to \eqref{eqn:int-GFF-cov} (e.g., the GFF in $d=2$), Kahane's theory gives rigorous meaning to the measure with density $\e^{ \gamma Y(x) } \d x$. A modern approach is to xregularize the field by mollification on scale $\eps$, and then take $\eps \to 0$ (along with a suitable renormalization) \cite{berestycki2017elementary}. One obtains a unique, non-trivial limit in the \emph{sub-critical phase}  $\gamma \in (0, 2 \sqrt{d})$.\footnote{Our convention differs from the GMC literature by a factor of $\sqrt{2}$ due to the random matrix scaling. For critical and supercritical $\gamma$ one can obtain a non-trivial measure through more delicate renormalization schemes.} The GMC is a fascinating object in its own right, being a random fractal measure with Hausdorf dimension $d- \frac{\gamma^2}{4}$. In particular it is both non-atomic and singular with respect to Lebesgue measure. The correspondence with random matrix theory is that since $\Phi_N(z)$ converges to a log-correlated Gaussian field, $|\det (X-z)|^\gamma$ would conjecturally converge to the GMC, up to a renormalization.

Given the ubiquity of log-correlated fields it is no surprise that the GMC has found numerous applications across mathematics. Famously, in $d=2$ it coincides with the volume form of Liouville Quantum Gravity \cite{ding2020tightness, gwynne2021existence, ding2023introduction}. It further arises in the context of random planar map scaling limits and the Brownian Map \cite{holden2023convergence, le2013uniqueness, miller2020liouville, duplantier2021liouville}, and in connection with Liouville Conformal Field Theory \cite{david2016liouville, kupiainen2020integrability}. We refer the reader to  surveys \cite{berestycki2024gaussian, rhodes2025two, sheffield2023random, rhodes2013gaussianmultiplicativechaosapplications} for further discussion in this direction, and the works \cite{powell2020criticalgaussianmultiplicativechaos,
bertacco2025does,bertacco2025uniqueness,duplantier2014critical} for various results and discussions of the critical and supercritical regimes.

 The paper \cite{webb2015characteristic} spurred further interest in connecting the GMC with random matrix theory. Indeed, Webb conjectured that this connection is ``quite universal.'' Webb's result for the Circular Unitary Ensemble (CUE) was extended to the entire subcritical range of $\gamma$ in \cite{nikula2020multiplicative}. The similar case of the Gaussian Unitary Ensemble (GUE) was handled in
\cite{berestycki2018random, claeys2021much}. Notably, there are only a few other cases in which convergence to the GMC is known. These include the   classical compact groups \cite{forkel2021classical}, the Circular $\beta$ Ensembles  ensembles \cite{chhaibi2025circle, lambert2024subcritical}, and a certain dynamical extension of the  CUE and GUE \cite{bourgade2206liouville,keles2025non}. In fact, convergence to the GMC is surprisingly not even known for the entire subcritical range of $\gamma$  for the classical Gaussian Orthogonal and Symplectic ensembles; only the partial results of \cite{kivimae2020gaussian} are available. Prior to ours, the only work treating non-Hermitian ensembles is the recent  \cite{bourgade2025fisher} that shows that the GMC measure appears in connection with general normal random matrices. This class of models includes the fundamental complex Ginibre ensemble, which is an i.i.d. matrix with complex Gaussian entries.

Despite the considerable effort since \cite{webb2015characteristic} the landscape of results is fragmentary. It is imperative to develop methods that can cover   large classes of generic examples. In particular, we emphasize that all results on convergence to the GMC concern only the so-called \emph{invariant ensembles.} These ensembles are rotationally invariant and  admit exact formulas for their joint eigenvalue distribution which are model-specific. A goal of our work is to address these challenges and develop probabilistic methods that can overcome the limitations of relying on invariance or exact formulas.

In Theorem \ref{theo:maingmcconv} we will show that the sequence of measures on $\dd$
\beq \label{eqn:int-det-measure}
\mu_N(z):=\nc\frac{ | \det (X -z ) |^\gamma}{ \ee[ | \det (X-z) |^\gamma]} \d z \d \bar{z}
\eeq
converges to the GMC in the entire subcritical phase for \emph{any} real or complex i.i.d. ensemble $X$. This is a considerable advance beyond any previous result, and our methods are a significant departure from prior works establishing convergence to the GMC.

Perhaps most importantly, the general class of ensembles we consider are non-invariant in the sense described above. At a technical level, all previous works on the GMC in random matrix theory   heavily leverage integrability or invariance. These properties allow one to import techniques specific to the models, such as the Riemann-Hilbert method,  orthogonal polynomials, 
or asymptotic analysis of explicit determinantal formulas.
These tools are unavailable for the ensembles we consider. Indeed, Webb \cite{webb2015characteristic} raised the question of ``what can one do in the two-dimensional case\footnote{I.e., in the non-Hermitian case.} or when $\beta\neq 2$ and a Riemann-Hilbert approach might not exist''. We develop a new dynamical approach to the GMC using  Dyson Brownian motion that we describe in Section~\ref{sec:methods} below.\footnote{While the models in \cite{bourgade2206liouville,keles2025non} are dynamical in nature, these works fundamentally rely on their specific determinantal structure} In particular, the only input we require from the exactly solveable models are the asymptotics  of $\ee[ |\det (X) |^\gamma]$ for the complex Ginibre ensemble which may be computed in an elementary way.

The strength of our dynamical approach allows us to prove our result for any i.i.d. matrix, invariant or not. This allows for an independent proof of the complex Ginibre case of \cite{bourgade2025fisher} (this being the only model in common). Moreover, we give the first GMC convergence for any i.i.d. matrix with real entries. This includes the real Ginibre ensemble (where the entries are real Gaussians) which, despite its invariance, has not yet been handled even by integrable methods. 

It is instructive to compare our result to the model considered in \cite{bourgade2025fisher}, being the closest precedent in the literature. The eigenvalue density of the unitarily invariant normal matrices considered there have a determinantal form, coinciding with the positions of free Fermions in the lowest Landau level. Convergence to the GMC is interpreted in \cite{bourgade2025fisher} as Liouville Quantum Gravity arising as a universal limit of a free Fermi gas with general potential. Our result shows that neither the Fermionic structure,  the invariance, or the complex symmetry class are essential features for the LQG limit to arise in random matrix theory. The matrix structure itself is the fundamental property leveraged by our methods.

The phenomenology of  real i.i.d. matrices is richer than the complex case, due to the symmetry of the spectrum about the real axis. This is exhibited by the fact that a second term arises in the analog of \eqref{eqn:int-GFF-cov}, introducing a second logarithmic singularity as $z, w \to \rr$. We investigate this phenomena by showing that the restriction of \eqref{eqn:int-det-measure} to the interval $(-1, 1)\subseteq \rr$ converges to a one-dimensional GMC. Indeed, the power of our approach is such that the changes between the complex and real, $2$-d and $1$-d cases are largely notational.

Aside from being an object of intrinsic interest, convergence to the GMC is partially motivated by studying the more delicate properties of the high points of $\Phi_N (z)$ (i.e., the predictions of \cite{fyodorov2014freezing,fyodorov2012freezing} as discussed above). In the next section we outline several applications, including determining the volume of high points as well as freezing phenomena. We also deduce pointwise central limit theorems for  $\Phi_N (z)$ from our results.


Our methods  allow us to compute the asymptotics of the $K$-point function of the characteristic polynomial evaluated at mesoscopically separated points and with fractional exponents (see Theorem \ref{thm:k-pt-intro}). Indeed, convergence to the GMC is more or less a by-product of computing the $2$-point function, due to the general framework of \cite{claeys2021much} which our work partially relies on. Our asymptotics are a significant generalization of several prior results in the literature which were only available for, e.g., matrices with Gaussian entries or with integer moments. In particular we prove (and generalize) a conjecture of Serebryakov and Simm \cite{serebryakov2025schur}. These aspects are discussed in further detail below.

\subsection{Definition of model and main results} 

In this section we will define our model and state our main results. The class of matrices we consider are the i.i.d. ensembles:

\bed[i.i.d. matrices] \label{def:iid} A real or complex i.i.d. matrix $X$ is an $N \times N$ random matrix whose entries $X_{ab}$ are i.i.d. random variables with distribution $X_{ab}\stackrel{\d}{=}N^{-1/2}\chi$, for some real or complex random variable $\chi$. 

We assume that the  random variable $\chi$ satisfies $\ee[\chi] = 0$ and $\ee[ |\chi|^2]=1$. In the complex case we further assume $\ee[ \chi^2]=0$. Additionally, we will always assume that for all $p \in \mathbb{N}$ there exists a $C_p >0$ so that,
\begin{equation}
\label{eq:moments}
\ee |\chi|^p \le C_p.
\end{equation}
Throughout, we will use the parameter $\beta$ to unify formulas that hold in the real and complex cases. Specifically, in the real case we set $\beta=1$ and in the complex case $\beta=2$. We  define the fourth cumulant $\kappa^{(\beta)}_4$ by,
\beq \label{eqn:kappa-def}
\kappa^{(\beta)}_4 := \ee[ N^2 |X_{ij}|^4]  - (4- \beta) .
\eeq 
To lighten the notation we will often write $\kappa_4 = \kappa_4^{(\beta)}$ when the context is clear. Note that if $\chi$ is Gaussian then $\kappa^{(\beta)}_4 = 0$. 
\eed
The primary object of interest in our work is the log-characteristic polynomial $\Phi_N (z)$ defined above in \eqref{eqn:phiN-int-def}. The random field $\Phi_N (z)$ behaves asymptotically like a Gaussian random field with covariance
\beq \label{eqn:int-cov}
\Kb (z, w ; \kappa_4) := - \frac{1}{2} \log |z-w| - \frac{ \1_{\beta=1}}{2} \log |z-\bar{w} | + \frac{ \kappa_4}{4} (1- |z|^2)(1-|w|^2).
\eeq
\noindent{\bf Notational convention.} \emph{In \eqref{eqn:int-cov}, the indicator function $\1_{\beta=1}$ means that the second term is not present in the complex i.i.d. case when $\beta=2$ but is present in the real i.i.d. case when $\beta=1$. We use this convention throughout the paper.}

\

If $\psi(z)$ is the Gaussian process on the unit disc with covariance kernel $\Kb (z, w ; \kappa_4)$, then the Gaussian multiplicative chaos is formally defined by
\beq
\d \mu_{\GMC, \kappa_4}^{\beta, \gamma} (z)  = \frac{ \e^{ \gamma \psi (z)}}{ \ee[ \e^{ \gamma \psi (z)}]} \d z \d \bar{z} .
\eeq
As we discussed earlier, it of course does not make sense to exponentiate a distribution. We postpone the rigorous definition of $\mu_{\GMC, \kappa_4}^{\beta, \gamma}$ to Section \ref{sec:GMC-def} below, following \cite{berestycki2017elementary}.

\bet\label{theo:maingmcconv} Let $X$ be an i.i.d. matrix. For $0 < \gamma < 2 \sqrt{2}$ we have that the random measures defined on the unit disc $\dd \subseteq \cc$ by 
\beq \label{eqn:maingmcconv}
\frac{ | \det (X -z )|^\gamma}{ \ee[ | \det (X-z)|^\gamma]} \d z \d \bar{z} 
\eeq
converge in distribution (with respect to the topology of weak convergence) to the measure $\mu_{\GMC, \kappa_4}^{\beta, \gamma }$

\eet

This result covers the entire subcritical phase of $\gamma$, and is new for real matrices even in the Gaussian case. In the complex case we recover a special case of \cite{bourgade2025fisher} for the complex Ginibre ensemble. Specifically, \cite{bourgade2025fisher} treats a class of invariant random normal matrices; one specific example coincides with the complex Ginibre ensemble. This is the only ensemble that our two works have in common.

In the case of real i.i.d. matrices, the symmetry of the spectrum across the real axis introduces the second term on the RHS of $K^\beta$ defined in \eqref{eqn:int-cov}. In particular, this causes a second logarithmic singularity as $\Im[z], \Im[w] \to 0$. However, the form of convergence in Theorem \ref{theo:maingmcconv} is insensitive to this behavior, due to the normalization by the expectation in \eqref{eqn:maingmcconv}. Let us explain this point. The expectation of the volume that the measure in \eqref{eqn:maingmcconv} assigns to any subset of the complex plane equals its Lebegue measure; therefore, one can remove a small neighbourhood of the real axis at an $o(1)$ error, thus removing the second logarithmic singularity on the RHS of \eqref{eqn:int-cov}. 

It is therefore of interest to find a different avenue by which to investigate the behavior near the real axis in the real i.i.d. case. There are many possible ways to approach this problem; we take perhaps the most natural and study the restriction of the measure in \eqref{eqn:maingmcconv} to the interval $(-1, 1)$. We will show it converges to a one-dimensional GMC. That is, if $\psi_\rr (x)$ is a Gaussian process on $(-1, 1)$ with covariance kernel $K^1 (x, y ; \kappa_4)$ then we formally define the GMC by $\d \mu^{\rr, \gamma}_{\GMC, \kappa_4} (x) = \frac{ \e^{ \gamma \psi_\rr (x) }}{ \ee[ \e^{ \gamma \psi_\rr (x) }]} \d x$ (again, see Section \ref{sec:GMC-def} below). 
\bet \label{thm:gmc-real} Let $X$ be a real i.i.d. matrix. For $0 < \gamma < \sqrt{2}$ we have that the random measures on the interval $(-1, 1)$ defined by
\beq
\frac{ | \det (X-x)|^\gamma}{ \ee[ | \det (X-x)|^\gamma]} \d x
\eeq
converge in distribution (with respect to the topology of weak convergence) to the measure $\mu^{\rr, \gamma}_{\GMC, \kappa_4}$. 
\eet
This corresponds again to the entire subcritical phase of $\gamma$. The same proof handles also the complex case which is of less interest and so is omitted.

As  stated in the introduction our proof  of convergence to the GMC relies on the framework of \cite{claeys2021much}, which requires one to compute (among other things) the asymptotics of the two point function of the characteristic polynomial. It is not any harder to handle the general case and 
the joint fractional moments are themselves of independent interest. 

This brings us to our second set of main results. In order to state them we first require some notation which will be used throughout the paper. The Barne's G-function $G(z)$ can be defined via the Weierstrass product
\beq
\label{eq:defG}
G(z+1)=(2\pi)^{z/2}\exp \left( \frac{-z+z^2(1+\gamma_{\mathrm{EM}})}{2} \right)\prod_{k\ge 1}\left(1+\frac{z}{k}\right)^k \e^{\frac{z^2}{2k}-z},
\eeq
where $\gamma_{\mathrm{EM}}$ is the Euler-Mascheroni constant. We will not use any specific properties of $G(z)$; it enters our work entirely through an exact computation for the Ginibre ensembles 
(see Proposition \ref{prop:ginibre} below).

We further require some notation for some specific subset of the unit disc. We denote the open disc of radius $r$ by $\dd_r := \{ z \in \cc : |z| <r \}$. To unify some notation we also introduce the sets $\ddb_r$, $\ddhb_r$ and $I_r$, for $\beta=1, 2$ and $r \in (0, 1)$, by,
\begin{align} \label{eqn:ddb-def}
\dd_r^2 := \dd_r, \quad \dd_r^1 := \dd_r \cap \{ z : \Im[z] > 1-r \} , \quad I_r := \{ z \in \rr : |z| <1 \} , \quad  \ddhb_r := \ddb_r \cup I_r.
\end{align}
Note that as $r\to 1$, the set $\dd^1_r$ exhausts the open half unit disc in the upper half plane. The significance of these sets will be clear momentarily.

\bet[$K$-point function] \label{thm:k-pt-intro} Fix any small $b >0$, any $L \geq 1$, and any $K \geq 1$. There are small $\mfc_1 = \mfc (b, K)$ and $\mfc_2 = \mfc (b, K, L)$ so that the following holds. 
Let $0 < r< 1$, let $z_1, \dots, z_K$ be points in $\ddhb_r$ such that
\beq
\min_{i \neq j} |z_i -z_j| \geq N^{-1/2+b},
\eeq
and let $\gamma_1, \dots , \gamma_k$ be points in $[0, L ] + \i [-\mfc_2, \mfc_2]$. Then,
\begin{align}
\label{eq:momasymp}
& \ee\left[ \prod_{i=1}^K | \det (X -z_i )|^{\gamma_i} \right] = (1 + \O ( N^{-\mfc_1} ) ) \prod_{i=1}^K \frac{\e^{\frac{\gamma_iN}{2}(|z_i|^2-1)}}{\e^{\frac{2\gamma_i-\gamma_i^2}{8}\kappa_4(|z_i|^2-1)^2}}N^{\frac{\gamma_i^2}{8}}\frac{(2\pi)^{\frac{\gamma_i}{4}}}{G\big(1+\frac{\gamma_i}{2}\big)} \notag\\
\times &\left( \prod_{1\leq i<j\leq K}\frac{\e^{\frac{\kappa_4\gamma_i\gamma_j}{4}(|z_i|^2-1)(|z_j|^2-1)}}{|z_i-z_j|^{\frac{\gamma_i\gamma_j}{2}}} \right) \times \e^{\1_{\beta=1}\EE} \times \prod_{i=1}^K  \left( \frac{ (2 \pi)^\frac{\gamma_i}{4} G(1/2)}{G \left( \frac{1+\gamma_i}{2} \right) } \right)^{\1_{\beta=1} \1_{ z_i \in \rr} },
\end{align}
where $G$ is the Barne's G function from \eqref{eq:defG} and
\beq
\EE := -\sum_{i, j=1}^K \frac{ \gamma_i \gamma_j}{8} \log ( |z_i  - \bar{z}_j |^2 \vee (4/N)) + \sum_{i=1}^K \frac{\gamma_i}{4} \log ( |z_i - \bar{z}_i|^2 \vee (2/N) ).
\eeq
The estimate is uniform over the choice of $z_i$ and $\gamma_i$ in the sets above.
\eet
We now discuss a few salient points about Theorem \ref{thm:k-pt-intro}. In our scaling, the interparticle distance is $N^{-1/2}$, and so our result applies to points separated by any mesoscopic distance. The contribution of the local statistics is the Barne's G function and arrives in product form. The correlations between points can be read off from \eqref{eqn:int-cov} and represent interactions at larger scales. Different  asymptotics  arise for microscopically separated $z_i$ \cite{deano2022characteristic}. 

Fractional moment asymptotics for i.i.d. matrices were  available for only a handful of cases prior to our work. The one-point function for the complex Ginibre ensemble was treated in \cite{webb2019moments}, and the general $K$-point case in \cite{bourgade2025fisher} (which also handled general normal matrices). The one-point function of the real Ginibre ensemble was available at the single point $z=0$ due to a connection with the Laguerre Orthogonal Ensemble \cite{serebryakov2025schur}. None of these results apply to the non-invariant ensembles we consider.

The application to real i.i.d. matrices of Theorem \ref{thm:k-pt-intro}  is particularly notable. The asymptotic behavior \eqref{eq:momasymp}  was not known even in the real Ginibre case. In fact, aside from the case $z=0$, no result for even the one-point function existed prior to our work. When $z_i$ is separated from the real axis  the local behavior is the same as the complex  case, reflecting the fact that the local statistics coincide with the complex Ginibre ensemble in this region. On the other hand, a different behavior arises for $z_i \in \rr$. The one-point formula  was conjectured by \cite[Corollary 4.9]{serebryakov2025schur} which is a consequence of Theorem \ref{thm:k-pt-intro}. Of course, Theorem \ref{thm:k-pt-intro} goes further, deriving joint asymptotics for mesoscopically separated points lying in $(-1, 1)$.

Let us comment now about the choice of domain 
in the real i.i.d. case. 
In this case the eigenvalues of $X$ are symmetric about the real axis and so the constraint that the $z_i$ lie in the upper half plane is artificial. The constraint that the $z_i$ must either be separated from the real line by a small constant or on the real axis is more stringent. We could significantly relax this constraint and allow $\Im[z_i] \to 0$ as $N \to \infty$   with a non-trivial amount of technical work; at a high level we expect our strategy  is applicable. 
However, these asymptotics  are not needed for our results on the GMC and in the interest of brevity we leave these various extensions for future work.

Various works treated the case of joint integer moments in the complementary regime $|z_i -z_j| \lesssim N^{-1/2}$. These appeared in \cite{maltsev2025bulk, osman2025bulk} for i.i.d. matrices with a Ginibre component and for general i.i.d. matrices in \cite{afanasiev2019correlation, afanasiev2022correlation}. Results on integer moments of the real Ginibre ensemble include  \cite{serebryakov2025schur,kivimae2024moments}. We remark here that integer moments are amenable to various techniques including supersymmetric methods; however, the fractional moments we handle require the  new approach of our paper. 

More generally, the moments of random characteristic polynomials for broad classes of random matrices have attracted significant interest. For example, integer moments have been studied in \cite{akemann2003characteristic, afanasiev2020correlation, borodin2006averages, forrester2009matrix, fyodorov2002negative, serebryakov2025schur} by algebraic and supersymmetric methods. Fractional exponents have been investigated in \cite{berestycki2018random, byun2025free, charlier2019asymptotics, charlier2025asymptotics, claeys2021much, deano2025asymptotics, its2007hankel, krasovsky2007correlations}  by Riemann-Hilbert methods. The work  \cite{bailey19}  studies moments of the moments of characteristic polynomials for CUE matrices. 
Motivated by connections to analytic number theory, moments of characteristic polynomials from the classical compact groups have been examined in \cite{assiotis2026exchangeable, assiotis2026joint, forkel2021classical}. Lastly, we mention  that the work \cite{keating2022critical}  considered ``critical moments'' of subcritical Gaussian multiplicative chaos for the CUE ensemble, and that  moments of characteristic polynomials beyond the typical GMC regime have also been recently considered for C$\beta$E ensembles \cite{assiotis2026moments}.

\subsubsection{Consequences}

In this subsection we list a few immediate consequences of our main results. The discussion largely follows \cite{bourgade2025fisher} and so we will be brief. From Theorem \ref{thm:k-pt-intro} we have the following pointwise convergence of the log-characteristic polynomial field. 
\bec[CLT for log-determinants]
Fix $b>0$, $0<r<1$, and fix any $z_1,\dots, z_K$, with  $z_i\in \ddhb_r$  and $\min_{i\ne j}|z_i-z_j|\ge N^{-1/2+b}$. Define
\beq
\Psi_N(z_i):=\frac{\log\big|\mathrm{det}(X-z_i)\big|-N\frac{1-|z_i|^2}{2} + \frac{\1_{z_i \in \rr} \1_{\beta=1}}{4} \log N}{ \sqrt{ \log N}} .
\eeq
Assume that the limits,
\begin{align}
C_{ij} &:= -\lim_{N \to \infty} \frac{1}{ \log N} \bigg\{  \frac{1}{4}\log (|z_i - z_j|^2+N^{-1}) + \frac{\1_{\beta=1}}{4}\log (|z_i - \bar{z}_j|^2+N^{-1})\bigg\}
\end{align}
exist. 
Then, $\big(\Psi_N(z_1),\dots, \Psi_N(z_K)\big)$ converges in distribution to a centered Gaussian vector with covariance matrix $C_{ij}$. 
\eec

We now state two immediate corollaries also relying on Theorem \ref{theo:maingmcconv}  (see also \cite[Corollaries 1.11 and 1.13]{bourgade2025fisher} for the complex Ginibre case). The argument in \cite[Proposition 3.8]{claeys2021much} yields the following.
\bec[Thick points]
\label{cor:thick}
For every compact $K\subset\mathbb{D}_1^\beta$ with non-empty interior, $\nu\in [0,1/\sqrt{2})$, and $\eps>0$,
we have
\beq
\label{eq:thickp}
\lim_{N\to +\infty}\mathbb{P}\left(N^{-2\nu^2-\eps}\le|\{z\in K :\Phi_N(z)\ge \nu\log N\}|\le N^{-2\nu^2+\eps}\right)=1.
\eeq
\eec
Using \eqref{eq:thickp} and following the same steps as in \cite[ Proof of Corollary 1.4]{arguin2017maximum}, we obtain the convergence for the free energy:
\bec[Freezing]
For every compact $K\subset\mathbb{D}_1^\beta$ with non-empty interior and $\gamma>0$, the free energy converges in probability as $N\to\infty$:
\beq
\lim_{N\to+\infty}\frac{\log\big(N\int_K \e^{\gamma\Phi_N(z)}\,\d z \d \bar{z} \big)}{\gamma\log N}=\begin{cases}
\frac{1}{\gamma}+\frac{\gamma}{8} &\mathrm{if}\quad \gamma\le 2\sqrt{2}, \\
\frac{1}{\sqrt{2}}&\mathrm{if}\quad \gamma> 2\sqrt{2}.
\end{cases}
\eeq
\eec

For simplicity we have excluded sets supported near the real axis in the above statements in the real i.i.d. case. Similar applications to subsets of $(-1, 1)$ are possible using Theorem \ref{thm:gmc-real}. With further work (see the discussion after Theorem \ref{thm:k-pt-intro}) one could extend the above results to subsets of the entire unit disc in the real case in an appropriate manner. See also the discussion of the real case in \cite{cipolloni2025maximum}. We leave such extensions for future investigation.

As mentioned in the introduction, convergence to the GMC is also  partially motivated by studying the maximum value of $\Phi_N (z)$. The complex Ginibre case was handled in \cite{lambert2020maximum} and the authors of the current article previously proved that
\beq
\lim_{N\to +\infty}\mathbb{P}\left[\max_{z\in \dd_r}\frac{\Phi_N(z)}{\log N}\in \left[\frac{1}{\sqrt{2}}-\eps,\frac{1}{\sqrt{2}}+\eps\right]\right]=1,
\eeq
for any $\eps >0$, for any i.i.d. matrix \cite{cipolloni2025maximum}. The upper bound is easier, being morally a union bound. In \cite{cipolloni2025maximum} we proved the lower bound using Kistler's modification of the second moment method \cite{kistler2014derrida}. Corollary \ref{cor:thick} gives an alternative proof of the more challenging lower bound. In the real case, Theorem \ref{thm:gmc-real} (together with the upper bound) yields that  the maximum value of $\Phi_N$ on $(-r, r)$ is $\frac{\log N}{\sqrt{2}} + o ( \log N)$. This was the one case not explicitly treated in \cite{cipolloni2025maximum} which had instead considered regimes where $z$ is mesoscopically separated from $\rr$, due to a connection with inhomogeneous branching random walks.

\subsection{Methods}
\label{sec:methods}

In this section we discuss the proofs of our main results, Theorems \ref{theo:maingmcconv} and \ref{thm:k-pt-intro}. For simplicity we assume $\beta=2$. Our proof of Theorem \ref{theo:maingmcconv} is based on a general framework for establishing the convergence of approximately Gaussian fields to the GMC developed in the work \cite{claeys2021much}. Here, the main inputs are a computation of the joint Laplace transform of $\Phi_N (z)$ evaluated at two mesoscopically separated points, together with a class of mesoscopically regularized processes that approximate the log-characteristic polynomial. 
The bulk of the technical work in proving Theorem \ref{theo:maingmcconv} is then the same as for proving Theorem \ref{thm:k-pt-intro}, and so we will mainly discuss the latter. 

The first observation going back to Girko is that
\beq
\label{eq:girkologdetnew}
\log | \det (X -z) | = \frac{1}{2} \log \det [ (X-z) (X^*-\bar{z} )] = \frac{1}{2} \Re \log \det ( H^z - \i 0 ),
\eeq
where $H^z$ is the Hermitization of $X$ defined by
\beq \label{eqn:Hermit-def}
H^z = H^z (X)= \left( \begin{matrix} 0 & X-z \\ X^* - \bar{z} & 0\end{matrix} \right).
\eeq
The spectrum of the $2N \times 2N$ Hermitian matrix $H^z$ is symmetric about the origin. We denote its non-negative eigenvalues by $\lambda_i^z$, for $i =1, \dots, N$, and $\lambda_{-i}^z = - \lambda_i^z$. To fix scales, it is helpful to keep in mind that the $i$-th eigenvalue $\lambda_i^z$ is located at (roughly) $\frac{i}{N}$ with fluctuations of order $N^{-1}$ (at least away from the spectral edges which are when $i = \pm N + o (N)$). 

This motivates the following regularization of the log characteristic polynomial
\begin{align}
    \Phi_N (z, \eta ) &=  \frac{1}{2} \Re \log \det ( H^z - \i \eta ) - ( \dots ) = \frac{1}{4} \sum_i \log ( (\lambda_i^z)^2 + \eta^2 ) - ( \dots )
\end{align}
for $\eta >0$. Here, $( \dots )$ is a leading order deterministic term coming from random matrix theory that approximately centers $\Phi_N (z, \eta)$ (see the precise definition in \eqref{eqn:Phi-def} below). We first discuss how $\Phi_N (z, \eta)$ helps us compute the one-point function (the special case of Theorem \ref{thm:k-pt-intro} where $K=1$) before turning to the general case later. 

The quantity $\Phi_N (z, \eta)$ is relatively insensitive to the behavior of the eigenvalues on scales much smaller than $ \eta$; in our scaling this corresponds to the eigenvalues $\lambda_i^z$ with $|i| \ll N \eta$. It therefore makes sense to decompose
\beq \label{eqn:intro-decomp}
\Phi_N (z, 0) = \left( \Phi_N (z, 0) - \Phi_N (z, \eta_*) \right) + \Phi_N (z, \eta_*),
\eeq
with $\eta_* = \frac{ ( \log N)^{C_*}}{N}$, for some sufficiently large $C_* >0$. The point is that the first term captures the local behavior of the small eigenvalues of $H^z$ and the second term the global behavior. The naive strategy is then to combine the two following guiding principles: (i) in random matrix theory the local behavior of the eigenvalues is known to be universal and therefore, the distribution of the first term in \eqref{eqn:intro-decomp} is universal and independent from the choice of $X$; (ii) the second term in \eqref{eqn:intro-decomp} lives on mesoscopic scales, and we therefore expect that we can compute its distribution using resolvent methods. 

While both of these facts are true as stated, they are alone insufficient to compute the one-point function (let alone the $K$-point function). This is due to the fact that there is still  correlation between the two terms in \eqref{eqn:intro-decomp} that is hard to compute directly, and even if the first term is universal, we still need to compute its Laplace transform. 

In order to overcome these basic obstacles we instead rely directly on the dynamical approach that is used to prove local universality in the first place. That is, we allow $X$ to evolve as a matrix-valued Brownian motion (Dyson Brownian motion or DBM) and promote $\Phi_N (z, \eta)$ to a three-parameter object $\Phi_N (z, \eta, t)$ (see Section \ref{sec:dynamics}). By tracking the evolution of $\Phi_N (z, \eta, t)$ in time and applying the homogenization theory of DBM (see \cite{landon2019fixed,che2019universality} and the more recent \cite{bourgade2021extreme, bourgade2024fluctuations}) we derive a \emph{local-global decomposition} (see Theorem \ref{thm:local-global} below) through the approximation
\begin{align} \label{eqn:int-lg-1}
\Phi_N (z, 0, t) &= ( \Phi_N (z, 0, t) - \Phi_N (z, \eta_* , t) ) + ( \Phi_N (z, \eta_*, t) - \Phi_N (z, \eta_0, 0) ) + \Phi_N (z, \eta_0, 0) \notag \\
& \approx L^{(1)} + L^{(2)} + G + o (1) =: L + G + o (1).
\end{align}
where $L$ captures the universal local fluctuations and $G$ the global fluctuations. The parameter $\eta_0 \approx N^{-1+c}$ is specifically chosen according to the limiting hydrodynamical equation describing DBM. Importantly, the two random variables $L$ and $G$ turn out to be independent and so
\beq \label{eqn:int-lg}
\ee[\e^{ \lambda \Phi_N (z, 0, t)}  ] = (1 + o (1) ) \ee[ \e^{ \lambda L} ] \ee[ \e^{ \lambda G} ].
\eeq
In particular, the factor $\ee[ \e^{ \lambda G} ]$ lives on the mesoscopic scale $\eta_0$ and is computed using resolvent methods; see the formulae in Section \ref{sec:meso-results}. This still leaves aside the problem of computing the local factor $\ee[ \e^{ \lambda L} ]$. The exact form of universality that we obtain is that $L$ is independent of not only the distribution of $X_{ij}$ but also the spectral parameter $z$ (or at least the dependence on $z$ is essentially by scaling). In particular, it coincides with the case where the LHS of \eqref{eqn:int-lg} comes from the complex Ginibre ensemble and $z=0$. In this case, the LHS of \eqref{eqn:int-lg} is readily computed using either Kostlan's theorem \cite{kostlan1992spectra} or relating it to the Laguerre Unitary Ensemble and applying Selberg's integral formula. This allows for the computation of the local factor $ \ee[ \e^{ \lambda L } ]$ which as we mentioned is universal. We remark that this is the only use of any exact formulas or integrability in our work. We do not rely on any computation for the complex Ginibre ensemble at nonzero $z$.

Our computation of the $K$-point function is based on generalizing the local-global decomposition of \eqref{eqn:int-lg-1} to
\beq \label{eqn:int-lg-2}
\Phi_N (z_j, 0, t) \approx L_j+G_j .
\eeq
The main issue now is to determine the joint distribution of all the local and global factors together. For this we rely on the generalization of the homogenization theory of \cite{cipolloni2021fluctuation, cipolloni2023central}. In these works the crucial observation was that if the $z_j$ are well-separated (i.e. $|z_j -z_k| \gg N^{-1/2}$) then the dynamics driving the smallest eigenvalues $\lambda^{z_j}_i (t)$ are approximately independent and can therefore be coupled to independent processes at the cost of a negligible error. We apply this fact and derive (in the $K$-point version of our local-global decomposition in Theorem~\ref{thm:local-global})  a representation \eqref{eqn:int-lg-2} where the $L_j$ are independent from each other and from the $G_j$. Therefore,
\beq
\ee \left[ \e^{ \sum_{j=1}^K \lambda_j \Phi_N (z_j, 0, t) } \right] =(1 + o (1) )  \left( \prod_{j=1}^K \ee[ \e^{ \lambda_j L_j }] \right) \ee \left[ \e^{ \sum_{j=1}^K \lambda_j G_j} \right].
\eeq
The local factors were already determined in our computation of the one-point function. The second factor concerns random variables on mesoscopic scales and so their distribution may be computed using resolvent methods. This computes the $K$-point function.

Neglected in the above discussion is that the DBM introduces a small Gaussian component to $X$, which must be removed. We achieve this through moment matching methods analogous to the Green's function comparison theorem or four moment theorem of Tao and Vu \cite{tao2011random}. While this is not a significant departure from previous results in the literature, there are a few technical subtleties regarding the $K$-point function which we omit from the present discussion, especially in terms of a certain sub-microscopic regularization related to what we will now briefly discuss (see in particular Proposition \ref{prop:gfct-kpt-apriori}).

Finally, we mention that as written the decompositions \eqref{eqn:int-lg-1}, \eqref{eqn:int-lg-2} are not quite what we prove. One must instead work with, on the LHS, the quantity $\Phi_N (z, \eta_m, t)$ where $\eta_m = N^{-1-\delta_m}$ for some small $\delta_m >0$, as $\Phi_N (z, 0, t)$ is hard to deal with directly. The approximation of $\Phi_N (z, 0, t)$ by $\Phi_N (z, \eta_m, t)$ is not in the very high probability sense and is  achieved in expectation only. That is, we directly show that the $K$-point function involving $\Phi_N (z, 0, t)$ is well approximated by the $K$-point function of $\Phi_N (z, \eta_m, t)$. This argument is based in part on level repulsion estimates (controlling the probability that $\lambda_1^z$ is small) as well as a simple deterministic observation Lemma \ref{lem:deterministic}. This is one source of the restriction that $\Re[\gamma] \geq 0$ in Theorem \ref{thm:k-pt-intro}. The extension to some range of $\gamma <0$ would require additional arguments that we leave for  future investigation (and would not be true for all matrices as written, as if the matrix entries are discrete, $| \det (X-z ) | = 0$ with positive probability for some values of $z$).

\subsection{Definition of GMC} \label{sec:GMC-def}

\subsubsection{Definition on the disc}
In this section we define the limiting objects of Theorem \ref{theo:maingmcconv}, the Gaussian Multiplicative Chaos (GMC). We will use the approach of \cite{berestycki2017elementary}. Recall the definition of $\kappa^{(\beta)}_4$ in \eqref{eqn:kappa-def} and the kernel $\Kb$ in \eqref{eqn:int-cov}.  Note that $\kappa^{(1)}_4 \geq -2$ and $\kappa^{(2)}_4 \geq -1$.
Throughout our work we will see that the kernel $\Kb ( \cdot ; \kapb )$ is the limiting covariance kernel of the process $\Phi_N (z)$. In Appendix~\ref{sec:posker} we show that $\Kb (z, w ; \kappa)$ is always non-negative definite for $\kappa \geq - \frac{4}{\beta}$.

For simplicity we extend $\Kb$ to a kernel on all of $\cc$ by setting it identically $0$ if either $z$ or $w$ are not in $\dd$. We first discuss the simpler case $\beta =2$. Let $\kappa \geq -1$. As in, e.g., \cite{berestycki2017elementary} one can construct a Gaussian process $\psi (z)$ living in the Sobolev space $H^{-s} ( \cc)$ for any $s >0$, with covariance kernel $K^{(2)} (z, w ; \kappa)$. Then, for any $\gamma \in [0, 2 \sqrt{2} )$ the work \cite{berestycki2017elementary} constructs the GMC measure $\d \mu_{\GMC, \kappa}^{(2), \gamma}$ as the limit (in the topology of weak convergence of measures on $\dd$) in probability of the measures,
\beq
\d \mu^\gamma_\eps (z) := \frac{ \e^{ \gamma \psi_\eps (z)}}{ \ee[ \e^{ \gamma \psi_\eps (z) } ]} \d z \d \bar{z},
\eeq
where $\psi_\eps$ is the convolution of $\psi$ with a smooth radial mollifier on scale $\eps$. In particular, the limit is unique and independent of the choice of the mollifier. 

We now discuss the case $\beta=1$.   Let $\kappa \geq -2$. Due to the additional singularity near the real axis, the kernel $K^{(1)}$ does not quite fit into the standard framework on the GMC (to our knowledge). However, due to the symmetry of the process and the fact that the real line is a small set (and is therefore negligible in the same sense that the boundary of the disc is) we can reduce its construction to that of \cite{berestycki2017elementary}. We now turn to a few details.

If $\tilde{\psi} (z)$ is the process constructed above with kernel $K^{(2)} (z, w ; \frac{\kappa}{2} )$ then $\psi (z) := \frac{ \tilde{\psi} (z) + \tilde{\psi} (\bar{z} ) }{ \sqrt{2}}$ has covariance kernel $K^{(1)} (z, w ; \kappa )$. Since $\psi (z) = \psi ( \bar{z} )$ it suffices to construct $ \d \mu^{(1), \gamma}_{ \GMC, \kappa}$ on the half disc $\dd^{(1)} := \dd \cap \{ z : \Im[z] >  0 \}$, setting $ \mu^{(1), \gamma}_{\GMC, \kappa} ( \{ z : z \in (-1, 1) \} ) = 0$ by definition. The measure $\d \mu^{(1), \gamma}_{\GMC, \kappa}$ can be constructed via the approach of \cite{berestycki2017elementary}, i.e., via the radial mollifiers. Let us briefly explain why the singularity near the real axis does not play a role in this construction. With $\mu^\gamma_\eps$ as above one has 
\beq
\ee[ \mu^\gamma_\eps ( \dd^{(1)} )] = \frac{\pi}{2} , \qquad\quad \ee[ \mu^\gamma_{\eps} ( \{ z : \mathrm{dist} (z, \del \dd^{(1)} ) < \delta \} )] \leq C \delta,
\eeq
for some $C >0$ independent of $\eps >0$. Due to the above estimates and the choice of topology of weak convergence of measures on $\dd^{(1)}$ it suffices (see in particular the argument in \cite[Section 6]{berestycki2017elementary}) to show that $\mu^\gamma_\eps (A)$ converges in probability for sets such that $\mathrm{dist} (A, \del \dd^{(1)} ) >0$. The singularity of $\Kb$ near the real line does not play a role here as the singularity remains bounded for $z \in A$. The arguments in \cite{berestycki2017elementary} go through without change for the convergence of $\mu^\gamma_\eps (A)$ for such sets $A$. 

\subsubsection{Definition on the interval}

In this case, we can consider the kernel $K^{(1)} (x, y ; \kappa)$ as a kernel on $(-1, 1)$. We set it to be $0$ outside the unit interval. As the restriction of  $K^{(1)} (z, w ; \kappa)$ to $\rr$, it is not hard to check that $K^{(1)} (x, y ; \kappa)$ is also non-negative definite. Therefore, the construction of $ \d \mu_{\GMC, \kappa}^{\rr, \gamma}$ as a limit of mollified measures is as in \cite{berestycki2017elementary} (see also the $\beta=2$ case above).

\subsection{Notation and conventions} \label{sec:notation} For integers $k\in\mathbb{N}$ we use the notation $[\![k]\!] := \{1,2,...,k\}$, and by $\rr_{\ge 0}$ we denote the non-negative real numbers $\{x\in \rr:x\ge 0\}$.
For $N$ dependent non-negative quantities $a_N$, $b_N$ we say that $a_N \lesssim b_N$ if $a_N \leq C b_N$ for an $N$-independent constant $C>0$, which depend only on the constants appearing in \eqref{eq:moments}. If the $a_N$ and $b_N$ depend on some parameter $i \in I$ we will specify that the constant $C$ is uniform over $i$. We write $a_N \asymp b_N$ if $a_N \lesssim b_N \lesssim a_N$.

We denote vectors by bold-faced lower case Roman letters ${\bf x}, {\bf y}\in \cc^d$, for some $d \in \mathbb{N}$, and their scalar product by
\[
\langle {\bf x}, {\bf y}\rangle := \sum_{i=1}^d \overline{x_i}y_i.
\]
For any $d\times d$ matrix $A$ we use the notation $\langle A\rangle := d^{-1}\mathrm{Tr}[A]$ to denote the normalized trace of $A$, and $A^\mathfrak{t}$ to denote its transpose.
We denote the $d$-dimensional identity matrix by $I = I_d$.

We say that an event holds “with overwhelming probability” if for any fixed $D>0$ the probability of the event is bigger than $1-N^{-D}$ if $N\ge N_0(D)$, with $N_0(D)$ possibly depending on the constants appearing in \eqref{eq:moments} of the definition of our model. 

We define the $2N \times 2N$ block matrices
\beq
E_1 := \left( \begin{matrix} 1 & 0 \\ 0 & 0 \end{matrix} \right), \qquad E_2 := \left( \begin{matrix} 0 & 0 \\ 0 & 1 \end{matrix} \right)
\eeq
We will often need to sum over certain combinations of $E_1$ and $E_2$, and so we introduce the notation,
\beq
\tilde{\sum}_{ij} A(i, j) := A(1, 2) + A(2, 1)
\eeq

For real-valued martingales $M_t, N_t$, we denote the covariation process by $\d [ M_t, N_t]$. For complex valued martingales $M_t = X_t+ \i Y_t, N_t = P_t + \i Q_t$ the covariation process is defined by, $\d [ M_t, N_t] := \d [X_t, P_t] - \d [Y_t, Q_t] + \i ( \d [Y_t, P_t] + \d [ X_t, Q_t])$. The total variation process of a real-valued martingale is denoted by $\d [M_t] := \d [ M_t, M_t]$. 

\subsection{Organization of paper}

Section \ref{sec:prelim} collects various preliminaries about the Hermitization of an i.i.d. matrix, the local laws as well as some basic estimates that will be used throughout the paper. In Section \ref{sec:meso-results} we stateasymptotics for the joint Laplace transform of the characteristic polynomial regularized on \emph{mesoscopic} scales; the proofs appear in Section \ref{sec:meso}. In Section \ref{sec:dynamics} we apply Dyson Brownian motion in order to establish our key local-global decomposition for the characteristic polynomial evaluated at mesoscopically separated $z_i$. Sections \ref{sec:1pt} and \ref{sec:kpt} contain important calculations concerning the $1$- and $K$-pt functions, respectively. Section \ref{sec:gmc-reg} is a short section concerning regularizing the observables which we later show converge to the GMC. In Section \ref{sec:claeys} we apply the framework of \cite{claeys2021much} to establish a preliminary convergence to the GMC (that of a regularized quantity for Gaussian divisible matrices). In Section \ref{sec:gfct} we establish Green's function comparison theorems for the various observables we consider which allow us to remove the Gaussian component introduced by DBM. It is convenient (for organizational purposes) to provide the proofs of our main results in this section. The proofs of Theorems \ref{theo:maingmcconv} and \ref{thm:gmc-real} are given at the end of Section \ref{sec:GMC-GFCT}. The proof of Theorem \ref{thm:k-pt-intro} in the case where $\gamma_i \in \rr$ is given at the end of Section \ref{sec:gfct-kpt}. The extension to complex $\gamma_i$ is discussed in Section \ref{sec:complex}.

This a long version of the paper submitted for publication and contains several longer proofs that were omitted from the submitted version in the interest of brevity. In the short version of our paper, this version is cited as \cite{long-version}.

\subsubsection{An organizational note on the proofs}

The bulk of our work concerns various estimates about the log-characteristic polynomial $\Phi_N (z)$ defined in \eqref{eqn:phiN-int-def}. There are three distinct cases that we deal with in our work: (i) $X$ is complex and $z \in \dd_r$; (ii) $X$ is real and $z \in \dd_r \cap \{ \Im[z] > c \}$ for some $c>0$; (iii) $X$ is real and $z \in (-r, r)$. We treat each of these cases using the same strategy; for the most part the differences between the cases is notational. The complex case is the simplest. In the second case there are various additional terms that arise in various calculations and estimates. We will use the indicator function $\1_{\beta=1}$ to unify certain formulas where these terms are present. If $\Im[z] \geq c$ then these terms are typically $\O (1)$ or vanishing. The third case when $\beta=1$ and $z \in \rr$ is the most complex; now the additional terms are often logarithmic. Recalling our notation of $\ddb_r$ and $\ddhb_r$ in \eqref{eqn:ddb-def} we will often distinguish between the third case and the other two by writing things like ``In the case that $z \in \ddb_r$, ...'' noting that if $\beta=1$, we have $\ddb_r \cap \rr = \emptyset$.

Finally, to understand the various sizes of terms in expressions, it is helpful to note that when $z \in \ddb_r$ we expect that the random variable $\Phi_N (z)$ is approximately a centered Gaussian with variance $\frac{1}{4} \log N$, up to $\O (1)$ terms. In the case that $\beta = 1$ and $z \in I_r$, then instead the expectation of $\Phi_N (z)$ is approximately $- \frac{1}{4} \log N$, and its variance is $\frac{1}{2} \log N$, again up to $\O (1)$ corrections.

\medskip

\noindent{\bf Acknowledgements.} The research of G. C. is supported by the Italian Ministry of University and Research (MUR) - Fondo Italiano per la Scienza (FIS3) - 2024 Call, project UBLOCO, CUP F53C25000940001, and also partially supported by the MUR Excellence Department Project MatMod@TOV awarded to the Department of Mathematics, University of Rome Tor Vergata, CUP E83C23000330006. Additionally, G.C. thanks INdAM (Istituto Nazionale di Alta Matematica “Francesco Severi”) and the group GNFM. The research of B.L. is partially supported by an NSERC Discovery Grant and a Connaught New Researcher award. We also thank Ievgenii Afanasiev for useful discussions on \cite{afanasiev2019correlation, afanasiev2022correlation}, Theo Assiotis for bringing to our attention the references \cite{assiotis2026exchangeable, assiotis2026joint, assiotis2026moments},  and Jon Keating for bringing to our attention the reference\cite{forkel2021classical}. B.L. warmly thanks Philippe Sosoe for helpful discussions.

\section{Preliminaries} \label{sec:prelim}

In this section we introduce notations and state several preliminary results that will be used throughout the paper. In particular, in Section~\ref{sec:MDE} we describe properties of the  Hermitization $H^z$ of $X$ and its limiting eigenvalue density. In Section~\ref{sec:locallaws} we state several local laws for the resolvent of $H^z$ and state a rigidity bound for the eigenvalues. Next, in Section~\ref{sec:ssv} we recall an estimate stating that the smallest singular value of $X-z$ is not unusually small with high probability. In Section~\ref{sec:extres}, we state a  result on the moments  of the characteristic polynomial of the complex Ginibre ensemble, as well as some a priori estimates on integer moments of characteristic polynomials of i.i.d. matrices. Finally, in Sections~\ref{sec:regquant}--\ref{sec:cutoff} and \ref{sec:detreg} we show that it is possible to regularize the characteristic polynomial on submicroscopic scales and that we can cut-off its high values at the price of a negligible error.

\subsection{Hermitization and MDE}
\label{sec:MDE}

 Girko's formula \cite{girko1985circular} is a fundamental tool relating the spectral statistics of a non-Hermitian matrix $X$ with a family of Hermitian matrices. More precisely, for any matrix $X\in\cc^{N\times N}$ and complex parameter $z \in \cc$ we recall the definition of $H^z= H^z (X)$ in \eqref{eqn:Hermit-def}.

As a consequence of the $2\times 2$ block structure, the spectrum of $H^z$ is symmetric with respect to $0$ (\emph{chiral symmetry}). In particular, we denote the non-negative eigenvalues of $H^z$ by $\{\lambda_i^z\}_{i\in [\![N]\!]}$ in increasing order and for negative indices define $\lambda_{-i}^z=-\lambda_i^z$.

In the context of log-determinants, Girko's Hermitization formula is given by \eqref{eq:girkologdetnew}. 
We define the resolvent and empirical Stieltjes transform of $H^z$ by,
\beq
G^z (w) = \frac{1}{H^z -w }, \qquad m_N^z (w) = \langle G_z (w) \rangle  = \frac{1}{2N} \tr \left( \frac{1}{H^z - w } \right) ,
\eeq
for $w \in \cc\setminus\rr$.
It is well known that in the large $N$ limit the resolvent $G^z$ becomes approximately deterministic. Its limit is given by the $2N \times 2N$ matrix $M^z (w)$, which has the $2\times 2$ block form (each block is $N \times N$)
\beq
\label{eq:bigM}
M^z (w) := \left( \begin{matrix} m^z (w) & - z u^z (w) \\ - \bar{z} \bar{u}^z (w) & m^z (w)\end{matrix}\right), \qquad\quad u^z(w):= \frac{m^z(w)}{w+m^z(w)}.
\eeq
Here $m^z(w)$ is the unique solution to
\beq
\label{eq:mde}
-\frac{1}{m^z(w)}=w+m^z(w)-\frac{|z|^2}{w+m^z(w)}, \qquad\quad \Im[ m^z (w)] \Im[w] > 0.
\eeq
We stress that on the imaginary axis $w=\i\eta$ the quantity $m^z(w)$ is purely imaginary, and consequently $u^z(w)$ is real. Given the solution to \eqref{eq:mde}, we can recover the limiting eigenvalue distribution of $H^z$ via
\beq
\label{eq:behrho}
\rho^z (x) := \lim_{\eta \to 0^+} \frac{ \Im[ m^z (x+ \i \eta)]}{\pi}.
\eeq
Next, we state some basic properties of $\rho^z$ which will be used throughout the paper (cf. \cite[Lemma 2.4]{cipolloni2025maximum}):
\bel
\label{lem:proprho} Fix $0 < r < 1$. Let $\rho^z(x)$ be the density defined in \eqref{eq:behrho}. Uniformly in $z$ satisfying $|z| \leq r$ we have,
\begin{enumerate}[label=(\roman*),font=\normalfont]

    \item The density $\rho^z$ is symmetric, and its support is given by $[-\mathfrak{e}_z, \mathfrak{e}_z]$ for an explicit $\mathfrak{e}_z>0$. In particular,  it consists of a single interval.

    \item The edge $\mathfrak{e}_z$ satisfies the bound $C^{-1}\le\mathfrak{e}_z\le C$, for some $C>0$.

    \item \label{it:rho-der} For any small $\delta>0$ we have that for $|x|\le \mathfrak{e}_z-\delta$ that $\rho^z(x)\asymp 1$ and $| \partial_x\rho^z(x)| \lesssim 1$.

    \item For any small $\delta >0$ we have for $|x|\leq \mathfrak{e}_z-\delta$ and $0<\eta\le \delta^{-1} $ that $\Im m^z(x+\i \eta)\asymp 1$.
    
\end{enumerate}
\eel
\proof All the statements of this lemma are already in \cite[Lemma 2.4]{cipolloni2025maximum}, with the exception of the second statement of \ref{it:rho-der}. This statement instead follows by solving \eqref{eq:mde} explicitly using Cardano's formula (see also the display above \cite[Eq. (3.14)]{cipolloni2024precise} and \cite[Proposition 3.2 (ii)]{erdHos2020cusp}, \cite[Remark 7.3]{alt2020dyson} for similar derivations).
\qed

\

Before concluding this section we also state an asymptotic expansion  for $u^z, m^z$ on the imaginary axis valid for small $\eta$ (see e.g. \cite[Eq. (6.15)]{cipolloni2025maximum}).
\bel
We have
\beq
\label{eq:exp}
u^z(\i\eta)=1-\frac{\eta}{\sqrt{1-|z|^2}}+\mathcal{O}(\eta^2), \qquad\quad m^z(\i\eta)=\i \sqrt{1-|z|^2}+\mathcal{O}(\eta).
\eeq
\eel
We will also use the Ward identity,
\beq \label{eqn:ward}
\sum_{i} | G^z_{ia} ( \i \eta )|^2 = \frac{ \Im[ G^z_{aa}]}{\eta} 
\eeq
which follows in elementary way from the spectral theorem.

\subsection{Local laws and rigidity}
\label{sec:locallaws}

From \cite[Theorem 3.1]{cipolloni2021fluctuation} we have,
\bet[Local laws I] \label{thm:ll-1} Let $0 < r < 1$. Denote $w = E + \i \eta$.  Uniformly in $|z| \leq r $, for any unit vectors $\bx, \by$, matrices $A$, and any $\eps >0$, we have
\beq
| \langle \bx , (G^z (w) - M^z (w) ) \by \rangle| \leq N^{\eps} \left( \sqrt{\frac{ \Im[ m^z(w)]}{N \eta}} + \frac{1}{N \eta} \right)
\eeq
and
\beq
| \langle (G^z (w) - M^z (w) A \rangle | \leq \| A \| \frac{N^{\eps}}{N \eta},
\eeq
with overwhelming probability, for all $|w| \leq \eps^{-1}$ and $\eta>0$. 
\eet

Next, we state a local law for the averaged trace of the resolvent, which improves on the polynomial error rate.

\bet[Local laws II] \label{thm:ll-2} Let $0 < r < 1$ and $\delta >0$. Uniformly in $|z| \leq r$ we have
\beq
| m_N^z (\i \eta ) - m^z ( \i \eta ) | \leq \frac{ ( \log N)^{1/2+\delta}}{N \eta} 
\eeq
with overwhelming probability, for all $0 < \eta < N^{-\delta}$. 
\eet
\proof The estimate for $\eta \geq ( \log N)^{1/2+\delta}/N$ follows from \cite[Proposition 3.2]{cipolloni2025maximum}. For smaller $\eta$ we use the fact that $m_N^z ( \i \eta ) = \Im [ m_N^z ( \i \eta)]$ and that $\eta \mapsto \eta \Im[ m_N^z( \i \eta)]$ is increasing. \qed

\

We define the semiclassical eigenvalue locations (or \emph{quantiles}) by
\beq \label{eqn:quantiles-def} 
\frac{i}{2N} = \int_0^{\gamma_i^z} \rho^z (x) \d x
\eeq
for $1 \leq  i \leq N$ and $\gamma_{-i}^z = - \gamma_i^z$.

\bet[Rigidity estimates] \label{thm:rigidity}
Let $0 < r < 1$ and $C$, $\delta >0$ be given. Uniformly in $|z| < r$ the following estimates hold with overwhelming probability. We have
\beq
\label{eq:complrig}
| \lambda_i^z - \gamma_i^z | \leq \frac{ ( \log N)^{1/2+\delta}}{N} + \1_{ \{ |i| > ( \log N)^C \}} \frac{ ( \log N)^{3/2+\delta}}{N} + \1_{ \{ |i| > N^{1-\delta} \}} \frac{N^{\delta}}{N^{2/3} (N+1 - |i| )^{1/3}}
\eeq
\eet
\proof The bound \eqref{eq:complrig} follows by \cite[Corollary 3.3]{cipolloni2025maximum} together with standard rigidity estimates (see e.g. \cite[Eq. (7.46)]{cipolloni2021fluctuation}). \qed

\subsection{Smallest singular value estimate}
\label{sec:ssv}

Here we recall a well known bound on the smallest singular value of matrices with i.i.d. entries:
\bel \label{lem:wegner} Let $X$ be an i.i.d. matrix.  There is a $c_W>0$ so that the following holds. Let $0 < r < 1$. There exists $C_r >0$ so that for $z \in \ddhb_r$ we have
\beq
\label{eq:ssvbound}
\pp\left[ \lambda_1^z \leq \frac{s}{N} \right] \leq C_r\left(  s^2 + s \1_{\beta=1} \1_{ z \in \rr}\right)+  N^{-c_W}. 
\eeq
\eel
\proof 
By \cite[Eq. (4a)]{cipolloni2020optimal}, for real or complex Ginibre matrices the bound \eqref{eq:ssvbound} follows (even without the $N^{-c_W}$ correction). Then, for general i.i.d. matrices \eqref{eq:ssvbound} follows by a combination of the Dyson Brownian motion (DBM) and the Green's function comparison argument (see e.g. \cite[Remark 4.2]{cipolloni2023central}; in the real case we use the DBM from \cite[Section 7]{cipolloni2021fluctuation} instead of \cite[Section 7]{cipolloni2023central}). We point out that an analogous bound also directly follows from \cite[Theorem 3.2]{tao2010randoms}.
\qed

\subsection{Definitions of regularized quantities}
\label{sec:regquant}

For any matrix $X$ we define,
\begin{align} \label{eqn:Phi-def}
\Phi_N (z, \eta ; X) := \frac{1}{4} \sum_{i=-N}^N \log ( ( \lambda_i^z)^2 + \eta^2 )  -\frac{N}{2} \int \log ( x^2 + \eta^2 ) \rho^z (x) \d x,
\end{align}
where $\lambda_i^z$ are the eigenvalues of the Hermitization $H^z(X)$ as in \eqref{eqn:Hermit-def}. 
In particular, 
\beq
\Phi_N (z, 0; X) = \log |\det (X-z) | - \frac{N}{2} (|z|^2 - 1 ),
\eeq
 and so $\Phi_N (z, \eta ; X)$ is a regularized version of the log-characteristic polynomial. Note that,
\beq \label{eqn:deterministic-cont}
\left| \del_\eta \int \log ( x^2 + \eta^2 ) \rho^z (x) \d x  \right| \lesssim 1,
\eeq
uniformly  for $|z| \leq r$ for any $r \in (0, 1)$. Indeed, the derivative $\partial_\eta$ is exactly $\Im[ m^z ( \i \eta)]$.

We warn the reader here that later we will need to allow $X$ to evolve dynamically and so $\Phi_N$ will gain a third time parameter $\Phi_N(z, \eta, t)$. In this case $\Phi_N(z, \eta, t ) \neq \Phi_N (z, \eta ; X_t)$ as they mean slightly different things (hence the use of the semicolon). 

\subsection{Cut-off functions}
\label{sec:cutoff}

As we often deal with Laplace transforms of various random variables we will need to introduce various cut-offs. Throughout the paper, we fix smooth monotonic functions $h_1, h_2 : \rr \to \rr$ such that
\beq \label{eqn:cut-off-def}
h_1 (x) = \begin{cases} x , & x \leq 1 \\
1.5 , & x > 2 \end{cases} , \qquad h_2 (x) = \begin{cases} x , & |x| \leq 1 \\ 1.5 \sgn (x) , & |x| >2 \end{cases} 
\eeq
We may assume that $|h_i(x)| \leq |x|$ for all $x \in \rr$ and $i=1, 2$. We will often consider quantities such as $(A \log N) h_1 \left( \frac{Z}{A \log N} \right)$ for random variables $Z$. Then this cut-off random variable coincides with $Z$ for $Z \leq A \log N$ and is above bounded by $1.5 A \log N$. In particular, taking the exponential of this quantity results in a random variable bounded above by $N^C$, which can be controlled using overwhelming probability bounds.  Similar considerations apply to cut-offs using $h_2$. 

\subsection{Integer moments of characteristic polynomials}
\label{sec:extres}

We state the following asymptotics for the fractional moments of the complex and real Ginibre ensemble. The proofs follow from elementary methods; in the complex case it can be deduced from Kostlan's theorem \cite{kostlan1992spectra} which states that the absolute values of the eigenvalues of the complex Ginibre ensemble have the same distribution as $N$ i.i.d. Gamma distributed random variables. Both cases may also be proven using the connections to the LOE/LUE and Selberg's integral formula. The complex case is \cite[Lemma 3.7]{bourgade2025fisher} and the real case is \cite[Theorem 4.8]{serebryakov2025schur}. We review the proofs in Appendix \ref{a:moments}. 
 
\bep \label{prop:ginibre}  Let $C>0$. Uniformly in $\gamma \in [0, C] + \i [-C, C]$ we have: 
\begin{enumerate}[label=(\roman*),font=\normalfont]
\item 
If $X$ is the complex Ginibre ensemble,
\beq \label{eqn:exact-complex}
\ee\left[ | \det (X) |^\gamma \right] = N^{-N \frac{\gamma}{2}} \frac{G(N+1+\frac{\gamma}{2} )}{G(N+1) G(1 + \frac{\gamma}{2} )} = \e^{ - \frac{\gamma}{2} N} N^{\gamma^2/8} \frac{ ( 2 \pi )^{\gamma/4}}{G(1 + \frac{\gamma}{2})} \left( 1+ \O ( N^{-1}  ) \right)
\eeq
\item If $X$ is the real Ginibre ensemble,
\begin{align} \label{eqn:exact-real}
\ee[ | \det (X)|^\gamma] &= N^{-\frac{\gamma}{2} N} 2^{ \frac{\gamma N}{2}} \frac{ G(1/2) G(N/2 + \gamma/2+1/2)G(N/2 + \gamma/2 + 1)}{G(1/2+\gamma) G(1+\gamma)G(N/2+1/2)G(N/2+1)} \notag\\
&= \e^{ - \frac{N \gamma}{2}} \left(\frac{N}{2}\right)^{\frac{\gamma^2}{4} - \frac{\gamma}{4}} \frac{ (2\pi)^{\gamma/2} G(1/2)}{G(1+ \gamma/2) G(1/2+\gamma/2)} (1 + \O ( N^{-1} ))
\end{align}
\end{enumerate}
\eep

\

In the light of the above we introduce the notation
\beq \label{eqn:GB-def}
\GB (\lambda ; z, \beta) :=  \frac{ G\left( 1 + \frac{\lambda}{2} \right)}{(2 \pi)^{\lambda/4}} \times \left(\frac{2^{\frac{\gamma(\gamma-1)}{4}} G(1/2+\lambda/2)}{(2 \pi )^{\lambda/4} G(1/2) } \right)^{\1_{\beta=1} \1_{z \in \rr}}.
\eeq
We now state several results on the even integer moments of the characteristic polynomial for general i.i.d. matrices available from the literature. For the most part these formulas follow from the  works of Afanasiev \cite{afanasiev2019correlation,afanasiev2020correlation,
afanasiev2022correlation}  who computed various correlation functions via supersymmetric methods. We detail how the following is deduced from the literature in Appendix \ref{a:moments}. 

\bet
\label{theo:momcomp} Fix $r \in (0, 1)$. 
For any  i.i.d. matrix $X$ and $m \in \nn$ the following holds, where $\kappa_4 = \kappa_4^{(2)}$ is as in \eqref{eqn:kappa-def}:
\begin{enumerate}[label=(\roman*), font=\normalfont] 
\item For $z \in \ddb_r$ we have,
\beq
\ee[ \e^{ 2m \Phi_N (z, 0)}] = \frac{ (2\pi)^{m/2}}{G(1+m)} \e^{ \frac{m^2-m}{2} (1 - |z|^2)^2 \kappa_{4}} N^{\frac{m^2}{2}} \e^{ \1_{\beta=1} m  \log ( |z - \bar{z} |^2)} (1 + o (1) ).
\eeq
\item For $z \in I_r$ and $\beta=1$ we have,
\beq
\label{eq:momentsreal}
\ee[ \e^{ 2m \Phi_N (z, 0)}] = \left( \frac{N}{2} \right)^{m^2 - \frac{m}{2}} \e^{\frac{m^2-m}{2} (1 - |z|^2 )^2 \kappa_4} \frac{(2 \pi)^{m} G(1/2)}{G(1+m) G(1/2+m)} (1 + o (1) ).
\eeq
\end{enumerate} 
\eet

\

The only use of the above result is in order to deduce the following a priori bound on the fractional moments. That is, by Theorem~\ref{theo:momcomp} and H{\"o}lder inequality we  obtain the following.
\bec \label{cor:char-poly-a-priori} Let $X$ be an i.i.d. matrix. Let $0 < r < 1$ and $L >0$. 
Uniformly for $z \in \ddhb_r$ and $\lambda \in [0, L]$ we have,
\beq
\ee[ \e^{ \lambda \Phi_N (z, 0) }] \lesssim N^{ \frac{\lambda(\lambda+2)}{8}} + \1_{ \beta=1} \1_{ \{ z \in I_r \}} N^{ \frac{ \lambda(\lambda+2)}{4} - \frac{\lambda}{4}}.
\eeq
\eec

\subsection{Basic estimates}
\label{sec:detreg}

We will use the following result to control $\Phi_N (z, \eta)$ on the event that $\lambda_1^z$ is small. 

\bel \label{lem:deterministic} Let $0 < r < 1$. There is a constant $C_r>0$ so that the following holds for all $ z\in \dd_r$. 
Let $0 < \eta_1 < 2^{-1/2} N^{-1}$. Suppose that the set $\{ i : | \lambda_i^z | \leq \eta_1 \} $ is nonempty. Then, for any $0 < \eta < \eta_1$ we have,
\beq
\Phi_N (z, \eta) \leq \log (N \eta_1 ) + C_r + \Phi_N (z, N^{-1} ) .
\eeq
\eel
\proof Call the set $S := \{ i : | \lambda_i^z| \leq \eta_1 \}$ and let $K= | S |$.  Then,
\begin{align}
4 \Phi_N (z , \eta ) =  \sum_{i \in S} \log ( (\lambda_i^z)^2 + \eta^2 ) +  \sum_{i \notin S} \log ( (\lambda_i^z)^2 + \eta^2) - 2 N \int \log (x^2 + \eta^2) \rho^z (x) \d x
\end{align}
Then, using that $\eta \leq \eta_1$,
\beq
  \sum_{i \in S} \log ( (\lambda^z_i)^2 + \eta^2 ) \leq K ( \log (\eta_1^2) + \log (2) )  =  K ( \log (N^2 \eta_1^2) + \log (2) ) - K \log (N^2) 
\eeq
and
\begin{align}
& \sum_{i \notin S} \log ( (\lambda_i^z)^2 + \eta^2) -2 N \int \log (x^2 + \eta^2) \rho^z (x) \d x \notag\\
&\qquad \leq  \sum_{i \notin S} \log ( (\lambda_i^z)^2 + N^{-2} ) -2 N \int \log (x^2 + N^{-2} ) \rhosc (x) \d x + C \notag\\
&\qquad= 4 \Phi_N (z, N^{-1}) - \sum_{i \in S} \log ( (\lambda_i^z)^2 + N^{-2} ) +C \leq 4 \Phi_N (z, N^{-1}) + K \log ( N^2 ) +C.
\end{align}
In the first inequality we used monotonicity of the first term in $\eta$ and also \eqref{eqn:deterministic-cont}. Adding the inequalities we have
\beq
\Phi_N (z, \eta) \leq \frac{K}{2} ( \log (N \eta_1 ) + \log (2^{1/2} ) ) + C + \Phi_N (z, N^{-1} ).
\eeq
The claim now follows by noting that $K \geq 2$
and the assumption that $N \eta_1 < 2^{-1/2}$ ensures that $\log (N \eta_1) + \log (2^{1/2} ) \leq 0$.  \qed

\

When $\lambda_1^z$ is not too small we have the following which will allow us to regularize $\Phi_N (z, 0)$ to $\Phi_N (z , N^{-1-\delta} )$ for an appropriate $\delta >0$ small. 

\bel \label{lem:rig-reg}
Let $\eta_1 < N^{-1}$ and let $\eta_2 \leq \eta_1$. Then, we have
\beq
 | \Phi_N (z, 0) - \Phi_N (z, \eta_2 ) |\lesssim N\eta_2 + \frac{\eta_2^2}{\eta_1^2} \log N
\eeq
with overwhelming probability on the event $\lambda_1^z > \eta_1$.
\eel
\proof Since $|\lambda_i^z| \geq \eta_1 \geq \eta_2$ for all $i$, we have
\beq
\left| \log ( ( \lambda_i^z)^2 + \eta_2^2 ) - \log ( ( \lambda_i^z )^2 ) \right| \lesssim  \frac{ \eta_2^2}{(\lambda_i^z)^2 }  \lesssim \frac{ \eta_2^2}{ (\lambda_i^z)^2 + \eta_1^2}.
\eeq
Therefore, with $\nu = ( \log N) / N$ and using \eqref{eqn:deterministic-cont}, we have
\beq
| \Phi_N (z, 0) - \Phi_N (z, \eta_2 ) | \lesssim N\eta_2 + \frac{\eta_2^2}{\eta_1^2} (N \nu ) \Im[ m_N^z ( \i \nu) ] 
\eeq
where we used the monotonicity of $y \mapsto y \Im[ m_N^z ( \i y ) ]$. The claim now follows from Theorem \ref{thm:ll-2}.  \qed

Finally we have  the following a priori bound; it is useful because $\Phi_N (z, \eta)$ with larger $\eta$ are  accessible via resolvent methods.

\bel \label{lem:a-priori-1}
Fix a $C_1 \geq 1$, $0  < r<1$, and $\delta \in (0, 1)$. For any two $0 < \eta_1 < \eta_2 < ( \log N)^{C_1} / N$ and $z \in \dd_r$, we have
\beq
\Phi_N (z, \eta_1 ) \leq \Phi_N (z, \eta_2 ) + \O \left(  (N \eta_2 ) \wedge ( \log N)^{\frac{1}{2}+\delta} \right) 
\eeq
with overwhelming probability.
\eel
\proof Let $\eta_3 = ( \log N)^{\frac{1}{2}+\delta} / N$. For $\eta_2 \leq \eta_3$ this is a simple consequence of monotonicity and \eqref{eqn:deterministic-cont}. If $\eta_2 > \eta_3$ it then suffices to assume $\eta_1 \geq \eta_3$ by first applying the inequality with the choice $\eta_2 = \eta_3$. Then, for $ \eta_3 \leq \eta_1 < \eta_2$ we have,
\beq
| \Phi_N (z, \eta_1 ) - \Phi_N (z, \eta_2 ) | \leq N \int_{\eta_3}^{\eta_2} | m_N^z ( \i  u ) - m^z ( \i u ) | \d u
\eeq
and so the claim follows from Theorem \ref{thm:ll-2}. \qed

\subsection{Results on mesoscopic scales} \label{sec:meso-results}

In this section we state various results about the regularized log-characteristic polynomial $\Phi_N (z, \eta)$ on mesoscopic scales; i.e., when $\eta \gg N^{-1}$. In order to state them we introduce the following quantities which we will see below are the asymptotic mean and covariance of $\Phi_N (z, \eta)$. We omit the $\beta$-dependence for notational simplicity.

\bed
For $z_i \in \cc,  \eta_i \geq 0$ and $\kappa_4 \in \rr$ we define the covariance functional:
\begin{align} \label{eqn:def-covar}
\cC (z_1, \eta_1, z_2, \eta_2 ) &:= \vV (z_1, \eta_1, z_2, \eta_2 ) + \1_{ \beta=1} \vV (z_1, \eta_1, \bar{z}_2, \eta_2 ) + \frac{\kappa_4}{4} ( m_1 m_2 )^2, \notag\\
\vV (z_1, \eta_1, z_2, \eta_2 ) &:= -\frac{1}{4}\log\big[1+(u_1u_2|z_1z_2|)^2-(m_1m_2)^2-2u_1u_2\Re ( z_1\overline{z_2} ) \big],
\end{align}
with $m_i := m^{z_i} ( \i \eta_i)$ and $u_i = u^{z_i} ( \i \eta_i )$. We also define the expectation correction:
\beq \label{eqn:bmE-def}
\bmE (z_i, \eta_i ) := - \frac{ \kappa_4}{4} ( m_i )^4 + \frac{\1_{\beta=1}}{4} \log \left[ 1 - u_i^2 + 2 u_i^3 |z_i|^2 - u_i^2 ( z_i^2 + \bar{z}_i^2 ) \right] 
\eeq
\eed
We will require the following. The proof is in Appendix \ref{app:tech}.

\bep \label{prop:parameter-properties}
 Let $ 0 < r < 1$. The functions $\cC$ and $\vV$ are jointly continuous on $\dd_r \times \rr_{\geq 0} \times \dd_r \times \rr_{\geq 0 }$ away from the set $(z, 0, z, 0)$. The function $\bmE$ is jointly continuous on $\dd_r \times \rr_{\geq 0}$. The following are uniform in $z, z_i \in \dd_r$ and $\eta,\eta_i \in [0, 1] $.

\begin{enumerate}[label=(\roman*), font=\normalfont]

\item We have
\begin{align}
\label{eqn:covar-order-0}
\mathcal{V}(z_1,\eta_1,z_2,\eta_2)&=-\frac{1}{4}\log\big[|z_1-z_2|^2+\sqrt{1-|z_1|^2}\eta_1+\sqrt{1-|z_2|^2}\eta_2\big]\notag\\
+ & \mathcal{O}\left(\eta_1+\eta_2+|z_1-z_2|^2\right).
\end{align}

\item We have, with $u = u^z ( \i \eta)$ and $\eta < |z- \bar{z}|^2$, 
\beq \label{eqn:meso-logE}
 \log \left[ 1 -u^2 + 2 u^3 |z|^2 - u^2 (z^2 + \bar{z}^2 ) \right] =  \log \left[ |z-\bar{z} |^2 \right] + \O \left( \frac{\eta}{ |z-\bar{z}|^2} \right).
\eeq 

\item If $z \in I_r$ we have,
\beq
 \log \left[ 1 -u^2 + 2 u^3 |z|^2 - u^2 (z^2 + \bar{z}^2 ) \right] = \log \eta + \log (2 \sqrt{1-z^2} ) + \O ( \eta ).
\eeq 
\end{enumerate}
\eep

The following is our result on the joint Laplace transform of $\Phi_N (z, \eta)$ on mesoscopic scales. The proof appears below in Section \ref{sec:mainproofs}.

\bep There is a constant $B>0$ so that the following holds.\footnote{In the proof, the constant $B$ will be seen to be the constant from Lemma~\ref{lem:rough} below.} Let $X$ be an i.i.d. matrix. 
\label{pro:jointlapdnew}
Let $\eps, \delta, L >0$ and $K, J\in\mathbb{N}$. There is a $C_1 >0$, depending only on $J, L$, so that the following holds. Let $\{z_i\}_{i\in[J+K]}$ be complex numbers, $\{\eta_i\}_{i\in [J]}$ so that  $\frac{N^\eps}{N}=:\eta_M \leq \eta_i \leq 1$, $\{ a_i \}_{i=1}^K$ such that $\delta \leq a_i \leq 1$, $\{\gamma_i\}_{i\in [J]}$ with $\gamma_i \in [-L,L]$, and let $\{ t_i \}_{i=1}^K$ be in $\rr$.
Denote ${\bm t}:=(\gamma_1, \dots, \gamma_J, t_1,\dots, t_K)$.
Then, for any sufficiently large $D>0$, in terms of $L, J, K, C_1, \{t_i\}_{i\in [K]}$,
we have
\beq
\begin{split}
\label{eq:biglap}
&\ee\left[ \exp \left( \sum_{i=1}^J\gamma_i\Phi_N (z_i, \eta_i)+ \sum_{i=1}^K t_i \Phi_N (z_{J+i}, a_i ) \right) \1_\F \right]  \\
&\quad=\left(1+\mathcal{O}_{K,\delta, L, \eps, J}\left(\frac{1}{(N\eta_M)^{1/4}}\right)\right)\exp\left(\frac{1}{2}\langle {\bm t},\mathcal{C}{\bm t}\rangle+\langle{\bm t}, {\bm E}\rangle\right) \\
+& \mathcal{O}_{K,t_i,\delta, L, \eps, J}\big(N^{JLC_1-C_1^2/B}+N^{-D}\big)
\end{split}
\eeq
uniformly in $|z_i|\leq r<1$, $\min_i \eta_i\in[ \eta_M,1]$, and $a_i\in [\delta,1]$. Here, for $\hat{\eta}_i = \eta_i \1_{ i \leq J} + a_{i-J} \1_{ i > J}$, we defined
\beq
\F = \bigcap_{i\in [\![J+K]\!]}\{ |\Phi_N (z_i, \hat{\eta}_i )| \leq C_1 \log N\}.
\eeq
Additionally, ${\bm E}_i={\bm E}(z_i, \hat{\eta}_i)$ 
denotes the correction to the expectation from \eqref{eqn:bmE-def}, and $\mathcal{C}$ denotes the $(J+K)\times (J+K)$ covariance matrix with entries given by $\mathcal{C}(z_i,\hat{\eta}_i,z_j,\hat{\eta}_j)$ from \eqref{eqn:def-covar}.
\eep

We also require the following basic estimates about a single $\Phi_N (z, \eta)$, whose proof is given in Section \ref{sec:proof-ben-uses}.

\bel \label{lem:meso-ben-uses}
Let $C_* \geq 100$ be sufficiently large, let $L \geq 10$, and let $ r\in (0, 1)$. For all $u \in [1, L]$,  $\eta \in [ ( \log N)^{C_*} / N, 1 ] $, and $ z \in \ddhb_r$, we have
\beq \label{eqn:meso-tail-bd} 
\pp\left[ \left| \Phi_N (z, \eta) - \frac{\1_{ \beta=1} \1_{z \in I_r}}{4} \log \eta \right| > u \log N \right] \lesssim \e^{ - 2 u^2 \log N } + \1_{\beta=1} \1_{ \{ z \in I_r \}} \e^{ - u^2 \log N }.
\eeq
If $A \geq 10 L + 100$ and $\lambda \in [0, L]$, then,
\beq \label{eqn:meso-1pt-mgf}
\ee\left[\e^{ \lambda \Phi_N (z, \eta) } \1_{ \{ \Phi_N (z, \eta) \leq A \log N \} } \right] = (1 + \O ( ( \log N)^{100} / (N \eta)^{1/4} ) ) \e^{ \frac{\lambda^2}{2} \cC (z, \eta, z, \eta) + \lambda \bmE (z, \eta) }
\eeq
For any $A \geq 1, L \geq 1$ and $\eps >0$, the following holds uniformly in $\eta \in [ N^{\eps-1}, 1 ]$ and $|\lambda| \leq L$:
\beq \label{eqn:meso-uni-bd}
\ee\left[ \e^{ \lambda \Phi_N (z, \eta ) } \1_{ \{ |\Phi_N (z, \eta) |\leq A \log N \} } \right] \lesssim \exp \left( - \frac{\lambda^2}{4(2-\1_{ \beta=1} \1_{ \{ z \in I_r \}} )} \log \eta + \1_{ \beta=1} \1_{ \{ z \in I_r \}} \frac{  \lambda}{4} \log \eta \right).
\eeq
\eel

 \section{Homogenization and application} \label{sec:dynamics}

\subsection{Definition of dynamical quantities} \label{sec:def-dynamics} 

In this section as well as Sections \ref{sec:1pt}, \ref{sec:kpt} and \ref{sec:claeys} we will consider a dynamically evolving matrix $X_t$. We begin by defining the dynamics we use as well as several dynamically evolving quantities that we will need. We define $X_t$ as the solution to,
\beq \label{eqn:Xt-def}
\d X_t = \frac{ \d B(t)}{ \sqrt{N}},
\eeq
where the initial data $X_0$ is an i.i.d. matrix. Here $B(t)$ is an $N \times N$ matrix of i.i.d. standard real or complex Brownian motions when $\beta=1$ or $2$ respectively. Our normalization is that $\ee[ |B_{ij} (t) |^2] =t$. 

Note that at time $t$ the variance of $X_{ij}(t)$ is $N \ee[ |X_{ij} (t)|^2] = 1 + t $. We therefore need to rescale the spectral measures for the Hermitization of $X_t$. Specifically, define $H_t^z := H^z (X_t)$ and denote the eigenvalues of $H_t^z$ by $\lambda_i^z (t)$. Let $G_t^z (w) = (H^z_t - w)^{-1}$. Let 
$
c_* (t) := \sqrt{1 + t} 
$
and define
\beq
\rho_t^z(x):= \frac{1}{c_*(t)} \rho^{z/c_*(t)} \left( \frac{x}{c_*(t)} \right) , \qquad m_t^z (w) := \frac{1}{c_*(t)} m^{z/c_*(t)} (w / c_* (t) ) .
\eeq
We also define the quantiles $\gamma_i^z (t)$ associated to $\rho_t^z$ in the same manner as \eqref{eqn:quantiles-def}. From the fact that $\del_x \rho^z (x) = \O (1)$ for $|x - \mathfrak{e}_z | > c$  which follows from Lemma~\ref{lem:proprho}\ref{it:rho-der}, as well as the second estimate of \eqref{eq:exp}, it is not hard to see that
\beq \label{eqn:gamma-t-expand}
\gamma_i^z (t) = \frac{i}{2N \rho^z (0) } (1 + \O (t) + \O (i/N) )
\eeq
for $i \in [\![N/2]\!]$. 
Therefore we define,
\beq \label{eqn:gamma-hat-def}
\hatgam_i^z := \frac{i}{2N \rho^z (0) }
\eeq
Associated with the flow \eqref{eqn:Xt-def} are the characteristics $(\eta_t, z_t)$ which solve
\beq \label{eqn:char-def}
\del_t \eta_t = - \Im[ m^{z_0}_t ( \i \eta_t ) ] = \i m^{z_0}_t ( \i \eta_t) , \qquad z_t = z_0 .
\eeq
Note that $z_t$ is constant in time. 
In particular, one can compute,
\beq
\del_t m^z_t (w) = m_t^z (w) \del_w m_t^z (w), \qquad \del_t m_t^{z_t} (\i \eta_t)  = 0. 
\eeq
We also define,
\begin{align} \label{eqn:phiNt-def}
\Phi_N (z, \eta, t ) &:= \frac{1}{4} \sum_{i=-N}^N \log ( ( \lambda_i^z (t))^2 + \eta^2) - \frac{N}{2} \int \log (x^2 +\eta^2) \rho_t^z (x) \d x
\end{align}
and 
\beq
m_{N, t}^z (w) := \frac{1}{2N} \sum_i \frac{1}{ \lambda_i^z (t) - w } = \langle G_t^z (w) \rangle.
\eeq
Note that $\Phi_N (z, \eta, t) \neq \Phi_N (z, \eta; X_t )$ due to the difference in the deterministic term. 

The scaling,
\begin{align} \label{eqn:t-scaling}
 &\Phi_N (z, \eta, t) = \frac{1}{2} \Re \log \det (H_t^z - \i \eta ) - \frac{N}{2} \int \log (x^2 +\eta^2) \rho_t^z (x) \d x \notag\\
 = & \frac{1}{2} \Re \log \det ( H^{z/c_*(t)} (X_t / c_*(t)) - \i \eta / c_*(t) ) - \frac{N}{2} \int \log (x^2 + (\eta/c_*(t) )^2 ) \rho^{z/c_* (t) } (x) \d x
\end{align}
is helpful to keep in mind (note that $X_t /c_*(t)$ has matrix elements with variance $1/N$ which coincides with the scaling in Definition \ref{def:iid}) and will be used elsewhere.

Under the evolution in \eqref{eqn:Xt-def} we have (see e.g. \cite[Appendix B]{cipolloni2023central},\cite[Appendix B]{cipolloni2021fluctuation} in the complex and real cases respectively) 
\beq
\label{eq:DBMz}
\d \lambda_i^z (t) = \frac{ \d b_i^z (t)}{ \sqrt{2N}} + \frac{1}{2N} \sum_{i \neq j } \frac{1 + \1_{\beta=1} \Lambda^z_{ij} (t) }{ \lambda^z_i (t) - \lambda_j^z (t)} \d t
\eeq
where the processes $\Lambda^z_{ij}(t)$ and Martingales $b_i^z(t)$ are as follows: when $\beta=2$, the $\{ b_i^z \}_{i=1}^N$ are i.i.d. standard Brownian motions; when $\beta=1$ the $b_i^z(t)$  have the  covariation process, 
\beq
\label{eq:corrbms}
\d [ b_i^z , b_j^z] = [\delta_{ij} - \delta_{i,-j} + \Lambda_{ij}^z], \quad \Lambda_{ij}^z(t) = 4 \Re \left[ \langle \bw_i^z (t) , E_1 \bw_j^{\bar{z}} (t) \rangle \langle \bw_j^{\bar{z}} (t) , E_2 \bw_i^z (t) \rangle \right].
\eeq
Here $\bw_i^z (t)$ are the orthonormal eigenvectors of $H^z (X_t)$.  Note that $\Lambda_{ij}^z(t)=\delta_{i,j}-\delta_{i,-j}$ for $z\in \rr$.  If $(\eta_t, z_t)$ is a characteristic then by It\^{o}'s lemma (see e.g. \cite[Eq. (5.19)]{cipolloni2025maximum}),
\begin{align} \label{eqn:psi-char}
    &\d \left[ \sum_i \log ( \lambda_i^z - \i \eta_t ) - 2 N \int \log (x - \i \eta_t) \rho_t^z (x) \d x \right] \\
    \notag
    &\qquad= \frac{1}{\sqrt{2N}} \sum_i \frac{ \d b_i^z}{ \lambda_i^z - \i \eta_t }  - N ( m_{N, t}^z ( \i \eta_t ) - m_t^z (\i \eta_t ) )^2   - \1_{  \beta=1 } \tilde{\sum}_{ij}\langle G_t^z(\i \eta_t )E_iG_t^{\overline{z}}(\i \eta_t )E_j\rangle \d t. 
\end{align}
The following gives a useful approximation for the characteristic $\eta_t$. 
\bel \label{lem:approx-char}
Let $0 < r <1$ and $T \in [0, 1]$. There is a $c_r >0$ so that for all $z \in \dd_r$ and all characteristics $\eta_t$ with $\eta_T \leq c_r$ and $T \leq c_r$ we have,
\beq
\eta_{T-s} = \eta_T + s \left( \sqrt{1- |z|^2} + \O ( \eta_T + T ) \right)
\eeq
\eel
\proof This follows from $\eta_{T-s} = \eta_T + s \Im[ m_T^z ( \i \eta_T ) ]$ and that $m_T^z ( \i \eta_T) = m_T^z ( \i 0) + \O ( \eta_T) = m^z ( \i 0) + \O (\eta_T + T)$, from \eqref{eq:exp} and the definition of $m_t^z$. \qed

\subsection{Statement of homogenization result}

We now show that the processes \eqref{eq:DBMz} for different $z$'s (sufficiently far away from each other) can be coupled to independent processes. The proof of this theorem is postponed to Appendix~\ref{app:tech}. Theorem \ref{thm:homog} is essentially a recapitulation of existing homogenization results for i.i.d. matrices. We state it  in order to explicitly summarize the salient aspects of the literature that we need for our methods.  Recall the notation $\ddhb_r$ defined in \eqref{eqn:ddb-def}. 

\bet \label{thm:homog}
Fix an integer $K \geq 2$, and $0 < r < 1$ and $0 < b < \frac{1}{100}$. There is a small constant $\frac{1}{100} > \mfb >0$, depending on $b >0$, so that the following holds.

Let $z_i$ be $K$ points in $\ddhb_r$ such that $|z_i -z_j| \geq N^{b-1/2}$ for all $i \neq j$. There are $K$ independent processes $\{ \mu_j^{(n)} (t) \}_{j=-N, \dots N}$ for $1 \leq n \leq K$ and $K$ independent families of independent  Brownian motions $\{ W_j^{(n)} \}_{j=1}^N$ with quadratic variation 
\beq  \label{eqn:W-qv}
\d [ W_i^{(n)} (t)] = \left( 1 + \1_{ \beta=1} \1_{  z_n \in I_r }  \right) \d t
\eeq
so that following holds.
\begin{enumerate}[label=(\roman*),font=\normalfont]
\item There are independent Ginibre matrices $G_1, \dots , G_K$ so that (here $H^0(G_n)$ is the Hermitization at $z=0$)
\beq
\mu_i ^{(n)} (0) = \frac{\rho^0(0)}{ \rho^{z_i} (0) } \lambda_i ( H^0 (G_n) ) .
\eeq
If $z_n \in \ddb_r$ then $G_n$ is a complex Ginibre matrix. If $z_n \in I_r$ and $\beta=1$, then $G_n$ is a real Ginibre matrix. Moreover,  (defining $W_{-i}^{(n)} := - W_i^{(n)}$) we have 
\beq \label{eqn:mu-n-i-dbm}
\d \mu^{(n)}_i (t) = \frac{ \d W_i^{(n)}}{ \sqrt{2 N}} + \frac{1}{2N} \sum_{j \neq i } \frac{1 - \delta_{i, -j} \1_{ \beta=1} \1_{  z_n \in I_r  } }{ \mu^{(n)}_i (t) - \mu_j^{(n)} (t) } \d t.
\eeq
\label{it:constrdbm}
\item For all $|i|, |j| \leq N^{\mfb}$ and $t \leq N^{\mfb-1}$ we have  
\beq \label{eqn:homog-new-bm}
\d \left| \left[ (b_i^{z_n} - W^{(n)}_i ) , (b_j^{z_n} -  W^{(n)}_j ) \right] \right| \leq N^{-\mfb} \d t.
\eeq
\label{it:asympind}
\item For any  $\omega \in (0, \mfb)$ there is an $\mfa >0$ so that for all $N^{\frac{\omega}{2}} \leq  Nt \leq  N^{\mfb}$ we have,
\beq \label{eqn:homog}
\sup_{ |i| \leq N^{\mfb}} | \lambda_i^{(z_n)} (t) - \mu_i^{(n)} (t) | \leq N^{-1-\mfa}
\eeq
with overwhelming probability. \label{it:inuniv} 
\item The Brownian motions $\{ W_j^{(n)} (t)\}_{j, n, t}$ are independent from the initial data $X_0$. \label{it:homog-ind} 
\end{enumerate}
\eet

\begin{remark} \emph{We emphasize above that when $z_n \notin I_r$ and $\beta=1$, the real i.i.d. matrices are  coupled to flows involving the complex Ginibre ensemble (evaluated at the origin) and not the real Ginibre ensemble. This reflects the fact that  the local statistics away from the real axis coincide with the complex Ginibre ensemble, even though $X$ is real.}
\end{remark} 

\subsection{Universal representation of characteristic polynomials} \label{sec:univ-rep}

In this section, let $b >0$ be the parameter from the previous subsection.  Let $\mfb>0$ be the corresponding parameter from Theorem \ref{thm:homog}, and then choose a small $\mfb > \omega_1 >0$, and let $\mfa$ be the corresponding parameters from Theorem \ref{thm:homog} so that \eqref{eqn:homog} holds. Fix further a small $\delta_m >0$  and a large $C_* \geq 1000$ satisfying,
\beq
\mfa > \delta_m , \quad \omega_1 < q_1 < \mfb
\eeq
and then define,
\beq
\eta_m := \frac{N^{-\delta_m}}{N} , \quad \eta_* := \frac{ ( \log N)^{C_*}}{N} \quad t_1 := \frac{N^{\omega_1}}{N}.
\eeq
Let $z_n \in \{ z_1, \dots z_K \}$ from the previous section.
Throughout this section we will adopt the shorthand notation $z= z_n$,  $\mu_i (t) = \mu_i^{(n)} (t)$ and $W_i = W_i^{(n)}$.  Let now
\beq \label{eqn:approx-char} 
\nu^z (s) := \eta_* + (t_1 -s ) \sqrt{1 - |z|^2}.
\eeq
According to Lemma \ref{lem:approx-char},  $\nu^z (s)$ is approximately the characteristic ending at $\eta_*$ at time $t_1$.  We write,
\begin{align} \label{eqn:lg-decomp}
\Phi_N (z, \eta_m, t_1) =& ( \Phi_N (z, \eta_m, t_1 ) - \Phi_N (z, \eta_*, t_1 ) ) + ( \Phi_N (z, \eta_*, t_1 ) )- \Phi_N (z, t_1 \sqrt{1-|z|^2}, 0) ) \notag\\
+& \Phi_N (z, t_1 \sqrt{1 - |z|^2}, 0).
\end{align}
In this section we will develop universal representations for the first two terms based on Theorem \ref{thm:homog}. The first term contains the local information about the characteristic polynomial required to go from the submicroscopic scale $N^{-1-\delta_m}$ to the mesoscopic scale $\eta_*$. Note that it only concerns the behavior at the fixed time $t_1$. The second term carries the effect of the dynamics; the evolution of $\Phi_N$ will be evaluated along the approximate characteristic $\nu^z(s)$. These first two terms will be correlated with each other but will be approximately independent from $\Phi_N (z, t_1 \sqrt{1-|z|^2}, 0)$ which carries the information from the initial data.

Our approximation for the first term on the RHS of \eqref{eqn:lg-decomp} is the following.

\bel \label{lem:univ-rep-1} 
In the above set-up we have,
\begin{align}  \label{eqn:univ-rep-1-1}
 & \Phi_N (z, \eta_m, t_1 ) - \Phi_N (z, \eta_*, t_1 ) = - \frac{1}{2} \int_{\eta_m}^{\eta_*} \sum_{ |i| < N^{\mfb}} \frac{ u}{ ( \mu_i (t_1))^2 + u^2} \d u \notag\\
 + & N \int_{\eta_m}^{\eta_*} \int_{|x| < \hatgam^z_{N^{\mfb}}} \frac{ u}{x^2+u^2} \rho^z (0) \d x \d u + \O(  ( \log N ) N^{\delta_m-\mfa} + N^{-\mfb} )
\end{align}
with overwhelming probability. For any $\delta >0$ we have with overwhelming probability that
\beq \label{eqn:univ-rep-1-2}
 - \frac{1}{2} \int_{\eta_m}^{\eta_*} \sum_{ |i| < N^{\mfb}} \frac{ u}{ ( \mu_i (t_1))^2 + u^2} \d u 
 +  N \int_{\eta_m}^{\eta_*} \int_{|x| < \hatgam^z_{N^{\mfb}}} \frac{ u}{x^2+u^2} \rho^z (0) \d x \d u \leq ( \log N)^{\frac{1}{2}+\delta} .
\eeq
\eel
\proof The cases $\beta=1, 2$ are identical. By definition,
\begin{align}
& \Phi_N (z, \eta_m, t_1 ) - \Phi_N (z, \eta_*, t_1 )
=  N \int_{\eta_m}^{ \eta_*} \left( \Im[ m^z_{t_1} ( \i u) ] - \Im[ m^z_{N, t_1}  ( \i u ) ] \right) \d u.
\end{align}
We have with overwhelming probability, using that $\mfa > \delta_m$, and \eqref{eqn:homog},
\begin{align}
& \int_{ \eta_m}^{\eta_*} \sum_{ |i| \leq N^{\mfb}} \left| \frac{ u}{ ( \lambda_i^z (t_1))^2 + u^2} - \frac{u}{ ( \mu_i (t_1) )^2 + u^2} \right| \d u
\lesssim  \frac{1}{N^{1+\mfa}} \int_{ \eta_m}^{\eta_*} \sum_{ |i| \leq N^{\mfb}} \frac{ u}{ | \lambda_i^z (t_1) |^3 + u^3} \d u.
\end{align}
For the quantity on the RHS we estimate first,
\beq
\frac{1}{N^{1+\mfa}} \int_{ \eta_m}^{\eta_*} \sum_{ |i| \leq \log N} \frac{ u }{ | \lambda_i^z (t_1) |^3 + u^3} \d u \lesssim \frac{ \log N}{N^{1 + \mfa}} \int_{\eta_m}^{\eta_*} \frac{1}{u^2} \d u \lesssim (\log N) N^{\delta_m - \mfa}  .
\eeq
On the other hand with overwhelming probability, with the notation $\eta_3 := ( \log N) / N$, 
\begin{align}
 & \frac{1}{N^{1+\mfa}} \int_{ \eta_m}^{\eta_*} \sum_{ |i| > \log N} \frac{ u  }{ | \lambda_i^z (t_1) |^3 + u^3} \d u \notag\\
\lesssim & N^{-\mfa}  \int_{\eta_m}^{\eta_*} \frac{u}{(u + \eta_3)^2} \Im[ m^z_{N, t_1}( \i (\eta_3 + u ) ] \d u \lesssim ( \log N) N^{-\mfa},
\end{align}
using Theorems \ref{thm:ll-2} and \ref{thm:rigidity}. Next, by Theorem \ref{thm:rigidity} and a standard calculation, with the notation $i_\mfb = N^{\mfb}$, we have for any $\eps >0$ with overwhelming probability
\begin{align}
& \left| N \int_{\eta_m}^{\eta_*} \sum_{ |i| > N^{\mfb}} \frac{1}{2N} \frac{u}{ ( \lambda_i^z (t_1) )^2 + u^2} \d u  - N \int_{\eta_m}^{\eta_*} \int_{ |x| > \gamma_{i_\mfb} ^z (t_1) } \frac{ u}{ u^2 +x^2} \rho^{(z)}_{t_1} (x) \d x  \d u \right|  \notag\\
\lesssim & N^{\eps/2} \int_{\eta_m}^{\eta_*} \frac{ N^2 u}{ N^{2 \mfb}} \d u \leq N^{\eps - 2 \mfb}.
\end{align}
Now, using $\rho_t^z (x) = \rho^z (0) + \O ( t+ |x|)$, we have
\begin{align}
  \int_{\eta_m}^{\eta_*} \int_{ |x| < \gamma_{i_\mfb} ^z (t_1) } \frac{N u}{u^2 + x^2} \rho^z_{t_1} (x) \d x \d u 
=   \int_{\eta_m}^{\eta_*} \int_{ |x| < \gamma_{i_\mfb} ^z (t_1) } \frac{ N u}{u^2 + x^2} \rho^z (0) \d x \d u+  \O ( N \eta_* ( t_1 + N^{\mfb-1} ) ).
\end{align}
Finally, using \eqref{eqn:gamma-t-expand} and $|\gamma_{i_\mfb} ^z (t_1) - \hatgam_{i_\mfb}^z| \ll  \gamma_{ i_\mfb}^z (t_1)$, we have
\begin{align}
 & \left| N \int_{\eta_m}^{\eta_*} \int_{ |x| < \gamma_{i_\mfb} ^z (t_1) } \frac{ u}{u^2 + x^2} \rho^z (0) \d x \d u - N \int_{\eta_m}^{\eta_*} \int_{ |x| < \hatgam_{i_\mfb} ^z  } \frac{ u}{u^2 + x^2} \rho^z (0) \d x \d u  \right| \notag\\
 \lesssim & N \eta_*^2 \frac{ | \gamma_{i_\mfb}^z (t_1) - \hatgam_{i_\mfb}^z|}{ ( \hatgam_{i_\mfb}^z)^2} \lesssim (N \eta_*^2) (1 + N^{\omega_1 - \mfb} )
\end{align}
In the last inequality we used that \eqref{eqn:gamma-t-expand} implies that $| \gamma_{i_\mfb}^z (t_1) - \hatgam_{i_\mfb}^z| \lesssim \hatgam_{i_\mfb}^z  ( t_1 + \hatgam_{i_\mfb}^z )$.
This completes the proof of the first estimate \eqref{eqn:univ-rep-1-1} of the lemma, where we simplified a few of the errors using $\mfb < \frac{1}{100}$. The second estimate \eqref{eqn:univ-rep-1-2} follows from \eqref{eqn:univ-rep-1-1} and Lemma \ref{lem:a-priori-1}. \qed

Our approximation for the second term on the RHS of \eqref{eqn:lg-decomp} is the following. Its proof appears at the end of this section, after a few preparatory remarks and intermediate lemmas. 

\bel \label{lem:univ-rep-2} With overwhelming probability, 
\begin{align} \label{eqn:univ-rep-2}
   & \Phi_N (z, \eta_*, t_1 ) - \Phi_N (z, t_1 \sqrt{1- |z|^2} , 0) + ( \log N)^4 \O ( N^{\omega_1/2 - \mfb/2}
   +  N^{-\mfa}  +  N^{-\omega_1/2 } ) \notag\\
    = & \frac{1}{2} \Re \frac{1}{ \sqrt{2 N}} \sum_{|i| < \Nfb}  \int_{0}^{t_1}  \frac{ \d W_i (s)}{ \mu_i (s) - \i \nu^z (s) } \notag\\
    + & \Re  \int_{t_1/2}^{t_1} \frac{N}{2}  \left( \frac{1}{2N} \sum_{ |i| < \Nfb} \frac{1}{ \mu_i (s) - \i \nu^z (s) }  - \int_{ |x| < \hatgam_{\Nfb}^z } \frac{\rho^z (0)  }{ x - \i \nu^z (s) } \d x \right)^2 \d s \\
    - & \1_{ \beta=1} \1_{  z \in I_r } \left\{ \Re \int_{t_1/2}^{t_1} \frac{1}{8N} \sum_{|i| \leq N^{\mfb}} \left( \frac{1}{( \mu_i (s) - \i \nu^z (s) )^2} + \frac{1}{ | \mu_i (s) - \i \nu^z (s) |^2}\right) \d s + \frac{1}{4} \log 2 \right\}
\end{align}
We  have that the second last line above is $\O ( ( \log N)^3 / ( \log N)^{C_*})$ with overwhelming probability. With overwhelming probability the last line equals $ \frac{1}{4} \log ( \sqrt{1-|z|^2} t_1/\eta_*) + \O ( ( \log N)^{-10} )$. 
\eel
In preparation for the proof of Lemma \ref{lem:univ-rep-2}, we let $\eta^z (s)$ be the characteristic ending at $\eta_*$ at time $t_1$ (defined by \eqref{eqn:char-def}). We write
\begin{align} \label{eqn:univ-rep-2-b1}
& \Phi_N (z, \eta_*, t_1 ) - \Phi_N (z, t_1 \sqrt{1-|z|^2} , 0) = \left( \Phi_N (z, \eta_*, t_1 ) - \Phi_N (z, \eta^z (0) , 0) \right) \notag\\
+ &\left( \Phi_N (z, \eta^z (0), 0) - \Phi_N (z, t_1 \sqrt{1-|z|^2}  , 0) \right)
\end{align}
We first bound the last line of \eqref{eqn:univ-rep-2-b1}. 
Since $\del_\eta \Phi_N (z, \eta, 0) = N \Im[ m_0^z ( \i \eta ) - m_{N, 0}^z ( \i \eta) ]$ and by Lemma \ref{lem:approx-char},  $|t_1 \sqrt{1-|z|^2}  - \eta^z (0) | \lesssim t_1 (t_1 + \eta_*/t_1 )$ we have by Theorem \ref{thm:ll-2}, 
\beq \label{eqn:univ-rep-2-d1} 
\left| \Phi_N (z, \eta^z (0), 0) - \Phi_N (z, t_1 \sqrt{1-|z|^2} , 0) \right| \leq ( \log N)^{2 +C_*} N^{-\omega_1},
\eeq
with overwhelming probability. Now by \eqref{eqn:psi-char} we have for the first term on the RHS of \eqref{eqn:univ-rep-2-b1}, 
\begin{align} \label{eqn:univ-rep-2-b2}
    & \Phi_N (z, \eta_*, t_1 ) - \Phi_N (z, \eta^z (0) , 0) \notag\\
    &=  \frac{1}{2} \Re \frac{1}{ \sqrt{2N}} \sum_i  \int_0^{t_1}  \frac{ \d b_i^z (s)}{ \lambda_i^z (s) - \i \eta^z (s) } +\Re  \int_0^{t_1} \frac{N}{2} ( m_{N, s}^z  ( \i \eta^z (s) ) - m_s^z ( \i \eta^z (s) ) )^2 \d s \notag\\
    &\quad - \frac{1}{2} \Re \int_0^{t_1} \1_{  \beta=1 } \tilde{\sum}_{ij}\langle G_s^z(\i \eta^z (s) )E_iG_s^{\overline{z}}(\i \eta^z (s) )E_j\rangle \d s .
\end{align}
The next three lemmas consider each of the terms on the RHS of \eqref{eqn:univ-rep-2-b2}. 
\bel \label{lem:univ-rep-2-1} 
In the set-up of Lemma \ref{lem:univ-rep-2} we have, with overwhelming probability, 
\begin{align} \label{eqn:univ-rep-2-b3} 
&\Re \frac{1}{ \sqrt{2N}} \sum_i  \int_0^{t_1}  \frac{ \d b_i^z (s)}{ \lambda_i^z (s) - \i \eta^z (s) } = \Re \frac{1}{ \sqrt{2N}} \sum_{|i| < \Nfb}  \int_{0}^{t_1}  \frac{ \d W_i^z (s)}{ \mu_i (s) - \i \nu^z (s) } \notag\\
+ &  ( \log N)^3 \O (  N^{\omega_1/2-\mfb/2} +  N^{-\mfa}  + N^{-\omega_1 } ). 
\end{align}
\eel
\proof We first consider the case $z \in \ddb_r$. That is, if $\beta=1$ we assume $z \notin I_r$, leaving the remaining case until the end of the proof.  We start by estimating the contribution of indices $i$ such that $|i| > N^{\mfb}$. The quadratic variation is, 
\beq \label{eqn:univ-rep-2-b9}
\begin{split} 
& \left[ \Re \frac{1}{ \sqrt{2N}} \sum_{|i| > N^{\mfb}}   \int_0^{t_1}  \frac{ \d b_i^z (s)}{ \lambda_i^z (s) - \i \eta^z (s) }  \right] 
=\int_0^{t_1}\frac{2}{N}\sum_{i>N^{\mfb}} \frac{\lambda_i^z(s)^2\,\d s}{|\lambda_i^z(s)-\i\eta^z(s)|^4} \notag\\
+&  \bm1_{\beta=1}\int_0^{t_1}\frac{2}{N}\sum_{i,j>N^{\mfb}} \frac{\lambda_i^z(s)\lambda_j^z(s) 4\Re[\langle {\bm w}_i^z,E_1{\bm w}_j^{\overline{z}}\rangle \langle {\bm w}_j^{\overline{z}},E_2{\bm w}_i^z\rangle]}{|\lambda_i^z(s)-\i\eta^z(s)|^2|\lambda_j^z(s)-\i\eta^z(s)|^2}\, \d s.
\end{split}
\eeq
Note that in the second term we dropped the time argument of the eigenvectors ${\bm w}^z_i = {\bm w}^z_i (s)$ for notational simplicity. For the first term on the RHS of \eqref{eqn:univ-rep-2-b9} we have with overwhelming probability,
\begin{align}
\frac{1}{N} \int_0^{t_1} \sum_{|i| > N^{\mfb} } \frac{(\lambda_i^z (s) )^2}{ | \lambda_i^z (s) - \i \eta_s |^4} \d s \lesssim N^{\omega_1 - \mfb}.
\end{align} 
We used Theorem \ref{thm:rigidity} and \eqref{eqn:gamma-t-expand} which imply that $\lambda_i^z (t_1) \asymp i/N$ for $i \geq \ell_1$ with overwhelming probability.

When $\beta=1$ we bound the additional term as follows. We have,
\beq
\begin{split}
\sum_{i,j>N^{\mfb}} &\frac{\lambda_i^z(s)\lambda_j^z(s)4\Re[\langle {\bm w}_i^z,E_1{\bm w}_j^{\overline{z}}\rangle \langle {\bm w}_j^{\overline{z}},E_2{\bm w}_i^z\rangle]}{|\lambda_i^z(s)-\i\eta^z(s)|^2|\lambda_j^z(s)-\i\eta^z(s)|^2}\leq 2 \sum_{i,j>N^{\mfb}} \frac{|\langle {\bm w}_i^z,E_1{\bm w}_j^{\overline{z}}\rangle|^2+ |\langle {\bm w}_j^{\overline{z}},E_2{\bm w}_i^z\rangle]|^2}{|\lambda_i^z(s)-\i\eta^z(s)||\lambda_j^z(s)-\i\eta^z(s)|}\\
= &\frac{1}{2}\sum_{|i|,|j|>N^{\mfb}} \frac{|\langle {\bm w}_i^z,E_1{\bm w}_j^{\overline{z}}\rangle|^2+ |\langle {\bm w}_j^{\overline{z}},E_2{\bm w}_i^z\rangle]|^2}{|\lambda_i^z(s)-\i\eta^z(s)||\lambda_j^z(s)-\i\eta^z(s)|} \leq \frac{1}{2} \sum_{k=1}^2\mathrm{Tr}\big[|G_s^z(\i\eta^{\overline{z}}(s))|E_k|G_s^z(\i\eta^z(s))|E_k\big].
\end{split}
\eeq
We point out that to go from the first to the second line we used ${\bm u}^z_{-i}={\bm u}^z_i$ and ${\bm v}^z_{-i}=-{\bm v}^z_i$, as a consequence of the chiral symmetry of $H^z$, where ${\bm w}^z_i = ( {\bm u}^z_i , {\bm v}^z_i )$. By a direct application of Lemma~\ref{lem:2grealb} we have with overwhelming probability
\beq
\frac{1}{N} \int_0^{t_1} \sum_{k=1}^2\mathrm{Tr}\big[|G_s^z(\i\eta^{\overline{z}}(s))|E_k|G_s^z(\i\eta^z(s))|E_k\big] \d s \leq N^{-1/10} ,
\eeq
Therefore, by Lemma \ref{lem:mart-bd}, we have with overwhelming probability
\beq \label{eqn:univ-rep-2-b4}
\Re \frac{1}{ \sqrt{2N}} \sum_i  \int_0^{t_1}  \frac{ \d b_i^z (s)}{ \lambda_i^z (s) - \i \eta^z (s) } = \Re \frac{1}{ \sqrt{2N}} \sum_{|i| < N^{\mfb}}  \int_0^{t_1}  \frac{ \d b_i^z (s)}{ \lambda_i^z (s) - \i \eta^z (s) } + \O ( ( \log N)^2 N^{\omega_1/2 - \mfb/2} ) .
\eeq
Now, using \eqref{eqn:homog-new-bm} we find, 
\begin{align} \label{eqn:univ-rep-2-qv1}
& \left[ \Re  \frac{1}{ \sqrt{2N}} \sum_{|i| <  \Nfb}   \int_0^{t_1}  \frac{ \d b_i^z (s)- \d W_i (s) }{ \lambda_i^z (s) - \i \eta^z (s) } \right] \leq \frac{1}{N} \sum_{ |i|, |j| \leq \Nfb} \int_0^{t_1} \frac{ \d \left| \left[ b_i^z -W_i , b_j^z - W_j \right] \right|}{ | \lambda_i^z - \i \eta^z (s) | | \lambda_j^z - \i \eta^z (s) | } \notag\\
\leq &\frac{N^{-\mfb}}{N}  \int_0^{t_1} \sum_{ |i|, |j| \leq \Nfb}  \frac{1}{ | \lambda_i^z - \i \eta^z (s) | | \lambda_j^z - \i \eta^z (s) |}\d s  \lesssim ( \log N)^2 N^{\omega_1-\mfb}.
\end{align}
Furthermore, by Lemma \ref{lem:approx-char} (recall the definition of $\nu^z(s)$ in \eqref{eqn:approx-char}) we have,
\beq \label{eqn:char-dif}
|\eta^z (s) - \nu^z (s) | \lesssim \eta^z (s) (t_1 + \eta_* ) .
\eeq
Therefore by applying Theorem \ref{thm:ll-2} we have with overwhelming probability,
\begin{align} \label{eqn:univ-rep-2-qv2}
& \left[ \Re \frac{1}{ \sqrt{2N}} \sum_{|i| < \Nfb}  \int_0^{t_1} \left(   \frac{1 }{ \lambda_i^z (s) - \i \eta^z (s) } - \frac{1 }{ \lambda_i^z (s) - \i \nu^z (s) }\right) \d W_i (s) \right] \notag\\ 
&\lesssim  \frac{1}{N} \int_0^{t_1} \sum_{ |i|  < \Nfb} \frac{ \eta^z (s)^2(t_1 + \eta_*)^2}{ | \lambda_i^z(s)  - \i \eta^z (s) |^4 } \d s \lesssim (t_1 + \eta_*)^2 \int_0^{t_1} \frac{1}{ \eta^z (s) } \d s \lesssim ( \log N) (t_1 + \eta_*)^2,
\end{align}
From \eqref{eqn:univ-rep-2-qv1} and \eqref{eqn:univ-rep-2-qv2} we have
\beq \label{eqn:univ-rep-2-b6}
\Re \frac{1}{ \sqrt{2N}} \sum_{|i| < \Nfb}  \int_0^{t_1}  \frac{ \d b_i^z (s)}{ \lambda_i^z (s) - \i \eta^z (s) } = \Re \frac{1}{ \sqrt{2N}} \sum_{|i| < \Nfb}  \int_0^{t_1}  \frac{ \d W_i (s)}{ \lambda_i^z (s) - \i \nu^z (s) } +\O ( ( \log N)^2 N^{\omega_1/2 - \mfb/2} ) ,
\eeq
with overwhelming probability.
Using \eqref{eqn:homog}, Theorem \ref{thm:rigidity},  \eqref{eqn:gamma-t-expand} and an analogous estimate for  the quantiles associated to $\mu_i (t)$ we see that for $|i| < \Nfb$ we have with overwhelming probability
\beq
| \mu_i (t) - \lambda_i^z (t) | \leq \frac{ N^{-\mfa} + \1_{ \{ t \leq t_1/2\}}( \log N)^2}{N}
\eeq
Therefore,
\begin{align} 
& \left[ \Re \frac{1}{ \sqrt{2N}} \sum_{|i| < \Nfb}  \int_{0}^{t_1}\left(   \frac{ 1}{ \lambda_i^z (s) - \i \nu^z (s) }  - \frac{1}{ \mu_i (s) - \i \nu^z (s) } \right) \d W_i (s) \right] \notag\\
\lesssim &\frac{1}{N^{3}} \sum_{ |i| < \Nfb} \int_{0}^{t_1} \frac{N^{-2\mfa} + \1_{ \{ s \leq t_1/2\}} ( \log N)^4}{ | \lambda_i^z (s) - \i \nu^z (s) |^4} \d s \lesssim \frac{1}{N^{2 \mfa} (N \eta_*)^2} + \frac{ ( \log N)^4}{(N t_1)^2}
\end{align}
and so 
\beq \label{eqn:univ-rep-2-b7}
 \frac{\Re}{ \sqrt{2N}} \sum_{|i| < \Nfb}  \int_{0}^{t_1}  \frac{ \d W_i (s)}{ \lambda_i^z (s) - \i \nu^z (s) } =  \frac{\Re}{ \sqrt{2N}} \sum_{|i| < \Nfb}  \int_0^{t_1}  \frac{ \d W_i (s)}{ \mu_i (s) - \i \nu^z (s) } + \O ( N^{-\mfa} + ( \log N)^3 N^{-\omega_1} ),
\eeq
with overwhelming probability.
The estimate \eqref{eqn:univ-rep-2-b3} now follows from \eqref{eqn:univ-rep-2-b4}, \eqref{eqn:univ-rep-2-b6} and \eqref{eqn:univ-rep-2-b7}. This completes the case that $z \in \ddb_r$. If $\beta=1$ and $z \in I_r$, then the proof is identical, except for the proof of \eqref{eqn:univ-rep-2-b4}. Now, we have $4\Re[\langle {\bm w}_i^z,E_1{\bm w}_j^{\overline{z}}\rangle \langle {\bm w}_j^{\overline{z}},E_2{\bm w}_i^z\rangle]=\delta_{ij}- \delta_{i, -j}$ and so the second term in \eqref{eqn:univ-rep-2-b9} is in fact equal to the first one, and so the proof of \eqref{eqn:univ-rep-2-b4} proceeds similarly to the $\beta=2$ case. The other estimates in the proof are the same. 
\qed

\bel \label{lem:univ-rep-2-2} In the set-up of Lemma \ref{lem:univ-rep-2} we have,  with overwhelming probability,
\begin{align} \label{eqn:univ-rep-2-c5} 
& \int_{0}^{t_1} N ( m_{N, s}^z ( \i \eta^z (s) ) - m_s^z ( \i \eta^z (s) ) )^2 \d s  \notag\\
= &  \int_{t_1/2}^{t_1} N  \left( \frac{1}{2N} \sum_{ |i| < \Nfb} \frac{1}{ \mu_i (s) - \i \nu^z (s) }  - \int_{ |x| < \hatgam_{\Nfb}^z } \frac{\rho^z (0)  }{ x - \i \nu^z (s) } \right)^2 \d s + ( \log N)^4 \O ( N^{-\omega_1 } +N^{-\mfa}) ,
\end{align}
The integrals on both sides of the above estimate are $\O ( ( \log N)^{3-C_*} )$ with overwhelming probability. 
\eel
\proof The cases $\beta=1, 2$ are identical. By Theorem \ref{thm:ll-2},
\beq \label{eqn:univ-rep-2-c6}
\int_0^{t_1/2} N | m_{N, s}^z ( \i \eta^z (s) ) - m_s^z ( \i \eta^z(s) ) |^2 \d s \lesssim ( \log N)^2 N^{-\omega_1},
\eeq
with overwhelming probability. For later use, we also record the inequality,
\beq \label{eqn:univ-a-1}
\int_{0}^{t_1} N | m_{N, s}^z ( \i \eta^z (s) ) - m_s^z ( \i \eta^z(s) ) |^2 \d s \lesssim \frac{ ( \log N)^2}{N \eta_*} ,
\eeq
with overwhelming probability, which is proven similarly to \eqref{eqn:univ-rep-2-c6}. 
We have, by Theorem \ref{thm:rigidity},
\beq \label{eqn:univ-rep-2-c2}
\left| \frac{1}{2N} \sum_{ |i| > \Nfb} \frac{1}{ \lambda_i^z (s) - \i \eta^z (s) }  - \int_{|x| > \gamma_{\Nfb}^z(s) } \frac{1}{ x - \i \eta^z (s) } \rho^z_u (x) \d x \right| \leq ( \log N)^2 N^{-\mfb}, 
\eeq
with overwhelming probability. Hence, with overwhelming probability 
\begin{align} \label{eqn:univ-rep-2-c1}
& \int_{t_1/2}^{t_1} N ( m_{N, s}^z ( \i \eta^z (s) ) - m_s^z ( \i \eta^z (s) ) )^2 \d s \notag\\
= &\int_{t_1/2}^{t_1} N  \left( \frac{1}{2 N} \sum_{ |i| < \Nfb} \frac{1}{ \lambda_i^z (s) - \i \eta^z (s) }  - \int_{ |x| < \gamma_{\Nfb}^z (s) } \frac{\rho^z_s (x) }{ x - \i \eta^z (s) } \d x+ \O  ( ( \log N)^2 N^{-\mfb}  ) \right)^2 \d u \notag\\
= &\int_{t_1/2}^{t_1} N  \left( \frac{1}{2 N} \sum_{ |i| < \Nfb} \frac{1}{ \lambda_i^z (s) - \i \eta^z (s) }  - \int_{ |x| < \gamma_{\Nfb}^z (s) } \frac{\rho^z_s (x) }{ x - \i \eta^z (s) } \d x \right)^2 \d s + ( \log N)^4 \O ( N^{-\mfb}).
\end{align}
Note that in the second inequality, in order to bound the cross-term when we expand out the integrand in the second line of \eqref{eqn:univ-rep-2-c1}, we used that the term in the parantheses on the last line of \eqref{eqn:univ-rep-2-c1}   is $( \log N)\O  ([N \eta^z (s)]^{-1} )$ with overwhelming probability (this being a consequence of Theorem \ref{thm:ll-2}, \eqref{eqn:univ-rep-2-c2} and that $\mfb > \omega_1$ by assumption). 

Using \eqref{eqn:homog} and \eqref{eqn:char-dif}  we have with overwhelming probability, for $s \in [t_1/2, t_1]$, 
\beq \label{eqn:univ-rep-3-est1}
\left| \frac{1}{N} \sum_{ |i| < \Nfb} \frac{1}{ \lambda_i^z (s) - \i \eta^z (s) }  - \frac{1}{N} \sum_{ |i| < \Nfb} \frac{1}{ \mu_i^z (s) - \i \nu^z (s) } \right| \lesssim \frac{ N^{-\mfa}}{N \eta^z (s)} + (t_1 + \eta_* ) 
\eeq
Next, using \eqref{eqn:gamma-t-expand}, Lemma \ref{lem:proprho}\ref{it:rho-der}, \eqref{eqn:char-dif}, and \eqref{eq:exp}, we have,
\begin{align} \label{eqn:univ-rep-3-est2}
\left| \int_{ |x| < \gamma_{\Nfb}^z (s) } \frac{\rho^z_s (x) }{ x - \i \eta^z (s) }\d x - \int_{ |x| < \hatgam_{\Nfb}^z } \frac{\rho^z (0)  }{ x - \i \nu^z (s) }\d x \right| \lesssim ( \log N) \left( t_1 + \gamma_{\ell_1}^z (t_1) + \eta_*  \right) 
\end{align}
Both of the error terms in \eqref{eqn:univ-rep-3-est1} and \eqref{eqn:univ-rep-3-est2} are $\O ( [ N \eta_z (s) ]^{-1} )$. Therefore by a similar argument to \eqref{eqn:univ-rep-2-c1} we have with overwhelming probability, 
\begin{align} \label{eqn:univ-rep-2-c3}
& \int_{t_1/2}^{t_1} N  \left( \frac{1}{N} \sum_{ |i| < \Nfb} \frac{1}{ \lambda_i^z (s) - \i \eta^z (s) }  - \int_{ |x| < \gamma_{\Nfb}^z (s) } \frac{\rho^z_s (x) }{ x - \i \eta^z (s) } \d x\right)^2 \d s \notag\\
= & \int_{t_1/2}^{t_1} N  \left( \frac{1}{N} \sum_{ |i| < \Nfb} \frac{1}{ \mu_i (s) - \i \nu^z (s) }  - \int_{ |x| < \hatgam_{\Nfb}^z  } \frac{\rho^z (0) }{ x - \i \eta^z (s) } \d x \right)^2 \d s 
+  ( \log N)^4 \O( N^{-\mfa}  ) 
\end{align}
The estimate \eqref{eqn:univ-rep-2-c5} now follows from \eqref{eqn:univ-rep-2-c6}, \eqref{eqn:univ-rep-2-c1} and \eqref{eqn:univ-rep-2-c3}. The other claim of the lemma is a consequence of \eqref{eqn:univ-rep-2-c5} and \eqref{eqn:univ-a-1}. \qed

\

Finally, we need to estimate the last term on the RHS of \eqref{eqn:univ-rep-2-b2} in the case $\beta=1$.

\bel \label{lem:univ-rep-3} In the case $\beta=1$ we have with overwhelming probability
\begin{align} \label{eqn:univ-rep-3-aa1}
&\int_0^{t_1}   \tilde{\sum}_{ij}\langle G_s^z(\i \eta^z(s) )E_iG_s^{\overline{z}}(\i \eta^z (s) )E_j\rangle  \d s + \O ( N^{-\mfa} + \O (N^{-\omega_1/2} ) ) \notag\\
= & \1_{ z \in I_r} \left\{  \int_{t_1/2}^{t_1} \frac{1}{4N} \sum_{|i| \leq N^{\mfb}} \left( \frac{1}{( \mu_i (s) - \i \nu^z (s) )^2} + \frac{1}{ | \mu_i (s) - \i \nu^z (s) |^2}\right) \d u + \frac{1}{2} \log 2 \right\}
\end{align}
Additionally when $ z\in I_r$ we have
\begin{align} \label{eqn:univ-rep-3-a2}
 & \int_0^{t_1}  \tilde{\sum}_{ij}\langle G_s^z(\i \eta^z(s) )E_iG_s^{\overline{z}}(\i \eta^z (s) )E_j\rangle  \d s
 = \frac{1}{2} \log ( \sqrt{1-|z|^2} t_1 / \eta_* ) + \O ( ( \log N)^{-10} )
\end{align}
In particular, with overwhelming probability the absolute value of the RHS of \eqref{eqn:univ-rep-3-aa1} is less than $2 \omega_1 \log N$. 
\eel
\proof The case when $z \notin I_r$ follows immediately from Lemma \ref{lem:2grealb}.  In the remainder of the proof we deal with the case $z \in I_r$. We have
\beq
\begin{split}
\tilde{\sum}_{ij} \langle G_s^z(\i \eta^z (s) )E_iG_s^{\overline{z}}(\i \eta^z(s) )E_j\rangle &= \frac{1}{2N}\sum_{i,j}\frac{\langle {\bm w}_i, E_1{\bm w}_j\rangle\langle {\bm w}_j, E_2{\bm w}_i\rangle+\langle {\bm w}_i, E_2{\bm w}_j\rangle\langle {\bm w}_j, E_1{\bm w}_i\rangle }{(\lambda_i^z(u)-\i\eta^z(s)(\lambda_j^z(s)-\i\eta^z(s))} \\
&= \frac{1}{4N} \sum_{i} \frac{1}{ ( \lambda_i^z(s) -\i \eta^z(s) )^2}  + \frac{1}{4N} \sum_i \frac{1}{ |\lambda_i^z(s)-\i \eta^z(s)|^2},
\end{split}
\eeq
where in the second equality we used that $2\langle {\bm w}_i, E_1{\bm w}_j\rangle=\delta_{ij}+\delta_{i,-j}$, $2\langle {\bm w}_i, E_2{\bm w}_j\rangle=\delta_{ij}-\delta_{i,-j}$, and that $\lambda_{i}^z (s) = - \lambda_{-i}^z(s)$. 

We first consider $s \in [0, t_1/2]$. By Theorem \ref{thm:ll-1} and the Cauchy integral formula we see that with overwhelming probability for any $\eps >0$
\begin{align} \label{eqn:univ-rep-3-aa2}
\int_0^{t_1/2} \frac{1}{2N} \sum_i \frac{1}{ (\lambda_i^z (s) - \i \eta^z(s) )^2} \d s &= \int_0^{t_1/2} (\del_w m_s^z ( \i \eta^z(s) ) ) \d s + \int_0^{t_1/2} \O \left( \frac{N^{\eps}}{N t_1^2 } \right) \d s \notag\\
&= \O ( t_1 + N^{\eps} N^{- \omega_1} ) = \O (N^{-\omega_1/2} ) .
\end{align}
Using Theorem \ref{thm:ll-1} as above, for the other term we have 
\begin{align} \label{eqn:univ-rep-3-a3}
& \int_0^{t_1/2} \frac{1}{4N} \sum_i \frac{1}{ | \lambda_i^z (s) - \i \eta^z(s)|^2} \d s =\frac{1}{2}  \int_0^{t_1/2}  \frac{\Im[ m_{N, s}^z (\i \eta^z(s) ) ]}{ \eta^z(s)} \d s \notag\\
&= \frac{1}{2} \int_0^{t_1/2} \frac{ \Im [ m_s^z ( \i \eta^z(s)) ]}{ \eta^z(s)} \d s + \O ( N^{-\omega_1/2} ) = \frac{1}{2} \left( \log ( \eta^z(0) ) - \log ( \eta^z(t_1/2))\right) + \O ( N^{-\omega_1/2} ) \notag\\
&= \frac{1}{2} \left( \log ( t_1 \sqrt{1-|z|^2} ) - \log \left( \frac{t_1}{2} \sqrt{1-|z|^2} \right) \right) + \O ( N^{-\omega_1/2} ) =   \frac{1}{2} \log 2 + \O ( N^{-\omega_1/2} ).
\end{align}
In the second line we used \eqref{eqn:char-def} and then in the last line used Lemma \ref{lem:approx-char}. For $s \in [t_1/2, t_1]$ we will apply \eqref{eqn:homog}. We first find,
\beq
\int_{t_1/2}^{t_1} \frac{1}{N} \sum_{ |i| > N^{\mfb}} \frac{1}{ | \lambda_i^z (s) - \i \eta^z(s) |^2} \d s \lesssim N^{\omega_1-\mfb}
\eeq
with overwhelming probability. Further, using \eqref{eqn:homog} as well as \eqref{eqn:char-dif} we have with overwhelming probability
\begin{align}
    & \int_{t_1/2}^{t_1} \frac{1}{4N} \sum_{|i| \leq N^{\mfb}} \left( \frac{1}{( \lambda_i^z (s) - \i \eta^z(s) )^2} + \frac{1}{ | \lambda_i^z (s) - \i \eta^z(s) |^2}\right) \d u \notag\\
    = &   \int_{t_1/2}^{t_1} \frac{1}{4N} \sum_{|i| \leq N^{\mfb}} \left( \frac{1}{( \mu_i (s) - \i \nu^z (s) )^2} + \frac{1}{ | \mu_i (s) - \i \nu^z (s) |^2}\right) \d s + \O ( N^{-\mfa} ) .
\end{align}
We used that the difference of the integrands  is $\O \left( \frac{N^{-\mfa}}{N( \eta^z (s))^2} + \frac{ t_1 + \eta_*}{\eta^z(s)} \right)$, which integrates to $\O ( N^{-\mfa})$. 
This completes the proof of \eqref{eqn:univ-rep-3-aa1}.  We now turn to \eqref{eqn:univ-rep-3-a2}. An argument similar to \eqref{eqn:univ-rep-3-aa2} using instead \cite[Proposition 3.2]{cipolloni2025maximum} and the Cauchy integral formula gives
\beq
\int_0^{t_1} \frac{1}{4 N} \sum_i \frac{1}{ ( \lambda_i^z (s) - \i \eta^z (s) )^2} \d s = \O ( ( \log N)^{-10} )
\eeq
with overwhelming probability. Moreover, a slight modification of \eqref{eqn:univ-rep-3-a3} yields
\beq
\int_0^{t_1} \frac{1}{4 N} \sum_i \frac{1}{ | \lambda_i^z(s) + \i \eta^z (s) |^2} \d s= \frac{1}{2} \log ( \sqrt{1- |z|^2} t_1 / \eta_* ) + \O ( ( \log N)^{-10} )
\eeq
with overwhelming probability. This completes the proof. \qed

\vspace{5 pt} 
\noindent{\bf Proof of Lemma \ref{lem:univ-rep-2}}. The estimate \eqref{eqn:univ-rep-2} follows from the decompositions \eqref{eqn:univ-rep-2-b1}, \eqref{eqn:univ-rep-2-b2} and the estimates \eqref{eqn:univ-rep-2-d1}, \eqref{eqn:univ-rep-2-b3}, \eqref{eqn:univ-rep-2-c5}, and Lemma \ref{lem:univ-rep-3}. The claim about the third line of \eqref{eqn:univ-rep-2} follows from the last claim of Lemma \ref{lem:univ-rep-2-2} and the claim about the fourth line from Lemma \ref{lem:univ-rep-3}. 
\qed

\subsection{Local-global decomposition}

In this section we will use the estimates of Section \ref{sec:univ-rep} to derive Theorem \ref{thm:local-global} below. This theorem gives a high probability representation of the characteristic polynomial (regularized at a scale $\eta_m = N^{-1-\delta_m}$) evaluated at $K$ different $z_i$ in terms of independent local variables and global variables which have some correlation. 

Eventually we will need to compute the contribution of the local component. In order to do so, we first need to compute the dependence of the local variable on the parameter $z$. This takes place in the next subsection.

\subsubsection{Rescaling the local variables} \label{sec:rescale}

In Section \ref{sec:univ-rep} we compared the process $\lambda_i^z (t)$ to a process $\mu_i^z (t)$, defined in the statement of Theorem \ref{thm:homog}. In this section we will show that the dependence on $z$ of the various quantities involving $\mu_i^z (t)$ that were introduced in Section \ref{sec:univ-rep} is through a simple re-scaling. 

Let us fix some $z \in \ddhb_r$ and let $B = B_z := \frac{ \rho^z (0)}{ \rho^0 (0)} = \sqrt{ 1 - |z|^2}$. Continuing with the notation of Section \ref{sec:univ-rep}, we  recall that $\mu_i^z (t) = \mu_i (t)$ solves  \eqref{eqn:mu-n-i-dbm} with initial data $B^{-1} \lambda_i (H^0 (G))$ with $G$ a complex or real Ginibre matrix. Let now $\tilmu_i (t) = B \mu_i (t /B^2)$. This process satisfies the equation, 
\beq \label{eqn:tilmu-eqn} 
\d \tilmu_i (t) = \frac{ \d \tilW_i (t)}{ \sqrt{2N}} + \frac{1}{2N} \sum_{j \neq i } \frac{1 - \delta_{i, -j} \1_{ \beta=1} \1_{ z \in \rr }}{ \tilmu_i (t) - \tilmu_j (t) } \d t
\eeq
for Brownian motions $\tilW_i (s) = B W_i (s /B^2)$ (here the $W_i = W_i^z$ are from \eqref{eqn:mu-n-i-dbm}) that have quadratic variation \eqref{eqn:W-qv} with initial data $\tilmu_i (0) = \lambda_i (H^0 (G))$. Similarly,  (recalling the definition of $\nu^z(s)$ in \eqref{eqn:approx-char}) we have
$
 B \nu^z (s / B^2)  = B \eta_* + B^2 t_1 - s .
$  The quantities on the RHS of \eqref{eqn:univ-rep-1-1} become
\beq
\frac{1}{2} \int_{\eta_m}^{\eta_*} \sum_{ |i| < N^{\mfb}} \frac{ u}{ ( \mu_i (t_1))^2 + u^2} \d u  = \frac{1}{2} \int_{B \eta_m}^{B \eta_*} \sum_{ |i| < N^{\mfb}} \frac{ u}{ ( \tilmu_i (B^2 t_1))^2 + u^2} \d u 
\eeq
and (recall the definition \eqref{eqn:gamma-hat-def})
\beq
N \int_{\eta_m}^{\eta_*} \int_{|x| < \hatgam^z_{N^{\mfb}}} \frac{ u}{x^2+u^2} \rho^z (0) \d x \d u =  N \int_{B\eta_m}^{B\eta_*} \int_{|x| < \hatgam^0_{N^{\mfb}}} \frac{ u}{x^2+u^2} \rho^0 (0) \d x \d u  .
\eeq
Similarly, the quantities on the RHS of \eqref{eqn:univ-rep-2} become
\beq
\frac{1}{2}\Re \frac{1}{ \sqrt{2N}} \sum_{|i| < \Nfb}  \int_{0}^{t_1}  \frac{ \d W_i^z (s)}{ \mu_i (s) - \i \nu^z (s) } = \frac{1}{2}\Re \frac{1}{ \sqrt{2N}} \sum_{|i| < \Nfb}  \int_{0}^{B^2 t_1}  \frac{ \d \tilW_i (s)}{ \tilmu_i (s) - \i (B \eta_* + B^2t_1 - s ) }
\eeq
and
\begin{align}
   &  \Re  \int_{t_1/2}^{t_1} \frac{N}{2}  \left( \frac{1}{2N} \sum_{ |i| < \Nfb} \frac{1}{ \mu_i (s) - \i \nu^z (s) }  - \int_{ |x| < \hatgam_{\Nfb}^z } \frac{\rho^z (0) \d x }{ x - \i \nu^z (s) } \right)^2 \d s \notag\\
     = & \Re \int_{ B^2 t_1/2}^{B^2 t_1} \frac{N}{2} \left( \frac{1}{2N} \sum_{ |i| < \Nfb} \frac{1}{ \tilmu_i (s) - \i ( B \eta_* + B^2 t_1 -s )} - \int_{ |x| < \hatgam_{\Nfb}^0} \frac{\rho^0 (0) \d x }{ x- \i ( B \eta_* + B^2 t_1 -s )}  \right)^2 \d s ,
    \end{align}
    and
    \begin{align}
        & \int_{t_1/2}^{t_1} \frac{1}{4N} \sum_{|i| \leq N^{\mfb}} \left( \frac{1}{( \mu_i (s) - \i \nu^z (s) )^2} + \frac{1}{ | \mu_i (s) - \i \nu^z (s) |^2}\right) \d s  \notag\\
        = &  \int_{B^2t_1/2}^{B^2t_1} \frac{1}{4N} \sum_{|i| \leq N^{\mfb}} \left( \frac{1}{( \tilmu_i (s) - \i(B \eta_* + B^2 t_1 -s ) )^2} + \frac{1}{ | \tilmu_i (s) - \i (B \eta_* + B^2 t_1 -s)) |^2}\right) \d s
    \end{align}
    The above motivates the introduction of the following \emph{local variables.}  For positive parameters $\eta_m < \eta_*,  \mfb, t_1$ and $z \in \cc$ and $\beta \in \{1, 2\}$ define the four random variables
\begin{align} \label{eqn:tilL-def}
\tilL_1 ( t_1, \eta_*, \eta_m,  \mfb; z, \beta ) &:= N \int_{\eta_m}^{\eta_*} \int_{|x| < \hatgam^0_{N^{\mfb}}} \frac{ u}{x^2+u^2} \rho^0 (0) \d x \d u  - \frac{1}{2} \int_{ \eta_m}^{ \eta_*} \sum_{ |i| < N^{\mfb}} \frac{ u}{ ( \tilmu_i (t_1))^2 + u^2} \d u,\notag\\
\tilL_2( t_1, \eta_*, \eta_m,  \mfb ; z, \beta ) &:= \frac{1}{2} \Re \frac{1}{ \sqrt{2N}} \sum_{|i| < \Nfb}  \int_{0}^{ t_1}  \frac{ \d \tilW_i (s)}{ \tilmu_i (s) - \i ( \eta_* + t_1 - s ) },\notag\\
 \tilL_3( t_1, \eta_*, \eta_m, \mfb ; z , \beta ) & \notag\\
:= \Re \int_{  t_1/2}^{ t_1} \frac{N}{2} & \left( \frac{1}{2N} \sum_{ |i| < \Nfb} \frac{1}{ \tilmu_i (s) - \i ( \eta_* +  t_1 -s )} - \int_{ |x| < \hatgam_{\Nfb}^0} \frac{\rho^0 (0) \d x }{ x- \i (  \eta_* + t_1 -s )}  \right)^2 \d s,
\end{align}
and 
\begin{align} \label{eqn:tilL-def-2}
\tilL_4 (t_1, \eta_*, \eta_m,  \mfb ; z , \beta )   \notag:&=  - \1_{ \beta=1} \1_{ \{ z \in I_r \}} \bigg\{ \Re  \int_{t_1/2}^{t_1} \frac{1}{8N} \sum_{|i| \leq N^{\mfb}} \frac{1}{( \tilmu_i (s) - \i( \eta_* +  t_1 -s ) )^2} \d s  \notag\\
&\quad + \int_{t_1/2}^{t_1} \frac{1}{8N} \sum_{|i| \leq N^{\mfb}} \frac{1}{| \tilmu_i (s) - \i( \eta_* +  t_1 -s ) |^2} \d s +  \frac{1}{4} \log 2 \bigg\}.
\end{align}
While not all of the  $\tilL_i$ depend on all of the parameters, we keep the arguments for a uniform notation. Note that the distribution depends on $z$ and $\beta$ in the following manner. First, it is the same for all $z$ and $\beta$ except if $\beta=1$ and $z \in I_r$. In the latter, case it is the same for all $z \in I_r$. 
 This will eventually allow us to compute (in Section \ref{sec:1pt}) their contribution in both the real and complex i.i.d. cases using only Proposition \ref{prop:ginibre}.

For additional parameters $A_1>0$ and $\gamma >0$,  we define the \emph{local factor}, 
\begin{align} \label{eqn:local-factor-def} 
\calL ( t_1, \eta_*, \eta_m, & \mfb , \gamma, A_1 ; z, \beta ) \notag\\
 &:= \ee \left[ \e^{ \gamma (\sum_{a=1}^4 \tilL_a ) } \1_{ \{ \tilL_1 \leq ( \log N)^{5/8} \} \cap  \{| \tilL_3 | \leq 1  \} \cap  \{ | \tilL_2 | \leq A_1 \log N \} } \1_{\{ | \tilL_4 |  \leq  \log N \} }\right]
\end{align}
We also define
\beq \label{eqn:tilZ-def}
\tilZ (t_1, \eta_*, \eta_m,  \mfb; z, \beta) := \frac{1}{2} \Re \frac{1}{ \sqrt{2N}} \sum_{ |i | < \Nfb} \int_0^{t_1} \frac{ \d \tilW_i (s) }{ \gamma_i^0 (s) - \i ( \eta_* +t_1 -s )}.
\eeq
This is a Gaussian random variable and for parameters as in Lemma \ref{lem:univ-rep-2} one has,
\beq  \label{eqn:tilZ-def-est} 
\Var (\tilZ ) = \frac{1 + \1_{ \beta=1}\1_{ z \in I_r}}{4} \log (t_1 / \eta_*)+ \O ( ( \log N)^{-1} ) , \qquad | \tilZ - \tilL_2 | \lesssim \frac{1}{( \log N)^2},
\eeq
with the latter estimate holding with overwhelming probability. We now explain the proof of \eqref{eqn:tilZ-def-est}. The second estimate is proven in the same manner as Lemma \ref{lem:univ-rep-2-1} (see also the similar proof of \cite[Lemma 5.2]{cipolloni2025maximum}). For the first estimate we have
\begin{align}
 & \frac{1}{2(1+\1_{ \beta=1} \1_{z \in I_r} ) N} \int_0^{t_1} \left[ \sum_{ |i| < \Nfb} \Re \frac{ \d \tilW_i (s)}{ \gamma_i^0 (s) - \i \nu^0 (s) } \right]  = \int_0^{t_1} \frac{1}{N} \sum_{ |i| < \Nfb} \frac{( \gamma_i^0 (s))^2}{  |\gamma_i^0 (s) - \i \nu^0 (s) |^4} \d s \notag \\
 &=  2 \int_0^{t_1} \frac{ (\pi x )^2}{ | (\pi x - \i \nu^0 (s) |^4} \d x \d s + \O ( ( \log N)^{-10} ) \notag \\
 & =  \int_0^{t_1} \int_{\rr} \frac{1}{ |x \pi - \i \nu^0 (s) |^2} \d x \d s + \Re \int_0^{t_1} \int_{\rr} \frac{1}{ (x \pi  - \i \nu^0 (s) )^2} \d x \d s + \O ( ( \log N)^{-10} ) \notag \\
 & =  \int_0^{t_1} \frac{ 1}{\eta_* + s } \d s + \O ( ( \log N)^{-10} )  =  \log \left( \frac{  t_1 }{\eta_*} \right) + \O ( ( \log N)^{-10} ) .
\end{align}

\subsubsection{Local-global decomposition}

We now combine the homogenization Theorem \ref{thm:homog}, the estimates of Section \ref{sec:univ-rep} and the rescaling done in Section \ref{sec:rescale} to derive the following.

\bet \label{thm:local-global}
Let $r >0$ and $b>0$ and let $z_1, \dots, z_K \in \ddb_r$ satisfy $|z_i -z_j| \geq N^{b-1/2}$ for $i \neq j$. There is a $\mfb >0$ so that the following holds. Choose $\omega_1 \in (0, \mfb)$ and let $\mfa >0$ be the exponent from Theorem \ref{thm:homog}. Choose exponents $0 < \delta_m < \mfa$, $\omega_1 <  \mfb$ and $C_* \geq 1000$. Define $\eta_m = N^{-1-\delta_m}$, $\ell_1 = N^{q_1}$ and $\eta_* = ( \log N)^{C_*}/N$. Define $B_i := \sqrt{1- |z_i|^2}$.

Let $X_0$ be a real or complex i.i.d. matrix, with $X_t$ as in \eqref{eqn:Xt-def}. Let the processes $\{ \lambda_j^{(z_i)} (t)\}_{j, i, t}$, $\{ b_j^{z_i}(t) \}_{j, i, t}$, and $\{ \Phi_N (z_i, \eta, t_1 ) \}_{i,\eta}$  be as in Section \ref{sec:def-dynamics}. 

There is a coupling of these processes and the initial data $X_0$ in which there are random variables $\{\tilL_a^{(i)}\}_{a \in [\![3]\!], i\in [\![K]\!]}$, $\{\tilZ^{(i)} \}_{ i \in [\![K]\!]}$, $\{\tilG^{(i)}  \}_{ i \in [\![K]\!]}$,  processes $\{ \tilW^{(i)}_j (t)\}_{i, j, t}$ and $\{ \tilmu_j^{(i)} (t) \}_{j, i, t}$ so that the following holds.

\begin{enumerate}[label=(\roman*),font=\normalfont]
    \item \label{it:Fi} The families $\tilF_i$ of random variables and processes defined by
    \beq
    \tilF_i := \left\{  \{ \tilW^{(i)}_j (t) , \tilmu_j^{(i)} (t) \}_{j, t}, \{ \tilL^{(i)}_a \}_{a \in [\![4]\!]}, \tilZ^{(i)} \right\}
    \eeq
    are mutually independent from each other as well as from the family $\{ \tilG^{(i)}\}_{i \in [\![K]\!]}$ and the matrix $X_0$. 
    \item For each $i$, the joint distribution of $\tilF_i$ is as follows. The $\tilW_j^{(i)}$ are independent  Brownian motions (with quadratic variation as in \eqref{eqn:W-qv}) and the $\tilmu^{(i)}_j$ are as in \eqref{eqn:tilmu-eqn}. In terms of the $\tilW_j^{(i)}$ and $\tilmu^{(i)}_j$, the random variables $\tilL_a^{(i)}$ and $\tilZ^{(i)}$ are defined by,
    \beq
    \tilL_a^{(i)} = \tilL_a (B_i^2 t_1, B_i \eta_*, B_i \eta_m, \ell_1 , \mfb; z_i, \beta) , \quad \tilZ^{(i)} = \tilZ (B_i^2 t_1, B_i \eta_*, B_i \eta_m, \ell_1 , \mfb  ;z_i,  \beta) 
    \eeq
    where the $\tilL_a$ are as in \eqref{eqn:tilL-def} and \eqref{eqn:tilL-def-2}, and the $\tilZ$ is as in \eqref{eqn:tilZ-def}. 
    \item \label{it:lg-G} We have that $\tilG^{(i)} = \Phi_N (z_i, t_1 B_i, 0)$. 
    \item \label{it:lg-Z} Each $\tilZ^{(i)}$ is a Gaussian with variance
    \beq
\Var ( \tilZ^{(i)} ) = \frac{1+ \1_{ \beta=1} \1_{ z\in I_r}}{4} \log ( B_i t_1 / \eta_*) + \O ( ( \log N)^{-1} )
    \eeq
    and $| \tilZ^{(i)} - \tilL_2^{(i)} | \leq ( \log N)^{-1}$ with overwhelming probability. Therefore, for any $D >0$ we have for all $u \geq 1$,
    \beq
\pp\left[ | \tilL_2^{(i)} | \geq u \log N  \right] \leq \e^{ - 10 u^2 \log N} + N^{-D} .
    \eeq

    Additionally for any $L \geq 1$ and $A \geq 1$ we have for all $ \lambda \in [0, L]$ that,
    \beq \label{eqn:lg-L2-mfg}
\ee\left[ \e^{ \lambda \tilL_2^{(i)}} \1_{ \{ \tilL_2^{(i)} \leq A \log N  \}} \right]  \lesssim N^{ \lambda^2/10}
    \eeq
    \item \label{it:lg-La} For any $\delta>0$ we have $\tilL_1^{(i)} \leq ( \log N)^{1/2+\delta}$ with overwhelming probability,  $| \tilL_3^{(i)} | \leq ( \log N)^{-1}$ with overwhelming probability and that  
    \beq
    \tilL_4^{(i)} = - \frac{\1_{\beta=1} \1_{ z \in I_r}}{4} \left( \log (B_i t_1 / \eta_* ) + \O ( ( \log N)^{-1} ) \right) 
    \eeq
    with overwhelming probability. 
    \item \label{it:lg-homog} We have with overwhelming probability,
    \beq  \label{eqn:local-global-homog}
| \lambda_1^{(z_i)} (t_1 ) - B_i^{-1} \tilmu^{(i)}_1 (B_i^2 t_1) | \leq N^{-1-\mfa} .
    \eeq
    \item With overwhelming probability, 
\begin{align} \label{eqn:lg-main-est}
& \Phi_N (z_i, \eta_m, t_1) = \sum_{a=1}^3 \tilL_a^{(i)} + \tilG^{(i)} 
+ ( \log N)^4 \O ( N^{\omega_1/2 - \mfb/2}  +  N^{\delta_m-\mfa}  +  N^{-\omega_1/2 } ) 
\end{align}
\item With overwhelming probability,
\beq \label{eqn:lg-main-bd} 
\Phi_N (z_i, \eta_*, t_1) = \tilZ^{(i)} + \tilG^{(i)}  - \1_{ \beta=1} \1_{ z \in I_r} \log (B_i t_1 / \eta_* )+ \O ( ( \log N)^{-1} ).
\eeq
\end{enumerate}
\eet
\proof  The processes $ \tilmu_j^{(i)} (t)$ and $\tilW^{(i)}_j (t)$ come from the set-up of Theorem \ref{thm:homog} and the rescalings done in Section \ref{sec:rescale}. The random variables $\tilG^{(i)}$ are defined by $\tilG^{(i)} := \Phi_N (z_i, t_1 B_i ; X_0)$ and so the first three items follow from Theorem \ref{thm:homog} and by definition. Item \ref{it:lg-Z} follows from \eqref{eqn:tilZ-def-est}. Item \ref{it:lg-La} follows from \eqref{eqn:univ-rep-1-2} and the claim immediately after \eqref{eqn:univ-rep-2}, translated into statements about the $\tilL_a^{(i)}$. Item \ref{it:lg-homog} follows from \eqref{eqn:homog}. The estimate \eqref{eqn:lg-main-est} follows from \eqref{eqn:lg-decomp} and Lemmas \ref{lem:univ-rep-1} and \ref{lem:univ-rep-2}. The estimate \eqref{eqn:lg-main-bd} follows from Lemma \ref{lem:univ-rep-2} and items \ref{it:lg-Z} and \ref{it:lg-La}. \qed

\subsection{A useful upper bound}

Later we will require an estimate for the characteristic polynomial when $z$ and $w$ are close. 

\bep \label{prop:useful-upper} 
Let $\eta_* = ( \log N)^{C_*}/ N, t_1 = N^{\omega_1-1}$ with $ \omega_1 \in (0, \frac{1}{100})$ and $C_* \geq 1000$.  Let $z, w \in \dd_r$. There are random variables $\hatZ_z$ and $\hatZ_w$  with the property that for any $D>0$,
\beq
\pp\left[ | \hatZ_u | > ( \omega_1)^{1/3} \log N \right] \lesssim \e^{ - \frac{ ( \omega_1)^{-1/3}}{10} \log N} + N^{-D}
\eeq
such that with overwhelming probability, for $u \in \{z, w\}$ and $ \eta \in [0, \eta_*]$ we have
\beq \label{eqn:useful-upper}
\Phi_N (u, \eta, t_1) \leq \O  ( ( \log N)^{1/2+\delta} ) + \Phi_N (u, B_u t_1 , 0) + \hat{Z}_u  - \frac{ \1_{\beta=1} \1_{z \in I_r} }{4} \omega_1 \log N 
\eeq
for any $\delta >0$. 
\eep
\begin{remark} \emph{Note that we make no claim about the joint Gaussianity of $\hatZ_z$ and $\hatZ_w$ or their independence from the initial data.}\end{remark}

\

\proof[Proof of Proposition~\ref{prop:useful-upper}] By Lemma \ref{lem:a-priori-1} it suffices to consider $\Phi_N (u, \eta_*, t_1)$. By \eqref{eqn:univ-rep-2-b1}, \eqref{eqn:univ-rep-2-b2}, Lemma \ref{lem:univ-rep-2-2} and Lemma \ref{lem:univ-rep-3} (note that these arguments do not depend on the coupling of Theorem \ref{thm:homog}) we see that
\begin{align}
\Phi_N (u, \eta_* , t_1) &= \left(\Re \frac{1}{\sqrt{8N}} \int_0^{t_1} \sum_i \frac{ \d b_i^u(s)}{ \lambda_i^u(s) - \i \eta^u (s) } \right) + \Phi_N (u, B_u t_1, 0) \notag\\
-& \frac{ \1_{\beta=1} \1_{z \in I_r} }{4}\log (B_u t_1/\eta_*)+ \O ( ( \log N)^{-1} ),
\end{align}
with overwhelming probability. We define the first random variable on the RHS to be $\hatZ_u$. By inspecting the proof of Lemma \ref{lem:univ-rep-2-1} we see that the quadratic variation of $\hatZ_u$ is bounded above by
\beq
[ \hatZ_u] \leq \int_0^{t_1} \frac{ \Im[m_{N, s}^u ( \i \eta^u (s) ]}{ \eta^u (s)} \d s + N^{-1/4} \leq \omega_1 \log N ,
\eeq
with overwhelming probability, applying Theorem \ref{thm:ll-2}. The claim now follows from this and Lemma \ref{lem:mart-bd}. \qed

\section{Preliminary analysis of one-point function}

In this section we will develop several computations concerning the one-point function $\ee[ | \det (X-z)|^\lambda]$. In particular, we will show how to regularize this quantity (i.e., compare it with $\Phi_N (z, \eta; X)$ for appropriate $\eta >0$) and derive a product representation when $X = X_t$, with $X_t$ as in Section \ref{sec:def-dynamics}. We will also compute the local factor $\calL$ defined in \eqref{eqn:local-factor-def} through the exact expressions available for the Ginibre ensemble of Proposition \ref{prop:ginibre}. Throughout this section we fix an $r \in (0, 1)$ and we will always work with $z \in \ddhb_r$. We define the following quantity

\beq \label{eqn:control-scalar}
\cD (z, \lambda) := \exp \left(  \frac{ \lambda^2}{8} \log N +  \1_{\beta=1} \1_{ z \in I_r}  \frac{ \lambda(\lambda-2)}{8} \log N \right) = \begin{cases} N^{\lambda(\lambda-1)/4} , & \beta =1, z \in I_r \\ N^{\lambda^2/8} , & \mathrm{else} \end{cases} ,
\eeq

which is the expected order of magnitude of the one-point function (compare with Proposition \ref{prop:ginibre}). The reader may find it convenient to ignore the case $\beta=1, z \in I_r$ and replace $\cD(z, \lambda)$ by $N^{\lambda^2/8}$ everywhere. The difference is largely notational.

\label{sec:1pt}

\subsection{Regularization of one-point function}

Let $X$ be a real or complex i.i.d. matrix and $\Phi_N (z, \eta ) = \Phi_N (z, \eta ; X)$. For $A>0$ we define the cut-off event
\beq \label{eqn:A-def} 
\A(z, \eta, A) := \{ \Phi_N (z, \eta) \leq A \log N \} .
\eeq
We first establish the following a-priori bounds for $\Phi_N (z, \eta)$. Note that when $\eta \lesssim N^{-1}$ we expect that the LHS of \eqref{eqn:1-pt-mgf-bd} is of order $  \cD(z, \lambda) $ for sufficiently large $A$.
\bel \label{lem:1-pt-apriori} 
Let $\eta_* := (\log N)^{C_*} /N$, for some $C_* \geq 1000$ sufficiently large and let $L \geq 100$ and $\delta >0$. Uniformly for $\lambda \in [0, L]$ and $A \in [0, L]$ and $ \eta \in [0, \eta_*]$ we have,
\beq \label{eqn:1-pt-mgf-bd} 
\ee\left[ \e^{ \lambda \Phi_N (z, \eta)} \1_{ \A (z, \eta, A) } \right] \lesssim \e^{ \lambda ( \log N)^{1/2+\delta}} \cD (z, \lambda ) .
\eeq
Uniformly in $\eta \in [0, 1]$, for $10 \leq u \leq L$, we have
\beq \label{eqn:1-pt-tail} 
\pp\left[ \Phi_N (z, \eta) > u \log N \right] \lesssim \e^{ - u^2 \log N } + \1_{ \beta=1} \1_{z \in I_r} \e^{ - \frac{u^2}{2} \log N } .
\eeq
\eel
\proof The estimate \eqref{eqn:1-pt-tail} follows from \eqref{eqn:meso-tail-bd} and Lemma \ref{lem:a-priori-1}. 

Turning to the estimate \eqref{eqn:1-pt-mgf-bd}, from \eqref{eqn:meso-1pt-mgf}  and Proposition \ref{prop:parameter-properties} we first see that
\beq \label{eqn:1-pt-mgf-bd-a1}
\ee\left[ \e^{ \lambda \Phi_N (z, \eta_*)} \1_{ \A (z, \eta_*, A) } \right] \lesssim \cD (z, \lambda) \exp\left( \1_{z \in I_r} \1_{ \beta=1} C_* \lambda \log \log N \right) .
\eeq
Note that the second factor on the RHS comes from the $\frac{\lambda}{4} \log \eta$ correction to the expectation in the $\beta=1$, $z \in I_r$ case. 

For smaller $\eta < \eta_*$, using \eqref{eqn:1-pt-tail}, we first have for any $C>0$ 
\beq
\ee\left[ \e^{ \lambda \Phi_N (z, \eta)} \1_{ \A (z, \eta, A)} \right] = \ee\left[ \e^{ \lambda \Phi_N (z, \eta)} \1_{\A (z, \eta, A) \cap \A (z, \eta_*, C)} \right] + \O ( N^{LA - C^2/2} ),
\eeq
We then estimate the term on the RHS using Lemma \ref{lem:a-priori-1} and  \eqref{eqn:1-pt-mgf-bd-a1}, yielding \eqref{eqn:1-pt-mgf-bd}. \qed

\

The following estimate allows us to control $\Phi_N (z, \eta)$ on the event that $\lambda_1^z$ is small. 

\bel \label{lem:1-pt-reg}
Let $\eta_w = N^{-1-\delta_w}$ with $0 < \delta_w \leq \frac{c_W}{2} \wedge \frac{1}{100}$ where $c_W$ is from Lemma \ref{lem:wegner}. Let $L\geq 1$. Uniformly for $\eta \in [0, \eta_w]$, $\lambda \in [0, L]$, and $A \in [0, L]$  we have for any $\delta \in (0, \frac{1}{2})$,
\begin{align} \label{eqn:1-pt-reg-a1}
 &\ee\left[ \e^{ \lambda \Phi_N (z,\eta ) } \1_{ \A (z, \eta, A) \cap \{ \lambda_1^z \leq \eta_w \} } \right] \lesssim \e^{ \lambda ( \log N)^{1/2+\delta} } \cD (z, \lambda) N^{-\lambda \delta_w - ( \frac{1}{2} \delta_w - \frac{\lambda^2}{4})_+ } \notag\\
 \lesssim & N^{- \frac{ \lambda \delta_w}{10} - \delta_w^{3/2} } \cD ( z, \lambda) .
\end{align}
\eel
\proof We will prove two upper bounds, optimizing at the end of the proof. For simplicity of notation, we define the event $\S := \{ \lambda_1^z \leq \eta_w \}$. By Lemma \ref{lem:1-pt-apriori} we have for $C>0$  large
\begin{align}
\ee\left[ \e^{ \lambda \Phi_N (z,\eta ) } \1_{ \A (z, \eta, A) \cap \S } \right] \leq & \ee\left[ \e^{ \lambda \Phi_N (z,\eta ) } \1_{ \A (z, \eta, A) \cap \A (z, N^{-1}, C )\cap \S } \right] + N^{-10} \notag\\
\lesssim & N^{-\lambda \delta_w} \ee\left[ \e^{ \lambda \Phi_N (z,N^{-1} ) } \1_{ \A (z, \eta, A) \cap \A (z, N^{-1}, C )\cap \S } \right] + N^{-10}
\end{align}
using Lemma \ref{lem:deterministic} in the second line. On the one hand, by Lemma \ref{lem:1-pt-apriori} we have,
\beq
\ee\left[ \e^{ \lambda \Phi_N (z,N^{-1} ) } \1_{ \A (z, \eta, A) \cap \A (z, N^{-1}, C )\cap \S } \right]  \lesssim \e^{ \lambda ( \log N)^{1/2+\delta} } \cD(z, \lambda).
\eeq
Additionally, by Cauchy-Schwarz, Lemma \ref{lem:1-pt-apriori} and Lemma \ref{lem:wegner}, we have
\begin{align}
& \ee\left[ \e^{ \lambda \Phi_N (z,N^{-1} ) } \1_{ \A (z, \eta, A) \cap \A (z, N^{-1}, C )\cap \S } \right] \leq \left( \ee \left[ \e^{ 2 \lambda \Phi_N (z, N^{-1} ) } \1_{ \A (z, N^{-1}, C) } \right] \pp[ \S ] \right)^{1/2}  \notag\\
\lesssim & N^{-\delta_w/2} \e^{ \lambda ( \log N)^{1/2+\delta} } \exp\left( \frac{ \lambda^2 + \1_{ z \in I_r} \1_{ \beta=1} \lambda(\lambda-1)}{4} \log N \right) \lesssim \e^{ \lambda ( \log N)^{1/2+\delta}} \cD (z, \lambda) N^{-\frac{\delta_w}{2} + \frac{\lambda^2}{4}}
\end{align}
The first estimate of \eqref{eqn:1-pt-reg-a1} follows from minimizing over the last two inequalities. The second inequality in \eqref{eqn:1-pt-reg-a1} follows from considering, e.g., whether $\lambda \geq \sqrt{1.5 \delta_w}$  and using $\delta_w < \frac{1}{100}$. \qed

\

The above lemma allows us to prove the following result on regularizing the characteristic polynomial from $\eta = 0$ to $\eta = \eta_m = N^{-1-\delta_m}$, for any sufficiently small $\delta_m >0$.  

\bel \label{lem:1-pt-reg-2}
Let $\eta_m = N^{-1-\delta_m}$ with $0 < \delta_m \leq \frac{c_W}{2} \wedge \frac{1}{100}$. Let $L \geq 10$. For any $ A_1 \geq 10 L^2 + 100$ we have,
\begin{align} \label{eqn:1-pt-reg-2}
 &  \ee\left| \e^{ \lambda \Phi_N (z, 0) } -  \e^{ \lambda \Phi_N (z, \eta_m) } \1_{ \A (z, \eta_m, A_1)}  \right|  \lesssim N^{- \frac{ \delta_m^{3/2}}{4}} \cD (z, \lambda) .
\end{align}
Let $t_1 = N^{\omega_1-1}$ for $\omega_1 \in (0, \frac{1}{100})$, and recall the random variables $\Phi_N (z, \eta, t_1)$ defined in \eqref{eqn:phiNt-def}. 
The same estimate also holds with each of the $\Phi_N (z, \eta)$ replaced by $\Phi_N (z, \eta, t_1)$ for $\eta \in \{ 0, \eta_m\}$ and the event $\A$ adjusted accordingly.
\eel
\proof By Corollary \ref{cor:char-poly-a-priori}, Lemma \ref{lem:1-pt-apriori} and Cauchy-Schwarz we have, 
\beq
 \ee[\e^{ \lambda \Phi_N (z, 0) } \1_{ [\A (z, 0, A_1) \cap \A (z, \eta_m, A_1)]^c } ] \leq N^{-1},
\eeq
as well as 
\beq
 \ee[ \e^{ \lambda \Phi_N (z, \eta_m) } \1_{ [\A (z, 0, A_1) ^c]  \cap \A (z, \eta_m, A_1)} ] \leq N^{-1}.
\eeq
Let now $\delta_w = \frac{\delta_m}{2}$, and let $\S := \{ \lambda_1^z \leq N^{-1-\delta_w} \}$. By Lemma \ref{lem:1-pt-reg} we have for $\eta \in \{ 0, \eta_m \}$,
\beq
\ee[ \e^{ \lambda \Phi_N (z, \eta) } \1_{ \A (z, 0, A_1)  \cap \A (z, \eta_m, A_1) \cap \S} ] \lesssim  N^{-\frac{ \lambda \delta_w}{10} - \delta_w^{3/2}} \cD (z, \lambda ).
\eeq
On the other hand by Lemma \ref{lem:rig-reg} we have with overwhelming probability,
\beq \label{eqn:1-pt-reg-2-a1}
 \1_{ \A (z, 0, A_1) \cap \A (z, \eta_m, A_1) } \1_{\S^c} \left| \e^{ \lambda \Phi_N (z, 0) }- \e^{ \lambda \Phi_N (z, \eta_m) }\right| \lesssim ( \log N)^2 N^{-\delta_m} \e^{ \lambda \Phi_N (z, \eta_m) } \1_{ \A (z, \eta_m, A_1)}.
\eeq
The expectation of the quantity on the RHS of \eqref{eqn:1-pt-reg-2-a1} is bounded above by the RHS of \eqref{eqn:1-pt-reg-2} by Lemma \ref{lem:1-pt-apriori}. The expectation of the LHS of \eqref{eqn:1-pt-reg-2-a1} on the complementary event where the inequality \eqref{eqn:1-pt-reg-2-a1} does not hold may be bounded by Cauchy-Schwarz and Lemma \ref{lem:1-pt-apriori}. 

The estimate \eqref{eqn:1-pt-reg-2} now follows from the above arguments and inequalities. The claim about $\Phi_N (z, \eta, t_1)$ is a simple consequence of \eqref{eqn:t-scaling}. \qed

\subsection{Product representation for one-point function}

In this section we will consider the one-point function of $X_t$ where $X_t$ is as in \eqref{eqn:Xt-def}.  Below we will apply Theorem \ref{thm:local-global} in order to derive a product representation of the Laplace transform of $\Phi_N (z, \eta, t)$. Note that we apply Theorem \ref{thm:local-global} in the case $K=1$ so we just let $ b >0$ be any small constant, and then the parameters $\mfb, \omega_1, \mfa$ as in Theorem \ref{thm:local-global}. Generalizing \eqref{eqn:A-def} slightly we introduce the event,
\beq \label{eqn:At-def}
\A(z, \eta, t, A) := \{ \Phi_N (z, \eta, t) \leq A \log N \} .
\eeq
We also recall the notation $B_z := \sqrt{1 - |z|^2}$. 
\bel \label{lem:1-pt-prod}
Let $b >0$. Let $\mfb >0$ be the corresponding exponent from Theorem \ref{thm:homog}. Choose $\omega_1 < \frac{\mfb}{10}$ and let $\mfa$ be the corresponding exponent from Theorem \ref{thm:homog}. Let $\delta_m \leq \frac{\mfa}{2} \wedge \frac{c_W}{2} \wedge \frac{1}{100}$ and $\eta_m = N^{-1-\delta_m}$.   Choose $\eta_* = ( \log N)^{C_*} / N$ with $C_* \geq 1000$ sufficiently large. Let $L >0$.  Let $A_1 \geq 10 L^2 + 100$. Uniformly in $ \lambda \in [0, L]$ we have, 
\begin{align} \label{eqn:1-pt-prod-c1}
& \ee[ \e^{ \lambda \Phi_N (z, \eta_m, t_1) } \1_{ \A (z, \eta_m, t_1, A_1)} ]  = \calL (B_z^2 t_1 , B_z \eta_*, B_z \eta_m,  \mfb, A_1 ; z , \beta ) \notag\\
\times  & \exp\left( - \frac{\lambda^2}{8} \log (2 B_z^2 t_1 ) + \frac{\lambda(\lambda-2)}{8} \kapb B_z^4  + \1_{\beta=1} \Eone (z, \lambda,  t_1 )\right)  +  \O ( \cD (z, \lambda) (N^{-\mfa/3} + N^{-\omega_1/5} ))
\end{align}
where $\calL$ is as in \eqref{eqn:local-factor-def}. Above, $\kapb$ is the fourth cumulant of the matrix entries of $X_0$ and
\beq \label{eqn:Eone-def}
\Eone (z, \lambda, \nu ) := - \frac{\lambda^2}{8} \log ( |z- \bar{z} |^2 \vee (2 B_z^2 \nu) \vee N^{-1}  )  + \frac{\lambda}{4} \log ( |z-\bar{z} |^2 \vee ( 2 B_z^2 \nu ) \vee N^{-1} ).
\eeq
\eel
\proof As stated above, we apply the coupling of Theorem \ref{thm:local-global} with $K=1$. Let $\{ \tilL_i\}_{i=1}^4$ and $\tilG$  be as in Theorem \ref{thm:local-global}. Let
\beq
\F := \{ \tilL_1 \leq ( \log N)^{5/8} \} \cap \{ | \tilL_3 | \leq 1 \} \cap \{ | \tilL_2 | \leq A_1 \log N \} \cap \{ | \tilL_4| \leq \log N \} \cap  \{ \tilG \leq A_1 \log N \}  
\eeq
and $\tilL := \sum_{i=1}^4 \tilL_i$. 
By Theorem \ref{thm:local-global}\ref{it:lg-Z},\ref{it:lg-La} and Lemma \ref{lem:1-pt-apriori} (applied to $\tilG$) we have $\pp[ \F^c] \lesssim \e^{ - \frac{A_1^2}{2} \log N }$. Therefore, by Cauchy-Schwarz and Lemma \ref{lem:1-pt-apriori}, we have
\beq
\ee[ \e^{ \lambda \Phi_N (z, \eta_m, t_1) } \1_{ \A (z, \eta_m, t_1, A_1)} ] = \ee[ \e^{ \lambda \Phi_N (z, \eta_m, t_1) } \1_{ \A (z, \eta_m, t_1, A_1) \cap \F} ] + \O ( N^{-1} ) .
\eeq
Using also \eqref{eqn:lg-L2-mfg}, by similar reasoning,
\beq
\ee\left[ \e^{\lambda \tilL} \e^{ \lambda \tilG} \1_{ \F}  \right] = \ee\left[ \e^{\lambda \tilL} \e^{ \lambda \tilG} \1_{\A (z, \eta_m, t_1, A_1) \cap \F}  \right] + \O (N^{-1} ).
\eeq
By the estimate \eqref{eqn:lg-main-est} of Theorem \ref{thm:local-global} (and similar arguments to the above to bound the complementary event where the estimate \eqref{eqn:lg-main-est} does not hold), we have
\begin{align}
&   \ee\left[ \left|  \e^{ \lambda \Phi_N (z, \eta_m, t_1) } - \e^{\lambda \tilL} \e^{ \lambda \tilG}\right|  \1_{\A (z, \eta_m, t_1, A_1) \cap \F}  \right]   
\lesssim   ( N^{-\mfa/2} + N^{-\omega_1/2} ) \e^{ ( \log N)^{3/4}} \cD (z, \lambda) 
\end{align}
Consequently we so far have that
\begin{align} \label{eqn:1-pt-prod-a1}
& \left| \ee[ \e^{ \lambda \Phi_N (z, \eta, t_1) } \1_{ \A (z, \eta_m, t_1 , A_1 } ] - \ee[ \e^{ \lambda ( \tilL + \tilG) } \1_\F] \right| 
\lesssim  ( N^{-\mfa/2} + N^{-\omega_1/2} ) \e^{ ( \log N)^{3/4}} \cD (z, \lambda )  .
\end{align}
On the other hand,
\beq
\ee\left[ \e^{\lambda \tilL} \e^{ \lambda \tilG} \1_{ \F}  \right] =  \calL (B_z^2 t_1 , B_z \eta_*, B_z \eta_m,  \mfb ,  A_1 ; z, \beta )   \times \ee\left[ \e^{ \lambda \tilG} \1_{\{ \tilG \leq A_1 \log N \}} \right] 
\eeq
By \eqref{eqn:meso-1pt-mgf}, \eqref{eqn:covar-order-0} and \eqref{eq:exp}, for $\beta=2$, we have
\begin{align} \label{eqn:1-pt-prod-a2}
    \ee\left[ \e^{ \lambda \tilG} \1_{\{ \tilG \leq A_1 \log N \}} \right] = (1 + \O (N^{\omega_1/100-\omega_1/4} ) ) \exp\left( - \frac{\lambda^2}{8} \log (2 B_z^2 t_1 ) + \frac{\lambda^2}{8} \kappa_4 B_z^4 - \frac{ \lambda \kappa_4}{4} B_z^4 \right)
\end{align}
When $\beta=1$ we get a similar equality with also the factor $\Eone$ on the RHS; we omit the details.
In order to bound the multiplicative error in the above we note that \eqref{eqn:1-pt-prod-a1} and Lemma \ref{lem:1-pt-apriori} imply that 
\beq
\ee\left[ \e^{\lambda \tilL} \e^{ \lambda \tilG} \1_{ \F}  \right] \lesssim \e^{  ( \log N)^{3/4}} \cD (z, \lambda) .
\eeq
With this the claim now follows from \eqref{eqn:1-pt-prod-a1} and \eqref{eqn:1-pt-prod-a2}. 
\qed

\

The following is an immediate corollary of Lemmas \ref{lem:1-pt-reg-2} and \ref{lem:1-pt-prod}. 
\bec \label{cor:product-1pt}
In the set-up of Lemma \ref{lem:1-pt-prod} we have with $\alpha = \min \{ \frac{\mfa}{3} , \frac{\omega_1}{5}, \frac{\delta_m^{3/2}}{4}\}$ ,
    \begin{align} \label{eqn:1-pt-prod-cor}
& \ee[ \e^{ \lambda \Phi_N (z, 0, t_1) } ]   =  \bigg\{ \calL (B_z^2 t_1 , B_z \eta_*, B_z \eta_m,  \mfb, A_1 ; z, \beta ) \notag\\
 \times & \exp\left( - \frac{\lambda^2}{8} \log (2 B_z^2 t_1 ) + \frac{\lambda^2}{8} \kappa_4 B_z^4 - \frac{ \lambda \kappa_4}{4} B_z^4 + \1_{\beta=1} \Eone (z, \lambda,  t_1) \right) \bigg\}  + \O (\cD (z, \lambda) N^{-\alpha} )
\end{align}

\eec

\subsubsection{Computation of local component and one-point result under dynamics}

We are now in position to compute the local factor $\calL$ defined in \eqref{eqn:local-factor-def}. 

\bet \label{thm:local-factor}
In the set-up of the previous section,
\beq \label{eqn:local-factor-thm-a1}
\calL ( t_1, \eta_*, \eta_m,  \mfb, A_1 ; z, \beta ) = \frac{ \exp \left( \frac{\lambda^2 + \1_{ z \in I_r} \1_{ \beta=1} \lambda(\lambda-2)}{8} \log (2 N t_1 ) \right)}{\GB ( \lambda; z , \beta)} (1 + \O (N^{-\alpha} ) ),
\eeq
where $\alpha = \min \{ \frac{\mfa}{3}, \frac{ \omega_1}{5} , \frac{ \delta_m^{3/2}}{4} \}$ and $\GB$ is as in \eqref{eqn:GB-def}. 
\eet
\begin{remark} \emph{In the above we take $t_1 = N^{\omega_1-1}$, $\eta_* = ( \log N)^{C_*} /N$ and $\eta_m = N^{-1-\delta_m}$. But it is clear from the proof that the above estimate also holds for any $t_1 \asymp N^{\omega_1-1}$, $\eta_* \asymp ( \log N)^{C_*} /N$ and $\eta_m \asymp N^{-1-\delta_m}$.}
\end{remark}

\

\proof
As noted in Section \ref{sec:rescale} the local factor $\calL$ can be computed from the case that $z=0$ and the initial data $X_0$ is either a real or complex Ginibre matrix.  Then $X_t$ is simply a rescaled Ginibre ensemble and so  by \eqref{eqn:t-scaling} we can apply Proposition \ref{prop:ginibre} to compute the quantity on the LHS of \eqref{eqn:1-pt-prod-cor}, up to an $\O (N^{-1})$ multiplicative error.  Proposition \ref{prop:ginibre} shows that in the Ginibre case the LHS of \eqref{eqn:1-pt-prod-cor} is $\asymp \cD (z, \lambda)$. Therefore, the additive error on the RHS of \eqref{eqn:1-pt-prod-cor} becomes a multiplicative error when moved to the LHS. The claim \eqref{eqn:local-factor-thm-a1} now follows from moving the factor on the second line of \eqref{eqn:1-pt-prod-cor} over to the LHS. \qed

\

We now compute the Laplace transform of the regularized characteristic polynomial for $\eta = \eta_m$, for Gaussian divisible i.i.d. matrices. Of course, the following together with Lemma \ref{lem:1-pt-reg-2} implies a result for the characteristic polynomial itself (i.e., without the regularization) but for organizational reasons we will deduce this later. In fact, we need the following intermediate result as the quantity $\Phi_N (z, \eta_m)$ is better suited for a \emph{Green's Function Comparison Theorem (GFCT)} argument than $\Phi_N (z, 0)$ itself when we remove the additive Gaussian component introduced by DBM \eqref{eqn:Xt-def}. 

\bec \label{cor:gde-1pt}
Let the parameters $t_1, \eta_m, \mfa, \mfb, \eta_*$ be as in Lemma \ref{lem:1-pt-prod}. Let $X$ be an i.i.d. matrix of the form $X = \sqrt{1-t_1} Y + \sqrt{t_1}G$ where $Y$ is an i.i.d. matrix and $G$ is an independent Ginibre matrix of the same symmetry class. Let $\alpha$ be as in Theorem \ref{thm:local-factor}. Let $L \geq 1$ and $A \geq 10 L^2 + 100$. Then, uniformly in $\lambda \in [0, L]$ and $z \in \ddhb_r$ we have,
\begin{align}
& \ee\left[ \e^{\lambda \Phi_N (z, \eta_m ; X)} \1_{\{  \Phi_N (z, \eta_m ; X) \leq A \log N \} } \right] \notag\\
=& (1 + \O ( N^{-\alpha} ) ) \frac{  N^{\lambda^2/8}}{\GB (\lambda; z , \beta)} \e^{ \frac{\lambda^2}{8} \kappa_4 (1 - |z|^2)^2 - \frac{\lambda \kappa_4}{4} (1- |z|^2)^2 + \1_{\beta=1} \Eone (z, \lambda, N^{-1})}.
\end{align}
\eec
\proof This follows immediately from Corollary \ref{cor:product-1pt}, Theorem \ref{thm:local-factor} (and the remark immediately thereafter), as well as \eqref{eqn:t-scaling}. \qed

\section{Preliminary analysis of the K-point function}

\label{sec:kpt}

Throughout this section we will work in the dynamical set-up of Section \ref{sec:dynamics}, and therefore consider $X_t$ as in \eqref{eqn:Xt-def} and $\Phi_N (z, \eta, t)$ as in Section \ref{sec:def-dynamics}. We also recall the definition of $\A (z, \eta,  t, A)$ from \eqref{eqn:At-def}. 

We therefore fix a $b>0$ , $0 <r < 1$, $K \geq 2$, let $\omega_1$, $\mfb$, $q_1$, $C_*\geq 1000$, and $\mfa$ be as in Theorem \ref{thm:local-global}. We may also assume $\omega_1  < b$. We let $\bz = (z_1, \dots, z_K) \in \cc^K$ such that $z_i \in \ddhb_r$ and $|z_i -z_j| \geq N^{b-1/2}$ for all $i \neq j$. Note that the assumption $\omega_1 < b$ implies that
$
|z_i - z_j|^2 \gtrsim N^{b} t_1
$
for all $i \neq j$. 
We fix an $L \geq 1$ and let $\blam = ( \lambda_1, \dots, \lambda_K ) \in \rr^K$ satisfy $\lambda_i \in [0, L]$. We need to generalize the control parameter $\cD (z, \lambda)$ defined in \eqref{eqn:control-scalar} to vector arguments. We therefore define,
\beq \label{eqn:D-def}
\cD ( \bz, \blam ) := \left( \prod_{i=1}^K \cD (z_i, \lambda_i ) \right) \exp\left(   \sum_{i \neq j } \frac{ \lambda_i \lambda_j}{8}\left[  \log ( |z_i -z_j|^{-2} ) + \1_{\beta=1} \log ( |z_i - \bar{z_j} |^{-2} )\right]  \right)
\eeq
and recall the notation $B_i = \sqrt{1- |z_i|^2}$. 
The quantity \eqref{eqn:D-def} is what we expect for the size of the $K$-point function and will appear in various upper bounds and estimates. Note that since $\lambda_i \geq 0$ we see that $\cD (\bz , \blam) \gtrsim N^{- \frac{K}{16}}$. Note also that the additional term in the $\beta=1$ case is constant order unless both $z_i$ and $z_j$ are real (recall $\dd^{1}_r \subseteq \{ z : \Im[z] > r \}$).

\subsection{Regularization}

We first require the following lemma which gives an a priori bound on the size of the $K$-point function of almost the correct order of magnitude. 

\bel \label{lem:k-pt-apriori} In the coupling of Theorem \ref{thm:local-global}, for any $\eta \in [0, \eta_*]$, we have for any $\delta >0$,
\beq \label{eqn:k-pt-upper}
\Phi_N (z_i, \eta, t_1) \leq ( \log N)^{1/2+\delta} + \tilZ^{(i)} + \tilG^{(i)} - \1_{\beta=1} \1_{z_i \in I_r} \frac{\omega_1}{4} \log N,
\eeq
with overwhelming probability. Moreover, for any $A_1 \geq 1$, and $\eta_i \in [0, \eta_*]$, we have
\beq \label{eqn:kpt-apriori-b1} 
\ee\left[ \e^{ \sum_{i=1}^K \lambda_i \Phi_N (z_i, \eta_i, t_1 )} \1_{ \bigcap_{i=1}^K \A ( z_i, \eta_i, t_1, A_1 )  }\right] \lesssim \e^{ \| \blam \|_1 ( \log N)^{1/2+\delta}} \cD ( \bz, \blam ).
\eeq
\eel
\proof The estimate \eqref{eqn:k-pt-upper} follows from \eqref{eqn:lg-main-bd} and Lemma \ref{lem:a-priori-1}. 

To prove the second estimate, we first choose $A_2 = 10 (KL A_1 + 1)^2$ to see from Lemma \ref{lem:1-pt-apriori} that
\begin{align}
 & \ee\left[ \e^{ \sum_{i=1}^K \lambda_i \Phi_N (z_i, \eta_i, t_i )} \1_{ \bigcap_{i=1}^K \A ( z_i, \eta_i,  t_1, A_1 )  }\right] \notag\\
= & \ee\left[ \e^{ \sum_{i=1}^K \lambda_i \Phi_N (z_i, \eta_i, t_i )} \1_{ \bigcap_{i=1}^K \A ( z_i, \eta_i, t_1, A_1 ) \cap \{ \tilG^{(i)} \leq A_2 \log N \}   }\right]+ \O (N^{-K} ). 
\end{align}
Applying \eqref{eqn:k-pt-upper}  we see that,
\begin{align}
    & \ee\left[ \e^{ \sum_{i=1}^K \lambda_i \Phi_N (z_i, \eta_i, t_1 )} \1_{ \bigcap_{i=1}^K \A ( z_i, \eta_i, A_1, t_1 ) \cap \{ \tilG^{(i)} \leq A_2 \log N \}   }\right] 
    \lesssim  \bigg\{ \e^{\| \blam \|_1 ( \log N)^{1/2+\delta} } \left(  \prod_{i=1}^K \ee[ \e^{ \lambda_i Z^{(i)} } ] \right) \notag\\
    & \times \ee\left[ \e^{ \sum_{i=1}^K \lambda_i \tilG^{(i)}} \1_{\bigcap_{i=1}^K  \{ \tilG^{(i)} \leq A_2 \log N  \} } \right]\e^{ - \1_{\beta=1} \sum_{i=1}^K \1_{z_i \in I_r} \lambda_i \frac{\omega_1}{4} \log N } \bigg\} + N^{-K}
\end{align}
where we used the independence of the $\tilZ^{(i)}$ from the $\tilG^{(i)}$ from Theorem \ref{thm:local-global}\ref{it:Fi}. 
By Theorem \ref{thm:local-global}\ref{it:lg-Z}, 
\beq \label{eqn:k-pt-apriori-Z}
 \prod_{i=1}^K \ee[ \e^{ \lambda_i Z^{(i)} } ] \lesssim \exp \left( \sum_{i=1}^K \frac{ \lambda_i^2(1+\1_{\beta=1} \1_{z_i \in I_r})}{8}  \log (N t_1 ) \right),
\eeq
and by Theorem \ref{thm:local-global}\ref{it:lg-G},  Proposition \ref{pro:jointlapdnew}, and Proposition \ref{prop:parameter-properties},  we have
\begin{align} \label{eqn:k-pt-apriori-G}
 & \ee\left[ \e^{ \sum_{i=1}^K \lambda_i \tilG^{(i)}} \1_{\bigcap_{i=1}^K  \{ \tilG^{(i)} \leq A_2 \log N  \} } \right] \lesssim  \exp \left( \sum_{i=1}^K \frac{ \lambda_i^2+ \1_{\beta=1} \1_{z_i \in I_r}\lambda_i(\lambda_i-2)}{8} \log (t_1^{-1} ) \right) \notag\\
 &\times \exp \left( \sum_{i \neq j } \frac{ \lambda_i \lambda_j}{8} \left[\log ( |z_i -z_j |^{-2} ) + \1_{\beta=1} \log ( |z_i - \bar{z}_j |^{-2} ) \right] \right). 
\end{align}
Here we used the fact that $|z_i-z_j|^2 \geq t_1$ by assumption to estimate the covariance terms. The claim now follows. \qed

\

We now wish to control the $K$-point function on the event that one of the $\lambda_1^{(z_i)} (t_1)$ are small. The following is the analog of Lemma \ref{lem:1-pt-reg}. However, here we rely on the dynamical set-up (in contrast to Lemma \ref{lem:1-pt-reg}) to gain from the approximate independence of the $\lambda_1^{(z_i)} (t_1)$. Essentially, the argument using Cauchy-Schwarz in the proof of Lemma \ref{lem:1-pt-reg} to handle the small $\lambda$ regime is too lossy if one of the $\lambda_i$ is small while the other $\lambda_j$'s are large. 

Fix now $\eta_w = N^{-1-\delta_w}$ with $0 < \delta_w \leq \frac{c_W}{2} \wedge \frac{1}{100} \wedge \mfa $. Let
\beq
\S_i := \{ \lambda_1^{(z_i)}   (t_1) \leq \eta_w \}, \qquad \S := \bigcup_{i=1}^K \S_i.
\eeq

\bel \label{lem:kpt-reg} 
Uniformly in $\eta_i \in [0, \eta_w]$, for any $A_1 \geq 1$, we have
\begin{align}
 \ee\left[ \e^{ \sum_{i=1}^K \lambda_i \Phi_N (z_i, \eta_i, t_1 ) } \1_{ \S \cap \bigcap_{i=1}^K \A (z_i, \eta_i,  t_1, A_1 ) } \right] \lesssim N^{-\alpha_w} \cD( \bz, \blam ),
\end{align}
where $\alpha_w = \min_i ( \frac{ \lambda_i \delta_w}{10} + \delta_w^{3/2})$.
\eel
\proof By a union bound it suffices to prove the bound with $\S$ replaced by $\S_i$. For notational simplicity we do it for $\S_1.$ By choosing $A_2 = 10 (K L A_1 +1 )^2$ we have
\begin{align}
 &\ee\left[ \e^{ \sum_{i=1}^K \lambda_i \Phi_N (z_i, \eta_i, t_1 ) } \1_{ \S_1 \cap \bigcap_{i=1}^K \A (z_i, \eta_i, A_1, t_1 ) } \right]  \notag\\
 = &  \ee\left[ \e^{ \sum_{i=1}^K \lambda_i \Phi_N (z_i, \eta_i, t_1 ) } \1_{ \S_1 \cap \bigcap_{i=1}^K \A (z_i, \eta_i, A_1, t_1 ) \cap \{ \tilG^{(i)} \leq A_2 \log N \} } \right] + \O (N^{-K} ) .
\end{align}
By Lemma \ref{lem:deterministic} and \eqref{eqn:k-pt-upper}, on the event $\S_1$, we have
\beq
\Phi_N (z_1, \eta_1, t_1) \leq - \delta_w \log N + ( \log N)^{3/4} + \tilZ^{(1)} + \tilG^{(1)} - \1_{\beta=1} \1_{z_1 \in I_r} \frac{\omega_1}{4} \log N
\eeq
with overwhelming probability. Therefore, also applying \eqref{eqn:k-pt-upper} to the terms with $i\neq 1$, we obtain
\begin{align}
    & \ee\left[ \e^{ \sum_{i=1}^K \lambda_i \Phi_N (z_i, \eta_i, t_1 ) } \1_{ \S_1 \cap \bigcap_{i=1}^K \A (z_i, \eta_i, A_1, t_1 ) \cap \{ \tilG^{(i)} \leq A_2 \log N \} } \right]  \notag\\
    \lesssim & N^{-\delta_w \lambda_1 } \e^{ \| \blam \|_1 ( \log N)^{3/4}} \ee\left[\e^{ \sum_{i=1}^K \lambda_i (\tilZ^{(i)} + \tilG^{(i)} ) }   \1_{ \S_1 \cap \bigcap_{i=1}^K  \{ \tilG^{(i)} \leq A_2 \log N \} }\right] \e^{ - \1_{\beta=1} \sum_{i=1}^K \1_{ z_i \in I_r} \frac{\omega_1}{4} \log N } .
\end{align}
Let now $\tilde{\S}_1 := \{ \tilmu_1^{(1)} (B_{1}^2 t_1) \leq 2 \eta_w \}$. Since $\mfa > \delta_w$ we have by Theorem \ref{thm:local-global}\ref{it:lg-homog} that $\1_{\S_1} \leq \1_{ \tilde{\S}_1}$ with overwhelming probability. Therefore,
\begin{align}
 & \ee\left[\e^{ \sum_{i=1}^K \lambda_i (\tilZ^{(i)} + \tilG^{(i)} )}   \1_{ \S_1 \cap \bigcap_{i=1}^K  \{ \tilG^{(i)} \leq A_2 \log N \} }\right]\notag \\
\leq &  \ee[ \e^{ \lambda_1 \tilZ^{(1)} } \1_{ \tilde{\S}_1} ] 
\left( \prod_{i=2}^K \ee[ \e^{ \lambda_i \tilZ^{(i)}}] \right) \ee\left[ \e^{ \sum_{i=1}^K \lambda_i \tilG^{(i)} } \1_{ \bigcap_{i=1}^K  \{ \tilG^{(i)} \leq A_2 \log N \} } \right]  + N^{-10KL} \label{eqn:k-pt-reg-a1},
\end{align}
where we used the independence properties in Theorem \ref{thm:local-global}\ref{it:Fi}. 

On the one hand, we see by \eqref{eqn:k-pt-apriori-Z} and \eqref{eqn:k-pt-apriori-G} that the RHS of \eqref{eqn:k-pt-reg-a1} is $\O ( \cD ( \bz, \blam ))$. On the other hand, we have by Cauchy-Schwarz and Lemma \ref{lem:wegner} that
\beq
\ee[ \e^{ \lambda_1 \tilZ^{(1)} } \1_{ \tilde{\S}_1} ]  \lesssim N^{-\frac{\delta_w}{2}} \e^{ \frac{ \lambda_1^2(1+ \1_{\beta=1} \1_{z_1 \in I_r})}{4} \log (N t_1) } \leq N^{-\frac{\delta_w}{2} + \frac{\lambda_1^2}{4}} \e^{ \frac{\lambda_1^2(1+ \1_{ \beta=1} \1_{ z \in I_r) }}{8} \log (N t_1) }. 
\eeq
We bound the remaining terms in the last line of \eqref{eqn:k-pt-reg-a1} as before. We conclude similarly to the end of the proof of Lemma \ref{lem:1-pt-reg}. \qed

\

\subsection{Computation of (regularized) K-point function}

We now derive a product representation for the $K$-point function of the regularized characteristic polynomial of  $X_t$. Let $\calL_i := \calL (B_{i}^2 t_1, B_{i} \eta_* , B_i \eta_m,  \mfb, A_1 ; z_i, \beta )$, recalling  \eqref{eqn:local-factor-def}.

\bep \label{prop:k-pt-prod} Let $\eta_m = N^{-1-\delta_m}$ with $\delta_m \leq \frac{\mfa}{2} \wedge \frac{c_W}{2} \wedge \frac{1}{100}$. For $A_1 \geq 100 K^2 L^2 + 100$, we have
\begin{align} \label{eqn:k-pt-prod}
& \bigg| \ee[ \e^{ \sum_{i=1}^K \lambda_i \Phi_N (z_i, \eta_m, t_1) } \1_{ \bigcap_{i=1}^K \A ( z_i, \eta_i, t_1, A_1 ) } ] - \left( \prod_{i=1}^K \calL_i \right) \e^{ \sum_{i=1}^K  \frac{\lambda_i^2}{8} \left( \kappa_4 B_{z_i}^4- \log (2 B_{z_i}^2 t_1 ) \right) - \frac{\lambda_i \kappa_4}{4} B_{z_i}^4 } \notag\\
&\times  \e^{  \sum_{i \neq j } \frac{ \lambda_i \lambda_j}{8} \left( \log ( |z_i -z_j |^{-2} ) + \kappa_4 B_{z_i}^2 B_{z_j}^2 \right)  + \1_{\beta=1} \Etwo ( \bz , \blam, t_1 )  }  \bigg| \lesssim  \cD( \bz , \blam ) ( N^{-\omega_1/5} + N^{-\mfa/3} ),
\end{align}
where
\beq \label{eqn:Etwo-def}
\Etwo ( \bz , \blam, t_1 ) :=- \sum_{i,j} \frac{ \lambda_i \lambda_j}{8}  \log( |z_i - \bar{z}_j |^2 \vee 2 B_i^2 t_1 \vee N^{-1}) + \sum_i \frac{ \lambda_i}{4} \log ( |z_i -\bar{z}_i |^2 \vee 2 B_i^2 t_1 \vee N^{-1} )
\eeq
\eep
\proof Let $\F = \bigcap_{i=1}^K \left(\F_i \cap  \{ \tilG^{(i)} \leq A_1 \log N \} \right)$ where
\begin{align}
\F_i &:=  \{ \tilL^{(i)}_1 \leq ( \log N)^{5/8} \} \cap \{ | \tilL_3^{(i)} | \leq 1 \} \cap \{ | \tilL_2^{(i)} | \leq A_1 \log N \}  \cap \{ | \tilL_4^{(i)} | \leq \log N \}
\end{align}
and define $\tilL^{(i)} = \sum_{a=1}^4 \tilL^{(i)}_a$. 
By Lemma \ref{lem:1-pt-apriori} and Theorem \ref{thm:local-global} we have $\pp[ \F^c] \lesssim N^{ -\frac{A_1^2}{2} }$.
By Lemma \ref{lem:k-pt-apriori} we have,
\beq
\ee[ \e^{ 2 \sum_{i=1}^K \lambda_i \Phi_N (z_i, \eta_m, t_1) } \1_{ \bigcap_{i=1}^K \A ( z_i, \eta_i, t_1, A_1 ) } ]^{1/2} \lesssim N^{ ( 2 KL)^2} .
\eeq
Therefore, 
\beq
\ee[ \e^{ \sum_{i=1}^K \lambda_i \Phi_N (z_i, \eta_m, t_1) } \1_{ \bigcap_{i=1}^K \A ( z_i, \eta_i, t_1,  A_1 ) } ] = \ee[ \e^{ \sum_{i=1}^K \lambda_i \Phi_N (z_i, \eta_m, t_1) } \1_{ \F \cap \bigcap_{i=1}^K \A ( z_i, \eta_i, t_1, A_1 ) } ] + \O ( N^{-K} ).
\eeq
By similar reasoning, we obtain
\beq
\ee\left[ \e^{ \sum_{i=1}^K \lambda_i (  \tilL^{(i)} + \tilG^{(i)} ) } \1_\F \right]  = \ee\left[ \e^{ \sum_{i=1}^K \lambda_i (\tilL^{(i)} + \tilG^{(i)} ) } \1_{ \F \cap \bigcap_{i=1}^K \A ( z_i, \eta_i, t_1, A_1 ) } \right]  + \O ( N^{-K} ). 
\eeq
By \eqref{eqn:lg-main-est} of Theorem \ref{thm:local-global} and Lemma \ref{lem:k-pt-apriori} we see that
\begin{align}
&   \ee\left[ \left| \e^{ \sum_{i=1}^K \lambda_i \Phi_N (z_i, \eta_m, t_1) } -\e^{ \sum_{i=1}^K \lambda_i (  \tilL^{(i)} + \tilG^{(i)} ) } \right| \1_{ \F \cap \bigcap_{i=1}^K \A ( z_i, \eta_i, A_1, t_1 ) } \right]  \notag\\
\lesssim & \e^{ ( \log N)^{3/4}} \left( N^{-\mfa/2} + N^{-\omega_1/2} \right) \cD ( \bz , \blam ) .
\end{align}
On the other hand,
\begin{align}
    \ee\left[ \e^{ \sum_{i=1}^K \lambda_i (  \tilL^{(i)} + \tilG^{(i)} ) } \1_\F \right]  &= \left( \prod_{i=1}^K \calL_i \right) \ee\left[ \e^{ \sum_{i=1}^K \lambda_i \Phi_N (z_i, t_1 B_i, 0) } \1_{ \bigcap_{i=1}^K \{ \Phi_N (z_i, t_1 B_i ) \leq A_1 \log N \}  } \right]  \label{eqn:k-pt-a1}.
\end{align}
By Lemma \ref{lem:1-pt-apriori} and our choice of $A_1$ we have for any $C_1 \geq A_1$ that
\begin{align}
&\ee\left[ \e^{ \sum_{i=1}^K \lambda_i \Phi_N (z_i, t_1 B_i, 0) } \1_{ \bigcap_{i=1}^K \{ \Phi_N (z_i, t_1 B_i ) \leq A_1 \log N \}  } \right] \notag\\
=& \ee\left[ \e^{ \sum_{i=1}^K \lambda_i \Phi_N (z_i, t_1 B_i, 0) } \1_{ \bigcap_{i=1}^K \{ \Phi_N (z_i, t_1 B_i ) \leq C_1 \log N \}  } \right]  + \O ( N^{-K} )
\end{align}
Now, when $\beta=2$, by Proposition \ref{pro:jointlapdnew},  \eqref{eqn:covar-order-0} and \eqref{eq:exp}, we have (taking $C_1$ sufficiently large)
\begin{align} \label{eqn:k-pt-prod-a1}
    & \ee\left[ \e^{ \sum_{i=1}^K \lambda_i \Phi_N (z_i, t_1 B_i, 0) } \1_{ \bigcap_{i=1}^K \{ \Phi_N (z_i, t_1 B_i ) \leq C_1 \log N \}  } \right]  (1 + \O ( N^{-\omega_1/4+\omega_1/100} ) ) \notag\\
    = & \exp\left( \sum_{i=1}^K  \frac{\lambda_i^2}{8} \left( \kappa_4 B_{i}^4- \log (2 B_{i}^2 t_1 ) \right) - \frac{\lambda_i \kappa_4}{4} B_{i}^4  + \sum_{i \neq j } \frac{ \lambda_i \lambda_j}{8} \left( \log ( |z_i -z_j |^{-2} ) + \kappa_4 B_{i}^2 B_{j}^2 \right)  \right) 
\end{align}
and a similar estimate with $\Etwo ( \bz , \blam, t_1)$ when $\beta=1$. Here in particular we used that \eqref{eqn:covar-order-0} implies that $\vV (z_i, t_1 B_i, z_j, t_1 B_j ) = - \frac{1}{4} \log ( |z_i -z_j|^2 ) + \O ( N^{-\omega_1} )$, for $i \neq j$, due to our assumption that $\omega_1 < b$, as well as that $\vV (z_i, t_1 B_i, z_i, t_1 B_i ) = - \frac{1}{4} \log ( 2 B_i^2 t_1) + \O ( N^{-1/2})$. In the case $\beta=1$ we also use \eqref{eqn:meso-logE}. 

The above estimates and Lemma \ref{lem:k-pt-apriori} imply that the  RHS  of \eqref{eqn:k-pt-a1} is $\O ( \e^{ \| \blam \|_1 ( \log N)^{3/4}} \cD ( \bz, \blam ))$, and so the multiplicative error in \eqref{eqn:k-pt-prod-a1} can be absorbed into the RHS of \eqref{eqn:k-pt-prod}. We therefore conclude the proof. \qed

\

\bec \label{cor:kpt-reg-gde}
Let $\alpha = \min \{ \frac{ \mfa}{3}, \frac{\omega_1}{5}, \frac{ ( \delta_m)^{3/2}}{4} \}$. Let $X$ be an i.i.d. matrix of the form $X = \sqrt{1-t_1} Y + \sqrt{t_1} G$ where $Y$ is an i.i.d. matrix and $G$ is an independent Ginibre matrix of the same symmetry class. We have,
\begin{align} \label{eqn:kpt-reg-gde}
    & \ee\left[ \e^{ \sum_{i=1}^K \lambda_i \Phi_N (z_i, \eta_m ; X)} \1_{ \bigcap_{i=1}^K \{ \Phi_N (z_i, \eta_m ; X) \leq A_1 \log N \} }  \right] 
    =   \prod_{i=1}^K \left( \frac{ N^{\lambda_i^2/8}}{ \GB ( \lambda_i; z_i , \beta)} \e^{ \frac{ \kappa_4 (1-|z_i|^2)^2 \lambda_i(\lambda_i-2)}{8}} \right) \notag\\
    \times & \e^{ \sum_{i \neq j} \frac{ \lambda_i \lambda_j}{8} \left( \kappa_4 (1 - |z_i|^2)(1-|z_j|^2) - \log ( |z_i -z_j|^2 ) + \1_{\beta=1} \Etwo ( \bz , \blam, 0 ) \right) } \left( 1 + \O ( N^{-\alpha} ) \right)
\end{align}
\eec
\proof From Proposition \ref{prop:k-pt-prod}, Theorem \ref{thm:local-factor} and \eqref{eqn:t-scaling} we obtain \eqref{eqn:kpt-reg-gde} but with all of the $z_i$ on the RHS replaced by $z_i c_*(t_1) = z_i\sqrt{1+t_1}$. Changing $z_i c_* (t_1)$ to $z_i$ everywhere induces a  negligible multiplicative error of the form $1 + \O ( N^{-1/4})$. \qed

\section{Submicroscopic regularization for the GMC} \label{sec:gmc-reg}

In this section we will develop a regularization for the GMC for general i.i.d. matrices, under an assumption about the asymptotics of the one-point function. In Section \ref{sec:claeys} we prove the GMC convergence for the regularized characteristic polynomial (at scale $\eta = \eta_m = N^{-1-\delta_m}$), for matrices with a Gaussian component. These matrices are called Gaussian divisible (or in short GDE) matrices. This latter regularized quantity is already well-suited for an application of a GFCT argument to remove the Gaussian component. The results in this section will then later be used to remove the regularization, \emph{after} we have proven the one-point asymptotics for general i.i.d. matrices (and so verify assumption \eqref{eqn:GMC-reg-assump} below).

Let $h_1$ be the cut-off function in \eqref{eqn:cut-off-def}. For any $A \geq 1$ and any i.i.d. matrix $X$ define,
\beq \label{eqn:yY-def}
Z_A (z, \eta ; X) = Z_A (z, \eta ) := (A \log N) h_1 \left( \frac{ \Phi_N (z, \eta ; X)}{A \log N } \right) .
\eeq

\bep \label{prop:gmc-reg-new}
Let $f$ be a bounded measureable function supported in $\ddb_r$. Fix $L \geq 100$, let $A \geq 100 L^2$, and let $\eta_m = N^{-1-\delta_m}$, with $\delta_m \leq \frac{c_W}{2} \wedge \frac{1}{100}$. Suppose that 
\beq \label{eqn:GMC-reg-assump}
\ee[ \e^{ \lambda \Phi_N (z, \eta_m ) } \1_{ \{ \Phi_N (z, \eta_m ) \leq A \log N \} } ] = (1 + \O (N^{-\alpha_1} ) )\frac{N^{\lambda^2/8}}{ \GB ( \lambda ; z, \beta ) } \e^{ \frac{ \lambda(\lambda-2)}{8} \kappa_4 (1 - |z|^2 )^2 + \1_{\beta=1} \Eone (z, \lambda, 0) } 
\eeq
uniformly in $z \in \ddb_r$ and $\lambda \in [0, L]$ for some $\alpha_1 >0$ (recalling the definition \eqref{eqn:Eone-def}).

Then,
\beq \label{eqn:GMC-reg-1}
\ee\left| \int f(z) \frac{ \e^{ \lambda \Phi_N (z, 0 )}}{\ee[ \e^{ \lambda \Phi_N (z, 0 )}]  } \dz - \int f(z) \frac{ \e^{ \lambda Z_A  (z, \eta_m) }}{\ee[ \e^{ \lambda Z_A (z, \eta_m) }]  }\dz \right| \lesssim N^{-\frac{ \delta_m^{3/2}}{4}},
\eeq
as well as,
\beq \label{eqn:GMC-reg-2}
\ee\left|  \int f(z) \frac{ \e^{ \lambda Z_A (z, \eta_m) }}{\ee[ \e^{ \lambda Z_A (z, \eta_m) }]  } \dz - \int f(z) \frac{ \e^{ \lambda Z_A (z, \eta_m) }}{N^{\lambda^2/8} \kB (z)  } \dz \right| \lesssim N^{-\alpha_1},
\eeq
where $\kB (z) := \GB (\lambda ; 0, 2) \e^{- \frac{ \lambda^2}{8} \kappa_4 (1 - |z|^2 )^2 + \frac{ \lambda \kappa_4}{4} (1 - |z|^2 )^2 - \1_{\beta=1} \Eone (z, \lambda, 0) }$.
\eep
\proof Due to Lemma \ref{lem:1-pt-apriori} and our choice of $A_1$ it is easy to see that
\beq \label{eqn:GMC-reg-b3}
\ee\left| \e^{ \lambda \Phi_N (z, \eta_m) } \1_{ \{ \Phi_N (z, \eta_m) \leq A_1 \log N  \}} - \e^{ \lambda Z_A (z, \eta_m) } \right| \leq N^{-10}. 
\eeq
From this and Lemma \ref{lem:1-pt-reg-2} we easily obtain
\beq \label{eqn:GMC-reg-b1}
\ee\left|  \e^{ \lambda \Phi_N (z, 0) }-  \e^{ \lambda Z_A (z, \eta_m)}  \right| \lesssim N^{- \frac{  \delta_m^{3/2}}{4}} N^{\lambda^2/8}.
\eeq
As a consequence of \eqref{eqn:GMC-reg-b1} and \eqref{eqn:GMC-reg-assump} we see that 
\beq \label{eqn:GMC-reg-b2}
\ee[ \e^{ \lambda \Phi_N (z, 0) } ] \asymp \ee[ \e^{ \lambda Z_A (z, \eta_m) } ] \asymp N^{\lambda^2/8}. 
\eeq
Then, the estimate \eqref{eqn:GMC-reg-1} follows in a straightforward way from \eqref{eqn:GMC-reg-b1} and \eqref{eqn:GMC-reg-b2}. The estimate \eqref{eqn:GMC-reg-2} follows from \eqref{eqn:GMC-reg-assump}, \eqref{eqn:GMC-reg-b3}, and \eqref{eqn:GMC-reg-b2}. \qed

The analog for GMC on $I_r$ is as follows. We omit the proof as it is the same as the proof of Proposition \ref{prop:gmc-reg-new}, requiring only notational differences. 

\bep \label{prop:GMC-reg-R}
Assume $\beta=1$. Let $f$ be a bounded measureable function in $I_r$. Fix $L \geq 100$ and $A_1 \geq 100L^2$. Let $Z_A(z, \eta_m)$ and $\delta_m$ be as in Proposition \ref{prop:gmc-reg-new}. Suppose for some $\alpha_1 >0$,
\beq
\ee[ \e^{ \lambda \Phi_N (z, \eta_m) } \1_{ \{ \Phi_N (z, \eta_m ) \leq A_1 \log N \}}] = (1 + \O (N^{-\alpha_1} ) ) \frac{N^{\frac{ \lambda(\lambda-1)}{4}} }{ \GB (\lambda ; z , 1)}  \e^{ \frac{ \lambda^2}{8} \kappa_4 (1 - |z|^2 )^2 - \frac{ \lambda \kappa_4}{4} (1 - |z|^2 )^2  }
\eeq
uniformly for $z \in I_r$ and $\lambda \in [0, L]$. Then,
\beq
\ee\left| \int_{I_r} f(x) \frac{ \e^{ \lambda \Phi_N (x, 0 )}}{\ee[ \e^{ \lambda \Phi_N (x, 0 )}]  } \d x - \int_{I_r} f(x) \frac{ \e^{ \lambda Z_A (x, \eta_m) }}{\ee[ \e^{ \lambda Z_A (x, \eta_m) }]  } \d x \right| \lesssim N^{-\frac{ \delta_m^{3/2}}{4}},
\eeq
as well as
\beq
\ee\left|  \int_{I_r} f(x) \frac{ \e^{ \lambda Z_A (x, \eta_m) }}{\ee[ \e^{ \lambda Z_A (x, \eta_m) }]  } \d x- \int f(x) \frac{ \e^{ \lambda Z_A (x, \eta_m) }}{N^{\frac{\lambda(\lambda-1)}{4}} \kB_\rr (x)  } \d x\right| \lesssim N^{-\alpha_1},
\eeq
where $\kB_\rr (x) := \GB (\lambda ; 0, 1) \e^{- \frac{ \lambda^2}{8} \kappa_4 (1 - |x|^2 )^2 + \frac{ \lambda \kappa_4}{4} (1 - |x|^2 )^2 }$.

\eep

\section{Convergence to the GMC}

\label{sec:claeys}

In this section we  use the methodology of \cite{claeys2021much} to deduce convergence to the GMC for Gaussian divisible ensembles on the disc. The changes for $(-1, 1)$ are outlined in Section \ref{sec:claeys-GMC-R}. The work  \cite{claeys2021much} provides a general framework for deducing GMC convergence; i.e., a list of assumptions (see \cite[Assumptions 2.5]{claeys2021much}) on an asymptotically Gaussian field $\xX(z)$ under which one has convergence to the GMC, \cite[Theorem 2.6]{claeys2021much}.

For the most part, the framework of \cite{claeys2021much}  applies without major changes to our set-up. However, there are a few differences which we now outline. We will be working with the matrices $X_t$ defined in \eqref{eqn:Xt-def} and quantities such as $\Phi_N (z, \eta, t)$ as defined in Section \ref{sec:def-dynamics}. For the discussion that follows we recall also the parameter $\eta_m = N^{-1-\delta_m}$, which has appeared previously in several sections.

\begin{enumerate}[label=(\roman*)]
\item The work \cite{claeys2021much}  regularizes the field $\xX(z)$ by convolving it with a smooth test function on some scale $\nu$, obtaining smoother fields $\xX_\nu (z)$. As can be seen by our methods, this is inconvenient for us, as our regularization parameter is the $\eta$ (actually in the language of \cite{claeys2021much} $\nu = \eta^{1/2}$) in $\Phi_N (z, \eta, t)$.

Instead, we will take (modulo some cut-offs and a rescaling, see Section \ref{sec:claeys-set-up} below) $\xX(z) =  \Phi_N (z, N^{-1-\delta_m} , t_1)$ and $\xX_\nu (z) =  \Phi_N (z, \nu^2, 0)$ for $\nu^2 \gg N^{-1/2}$. At this point note that $\xX(z)$ is evaluated at $\eta = \eta_m$ instead of $\eta = 0$ and at $t= t_1$. On the other hand, $\xX_\nu (z)$ is evaluated at $t=0$. The proof would also work with $\eta=0$ and $t=t_1$ everywhere, but the various computations would be more tedious. 

We will prove asymptotics for the joint Laplace transforms of these fields which are sufficient to apply the method of \cite{claeys2021much}. For the most part, \cite{claeys2021much}  makes no explicit use of the fact that the regularizations come from convolutions. In Section \ref{sec:claeys-application} below we explain why the parts of \cite{claeys2021much}  that do rely on this are still applicable to our set-up.

\item We cannot quite verify every assumption of \cite[Assumptions 2.5]{claeys2021much}, as they were sometimes written in a notationally convenient way, that turns out to be slightly stronger than what was necessary. Below we present our joint Laplace transform asymptotics in Proposition \ref{prop:claeys-assumption-1-new}. By reading the proof of \cite[Theorem 2.6]{claeys2021much}, one can easily see that these asymptotics are sufficient for the proofs to work. Details are given below in Section \ref{sec:claeys-application}.
\item Part of the method of \cite{claeys2021much}  is based on a modified $L^2$ argument, which leads to controlling the joint moments of $\xX(z)$ and $\xX(w)$. In particular, one needs to cut-off a contribution where $|z-w | \leq N^{-1/2+b}$; the assumptions in \cite{claeys2021much}  to control this contribution are stated in way that we cannot quite verify them. Instead, we provide an upper bound (Lemma \ref{lem:claeys-assumption-2}) and a further argument (Lemma \ref{lem:new-short-scale}) that shows this upper bound is sufficient to remove this contribution. Essentially, we are losing a small $N^{\eps}$ factor over \cite[Eq. (2.15)]{claeys2021much} which can instead be absorbed into \cite[Eq. (2.16)]{claeys2021much}.
\end{enumerate}

\subsection{Set-up} \label{sec:claeys-set-up}

Let $\gamma \in (0, 2)$, and define
\beq
c_\gamma := \min \left\{ \gamma, 2 - \gamma , \frac{1}{100} \right\} , \quad b = \frac{c_\gamma^2}{10^6} .
\eeq
Let $\mfb >0$ be the corresponding parameter from Theorem \ref{thm:local-global}. Let us choose
\beq \label{eqn:claeys-setup-1}
\omega_1 = \min \left\{ \frac{ c_\gamma^6}{10^6}, 10^{-10^6}, \frac{\mfb}{100} \right\},
\eeq
and let $\mfa >0$ be the parameter corresponding to this choice of $\omega_1 >0$ as in Theorem \ref{thm:local-global}. Finally, let
\beq \label{eqn:claeys-setup-2}
0 < \delta_m < \min\left\{ \frac{\mfa}{2} , \frac{c_W}{2} , \frac{1}{100}\right\}  , \qquad \eeps_{N} := \frac{N^{b}}{N^{1/2}}.
\eeq
Let $C_* \geq 1000$ be sufficiently large. We then recall the notation $\eta_m = N^{-1-\delta_m}$,  $t_1 = N^{\omega_1-1}$, and $\eta_* = ( \log N)^{C_*} / N$. These parameters will all play the same roles as they did in previous sections. On the other hand, the parameter $\eeps_N$ is new and corresponds to the $\eeps_N$ of \cite[Section 2]{claeys2021much}.  In order to have a definite limit, we will assume the following throughout Section \ref{sec:claeys}.
\begin{ass} \label{ass:claeys-kappa}
The distribution of the matrix entries of $X_0$ may depend on $N$, but in a way so that so that their fourth cumulants are $\kappa_4 + \O ( N^{-1/100} )$ for some fixed $\kappa_4 \geq - \frac{2}{\beta}$. 
\end{ass}
We require this as the limiting GMC measure involves the parameter $\kappa_4$. The $\O (N^{-1/100})$ error could be changed to $o(1)$ at the cost of adjusting the various estimates in this section appropriately. 

Fix now,
\beq
A_1 := 10^{10^{100}} 
\eeq
and let us define the fields,
\begin{align} \label{eqn:claeys-setup-3}
\xX(z) &:= (\sqrt{2} A_1 \log N) h_1 \left( \frac{  \Phi_N (z, \eta_m, t_1)}{A_1 \log N } \right) - \sqrt{2} \bmE (z, 0),
\notag\\
 \xX_\nu ( z) &:= ( \sqrt{2} A_1 \log N) h_2 \left( \frac{ \Phi_N (z, \nu^2, 0)}{A_1 \log N } \right) -\sqrt{2} \bmE (z, \nu^2),
\end{align} 
where $h_1$ and $h_2$ are the cut-off functions in Section \ref{sec:cutoff} and $\bmE$ is as in \eqref{eqn:bmE-def}. The re-scaling by $\sqrt{2}$ is to match the set-up of \cite{claeys2021much}. The choice of the different cut-off functions $h_1$ and $h_2$ is deliberate. Laplace transforms of $\xX(z)$ will be evaluated at only non-negative exponents and so we only need to cut-off the large positive values (we also do not have much control on the lower tail of $\xX(z)$ so it is convenient to not cut-off the large negative values). We will need to evaluate the Laplace transform of $\xX_\nu (z)$ with both positive and negative parameters, so we cut-off $\xX_\nu (z)$ in both directions.

The main result of this section is the following:
\bet \label{thm:GMC-gde} Let $0 < r < 1$. 
Let $f$ be a bounded continuous function supported in $\ddb_r$. Let $\gamma \in (0, 2)$ and let $\xX (z)$ be as above.  Then,
\beq
\lim_{N \to \infty} \int f(z) \frac{ \e^{ \gamma \xX (z)}}{ \ee[ \e^{ \gamma \xX (z)}]}  \dz = \int f(z) \d \mu_{\GMC, \kapb}^{\beta,\sqrt{2}\gamma}
\eeq
in distribution.  Here $\mu_{\GMC, \kapb}^{\beta,\sqrt{2}\gamma}$ is the measure defined in Section \ref{sec:GMC-def}. 
\eet

We require some estimates on the quantities $\xX(z)$. 
\bep \label{prop:claeys-apriori} Let the fields $\xX(z)$ and $\xX_\nu (z)$ be as above. The following holds uniformly for $z \in \ddb_r$. 
\begin{enumerate}[label=(\roman*),font=\normalfont]
\item \label{it:claeys-tail}
For all $u \geq 10$, and any $\nu \in [ N^{\omega_1/10-1/2}, 1]$, we have
\beq \label{eqn:claeys-tail} 
\pp\left[ \xX (z) \geq u \log N \right] \lesssim \e^{ - \frac{u^2}{3} \log N }, \quad \pp\left[ |\xX_\nu (z) | > u \log N \right] \lesssim \e^{ - \frac{u^2}{3} \log N }.
\eeq
\item 
Uniformly in $\lambda \in [0, 10^4]$, we have
\beq \label{eqn:claeys-mgf-bd-1}
\ee[\e^{ \lambda \xX(z)} ] \lesssim N^{\lambda^2} .
\eeq
\item 
Uniformly in $\lambda \in [-10^4, 10^4]$  and $\nu \in [ N^{\omega_1/10-1/2} , 1 ]$, we have
\beq \label{eqn:claeys-mgf-bd-2}
\ee[ \e^{ \lambda \xX_\nu (z) } ] \lesssim N^{\lambda^2} .
\eeq
\item Let $L>0$ and $a \in (0, 1)$  be given. For $N$ sufficiently large, uniformly in $\nu \in [a, 1]$ and $|\lambda| \leq L$, we have
\beq \label{eqn:claeys-mgf-bd-3}
\ee[ \e^{ \lambda \xX_\nu (z) } ] \lesssim 1.
\eeq
\end{enumerate}
\eep
\proof The estimates for $\xX (z)$ in \eqref{eqn:claeys-tail} and \eqref{eqn:claeys-mgf-bd-1} follow in a straightforward manner from Lemma \ref{lem:1-pt-apriori}. The second estimate of \eqref{eqn:claeys-tail} follows from \eqref{eqn:meso-tail-bd}, and \eqref{eqn:claeys-mgf-bd-2} follows in a straightforward manner from Proposition \ref{pro:jointlapdnew}. 

For the last item, we let $A_2 \geq A_1 + 10$. We assume that $\lambda \geq 0$, the other case is analogous. Then,
\begin{align}
    \ee[ \e^{ \lambda \xX_\nu (z) } ] &\lesssim 1 + \ee[ \e^{ \lambda \sqrt{2} \Phi_N (z, \nu^2, 0)} \1_{ \{ \Phi_N (z, \nu^2, 0) \leq A_2 \log N \} }] \notag\\
    &\quad +  N^{2 L A_1} \pp\left[ \Phi_N (z, \nu^2, 0) \geq A_2 \log N \right]
\end{align}
We also used the fact that $h_2(x) \leq x$ by assumption. By taking $A_2$ sufficiently large we see by Lemma \ref{lem:rough} that the last term is $\O (N^{-1})$. The second term is $\O (1)$ by Proposition \ref{pro:jointlapdnew}. \qed

\

We also define the covariance functional,
\begin{align}
\hC (z_1, z_2, \nu_1, \nu_2 ) &:= 2 \vV ( \nu_1^2, z_1, \nu_2^2, z_2) + \frac{\kappa_4}{2} (m^{z_1} ( \i \nu_1^2) m^{z_2} ( \i \nu_2^2 ) )^2 + \1_{\beta=1} 2 \vV (\nu_1^2, z_1, \nu_2^2 , \bar{z}_2 ),
\end{align}
where $\vV$ is as in \eqref{eqn:def-covar}. For $z \in \ddb_r$, by Proposition \ref{prop:parameter-properties}, we have
\beq \label{eqn:claeys-var}
\hC (z, z, \nu, \nu ) = - \frac{1}{2} \log (2 \sqrt{1-|z|^2} \nu^2) + \frac{\kappa_4}{2} (1-|z|^2)^2 -\frac{\1_{\beta=1}}{2} \log( |z-\bar{z}|^2) + \O ( \nu^2 ).
\eeq
Additionally, if $\nu_2 + |z_1 -z_2 | \geq \nu_1 $ then, for $z_1, z_2 \in \ddb_r$, by \eqref{eqn:covar-order-0}, we have
\beq \label{eqn:claeys-cov}
\hC(z_1, z_2, \nu_1, \nu_2 ) = \hC(z_1, z_2, 0, \nu_2 ) +\O \left( \frac{ \nu_1}{ |z_1-z_2| + \nu_2} \right) .
\eeq
By \eqref{eqn:covar-order-0} we have
\beq \label{eqn:hC-asymp}
\hC (z_1, z_2, \nu_1, \nu_2 ) = - \log \left( \nu_1 + \nu_2 + |z_1 - z_2 | \right) + \O (1).
\eeq
Note also that,
\beq
\hC(z_1, z_2, 0, 0) = \log ( |z_1-z_2|^{-1} ) + \frac{\kappa_4}{2}  (1-|z_1|^2)(1-|z_2|^2) +\1_{\beta=1} \log ( |z_1 - \bar{z_2} |^{-1} ).
\eeq

\subsection{Estimates} \label{sec:claeys-estimates}

In this section we prove the  estimates needed for the application of \cite{claeys2021much} . We begin with some preliminaries, based on Theorem \ref{thm:local-global}. We will use this result in the case $K=2$. Denote the spectral parameters by $z_1, z_2$ such that $|z_1 -z_2| \geq N^{b-1/2} = \eeps_N$. Denote the corresponding random variables from Theorem \ref{thm:local-global} by $\tilL^{(i)}_a$ and recall the notation $B_i := \sqrt{1 -|z_i|^2}$. Define,
\begin{align} \label{eqn:claeys-quantities-def}
    \hatL_1^{(i)} &:= \sqrt{2} ( \log N)^{3/4} h_1 \left( \frac{ \tilL^{(i)}_1 }{ ( \log N)^{3/4}} \right) , \qquad 
    \hatL_2^{(i)} :=\sqrt{2} ( A_1\log N) h_2 \left( \frac{ \tilL^{(i)}_2 }{ A_1 \log N } \right) \notag\\
    \hatL_3^{(i)} &:= \sqrt{2}  h_2 \left(  \tilL^{(i)}_3  \right) , \qquad \hatL^{(i)} := \sum_{a=1}^3 \hatL^{(i)}_a, \qquad 
    \hnu_i := \sqrt{B_i t_1} .
\end{align}
For $\lambda \in [0, \gamma]$ define
\beq \label{eqn:claeys-local-def}
\hatcalL (z_i, \lambda) := \ee\left[ \e^{\lambda \hatL^{(i)}} \right].
\eeq

Applying Theorem \ref{thm:local-global} we obtain the following.
\bel In the set-up of Theorem~\ref{thm:local-global} we have that there is a $\mfc_1 >0$ so that, (recalling the notation in \eqref{eqn:claeys-quantities-def})
\beq \label{eqn:claeys-homog}
\left| \xX (z_i) - \left( \hatL^{(i)} + \xX_{\hnu_i} (z_i) \right)\right| \lesssim N^{-\mfc_1},
\eeq
with probability at least $ 1- \O ( N^{-A_1^2/4})$, as well as
\beq \label{eqn:claeys-local}
\hatcalL (z_i, \lambda )  = \frac{  \e^{ \frac{\lambda^2}{4} \log (2 N B_i^2 t_1) }}{\GB ( \sqrt{2} \lambda ; 0 , 2 )} (1 + \O( N^{-\mfc_1} ) ).
\eeq
Moreover, uniformly for $\lambda \in [0, 10^4]$ and $a =1, 2, 3$ and $i=1, 2$, we have,
\beq \label{eqn:La-mgf-bd} 
\ee\left[ \e^{ \lambda \hatL^{(i)}_a  } \right] \lesssim N^{\lambda^2}. 
\eeq
\eel
\proof We start with \eqref{eqn:claeys-homog}. By Proposition \ref{prop:claeys-apriori} we have $\xX (z_i) = \sqrt{2} \Phi_N (z_i , \eta_m, t_1) - \sqrt{2} \bmE (z_i, 0)$, with probability at least $ 1- N^{-A_1^2/4}$, and a similar statement for $\xX_{\hnu_i} (z_i)$. Note also that $| \bmE (z_i , 0) - \bmE(z_i, \hnu_i^2) | \lesssim N^{-1/2}$. Similarly, by Theorem \ref{thm:local-global}\ref{it:lg-Z},\ref{it:lg-La} we have that $\hatL_a^{(i)} = \sqrt{2} \tilL_a^{(i)}$ with probability at least $1 - N^{-A_1^2}$, and so \eqref{eqn:claeys-homog} follows from \eqref{eqn:lg-main-est}.

For \eqref{eqn:claeys-local}, it is easy to see from Theorem \ref{thm:local-global} that, 
\beq
\left| \hatcalL(z_i, \lambda)  - \calL (B_i^2 t_1, B_i \eta_*, B_i \eta_m, \mfb, \sqrt{2} \lambda, A_1 ; z_i , \beta ) \right| \leq N^{-100} .
\eeq
The estimate \eqref{eqn:claeys-local} then follows from  Theorem \ref{thm:local-factor}, after possibly decreasing the value of $\mfc_1 >0$. 

The  estimate \eqref{eqn:La-mgf-bd} follows in a straightforward manner from their definitions in the case of $a=1, 3$ and from Theorem \ref{thm:local-global} \ref{it:lg-Z} for $a=2$.
\qed

\bep \label{prop:claeys-assumption-1-new} There is a small $\mfc_2 >0$ so that the following holds. 
Let $L \geq 1$, $K \geq 1$, $a \in (0, 1)$, and $0 < r< 1$. Choose the following parameters:
\begin{enumerate}[label=(\roman*),font=\normalfont]
\item Let $z_1, z_2 \in \ddb_r$ satisfy $|z_1 -z_2| \geq \eeps_N$ and $\gamma_1, \gamma_2 \in [0, \gamma]$. 
\item Let $\blam=  ( \lambda_1, \dots \lambda_{K+3} ) \in \rr^{K+3}$ satisfy $|\lambda_i| \leq 10$, for $i=1, 2, 3$, and $|\lambda_i| \leq L$, for $4 \leq i \leq K+3$. 
\item Let $\bnu = ( \nu_1, \dots, \nu_{K+3} ) \in \rr^{K+3}$ satisfy $\nu_i \in [\eeps_N, 1]$, for $i=1, 2, 3$, and $\nu_i \in [a, 1]$, for $4 \leq i \leq K+3$.
\item Let $\bw = (w_1, \dots, w_{K+3} ) \in (\ddb_r)^{K+3}$. 
\end{enumerate}
Uniformly in all choices of the above parameters we have,
\begin{align} \label{eqn:claeys-laplace-new}
   & \ee\left[ \e^{ \gamma_1 \xX (z_1) + \gamma_2 \xX (z_2) + \sum_{i=1}^{K+3} \lambda_i \xX_{\nu_i} (w_i )} \right] \notag\\
   = & ( 1 + \O ( N^{-\mfc_2} ) ) \frac{ N^{\frac{ \gamma_1^2 + \gamma_2^2}{4}} |z_1-z_w|^{ - \gamma_1 \gamma_2}}{\GB ( \sqrt{ 2} \gamma_1  ; 0, 2 ) \GB ( \sqrt{2} \gamma_2 ; 0 , 2 ) } \e^{ \frac{\kappa_4}{4} \left( \sum_{i=1}^2 \gamma_i (1 - |z_i|^2 ) \right)^2} 
   \notag\\
   \times & \e^{ \sum_{i=1}^2 \sum_{j=1}^{K+3} \gamma_i \lambda_j \hC ( z_i, w_j, 0, \nu_j ) + \frac{1}{2} \sum_{i, j=1}^{K+3} \lambda_i \lambda_j \hC (w_i, w_j, \nu_i, \nu_j )   - \frac{\1_{\beta=1}}{2} \sum_{i,j=1}^2 \gamma_i \gamma_j \log ( |z_i - \bar{z_j} |)} .
\end{align}
\eep
\proof We prove only the case $\beta=2$, the case $\beta=1$ being almost identical. Let $\F$ be some event. We first observe the following estimate, which is uniform over the choices of parameters in the statement of the Proposition:
\begin{align} \label{eqn:claeys-general}
& \ee\left[ \1_\F \e^{ \sum_{i=1}^2 \gamma_i \xX (z_i) + \sum_{j=1}^{K+3} \lambda_j \xX_{\nu_j} (w_j) } \right] \notag\\
\leq & \pp[ \F]^{1/2} \left( \ee[ \e^{ 4 \sum_{i=1}^2 \gamma_i \xX (z_i) + 4\sum_{j=1}^{3} \lambda_j \xX_{\nu_j} (w_j) }] \right)^{1/4} \left( \ee[ \e^{ 4  \sum_{j=4}^{K+3} \lambda_j \xX_{\nu_j} (w_j) } ] \right)^{1/4} \notag\\
\leq &\pp[ \F]^{1/2} N^{10^5},
\end{align}
for all $N$ sufficiently large. The last estimate follows from Proposition \ref{prop:claeys-apriori}. In particular, the exponent of $N$ in the last line is fixed independent of $K$ and $L$, so as long as $\pp[\F]$ tends to $0$ sufficiently fast (for example if $\pp[\F] \leq N^{-A_1^2/100}$) then the above tends to $0$ as $N\to \infty$. It is also easy to see that the same estimate also holds if we replace the $\xX (z_i)$ by $\hatL^{(i)} + \xX_{\hnu_i} (z_i)$, using also \eqref{eqn:La-mgf-bd}. 

By \eqref{eqn:claeys-homog} we have,
\begin{align}
& \ee\left[ \e^{ \sum_{i=1}^2 \gamma_i \xX (z_i) + \sum_{j=1}^{K+3} \lambda_j \xX_{\nu_j} (w_j) }\right] \notag\\
&= (1 + \O ( N^{-\mfc_1} ) )  \ee\left[ \e^{ \sum_{i=1}^2 \gamma_i (\xX_{\hnu_i} (z_i) + \hatL^{(i)} )))+ \sum_{j=1}^{K+3} \lambda_j \xX_{\nu_j} (w_j) }\right] + \O (N^{-A_1^2/5} )
\end{align}
We used also \eqref{eqn:claeys-general} to bound the complementary event where \eqref{eqn:claeys-homog} does not hold.
Now by definition,
\begin{align} \label{eqn:claeys-laplace-a1}
     \ee\left[ \e^{ \sum_{i=1}^2 \gamma_i (\xX_{\hnu_i} (z_i) + \hatL^{(i)} )))+ \sum_{j=1}^{K+3} \lambda_j \xX_{\nu_j} (w_j) }\right] 
    = \left( \prod_{i=1}^2  \hatcalL (z_i, \gamma_i ) \right) \ee\left[ \e^{ \sum_{i=1}^2 \gamma_i \xX_{\hnu_i} (z_i)+ \sum_{j=1}^{K+3} \lambda_j \xX_{\nu_j} (w_j) }\right]
\end{align}
For any $A \geq 1$ define the event,
\beq
\F(A) := \left( \bigcap_{i=1}^2  \left\{ | \Phi_N (z_i, \hnu_i^2 , 0)| \leq A \log N \right\} \right) \cap \left( \bigcap_{j=1}^{K+3} \left\{  |\Phi_N (w_j, \nu_j^2, 0) | \leq A \log N \right\} \right).
\eeq
From \eqref{eqn:meso-tail-bd} we have,
\beq
\pp\left[ \F(A)^c\right] \lesssim N^{-A^2/2}
\eeq
with the implicit constant depending on $A$ and $K$. In particular, by \eqref{eqn:claeys-mgf-bd-2} and \eqref{eqn:claeys-mgf-bd-3} one obtains (in a manner analogous to \eqref{eqn:claeys-general}) that,
\begin{align}
& \ee\left[ \e^{ \sum_{i=1}^2 \gamma_i \xX_{\hnu_i} (z_i)+ \sum_{j=1}^{K+3} \lambda_j \xX_{\nu_j} (w_j) }\right]  \notag\\
= & \ee\left[ \e^{ \sqrt{2} \sum_{i=1}^2 \gamma_i \Phi_N (z_i, \hnu_i^2, 0)+ \sqrt{2} \sum_{j=1}^{K+3} \lambda_j \Phi_N (w_j, \nu_j^2, 0) } \1_{ \F (\frac{A_1}{2} ) } \right] Q_E + \O (N^{-A_1^2/10} ).
\end{align}
Here we denoted $Q_E := \e^{ - \sqrt{2} \sum_{i=1}^2 \gamma_i \bmE (z_i, \hnu_i^2 ) - \sqrt{2} \sum_{j=1}^{K+3} \lambda_j \bmE (w_j, \nu_j^2 )}$. 
In a similar manner, for any $A_2 \geq 10 A_1$, by \eqref{eqn:meso-uni-bd}, we have
\begin{align}
    & \ee\left[ \e^{ \sqrt{2} \sum_{i=1}^2 \gamma_i \Phi_N (z_i, \hnu_i^2, 0)+ \sqrt{2} \sum_{j=1}^{K+3} \lambda_j \Phi_N (z_j, \nu_j^2, 0) } \1_{ \F (\frac{A_1}{2} ) } \right] \notag\\
    = & \ee\left[ \e^{ \sqrt{2} \sum_{i=1}^2 \gamma_i \Phi_N (z_i, \hnu_i^2, 0)+ \sqrt{2} \sum_{j=1}^{K+3} \lambda_j \Phi_N (w_j, \nu_j^2, 0) } \1_{ \F (A_2 ) } \right] + \O (N^{-A_1^2/10} ).
\end{align}
If $A_2$ is sufficiently large, then by Proposition \ref{pro:jointlapdnew} we have,
\begin{align}
    & \ee\left[ \e^{ \sqrt{2} \sum_{i=1}^2 \gamma_i \Phi_N (z_i, \hnu_i^2, 0)+ \sqrt{2} \sum_{j=1}^{K+3} \lambda_j \Phi_N (w_j, \nu_j^2, 0) } \1_{ \F (A_2 ) } \right] Q_E 
    =  (1 + \O ( N^{-\mfc_2} )) \notag\\
    \times &\e^{  \frac{1}{2} \sum_{i,j=1}^2 \gamma_i \gamma_j \hC ( z_i,z_j, \hnu_i,\hnu_j) + \sum_{i=1}^2 \sum_{ j=1}^{K+3} \gamma_i \lambda_j \hC (z_i,w_j,\hnu_i,\nu_j ) + \frac{1}{2} \sum_{i,j=1}^{K+3} \lambda_i \lambda_j \hC (w_i,w_j,\nu_i,\nu_j ) } .
\end{align}
The above estimate, together with \eqref{eqn:claeys-local} (and also \eqref{eqn:claeys-var} and \eqref{eqn:claeys-cov} to rewrite the $\hC$ terms involving a $z_i$) proves the estimate \eqref{eqn:claeys-laplace-new}, up to an additional error of $\O (N^{-A_1^2/10} )$ on the RHS. We just need to check that this additive error can be absorbed into the multiplicative one on the RHS of \eqref{eqn:claeys-laplace-new}. For this, we note that by \eqref{eqn:covar-order-0} all the $\hC$ terms involving $w_j$ for $j \geq 4$ satisfy $|\hC | \lesssim 1$. For the other 15 appearances of $\hC$ we have the general estimate $|\hC (u_1, u_2, \nu_1, \nu_2 )| \leq |\log (\nu_1+\nu_2) | + \O (1)$, by \eqref{eqn:hC-asymp}, and so the RHS of \eqref{eqn:claeys-laplace-new} is bounded below by $N^{-10^5}$, which is sufficient. This completes the proof. 
\qed

\

Proposition \ref{prop:claeys-assumption-1-new} holds only when $|z_1-z_2| \gtrsim N^{-1/2}$. We prove the following estimate valid for all choices of $z_1, z_2$ which loses the small factor $N^{(\gamma_1+\gamma_2)(\omega_1)^{1/3}}$ over the RHS of \eqref{eqn:claeys-laplace-new}.

\bel \label{lem:claeys-assumption-2}
Let $z_1, z_2, z_3 \in \ddb_r$, $\nu \in [\eeps_N, 1]$, $|\lambda| \leq 10$, and $\gamma_1, \gamma_2 \in [0,\gamma]$. Then,
\begin{align} \label{eqn:claeys-upper}
    & \ee[ \e^{\gamma_1 \xX (z_1) + \gamma_2 \xX (z_2) + \lambda \xX_\nu (z_3) } ] \notag\\
    & \lesssim \e^{(\gamma_1+\gamma_2) \omega_1^{1/3} \log N } \frac{ N^{\frac{\gamma_1^2+\gamma_2^2}{4} }}{(|z_1-z_2| + N^{-1/2} )^{\gamma_1\gamma_2}} \e^{ \frac{\lambda^2}{2} \log ( \nu^{-1} ) - \lambda ( \gamma_1 \log ( \nu + |z_1-z_3| ) + \gamma_2 \log ( \nu + |z_2 -z_3| ) ) }.
\end{align}
\eel
\proof By Proposition \ref{prop:useful-upper} and the choice of $\omega_1 \leq 10^{-10^6}$, for $i=1, 2$, we have
\beq \label{eqn:claeys-change-1}
\xX (z_i ) \leq \xX_{\hnu_i} (z_i) + \omega_1^{1/3} \log N
\eeq
with probability at least $1 - \e^{-\frac{1}{100}  \omega_1^{-1/3} \log N} - N^{-A_1^2/10}\geq 1- N^{-10^{100} }$. Therefore, using Proposition \ref{prop:claeys-apriori} on the complementary event where the above estimate does not hold we find,
\beq \label{eqn:claeys-upper-a1}
 \ee[ \e^{\gamma_1 \xX (z_1) + \gamma_2 \xX (z_2) + \lambda \xX_\nu (z_3) } ]  \leq \e^{ ( \gamma_1 + \gamma_2) \omega_1^{1/3} \log N } \ee[ \e^{\gamma_1 \xX_{\hnu_1} (z_1) + \gamma_2 \xX_{\hnu_2} (z_2) + \lambda \xX_\nu (z_3) } ]  + \O (N^{-10^{50}} ).
\eeq
Note that the error on the RHS can be absorbed into the RHS of \eqref{eqn:claeys-upper}, which can be seen to be lower bounded by $N^{-1000}$. 

The computation of the first term on the RHS of \eqref{eqn:claeys-upper-a1} follows from Proposition \ref{pro:jointlapdnew} the same way as of \eqref{eqn:claeys-laplace-a1} in the proof of Proposition \ref{prop:claeys-assumption-1-new}. The various $\hC$ terms that arise are bounded by \eqref{eqn:claeys-var} and \eqref{eqn:hC-asymp}.  \qed

\subsection{Application of the framework of Claeys, Fahs, Lambert and Webb} \label{sec:claeys-application}

The purpose of this section is to explain how to adapt the method of \cite{claeys2021much} to prove Theorem \ref{thm:GMC-gde}. At this point we do first the case $\beta=2$. The changes for $\beta=1$ are discussed in Section \ref{sec:claeys-real} below. We first introduce some notation, following closely 
 \cite{claeys2021much}. Let $\mfx >0$ be chosen as,
\beq
\mfx = \frac{c_\gamma}{10^6},
\eeq
which should be compared with \cite[Eq. (2.27)]{claeys2021much}. Note that the parameter $\tau >0$ is not present in our set-up, and in fact the lower bound $\mfx > \sqrt{2 \tau}$ of \cite{claeys2021much} is not actually needed in that work (the constraint is only used to bound quantities that get \emph{larger} as $\mfx$ increases - specifically \cite[Eq. (2.30)]{claeys2021much} -  so it is no loss of generality to dispense with the assumption).

For any $\ell < L$ define the barrier events and parameters, 
\beq
A_{\ell, L} (z) := \bigcap_{\ell \leq k \leq L} \{ \xX_{\e^k} (z) \leq (\gamma+ \mfx ) k \}, \qquad \nu_L := \e^{-L}, \qquad M_N := \log ( \eeps_N^{-1} ) .
\eeq
Let $f$ be a bounded continuous function supported in $\dd_r$. By \cite[Lemma 2.7]{claeys2021much} we have,
\begin{align} \label{eqn:barrier}
\limsup_{N \to \infty} \ee\left[ \int |f(z)| \1_{A_{\ell, M_N } (z)^c} \frac{ \e^{ \gamma \xX(z)}}{\ee[ \e^{ \gamma \xX(z)}]} \dz \right] \lesssim \e^{ - \mfx^2\ell/2}, \notag\\
\limsup_{N \to \infty} \ee \left[ \int |f(z)| \1_{A_{\ell, L} (z)^c} \frac{ \e^{\gamma  \xX_{\nu_L} (z)}}{\ee[ \e^{ \gamma \xX_{\nu_L} (z)}]} \dz \right] \lesssim \e^{ - \mfx^2\ell/2},
\end{align}
uniformly for any choice of $\ell$ and  $L$. Indeed, the proof of the first estimate, given in \cite[Lemma 2.7]{claeys2021much}  uses only \eqref{eqn:claeys-laplace-new} with: (i) the choices $\gamma_1 = \gamma$, $z_1=z$, and the other parameters equal to $0$ to estimate the expectation in the denominator; (ii) the choice $\gamma_1 = \gamma$, $\lambda_1 = \mfx$,  $\nu_1 \in [ \eeps_N, 1]$, $z_1 = w_1 = z$, and other parameters equal to $0$ to estimate the numerator via Markov's inequality. Similar considerations apply for the second estimate (in which case $\gamma_1 = \gamma_2 = 0$ and only $\lambda_1, \lambda_2$ are ever non-zero). 

We restate \cite[Lemma 2.8]{claeys2021much}  as the following.
\bel \label{lem:2.8} 
With $\mfx >0$ as above, we have 
\begin{align}
    \lim_{ \ell \to \infty} \limsup_{L \to \infty} \limsup_{N \to \infty} \ee \left| \int f(z) \1_{A_{\ell, M_N} (z)} \frac{\e^{ \gamma \xX(z)}}{\ee[\e^{ \gamma \xX(z)}]} \dz -\int f(z) \1_{A_{\ell, L} (z)} \frac{\e^{ \gamma \xX_{\nu_L} (z)}}{\ee[\e^{ \gamma \xX_{\nu_L} (z)}]} \dz \right|^2  = 0.
\end{align}
\eel
The proof of Lemma \ref{lem:2.8} is essentially the same as in \cite{claeys2021much}. In Section \ref{sec:2.8} below we summarize the minimal changes (in particular the cut-off of the short-range contribution where $|z-w| \leq \eeps_N$; see Lemma \ref{lem:new-short-scale}).

We now have the ingredients to complete the proof of Theorem \ref{thm:GMC-gde}.

\vspace{3 pt}

\noindent{\bf Proof of Theorem \ref{thm:GMC-gde}}. Naturally, we follow closely the proof of \cite[Theorem 2.6]{claeys2021much}, encouraging the reader to refer there where necessary. Let $F : \rr \to \rr$ be a Lipschitz function. Following exactly as in   \cite[Proof of Theorem 2.6]{claeys2021much} (specifically, the proof of \cite[Eq. (2.21)]{claeys2021much}), we see that Lemma \ref{lem:2.8} and \eqref{eqn:barrier} imply
\beq
\lim_{N \to \infty} \ee F\left( \int f(z) \frac{ \e^{\gamma \xX(z)}}{\ee[\e^{ \gamma \xX(z)}]} \d^2 z \right)  = \lim_{L \to \infty} \lim_{N \to \infty} \ee F\left( \int f(z) \frac{ \e^{\gamma \xX_{\nu_L} (z)}}{\ee[\e^{ \gamma \xX_{\nu_L} (z)}]} \d^2 z \right),
\eeq
as long as the limit on the RHS exists. We therefore turn to the computation of the RHS.

For fixed $L$, we have from \eqref{eqn:claeys-laplace-new} (with $\gamma_i = \lambda_j = 0$, for $i =1, 2$ and $j=1, 2, 3$, and $\lambda_j$ arbitrary for $4 \leq j \leq K+3$) that, as $N\to \infty$, $\{\xX_{\nu_L} (z)\}_{z \in \dd_r}$  converges jointly (in the finite dimensional marginals sense) to a Gaussian process with covariance,
\beq
K_{\nu_L} (z, w) :=  \hC (z,  w, \nu_L,\nu_L ) = - \log ( |z-w| + \nu_L ) + \O (1) .
\eeq
With this convergence as an input, the proof
 of \cite[Lemma 2.2]{lambert2018subcritical}  applies without major changes (in particular, see the discussion immediately after \cite[Lemma 2.3]{claeys2021much} which is just a restatement of \cite[Lemma 2.2]{lambert2018subcritical})  and shows that,
\beq \label{eqn:claeys-thm-a1}
\lim_{N \to \infty} \ee F\left( \int f(z) \frac{ \e^{\gamma \xX_{\nu_L} (z)}}{\ee[\e^{ \gamma \xX_{\nu_L} (z)}]} \d^2 z \right)   = \ee F\left( \int f(z) \frac{ \e^{\gamma \yY_{\nu_L} (z)}}{\ee[\e^{ \gamma \yY_{\nu_L} (z)}]} \d^2 z \right)  
\eeq
where $\yY_{\nu_L}(z)$ is the Gaussian process with kernel $K_{\nu_L} (z, w)$. We now claim that as $L \to \infty$, the RHS converges to the GMC. Indeed, this is in some sense the definition of GMC, as for fixed $\nu_L >0$, the kernel $K_{\nu_L} (z, w)$ is a regularization of the kernel $K_0 (z, w) = - \log |z-w| + \O (1)$. 

However, the GMC is usually defined via a different regularization, i.e., via mollification of the process $\yY_0 (z)$ with a radial test function. In our setting  it is not guaranteed that for fixed $\nu_L > 0$ the field $\yY_{\nu_L} (z)$ has the distribution of a mollification. In the framework of \cite{claeys2021much}, $\yY_{\nu_L}$ is literally a mollification, so the convergence of the RHS of \eqref{eqn:claeys-thm-a1} to the GMC is by definition (and that is how the proof of \cite[Theorem 2.6]{claeys2021much} wraps up).

In our case, we can deduce the convergence of the RHS of \eqref{eqn:claeys-thm-a1} to the GMC  by simply applying \cite[Theorem 2.6]{claeys2021much}  to the process $\yY_{\nu_L}$ itself. I.e., regard $n = \nu_L^{-2}$ as the fundamental large parameter, set $X(z) = \yY_{\nu_L} (z)$ and $X_\eta (z)$ to be the convolution of $\yY_{\nu_L}$ with a radial smooth mollifier. As $X(z)$ is a Gaussian process, the assumptions of \cite[Theorem 2.6]{claeys2021much} trivially hold. 
 Therefore, the RHS of \eqref{eqn:claeys-thm-a1} converges to the GMC, completing the proof. \qed

\subsubsection{Sketch of proof of Lemma \ref{lem:2.8}}  \label{sec:2.8}

The proof of Lemma \ref{lem:2.8} follows (very) closely that of \cite[Lemma 2.8]{claeys2021much} . The proof there is based on three lemmas, \cite[Lemmas 2.10--2.12]{claeys2021much}. The proof of \cite[Lemma 2.11]{claeys2021much}  is unchanged and relies  only on various cases of the estimate \eqref{eqn:claeys-laplace-new}. We therefore have no need to comment further on its proof.

The proof of \cite[Lemma 2.10]{claeys2021much} almost goes through without change. First, the proof of \cite[(2.29)]{claeys2021much} is identical. The proof of \cite[(2.30)]{claeys2021much} needs a minor change in order to estimate the contribution of
\begin{align} \label{eqn:claeys-2.8-a1}
    \int_{ |z-w| < \eeps_N, (z, w) \in \dd_r^2 } \frac{ \ee[ \e^{ \gamma (\xX(z) + \xX(w))} \1_{ A_{\ell, M_N} (z) \cap A_{\ell, M_N} (w)} ] }{\ee[\e^{\xX(z)}] \ee[ \e^{\xX(w)}]} \dz  \dw.
\end{align}
That is, we need a replacement for \cite[Eq. (2.38)]{claeys2021much} which will show that the above term tends to $0$. This will be seen to follow from the following lemma. Note that the proof is still similar to \cite[Lemma 2.10]{claeys2021much}, substituting \eqref{eqn:claeys-upper} for \cite[Eqs. (2.15)--(2.16)]{claeys2021much}.  
\bel \label{lem:new-short-scale} 
For any choice of $\mfx \in (0, \gamma)$ we have,
\beq \label{eqn:new-short-scale}
 \int_{ |z-w| < \eeps_N, (z, w) \in \dd_r^2 } \frac{ \ee[ \e^{ \gamma (\xX(z) + \xX(w))} \1_{ \{\xX_{\eeps_N} (z) \leq (\gamma+\mfx) \log ( \eeps_N^{-1} ) \} } ] }{\ee[\e^{\xX(z)}] \ee[ \e^{\xX(w)}]} \dz \dw \lesssim N^{4 \omega_1^{1/3}} N^{4 b + \mfx - (1- \gamma^2/4) }. 
\eeq
\eel
\proof Let $z, w$ be as in the integrand above. We have, by Markov's inequality (note that $\mfx < \gamma$) and \eqref{eqn:claeys-upper} (using that $\eeps_N + |z-w| \asymp \eeps_N$)
\begin{align}
    & N^{-\gamma^2/2} \ee[ \e^{ \gamma (\xX(z) + \xX(w))} \1_{ \{\xX_{\eeps_N} (z) \leq (\gamma+\mfx) \log ( \eeps_N^{-1} ) \} } ] \notag\\
    \leq & N^{-\gamma^2/2} \e^{ (\gamma^2-\mfx^2) \log (\eeps_N^{-1} ) }  \ee[ \e^{ \gamma (\xX(z) + \xX(w)) + (\mfx-\gamma) \xX_{\eeps_N} (z) } ] \notag\\
    \lesssim  & N^{-\gamma^2/2} \e^{ (\gamma^2-\mfx^2) \log (\eeps_N^{-1} )  + 2 \gamma \omega_1^{1/3} \log N} N^{\gamma^2} \e^{ \frac{ ( \gamma - \mfx)^2}{2} \log \eeps_N^{-1} + (\mfx - \gamma) 2 \gamma \log \eeps_N^{-1} } \notag\\
    = & N^{\gamma^2/2} N^{2 \gamma \omega_1^{1/3}} ( \eeps_N)^{\frac{ (\gamma - \mfx)^2}{2}}.
\end{align}
Therefore, since by \eqref{eqn:claeys-laplace-new} we have $\ee[ \e^{ \xX (z) } ] \asymp N^{\gamma^2/4}$, we obtain
\begin{align}
& \int_{ |z-w| < \eeps_N, (z, w) \in \dd_r^2 } \frac{ \ee[ \e^{ \gamma (\xX(z) + \xX(w))} \1_{ \{\xX_{\eeps_N} (z) \leq (\gamma+\mfx) \log ( \eeps_N^{-1} ) \} } ] }{\ee[\e^{\xX(z)}] \ee[ \e^{\xX(w)}]} \dz \dw  \notag\\
\lesssim & N^{2 \gamma \omega_1^{1/3}} N^{\gamma^2/2} (\eeps_N)^{2 + \frac{( \gamma-x)^2}{2}} \notag\\
= & N^{2 \gamma \omega_1^{1/3}} N^{\gamma^2/2 -1+2b - \frac{ ( \gamma-\mfx)^2}{2} (1/2-b)} \notag\\
= & N^{2 \gamma \omega_1^{1/3}} N^{ 2b - \frac{\mfx^2}{2} (1/2-b)  + \gamma \mfx (1/2-b) + \frac{\gamma^2}{2} b - (1- \gamma^2/4)} \leq N^{4 \omega_1^{1/3}} N^{4 b + \mfx - (1- \gamma^2/4) },
\end{align}
which proves the claim. Note that the last inequality follows from the fact that $\gamma \leq 2$. 
\qed

\

By our choice of $\omega_1$,  $b$ and $\mfx$ we have,
\beq
4 \omega_1^{1/3} + 4b + \mfx < \frac{1}{4} \left( 1- \frac{\gamma^2}{4} \right),
\eeq
and so the contribution \eqref{eqn:claeys-2.8-a1} tends to $0$ as $N \to \infty$. 
The other parts of the proof of \cite[Lemma 2.10]{claeys2021much}  go through without change.

Finally, we turn to the proof of \cite[Lemma 2.12]{claeys2021much} . The first part of the proof is based on computing the limits on the LHS of \cite[Eqs. (2.33)--(2.34)]{claeys2021much}. In our setting, the estimates in \eqref{eqn:claeys-laplace-new} allow us to do so, where we  leverage the full generality of our estimate; i.e., allowing the $\{\nu_i \}_{i=4}^{K+3}$ and $\{\lambda_i \}_{i=4}^{K+3}$ to be arbitrary.

A priori, the second part of the proof of \cite[Lemma 2.12]{claeys2021much} uses the fact that the regularized quantities $\xX_\nu (z)$ come from averages of the field $\xX(z)$. However, the only specific property that is used here is the continuity of the limiting covariance function (in our case $\hC(z_1, z_2, \nu_1, \nu_2)$) in the parameters $\nu_i, z_i$ away from the diagonal $\{ (z, z, 0, 0) \}$. In our case the continuity follows directly from the exact formula we have for $\hC$, and was previously verified in Proposition \ref{prop:parameter-properties}.

With the analogs of \cite[Lemmas 2.10--2.12]{claeys2021much} established, the proof of Lemma \ref{lem:2.8} follows as in \cite{claeys2021much}. \qed

\subsubsection{The real half-disc case} \label{sec:claeys-real}

The proof for the case $\beta=1$ is almost identical. The main difference is that we are now working with functions supported in the half disc $\ddtwo_r$ (recall the definition \eqref{eqn:ddb-def}) and the covariance functional (and expectation function $\bmE (z, \nu)$) have a singularity near the real axis.

There is no difference between working in the full and half-disc, as the method of \cite{claeys2021much} is formulated for any subset of $\cc$. Additionally, since we are working with test functions $f$ that are supported away from the real axis the singularity in the covariance functional does not cause any issue. For $z, w \in \ddtwo_r$ we still have
\beq
\hC (z, w, 0, 0) = - \log( |z-w| ) + \O ( 1)
\eeq
which is the main requirement for the method; see \cite[(2.2)]{claeys2021much}. With this in mind, the proof of Theorem \ref{thm:GMC-gde} given above goes through with only notational changes.

\subsection{The real case on the interval} \label{sec:claeys-GMC-R}

We now outline the changes to the previous arguments in order to prove convergence to the GMC on the interval $(-1, 1)$ in the case that $\beta=1$. The framework in \cite{claeys2021much} does not depend in a nontrivial manner on the dimension of the ambient space, and the methodology by which we compute the various joint Laplace transforms is largely insensitive as to whether the $z_i \in \rr$ or not. The changes required are therefore largely notational and we will only sketch the analogous arguments.

In the one-dimensional case, the parameter $\gamma$ will vary in the range $\gamma \in (0, \sqrt{2})$. Accordingly, we define $\hat{c}_\gamma = \min \{ \gamma, \sqrt{2} - \gamma, \frac{1}{100}\}$. We let $b = \frac{ ( \hat{c}_\gamma)^2}{10^6}$, and let $\mfb >0$ be the corresponding parameter from Theorem \ref{thm:local-global}. We let $\omega_1$, $\eeps_N$ and $\delta_m$ be as in \eqref{eqn:claeys-setup-1} and \eqref{eqn:claeys-setup-2} with the $c_\gamma$ there replaced by $\hat{c}_\gamma$ here. We let $\eta_m$, $t_1$ and $\eta_*$ be as in Section \ref{sec:claeys-set-up} and continue to assume Assumption \ref{ass:claeys-kappa}. Instead of \eqref{eqn:claeys-setup-3} we define,
\begin{align}
\hxX (x) &:= (A_1 \log N) h_1 \left( \frac{ \Phi_N (x, \eta_m, t_1) }{ A_1 \log N} \right) +\frac{\kappa_4}{4} (1 - |x|^2 )^2 +\frac{1}{4} \log N \notag\\
\hxX_\nu (x) &:=( A_1 \log N) h_2 \left( \frac{ \Phi_N (x, \nu^2, 0)}{ A_1 \log N } \right) - \bmE (x, \nu^2) .
\end{align}
The analog of Theorem \ref{thm:GMC-gde} is the following.
\bet \label{thm:GMC-gde-line}
Let $0 < r < 1$ and let $f$ be a bounded continuous function supported in $I_r$. Let $\gamma \in (0, \sqrt{2} )$ and $\hxX (x)$ be as above. Then,
\beq
\lim_{N \to \infty} \int f(x) \frac{ \e^{ \gamma \hxX (x)}}{ \ee[ \e^{ \gamma \hxX (x)}]} \d x = \int f(x) \d \mu^{\rr, \gamma}_{\GMC, \kappa_4} (x)
\eeq
\eet
In the remainder of the section we  outline how the necessary modifications to Sections \ref{sec:claeys-estimates} and \ref{sec:claeys-application} in order to prove the above theorem. The changes are largely notational and so we will be brief. 

The fields $\hxX (x)$ are approximately centered with limiting covariance kernel
\beq
\rcC (x_1, x_2, \nu_1, \nu_2 ) := 2 \vV (\nu_1^2 , x_1, \nu_2^2, x_2 ) + \frac{\kappa_4}{2} \left( m^{x_1} ( \i \nu_1^2) m^{x_2} ( \i \nu_2^2 ) \right) 
\eeq
In order to prove the analog of the joint Laplace transform asymptotics in Proposition \ref{prop:claeys-assumption-1-new} one must modify the definitions of the local quantities \eqref{eqn:claeys-quantities-def} in the obvious manner: one omits the $\sqrt{2}$ factors, introduces $\hatL^{(i)}_4 := (\log N) h_2 \left( \frac{ \tilL^{(i)}_4}{ \log N } \right)$, and then includes the additional local factor $\hatL^{(i)}_4$ in the definition of the local factor \eqref{eqn:claeys-local-def}.

With these modifications one finds by a similar proof to that of Proposition \ref{prop:claeys-assumption-1-new} the following:
\bep
There is a small $\mfc_3 >0$ so that the following holds. Let $L \geq 1$, $K \geq 1$, $a \in (0, 1)$ and $0 < r < 1$. 
 Let $x_1, x_2 \in I_r$ satisfy $|x_1-x_2| \geq \eeps_N$ and $\gamma_1, \gamma_2 \in [0, \gamma]$. Let $\blam, \bnu \in \rr^{K+3}$ be as in Proposition \ref{prop:claeys-assumption-1-new} and $(y_1, \dots, y_{K+3} ) \in (I_r)^{K+3}$.

 Uniformly in all of these parameters we have,
 \begin{align}
 & \ee\left[ \e^{ \gamma_1 \hxX (x_1 ) + \gamma_2 \hxX (x_2) + \sum_{i=1}^{K+3} \lambda_i \hxX_{\nu_i} (y_i ) } \right] 
=  \frac{N^{\frac{ \gamma_1^2+\gamma_2^2}{4}} |x_1 -x_2 |^{-\gamma_1 \gamma_2}}{ \GB (\gamma_1 ; 0, 1) \GB ( \gamma_2 ; 0 , 1 ) } \e^{ \frac{\kappa_4}{4}  \left( \sum_{i=1}^2 \gamma_i (1- |x_i |^2 ) \right)^2} \notag\\
\times &  \e^{  \sum_{i=1}^2 \sum_{j=1}^{K+3} \gamma_i \lambda_j \rcC ( x_i, y_j, 0, \nu_j ) + \frac{1}{2} \sum_{i, j=1}^{K+3} \lambda_i \lambda_j \rcC (y_i, y_j, \nu_i, \nu_j )  } \times \left(1 + \O (N^{-\mfc_2} ) \right)
 \end{align}
\eep

With the above asymptotics the arguments of Section \ref{sec:claeys-application} go through with only one minor change. That is, the estimate of the short-range part in Lemma \ref{lem:new-short-scale} is slightly different. In order to prove it, one notes that instead of \eqref{eqn:claeys-change-1} we have
\beq \label{eqn:claeys-change-1a}
\hxX (x_i ) \leq \hxX_{ \hat{\nu}_i} ( x_i ) + ( \omega_1)^{1/3} \log N
\eeq
with probability at least $ 1- N^{-10^{100}}$ by writing the last term on RHS of \eqref{eqn:useful-upper} as
\beq
\frac{\omega_1}{4} \log N = \bmE (x_i, \hat{\nu}_i ) + \frac{1}{4} \log N + \O ( 1) .
\eeq
Using \eqref{eqn:claeys-change-1a} as input, the estimate \eqref{eqn:new-short-scale} becomes
\beq
\int_{ |x-y| < \eeps_N } \frac{ \ee[\e^{ \gamma ( \hxX (x) + \hxX (y) )} \1_{ \{ \hxX_{\eeps_N} (x) \leq ( \gamma + \mfx) \log ( \eeps_N^{-1} ) \}}]}{ \ee[ \e^{ \gamma \hxX (x)}] \ee[ \e^{ \gamma \hxX (y)}]} \d x \d y \lesssim N^{4 \omega_1^{1/3}} N^{4b + \mfx - (1/2 - \gamma^2/4)}
\eeq
by the same proof. By the choice of parameters $b, \omega_1$ and $\mfx$, the RHS indeed goes to $0$ (recall $\gamma^2 <2$ now in the one-dimensional case). The remainder of the arguments in Section \ref{sec:claeys-application} go through with only the appropriate notational changes. This concludes our sketch of the proof of Theorem \ref{thm:GMC-gde-line}. 

\section{Proof of mesoscopic estimates} \label{sec:meso}

In this section we present the proofs of Proposition \ref{pro:jointlapdnew} and Lemma \ref{lem:meso-ben-uses}. In order to prove them we will first require the following a priori bound on the tail of $\Phi_N (z, \eta)$.  The constant $B>0$ appearing below is the same constant as in the statement of Proposition \ref{pro:jointlapdnew}.

\bel \label{lem:rough}
There are constants $B >0$ and $C_*>0$ so that for any i.i.d. matrix for any $\eps >0$ and $D>0$, we have
\beq
\label{eq:aprioritail}
\pp\left[ | \Phi_N (z, \eta) | > u \log N \right] \lesssim \e^{ - \frac{ u^2 \log N}{B}} + N^{-D}
\eeq
for all $u \geq 1$, $\eta$'s with $ ( \log N)^{C_*} / N =:\eta_*\leq \eta \leq 1$, and all $N$ large enough. 
\eel
The proof of the above and Proposition \ref{pro:jointlapdnew} are given in the following Section \ref{sec:mainproofs}.

\subsection{Proof of tail bound and joint Laplace asymptotic.} \label{sec:mainproofs} \, We start with the main body of the proof of Lemma~\ref{lem:rough} (tail bound) and then explain the differences to obtain Proposition~\ref{pro:jointlapdnew} (Laplace asymptotic) in Section~\ref{sec:sub1}. Then, in Section~\ref{sec:sub2} we present the proof of some technicalities used in the proofs of Section~\ref{sec:sub1}.

\subsubsection{Proof of Lemma~\ref{lem:rough} and Proposition~\ref{pro:jointlapdnew}.} \label{sec:sub1}

In this section we present the main body of the proofs of Lemma~\ref{lem:rough} and Proposition~\ref{pro:jointlapdnew}, in this order. A key technical input for both proofs will be a combination of Stein's method (to compute the Laplace transform) and cumulant expansion. Due to its length this part is postponed to Section~\ref{sec:sub2}.

\

\proof[Proof of Lemma~\ref{lem:rough}] Let $h(x)$ be an even smooth function function such that $h(x) = 1$ for $|x| \leq 1$, $h(x) = 1/|x|$ for $|x| > 2$, and for $|x| \in [1, 2]$ take $1/|x| \leq h(x) \leq 1$. Let $A = C_1 \log N$ for some large $C_1 \geq 10$, and let 
\beq \label{eqn:meso-d1}
Y := -\int_{\eta}^{N^{10}} N ( \Im [ m_N ( \i u ) - m^z ( \i u ) ] ) \d u = \Phi_N (z, \eta)  + \O ( N^{-5} ) , \qquad X := Y h ( Y/A).
\eeq
where the second equality follows by a simple integration in $u$. Additionally, we define
\beq
\label{eq:defpsi1}
\psi (t) := \ee[ \e^{ t X} ].
\eeq
We consider $t \in [0, L]$ for some fixed $L>0$. We will then use the following bound, whose proof, relying on Stein's method and cumulant expansion, is postponed to Section~\ref{sec:sub2}. For an explicit $B$ (see below \eqref{eq:newasymtot}), we have
\beq
\label{eq:steinbound}
\psi'(t) \leq B \log N ( t + 1 ) \psi (t) + N^{-D}.
\eeq
We point out that the precise value of $B$ does not matter for the proof, and $D>0$ is arbitrarily large. Integrating \eqref{eq:steinbound} in $t$, we thus obtain
\beq
\psi(t) \lesssim \e^{B (t^2 \log N + t \log N ) }.
\eeq
The tail bound \eqref{eq:aprioritail} now follows from Markov's inequality.
More precisely, fix some constant $Q \geq 1$. We want to bound $\pp[ \Phi_N (z, \eta) \geq Q \log N ]$. The key point is that if $\Phi_N (z, \eta) \geq Q \log N$, then we have that (we ignore the difference between $Y$ and $\Phi_N$ here for simplicity)
\beq
X := \Phi_N (z, \eta) h \left( \frac{ \Phi_N (z, \eta) }{C_1 \log N } \right) \geq Q \log N
\eeq
as long as $C_1 > Q$
because the function $x \to x h(x/C_1 \log N)$ is equal to $x$ for $0<x \leq C_1 \log N$ and then if $x \geq C_1 \log N$ we have that $x h(x/{C_1 \log N } ) \geq C_1 \log N$ (the function $s \to s h(s)$ is equal to $s$ for $|s| < 1$ and then smoothly goes to a constant function $s h(s) = 1$ for $s > 2$). Therefore, taking $t \geq 1$,
\beq
\pp\left[ \Phi_N (z, \eta ) \geq Q \log N \right] \leq \pp\left[ X \geq Q \log N \right] \leq \e^{ -t Q \log N } \ee[ \e^{ t X} ] \leq C \e^{- t Q \log N } \e^{ 2 B t^2 \log N}
\eeq
We choose $t = Q / (4B)$ and possibly take $C_1$ even larger so that the computation applies for $t$ this large, showing the desired bound for $\pp\left[ \Phi_N (z, \eta ) \geq Q \log N \right]$. A similar argument applies to $\pp\left[ \Phi_N (z, \eta ) \leq -Q \log N \right]$, concluding the proof. \qed

\

\proof[Proof of Proposition~\ref{pro:jointlapdnew}] Recall the definition $\hat{\eta}_i:=\eta\bm1_{i\le J}+a_{i-J}\bm1_{i>J}$. Consider the $J+K$-dimensional vector ${\bm Y}$ whose entries are defined by
\beq
\label{eq:approx}
Y_i:=\i\int_{\hat{\eta}_i}^{N^{10}} N\langle G^{z_i}(\i u_1)-m^{z_i}(\i u_1)\rangle \, \d u_1=\Phi_N(z_i,\hat{\eta}_i)+\mathcal{O}(N^{-5}) ,
\eeq
where the second equality follows as below \eqref{eqn:meso-d1}.
By \eqref{eq:approx} it follows that it is enough to compute the Laplace transform of ${\bm Y}$ to get the desired result.

Let $\rho(x)$ be a smooth cut-off function such that $\rho(x)=1$ for $|x|\le C_1\log N$ and $\rho(x)=0$ for $|x|\ge C_1\log N+1$. Additionally, in the following we often use the short-hand notation $\rho({\bm Y}):=\prod_i\rho(Y_i)$. We now compute 
\beq
\psi({\bm t}):=\ee\left[e^{\langle {\bm t}, {\bm Y}\rangle}\rho({\bm Y})\right],
\eeq
using Stein's method. For this purpose we need to ensure that the $\Phi_N(z_i,\hat{\eta}_i)$'s, and thus the $Y_i$'s, are not too large with high probability. In fact, we now show that for $j\in [J]$ and $k\in [K]$, we have
\beq
\label{eq:tailbneededa}
\pp\big(\Phi_N(z_j,\eta_j)\ge C_1\log N\big)\le N^{-C_1^2/B}, \qquad\quad \pp\big(\Phi_N(z_{J+k},a_k)\ge C_1\log N\big)\le N^{-D}.
\eeq
The first bound immediately follows from Lemma~\ref{lem:rough}. The second bound follows by noticing
\beq
\label{eq:reqblapnewa}
\ee\left[\e^{t_{J+k} \Phi_N (z_{J+k}, a_k)}\widetilde{\rho}(\Phi_N (z_{J+k}, a_k))\right]\le C_{t_{J+k},\delta,\eps}+N^{-C_2/2},
\eeq
together with a Markov's inequality. Here $\widetilde{\rho}$ is a cut-off function defined analogously to $\rho$ but with $C_1$ replaced with $C_2$, for a sufficiently large constant $C_2$ depending on $C_1$, $t_{J+k}$, and $D$. The bound \eqref{eq:reqblapnewa} readily follows by another application of Stein's method and cumulant expansion, similarly to the proof of \eqref{eq:steinbound}.

Given the a priori bounds in \eqref{eq:tailbneededa}, using again Stein's method together with a cumulant expansion, for any $i\in [J+K]$, we conclude (the proof of this equality is similar to the proof of \eqref{eq:steinbound} and it is postponed to Section~\ref{sec:sub2})
\beq
\label{eq:explforma}
\partial_{t_i}\psi({\bm t})=\left[\sum_i\mathcal{C}_{ji}t_i+{\bm E}_i+
\mathcal{O}_{K,\delta, L, \eps, J}\left(\frac{1}{\sqrt{N\eta_m}}\right)\right]\psi({\bm t})+\mathcal{O}_{K,t_i,\delta, L, \eps, J}(N^{-D}),
\eeq
where ${\bm E}_i=E(z_i,\hat{\eta}_i)$ is from \eqref{eqn:bmE-def} and $\mathcal{C}_{ji}=\mathcal{C}(z_i,\hat{\eta}_i,z_j,\hat{\eta}_j)$ is from \eqref{eqn:def-covar}. Finally, integrating \eqref{eq:explforma} we readily conclude the proof of \eqref{eq:biglap}. \qed

\subsubsection{Proof of technical results from Section~\ref{sec:sub1}: Stein's method}
\label{sec:sub2}

\proof[Proof of \eqref{eq:steinbound}] Recall the definitions of $X, Y$ from \eqref{eqn:meso-d1}, the one of $\psi(t)$ from \eqref{eq:defpsi1}, and that $t\in [0,L]$. We start computing (recall that $A=C_1\log N$)
\beq
\label{eq:steina}
\partial_t\psi(t)=\ee Yh(Y/A) e^{tX} =\i \ee \int_\eta^{N^{10}} Nh(Y/A)\langle G^z(\i u)-m^z(\i u)\rangle \, e^{tX}\, \d u,
\eeq
where in the second equality we used that $\Im m_N=\Im \langle G^z\rangle=-\i \langle G^z\rangle$ (and a similar relation for its deterministic approximation). In the following we will often use the local laws (for deterministic unit vectors ${\bm x}, {\bm y}$ and deterministic matrices $A$ with $\lVert A\rVert\lesssim 1$)
\beq
\begin{split}
\label{eq:llawlogscale}
| \langle \bx , (G^z (\i \eta) - M^z (\i \eta) ) \by \rangle| &\le \frac{ (\log N)^2+\bm1_{\{\eta\ge N^{-10\xi}\}}N^\xi}{\sqrt{N\eta}}, \\
| \langle (G^z (\i \eta) - M^z (\i \eta)) A \rangle | &\leq \frac{(\log N)^2+\bm1_{\{\eta\ge N^{-10\xi}\}}N^\xi}{N \eta},
\end{split}
\eeq
with overwhelming probability uniformly in $\eta\ge \eta_*$, for
any small $\xi>0$. The local laws in \eqref{eq:llawlogscale} follows directly from \cite[Theorem 3.1]{cipolloni2025optimal2}.

Denote $G:=G^z(\i u)$ and $M=M^z(\i u)$ (recall its definition from \eqref{eq:bigM}). To compute the second line of \eqref{eq:steina}, we use the equation
\beq
\label{eq:useqa}
\langle G-M\rangle=-\langle \underline{WG}\mathcal{A}\rangle+\langle G-M\rangle\langle (G-M)\mathcal{A}\rangle, \qquad\quad \mathcal{A}=\mathcal{A}(z,\i u):=\frac{M^z(\i u)}{1-\langle [M^z(\i u)]^2\rangle},
\eeq
which follows by the definitions of $M$ and of the resolvent $G$. Here $W:=H^0(X)$ (see its definition in \eqref{eqn:Hermit-def}) and $\underline{WG}:=WG+\langle G\rangle G$. We also point out that $\lVert \mathcal{A}(z,\i u)\rVert\lesssim 1$, since we are in the bulk of the spectrum, i.e. $|z|\le r<1$. Using \eqref{eq:useqa} in \eqref{eq:steina}, we thus obtain
\beq
\label{eq:stein2a}
\partial_t\psi(t)=-\i \ee \int_\eta^{N^{10}} Nh(Y/A)\langle \underline{WG}\mathcal{A}\rangle \, e^{tX}\, \d u+\mathcal{O}\left(\frac{(\log N)^4}{N\eta_*}\right)\psi(t)+\mathcal{O}(N^{-D}),
\eeq
for any sufficiently large $D>0$, where “sufficiently large" depends on $C_1, L$, as a consequence of $X\le 2C_1\log N$ deterministically.
Here and throughout this proof $\xi>0$ is an arbitrary small constant.
The second term in \eqref{eq:stein2a} comes from $\langle G-M\rangle\langle (G-M)\mathcal{A}\rangle$ after using the averaged local law from \eqref{eq:llawlogscale} and that $h$ is uniformly bounded in $N$.

In the following we will often use the notation
\beq
\sum_{ab}:=\sum_{1\le a\le N,  N+1\le b\le 2N}+\sum_{N+1\le a\le 2N,  1\le b\le N}.
\eeq
To compute the first term in the right-hand-side of \eqref{eq:stein2a} we now perform cumulant expansion in $W$:
\beq
\begin{split}
\label{eq:cumexpa}
\ee \int_\eta^{N^{10}} Nh(Y/A)\langle \underline{WG}\mathcal{A}\rangle e^{tX}\, \d u&=\ee \int_\eta^{N^{10}} N\ee\big[e^{tX}h(Y/A)\big]\ee\langle \underline{WG}\mathcal{A}\rangle \, \d u \\
&\quad+ \ee \int_\eta^{N^{10}} \sum_{ab}\langle \Delta^{ab}G\mathcal{A}\rangle \partial_{ba}\big(e^{tX}h(Y/A)\big) \, \d u \\
&\quad+N\ee \int_\eta^{N^{10}} \sum_{k=2}^R \sum_{{\bm \alpha}\in\{ab,ba\}^k} \frac{\kappa(ab,{\bm \alpha})}{k!}\partial_{{\bm \alpha}} \big(\langle \Delta^{ab}G\mathcal{A}\rangle e^{tX}h(Y/A)\big)\, \d u \\
&\quad+\mathcal{O}(N^{-D}),
\end{split}
\eeq
for a very large $R>0$ depending on the desired $D>0$, and thus on $C_1, L$. In \eqref{eq:cumexpa}, $\partial_{ab}:=\partial_{W_{ab}}$ denotes the directional derivative in the direction $W_{ab}$, $\partial_{{\bm \alpha}}:=\partial_{\alpha_1}\dots  \partial_{\alpha_k}$, with $\alpha_i\in \{ab,ba\}$, $\Delta^{ab}$ is the matrix with the $(a,b)$-entry being the only non-zero one, and $\kappa(ab,{\bm \alpha})$ denotes the $k+1$-th cumulant of the random variables $W_{ab}, W_{\alpha_1},\dots, W_{\alpha_k}$, with ${\bm \alpha}:=(\alpha_1,\dots, \alpha_k)$. Notice that in \eqref{eq:cumexpa} we truncated the cumulant expansion at the level $R$ for a very large $R>0$ (see e.g. \cite[Proposition 3.2]{erdHos2019random}). For this we used that $|h^{(l)}|\lesssim 1$ for any $l\ge 0$, that $|\kappa(ab,{\bm \alpha})|\lesssim N^{-(R+2)/2}$, for $k\ge R+1$, that $\max_{a,b}|G_{ab}|\lesssim 1$ and $|Y|\lesssim  N^\xi$ with overwhelming probability by \eqref{eq:llawlogscale}, and that $e^{tX}\le N^{2C_1L}$ deterministically. Note that using the convention $\partial_{ab}G=-G\Delta^{ab}G-\bm1_{\beta=1}G\Delta^{ba}G$, the expansion \eqref{eq:cumexpa} holds both for the real and complex cases.

We now start computing the various terms in the RHS of \eqref{eq:cumexpa}. We start with the first line. For such term we have (this follows directly from the proof of \cite[Proposition 3.3]{cipolloni2023central} and \cite[Proposition 3.3]{cipolloni2021fluctuation} in the complex and real case, respectively)
\beq
\label{eq:expcompnew}
-\i \ee \int_\eta^{N^{10}} N\langle \underline{WG}\mathcal{A}\rangle={\bm E}+\mathcal{O}_L\left(\frac{(\log N)^{10}}{\sqrt{N\eta_*}}\right)\psi(t),
\eeq
where ${\bm E}={\bm E}(z,\eta)$ is from \eqref{eqn:bmE-def}. We point out that strictly speaking the results from \cite{cipolloni2021fluctuation, cipolloni2023central} used local laws with $N^\xi$-error, as in Theorem~\ref{thm:ll-1}, but the same proof applies verbatim using the local laws in \eqref{eq:llawlogscale}. We now compute the term in the second line of \eqref{eq:cumexpa} (here we neglect $Y/A$ from the argument of $h$ to shorten the notation, and we denote $G_i^z:=G^z(\i u_i)$, $\mathcal{A}_1:=\mathcal{A}(z,\i u_1)$):
\beq
\begin{split}
\label{eq:intstepn}
&\ee \int_\eta^{N^{10}} \sum_{ab}\langle \Delta^{ab}G\mathcal{A}\rangle \partial_{ba}\big(e^{tX}h\big) \, \d u \\
&=-\tilde{\sum}_{kl}\ee \int\int_\eta^{N^{10}} \bigg\{ \left(th^2+\frac{tYhh'+h'}{A}\right)e^{tX}\big[\langle G_1^z\mathcal{A}_1 E_k (G_2^z)^2 E_l\rangle\\&\quad \qquad+\bm1_{\beta=1}\langle G_1^z\mathcal{A}_1 E_k [(G_2^z)^2]^\mathfrak{t} E_l\rangle\big] \bigg\}\,\d u_1\d u_2,
\end{split}
\eeq
where we recall that $\mathfrak{t}$ denotes the transpose of a matrix. To compute the leading order asymptotic of \eqref{eq:intstepn}, we use the local laws of Lemma \ref{lem:2G-large} and Proposition \ref{prop:2G-law}. In particular, the quantity $M(z_1,\eta_1,B_1, z_2,\eta_2)$ appearing there can be written as $-\partial_{\eta_1}\mathcal{V}(z_1,\eta_1,z_2,\eta_2)$ by explicit computations (see \cite[Eqs. (7.12) and (7.27)]{cipolloni2026eigenvalues} for analogous computations) in the special case $B_1=\mathcal{A}_1E_k$, $B_2=E_l$. 
 Using this observation, the facts that $\partial_{\i u_2}G_2^z=(G_2^z)^2$, $(G^z)^\mathfrak{t}=G^{\overline{z}}$, and performing explicitly the $u_2$-integration before using Lemma \ref{lem:2G-large} and Proposition \ref{prop:2G-law}, we thus obtain
\beq
\begin{split}
\label{eq:varcompnew}
&\ee \int_\eta^{N^{10}} \sum_{ab}\langle \Delta^{ab}G\mathcal{A}\rangle \partial_{ba}\big(e^{tX}h\big) \, \d u \\
&=\i\ee \left(th^2+\frac{tYhh'+h'}{A}\right)e^{tX}\big[\mathcal{V}(z,\eta,z, \eta)+\bm1_{\beta=1}\mathcal{V}(z,\eta,\overline{z}, \eta)\big]+\mathcal{O}_L\left(\frac{(\log N)^C}{N\eta_*}\right)\psi(t),
\end{split}
\eeq
where in the estimate of the error term we used $|h|+|h'|\le 1$ and $|Yh|\le 2A$ deterministically.

We are now left with the terms in the third line of \eqref{eq:cumexpa}. In the estimates of these terms the following bounds will be used repeatedly
\beq
\label{eq:usb}
\big| Y \partial_{ab}^lh(Y/A)\big|\lesssim (\log N)^2, \qquad\quad l\ge 0,
\eeq
which holds with overwhelming probability, where we used $|Y|\lesssim (\log N)^2$, by \eqref{eq:llawlogscale}, that $|\partial_{ab}^k Y|\lesssim 1$ as a consequence that for $k\ge 1$ we can express $\partial_{ab}^k Y$ as products of various entries of $G$ evaluated at $u=\eta$ (i.e. the integral can be performed explicitly), and that $h$ and all its derivatives are bounded uniformly in $N$. Here by $\partial_{ab}^l$ we denote any combination of $l$-derivatives $\partial_{ab}, \partial_{ba}$. By a simple power counting for $k\ge 4$, using \eqref{eq:usb}, we see that these terms give a contribution
\beq
\label{eq:powerc}
N N^2 N^{-(k+1)/2} N^{-1} (\log N)^{2k}\psi(t)+N^{-D}\le N^{-1/2}(\log N)^{2k}\psi(t)+N^{-D},
\eeq
where the first $N$ comes from the pre-factor in the third line of \eqref{eq:cumexpa}, the $N^2$ comes from the $ab$-summation, the factors $N^{-(k+1)/2}$ and $ N^{-1}$ come from the bound on the cumulants and the normalization of the trace, respectively, and $(\log N)^{2k}$ comes from \eqref{eq:usb}. We point out that the remaining $\eta$-integral in the regime $\eta\gg 1$ does not cause any problem since $\lVert \mathcal{A}\rVert\lesssim (\eta+1)^{-1}$ and hence it can always be performed at the price of a negligible $\log N$-factor which we omit in the bounds.

Next, we consider the third order cumulants, i.e. with $k=2$, which we will also show to be negligible. For brevity we present the estimate of two representative terms, all the other ones can be estimated analogously. Notice that the more off-diagonal resolvent entries are present the smaller is the contribution (as a consequence of of \eqref{eq:llawlogscale}). For this reason, we now consider terms where there is only one diagonal resolvent entry.  When both derivatives hit $G$ we have (up to constants)
\beq
\begin{split}
\label{eq:third1a}
\frac{1}{N^{3/2}}\ee \int_\eta^{N^{10}} \sum_{ab} G_{aa}G_{bb}(G\mathcal{A})_{ab} \, e^{tX}\, \d u &=\frac{1}{N^{3/2}}\ee \int_\eta^{N^{10}} m^2 \langle \bm1, G\mathcal{A} \bm1\rangle \, e^{tX}\, \d u \\
&\quad+\mathcal{O}_L\left(\frac{N^\xi}{\sqrt{N}}+\frac{(\log N)^2}{N\sqrt{\eta_*}}\right)\psi(t)+\mathcal{O}_L(N^{-D}) \\
&=\mathcal{O}_L\left(\frac{N^\xi}{\sqrt{N}}+\frac{(\log N)^2}{N\sqrt{\eta_*}}\right)\psi({\bm t})+\mathcal{O}_L(N^{-D}),
\end{split}
\eeq
where $\bm1:=(1,\dots, 1)\in\cc^{2N}$, $D$ is any sufficiently large constant depending on $C_1, L$, and to estimate the errors we used the local law in the first line of \eqref{eq:llawlogscale}. When, instead, both derivatives hit the same $Y$ we have (up to constants, and using the same notation as described in the lines above \eqref{eq:intstepn})
\beq
\begin{split}
\label{eq:third2a}
\frac{1}{N^{3/2}}\ee \int\int_\eta^{N^{10}} \sum_{ab} (G_1\mathcal{A}_1)_{ab}(G_2)_{bb}(G_2^2)_{aa} \, e^{tX}\, \d u&=\frac{\int_\eta^{N^{10}} mm'}{N^{3/2}}\ee \int_\eta^{N^{10}} \langle \bm1, G\mathcal{A} \bm1\rangle \, e^{tX}\, \d u \\
&\quad+\mathcal{O}_L\left(\frac{N^\xi}{\sqrt{N}}+\frac{(\log N)^2}{N\sqrt{\eta_*}}\right)\psi(t)+\mathcal{O}_L(N^{-D}) \\
&=\mathcal{O}_L\left(\frac{N^\xi}{\sqrt{N}}+\frac{(\log N)^2}{N\sqrt{\eta_*}}\right)\psi(t)+\mathcal{O}_L(N^{-D}),
\end{split}
\eeq
This shows that the third order terms ($k=2$) are also negligible. Similarly, it is possible to show that the only non-negligible term for $k=3$ is the term when one derivative $\partial_{ba}$ hits $G$ and two derivatives $\partial_{ab}\partial_{ba}$ hit the same $Y$. This term gives a contribution (here we omit the $\kappa_4$-factor)
\beq
\label{eq:4thorderterm}
\begin{split}
&\ee\int\int_{\eta}^{N^{10}}\left(th^2+\frac{tYhh'+h'}{A}\right)\frac{e^{tX}}{N^2}\sum_{ab}(G_1)_{aa}(G_1\mathcal{A}_1)_{bb} (G_2^2)_{aa} (G_2)_{bb}\, \d u_1\d u_2 \\
&=\frac{m^4}{4}\ee\left[\left(th^2+\frac{tYhh'+h'}{A}\right)e^{tX}\right]+\mathcal{O}_L\left(\frac{(\log N)^2}{\sqrt{N\eta_*}}\psi(t)+N^{-D}\right),
\end{split}
\eeq
where we used that $(M_1\mathcal{A}_1)_{bb}=\partial_{\i u_1}m(\i u_1)$, which immediately follows from \eqref{eq:bigM}--\eqref{eq:mde}, and the local law in the first line of \eqref{eq:llawlogscale}. All the other terms for $k=3$ give a contribution $(N\eta_*)^{-1/2}\psi(t)+N^{-D}$.

Combining \eqref{eq:expcompnew}, \eqref{eq:varcompnew}, \eqref{eq:powerc}--\eqref{eq:4thorderterm}, we thus conclude 
\beq
\begin{split}
\label{eq:newasymtot}
\psi'(t)&= {\bm E}\ee\big[he^{tX}\big]+\big[\mathcal{V}(z,\eta,z, \eta)+\bm1_{\beta=1}\mathcal{V}(z,\eta,\overline{z}, \eta)\big]\ee\left[\left(th^2+\frac{tYhh'+h'}{A}\right)e^{tX}\right]\\
&\quad+\frac{m^4\kappa_4}{4}\ee\left[\left(th^2+\frac{tYhh'+h'}{A}\right)e^{tX}\right]+\mathcal{O}_L\left(\frac{(\log N)^2\psi(t)}{\sqrt{N\eta_*}}+N^{-D}\right),
\end{split}
\eeq

Finally, using that for some constant $\tilde{C}>0$ we have $|\mathcal{V}(z,\eta,z,\eta)|\le (1+\bm1_{\beta=1})\log N+\tilde{C}$ from Proposition~\ref{prop:parameter-properties}, a similar estimate for $\mathcal{V}(z,\eta,\overline{z},\eta)$, $|{\bm E}|\le \tilde{C}+\bm1_{\beta=1}\log N$, and $|h|+|h'|\le 1$ and $|Yh|\le 2A$, we estimate the various terms in \eqref{eq:newasymtot} and obtain \eqref{eq:steinbound}. The choice $B=6+12\bm1_{\beta=1}$ does the job. \qed

\

\proof[Proof of \eqref{eq:explforma}] The proof of this equality is very similar to the proof of \eqref{eq:newasymtot} (and hence of \eqref{eq:steinbound}), we thus explain more carefully only the main differences and omit very similar details. We point out that for the purpose of this proof we can use directly the local laws in Theorem~\ref{thm:ll-1} instead of \eqref{eq:llawlogscale}, as in this proof we have $\eta_M\ge N^{-1+\eps}$. 

We start by computing the derivative with respect to $t_1$ as the other cases are completely analogous. Denote $G_i:=G^{z_i}(\i u_i)$, and compute
\beq
\label{eq:stein}
\partial_{t_1}\psi({\bm t})=\ee Y_1 e^{\langle {\bm t}, {\bm Y}\rangle}\rho({\bm Y}) =\i \ee \int_{\eta_1}^{N^{10}} N\langle G_1-m_1\rangle \, e^{\langle {\bm t}, {\bm y}\rangle}\rho({\bm Y})\, \d u_1.
\eeq
To compute the RHS of \eqref{eq:stein}, we proceed similarly to \eqref{eq:useqa}--\eqref{eq:stein2a} and obtain (here $\mathcal{A}_1:=\mathcal{A}(z_1,\eta_1)$ is from \eqref{eq:useqa})
\beq
\label{eq:stein2}
\partial_{t_1}\psi({\bm t})=-\i \ee \int_{\eta_1}^{N^{10}} N\langle \underline{WG_1}\mathcal{A}_1\rangle \, e^{\langle {\bm t}, {\bm Y}\rangle}\rho({\bm Y})\, \d u_1+\mathcal{O}\left(\frac{N^\xi}{N\eta_M}\right)\psi({\bm t})+\mathcal{O}(N^{-D}),
\eeq
for any sufficiently large $D>0$, where “sufficiently large" depends on $L, J, K$, and $C_1$. Here and throughout this proof $\xi>0$ is an arbitrary small constant.

To compute the first term in the RHS of \eqref{eq:stein2} we perform again a cumulant expansion in $W$ (see the paragraph below \eqref{eq:cumexpa} for a detailed explanation and for the defintion of the various notations):
\beq
\begin{split}
\label{eq:cumexp}
\ee \int_{\eta_1}^{N^{10}} N\langle \underline{WG_1}\mathcal{A}_1\rangle e^{\langle {\bm t}, {\bm y}\rangle} & \rho({\bm Y})\, \d u_1=\psi({\bm t})\ee \int_{\eta_1}^{N^{10}} N\langle \underline{WG_1}\mathcal{A}_1\rangle\, \d u_1  \\
&\quad+ \ee \int_{\eta_1}^{N^{10}} \sum_{ab}\langle \Delta^{ab}G_1A_1\rangle \partial_{ba}\big(e^{\langle {\bm t}, {\bm y}\rangle}\rho({\bm Y})\big)\, \d u_1 \\
&\quad+N\ee \int_{\eta_1}^{N^{10}} \sum_{k=2}^R \sum_{{\bm \alpha}\in\{ab,ba\}^k} \frac{\kappa(ab,{\bm \alpha})}{k!}\partial_{{\bm \alpha}} \big(\langle \Delta^{ab}G_1\mathcal{A}_1\rangle e^{\langle {\bm t},{\bm Y}\rangle}\rho({\bm Y})\big)\, \d u_1 \\
&\quad+\mathcal{O}(N^{-D}).
\end{split}
\eeq

In the following we neglect all the terms when one or more derivatives hit $\rho$ as they give a negligible contribution. More precisely, if $K=0$ then these terms are all of order $\mathcal{O}(N^{JLC_1-C_1^2/B})$ as a consequence of the first bound in \eqref{eq:tailbneededa}. This error is negligible by choosing $C_1$ large in terms of $L$ and $J$. If instead $K\ge 1$, then there are two cases: i) a derivative hits $\rho(Y_i)$ for $i\in [J]$; ii) a derivative hits $\rho(Y_i)$ for $i\in [J+1, J+K]$. In case ii), this results to a negligible contribution $N^{-D}$, as a consequence of the second bound in \eqref{eq:tailbneededa}. In case i), we get a bound 
\beq
\label{eq:somerrb}
N^{JLC_1-C_1^2/B}\prod_{i=1}^K\ee\left[\e^{Kt_i \Phi_N (z_{J+i}, a_i)}\rho(\Phi_N (z_{J+i}, a_i))\right]=\mathcal{O}_{t_i,K,\delta,\eps}\left(N^{JLC_1-C_1^2/B}\right),
\eeq
which is also negligible. We point out that in \eqref{eq:somerrb} we also used
\begin{equation}
\label{eq:addresneed}
\ee\left[\e^{t_i \Phi_N (z_{J+i}, a_i)}\rho(\Phi_N (z_{J+i}, a_i))\right]\le C_{t_i,\delta,\eps}+N^{-D},
\end{equation}
which follows analogously to \eqref{eq:reqblapnewa}.

We now compute the right-hand-side of \eqref{eq:cumexp}. Proceeding analogously to the proof of \eqref{eq:steinbound} we see that all the terms for $k=2$ and $k\ge 4$ are negligible. The only terms contributing from \eqref{eq:cumexp} are those in the first two lines and the one in the third line when $k=3$ and when one derivative $\partial_{ba}$ hits $G_1$ and two derivatives $\partial_{ab}\partial_{ba}$ hit the same $Y_i$ (as in the proof of \eqref{eq:steinbound}). Hence, following the computations in \eqref{eq:expcompnew}--\eqref{eq:newasymtot} verbatim we obtain \eqref{eq:explforma}, concluding the proof. \qed

\subsection{Consequences of previous arguments} 

\label{sec:impl}

By examining the arguments in the previous section we are able to conclude the following. Note that for a single $\Phi_N (z, \eta)$ it improves over the statement of Proposition \ref{pro:jointlapdnew} by allowing $\eta \geq ( \log N)^{C_*} / N$. 
\bel
\label{lem:impleas}
Let $\eta_*$ from Lemma~\ref{lem:rough}. There is a $B>0$ so that the following holds. For any $C_1 \geq 100,  D \geq 1, L>1$, we have that 
\beq
\begin{split}
\label{eq:rolap}
\ee[ \e^{ \lambda \Phi_N (z, \eta ) }  \1_{ \{ | \Phi_N (z, \eta) | \leq C_1 \log N \} } ] &= \e^{ \frac{ \lambda^2}{2} \cC (z, \eta, z, \eta)+ \lambda \bmE(z, \eta ) } \left[ 1 +\O ( ( \log N)^{100} / (N \eta)^{1/4} ) \right] \\
&\quad+ \O (N^{10+ B L^2+ L C_1 - C_1^2/B }+N^{-D}),
\end{split}
\eeq
uniformly in $\eta_* \leq \eta \leq 1$ and $|\lambda| \leq L$.
\eel
\proof
This immediately follows in a completely analogous way to the proofs of Lemma~\ref{lem:rough} and Proposition~\ref{pro:jointlapdnew}, for $J=1$ and $K=0$, as a consequence of the fact that on the imaginary axis we use local laws with a logarithmic precision $(\log N)^2$, see e.g. \eqref{eq:llawlogscale}.
We omit the details for brevity.
\qed

\subsubsection{Proof of Lemma \ref{lem:meso-ben-uses}} 
\label{sec:proof-ben-uses}
We start with \eqref{eqn:meso-tail-bd}. Let $B$ be the constant from Lemma \ref{lem:rough} and let $A \geq 100 (B+1) (L+10) + 100$. Then,
\beq
\pp\left[  \Phi_N (z, \eta)  > u \log N \right] \leq \pp\left[  \Phi_N (z, \eta)  > u \log N , |\Phi_N (z, \eta) | \leq A \log N \right] + \e^{ -10 L^2 \log N }
\eeq
The first term on the RHS is then bounded using Markov's inequality and Lemma \ref{lem:impleas}. A similar argument holds for $\pp [ \Phi_N (z, \eta) < - u \log N ]$, completing the proof of \eqref{eqn:meso-tail-bd}. The estimate \eqref{eqn:meso-1pt-mgf} follows from Lemma \ref{lem:impleas} and using \eqref{eqn:meso-tail-bd} to compare the LHS of \eqref{eqn:meso-1pt-mgf} to the same quantity with $A$ replaced by $C_1$ for some sufficiently large $C_1$. The proof of \eqref{eqn:meso-uni-bd} is similar. \qed

\section{GFCT results} \label{sec:gfct}

In this section we will prove various \emph{Green's function comparison theorems} (GFCT) for the spectral statistics we are interested in. For an extensive discussion of such theorems see \cite{erdHos2017dynamical}. We will then use the GFCT and the previous results of the paper to prove most of the main results.

\subsection{Preliminaries} \label{sec:gfct-prelim}

We begin by collecting some notation. For a function $F$ on the space of $N \times N$ matrices with complex entries $X$ we will denote the usual Wirtinger derivatives by,
\beq
\del_{X_{ab}} F(X) := \frac{1}{2} \left( \del_{ \Re[X_{ab}]} - \i \del_{\Im[X_{ab}]} \right) F(X) , \quad \del_{\bar{X}_{ab}} F(X) := \frac{1}{2} \left( \del_{ \Re[X_{ab}]} + \i \del_{\Im[X_{ab}]} \right) F(X)
\eeq
In the reminder of this section we will use the following minor abuse of notation: when $X$ is a random matrix will write,
\beq
\del_{X_{ab}} F(X) := ( \del_{M_{ab}} F (M) ) \big\vert_{M = X} .
\eeq
The following matrix will arise in various Taylor expansions. For a parameter $\theta \in \cc$ and indices $a,b \in [\![1, N ]\!]$ and an $N \times N$ matrix $X$ we define the matrix $X(\theta_{ab})$ by,
\beq \label{eqn:X-theta-def}
X(\theta_{ab} )_{ij} := \begin{cases} X_{ij}  & (i, j) \neq (a, b), \\ \theta & (i,j) = (a, b). \end{cases}
\eeq
For a matrix $X$, recall the Hermitization $H^z (X)$ defined in \eqref{eqn:Hermit-def}. We recall also the notation $\Phi_N (z, \eta ; X)$ in \eqref{eqn:Phi-def}. The following allows us to compare various quantities involving $X(\theta_{ab} )$ to those of $X$ itself.   Similar arguments have appeared numerous times in the literature.
\bep
Let $\delta \in (0, \frac{1}{100} )$ and $0 < r < 1$. Let $X$ be an i.i.d. matrix. For $\theta \in \cc$ and $a, b \in [\![1, N]\!]$ define the resolvent by, 
\beq 
G^z (\theta_{ab} ,  \i \eta)  := (H^z (X(\theta_{ab} ) ) - \i \eta)^{-1}.
\eeq
For any $\eps >0$ we have uniformly for $ z\in \dd_r$ and $\eta \in [N^{-1-\delta}, 1]$ that 
\beq \label{eqn:theta-g-est}
\sup_{ |\theta| \leq N^{1/10 -1/2} } | G^z_{ij} (\theta_{ab} , \i \eta ) - G^z_{ij} (X_{ab}, \i \eta ) | \leq N^{-1/4}, \quad \sup_{ |\theta| \leq N^{1/10 -1/2} } | G^z_{ij} (\theta_{ab} , \i \eta ) | \leq N^{\eps+\delta},
\eeq
and 
\beq \label{eqn:theta-phi-est}
\sup_{ |\theta| \leq N^{1/10-1/2} } 
| \Phi_N ( z , \eta ; X(\theta_{ab} ) )  - \Phi_N ( z, \eta ; X ) | \leq N^{-1/4},
\eeq
and
\beq \label{eqn:theta-mN-est}
\sup_{ | \theta| \leq N^{1/10-1/2}} \left| (N \eta ) \langle G^z ( \theta_{ab} , \i \eta ) - G^z ( \i \eta ) \rangle \right| \leq N^{-1/4},
\eeq
with overwhelming probability.
\eep
\proof We write $G^z(\i \eta) := (H^z(X) -\i \eta)^{-1}$ and let $\eta_\delta := N^{-1-\delta}$.  By \cite[Lemma 2.1]{bauerschmidt2017local} and Theorem \ref{thm:ll-1} we have,
\beq \label{eqn:theta-est-a1}
\sup_{i,j} \sup_{ \eta \in [ \eta_\delta,1 ]} |G_{ij}^z (\i \eta) | \leq N^{\eps+\delta} 
\eeq
with overwhelming probability. Let $\Delta$ be the $2 N \times 2N$ self-adjoint matrix whose only non-zero entry in the upper triangular part is $\Delta_{a,b+N} := (\theta - X_{ab} )$. By iterating the resolvent identity and using the deterministic estimate $|G^z_{ij}( \theta_{ab}, \i \eta ) | \leq \eta^{-1} $ we see that for $|\theta| \leq N^{1/10-1/2}$ and $\eta \geq \eta_\delta$, we have
\begin{align}
    G^z_{ij} ( \theta_{ab}, \i \eta  ) - G^z_{ij} ( \i \eta ) = \sum_{n=1}^{10} \langle \delta_i , [G^z ( \i \eta ) \Delta ]^n G^z ( \theta_{ab} , \i \eta ) \delta_j \rangle + \O ( N^{-1/2} ),
\end{align}
with overwhelming probability. The first estimate of \eqref{eqn:theta-g-est} follows from the above and \eqref{eqn:theta-est-a1}. The second estimate of \eqref{eqn:theta-g-est} now follows immediately from \eqref{eqn:theta-est-a1}. 

For a general non-singular matrix $A(s)$ depending in a differentiable way on a parameter $s \in \rr$ we have
\beq \label{eqn:det-deriv} 
\frac{\d}{\d s} \log \det (A(s) ) = \tr \left( A(s)^{-1} \del_s A(s) \right) .
\eeq
The estimate \eqref{eqn:theta-phi-est} then follows from integrating the difference along an appropriate contour connecting $\theta$ and $X_{ab}$ in the complex plane, and using \eqref{eqn:theta-g-est} to estimate the derivative. In particular, the derivative involves only two resolvent entries $G^z_{a,b+N} ( \varphi_{ab} , \i \eta)$ and $G^z_{b+N, a} (\varphi_{ab} , \i \eta )$ for appropriate $\varphi \in \cc$. 

The estimate \eqref{eqn:theta-mN-est} follows in a similar fashion to \eqref{eqn:theta-phi-est} using
\beq
\frac{\d }{\d s} \tr (A^{-1} (s) ) = \tr \left(  A^{-2} (s) \del_s A(s) \right)
\eeq
as well as the Ward identity \eqref{eqn:ward} and Cauchy-Schwarz. 
\qed

\

We also require the following deterministic lemmas which allow us to smooth out the indicator function of the event  $\{ \lambda^z_1 \leq \eta \}$ for appropriate $\eta$. 

\bel \label{lem:deterministic-small} Let $X$ be any matrix and $\lambda_i^z$ be the eigenvalues of $H^z (X)$. 
For any $\eta >0$, we have that if $\lambda_1^z \leq \eta$, then 
\beq
N \eta \Im[ m_N^z ( \i \eta ) ] \geq \frac{1}{4} .
\eeq
For $N \geq 3$ the following statement holds. 
Let $\eta \leq N^{-1}$, $\eta_1 = \eta / ( \log N)^3$, and let $\eta_2 = \frac{ \log N}{N}$. If $\lambda_1^z > \eta$ then, 
\beq
(N \eta_1) \Im[ m^z_N ( \i \eta_1) ] \leq \frac{ \Im[ m_N^z ( \i \eta_2) ]}{( \log N)^3}.
\eeq
Therefore, the following statement holds for $N \geq 3$. Let $0 <\nu \leq N^{-1} ( \log N)^{-4}$ and $\eta_2 = \frac{ \log N}{N}$. Then,
\beq
\Im[ m_N^z ( \i \eta_2) ] \leq  \log N  \mbox{ and } (N \nu ) \Im[ m_N^z ( \i \nu) ] \geq \frac{1}{ \log N} \quad \implies \quad \lambda_1^z \leq ( \log N)^3 \nu .
\eeq
\eel
\proof For the first estimate,
\beq
2 N \eta \Im[ m_N^z ( \i \eta ) ] = \sum_i \frac{ \eta^2}{ ( \lambda_i^z)^2 + \eta^2 } \geq 2 \frac{\eta^2}{ ( \lambda_1^z)^2 + \eta^2} \geq 1,
\eeq
using $\lambda_1^z \leq \eta$ in the last inequality. For the second estimate,
\begin{align}
    (2 N \eta_1) \Im[ m_N^z ( \i \eta_1) ] &= \sum_i \frac{ ( \eta_1)^2}{ ( \lambda_i^z)^2 + \eta_1^2} 
   \leq  2 \sum_i \frac{ ( \eta_1)^2}{ ( \lambda_i^z)^2 + \eta^2} = 4 \frac{ (\eta_1)^2}{\eta^2} (N \eta) \Im[ m_N^z ( \i \eta)] \notag\\
   &\leq  4 \frac{(\eta_1)^2}{(\eta)^2} (N \eta_2 ) \Im[ m_N^z ( \i \eta_2)] \leq \frac{ \Im[ m_N^z ( \i \eta_2)]}{( \log N)^3}.
\end{align}
In the first inequality we used $2 ( \lambda_i^z)^2 \geq ( \lambda_i^z)^2 + \eta^2$. In the first inequality on the second line we used that $y \to y \Im[m_N^z ( \i y)]$ is an increasing function. The final statement of the lemma is clear. \qed

\bel \label{lem:smoothed-indicator} In the set-up of Lemma \ref{lem:deterministic-small} the following holds. 
Let $\eta = N^{-1-\delta}$ for some $\delta >0$. Let $f : \rr \to \rr$ be a smooth function such that $f(x) =0$ for $x < \frac{1}{20}$ and $f(x) = 1$ for $x > \frac{1}{10}$. Then, 
\beq
\label{eq:boundsssv}
\1_{ \{ \lambda_1^z \leq \eta \}} \leq f (N \eta \Im[ m_N^z ( \i \eta) ] ) \leq \1_{ \{ \lambda_1^z \leq ( \log N)^3 \eta\}} + \1_{ \{ \Im[m_N^z ( \i N^{-1} \log N )] > \log N \} }.
\eeq
\eel
\proof The bounds in \eqref{eq:boundsssv} immediately follow from Lemma \ref{lem:deterministic-small}. \qed

\

Finally we collect here estimates for the derivatives with respect to matrix entries of the various quantities we will need to consider. 
\bel
Fix $\delta \in (0, \frac{1}{100} )$ and $ 0 < r < 1$. Let $X$ be an i.i.d. matrix and $X(\theta_{ab})$ as above. For any $k, l \geq 0$ we have with overwhelming probability uniformly for $z \in \dd_r$ and $\eta \in [ N^{-1-\delta}, 1]$
\begin{align} \label{eqn:del-phi-est}
\sup_{|\theta| \leq N^{1/10-1/2}} \left|  \del_{X_{ab}}^k \del_{\bar{X}_{ab}}^l \Phi_N (X(\theta_{ab} ) ,z, \eta) \right| \leq N^{\eps+(k+l) \delta } ,
\end{align}
and
\begin{align} \label{eqn:del-G-est} 
    \sup_{|\theta| \leq N^{1/10-1/2}} \left|  \del_{X_{ab}}^k \del_{\bar{X}_{ab}}^l N \eta \langle G^z ( \theta_{ab} , \i \eta ) \rangle  \right| \leq N^{\eps + (k+l) \delta}.
\end{align}
\eel
\proof
As functions on the space of $N \times N$ matrices we have the formulas
\begin{align}
\del_{X_{ab}} \Phi_N (z, \eta ; X ) &= G^{z}_{a,b+N} (\i \eta) , \quad \del_{\bar{X}_{ab}} \Phi_N (z, \eta ; X) = G^z_{b+N,a} ( \i \eta) \notag\\
\del_{X_{ab}} G^z_{ij}  &= G^z_{ia} G^z_{b+N,j} , \quad \del_{\bar{X}_{ab}} G^z_{ij} = G^z_{i,b+N} G^z_{aj}
\end{align}
The estimate \eqref{eqn:del-phi-est} follows from this and the second estimate of \eqref{eqn:theta-g-est}. For the second estimate we see from the above estimates that,
\begin{align}
     &  \left|  \del_{X_{ab}}^k \del_{\bar{X}_{ab}}^l N \eta \langle G^z ( \theta_{ab} , \i \eta ) \rangle  \right|  \notag\\
     \lesssim & N^{\eps+\delta(k+l-1)} \eta \sum_i ( |G_{ia} ^z(\theta_{ab} , \i \eta )| + |G_{ib+N} ^z(\theta_{ab} , \i \eta )| )( |G_{ai} ^z(\theta_{ab} , \i \eta )| + |G_{b+N,i} ^z(\theta_{ab} , \i \eta )| ) \notag\\
      \lesssim & N^{2 \eps + \delta(k+l)}, 
\end{align}
with the last inequality following from the Ward identity \eqref{eqn:ward} and Cauchy-Schwarz. \qed

\subsection{GFCT for one-point function} \label{sec:gfct-1pt}

Recall the definition of $Z_A (z, \eta ; X)$ in \eqref{eqn:yY-def} and let $A \geq 10$. An elementary argument using Lemma \ref{lem:1-pt-apriori} yields the estimate,
\beq \label{eqn:1-pt-gfct-apriori}
\ee[ \e^{ \lambda Z_A(z, \eta ; X) } ] \lesssim \e^{ \lambda  ( \log N)^{3/4}} \cD (z, \lambda )
\eeq
for $\eta \in [0, (\log N)^{C_*} / N]$, for any $C_* \geq 100$ and $ z\in \ddhb_r$. 
The following defines what it means for two i.i.d. matrices to match moments. 
\bed \label{def:T-match} We say that two i.i.d. matrices $X$ and $Y$ $T$-match, for some $T >0$, if for all non-negative integers $a, b$ with $a+b \leq 4$,
\beq
\left| \ee[ X_{ij}^a \bar{X}_{ij}^b] - \ee[ Y_{ij}^a \bar{Y}_{ij}^b] \right| \lesssim \1_{ \{ a+b=4 \}} T N^{-2}.
\eeq
\eed
\begin{remark}
\emph{By the construction of  \cite[Lemma 3.4]{erdos2010universality}, given any matrix $X$ one can always find a matrix of the form $Y = \sqrt{T} G + \sqrt{1-T} Y_1$ where $G$ is an independent Ginibre matrix that $T$-matches $X$. We will use this repeatedly in what follows without further comment.}

\emph{We also note that if $X$ and $Y$ $T$-match then $| \kappa_4 (X) - \kappa_4 (Y) | \lesssim T$. Therefore in  various appearances of $\kappa_4$ it will not matter if we use $\kappa_4 (X)$ or $\kappa_4 (Y)$ as the difference can be absorbed into the error in the various estimates below and so we will not distinguish between the two if it is not necessary.}
\end{remark}

\

The following compares the one-point function of two matrices with matching moments. 

\bep \label{prop:gfct-1pt}
Let $\delta \in (0, \frac{1}{100} )$, $C_* \geq 1000$, $L \geq 10$, $A \geq 10$ and $r \in (0, 1)$ be given. Let $\eta \in [N^{-1-\delta}, ( \log N)^{C_*} / N]$. Let $X$ and $Y$ be two i.i.d. matrices that $T$-match for some $T = N^{-\omega_T}$. Let $L \geq 1$ be given. Uniformly for $\lambda \in [0, L]$ and $z \in \ddhb_r$ we have,
\beq \label{eqn:gfct-1pt-c1}
\left| \ee [ \e^{ \lambda Z_A ( z, \eta; X) }] - \ee [ \e^{\lambda Z_A (z, \eta; Y) }] \right| \lesssim \cD ( z, \lambda ) N^{6 \delta } ( T + N^{-1/2} ).
\eeq
\eep
\proof Consider for the moment the  matrix $W$ that equals $X$ except for the $(a, b)$-th entry, where it equals $Y_{ab}$. Let $V$ be the matrix that equals $X$ except its $(a, b)$-th entry is $0$. Recall the notation $X(\theta_{ab})$ defined in \eqref{eqn:X-theta-def}.

Due to \eqref{eqn:theta-phi-est} and \eqref{eqn:1-pt-gfct-apriori} we have
\beq \label{eqn:1-pt-gfct-apriori-2} 
\ee\left[ \sup_{ |\theta| \leq N^{1/10-1/2}} \e^{ \lambda Z_A (  z ,\eta ; X(\theta_{ab} )) } \right] \lesssim \e^{ \lambda ( \log N)^{3/4}} \cD (z, \lambda ) .
\eeq
Using this, the fact that $|X_{ab}| \leq N^{\eps-1/2}$ with overwhelming probability, the estimates \eqref{eqn:del-phi-est}, and that $Z_A(z, \eta; M) \leq 2 A \log N$ for any matrix $M$, we have by a Taylor expansion to fifth order that
\begin{align}
\ee[\e^{ \lambda Z_A (z, \eta; X)  }] =& \sum_{n=0}^4 \sum_{i + j = n} c_{ijn} \ee[ X_{ab}^i \bar{X}_{ab}^j] \ee\left[  \left( \del_{M_{ab}}^i \del_{\bar{M}_{ab}}^j \e^{ \lambda Z_A ( z, \eta; M ) } \right)  \bigg\vert_{M=V} \right] \notag\\
+& \O \left(  N^{\eps+5 \delta-5/2} \cD (z, \lambda) \right) ,
\end{align}
for some combinatorial factors $c_{ijn}$. 
By \eqref{eqn:1-pt-gfct-apriori-2} and \eqref{eqn:del-phi-est} we have for the coefficients above,
\beq
\left| \ee\left[  \left( \del_{M_{ab}}^i \del_{\bar{M}_{ab}}^j \e^{ \lambda Z_A ( z, \eta; M ) } \right)  \bigg\vert_{M=V} \right] \right| \leq N^{\eps+\delta(i+j) } \cD (z, \lambda ) .
\eeq
Therefore, repeating the argument with $W$ in place of $X$ and comparing the Taylor expansions we see that,
\beq
\left| \ee\left[ \e^{ \lambda Z_A (X, z, \eta) }\right] - \ee\left[ \e^{ \lambda Z_A (W, z, \eta) }\right] \right| \leq N^{2 \eps+5 \delta }(N^{-1/2} + T ) N^{-2} \cD (z, \lambda ) .
\eeq
Repeating this argument by replacing the entries of $X$ by $Y$ one by one in a telescoping sum yields the claim (i.e., we have just estimated one step in the $N^2$ replacements). \qed

\subsubsection{One point function for general i.i.d. matrices}

We can now prove the one-point asymptotics for a general i.i.d. matrix. Recall the definition of $\Eone$ in \eqref{eqn:Eone-def}. 
\bet \label{thm:general-1pt} 
Let $X$ be an i.i.d. matrix and let $r \in (0, 1)$.  There is an $\alpha >0$ so that the following holds. Let $L \geq 10$, $\lambda \in [0, L]$, and $z \in \ddhb_r$. Then,
\beq \label{eqn:general-1pt}
\frac{ \ee[ | \det (X-z)|^\lambda]}{\e^{ \frac{ \lambda}{2} N (|z|^2-1)}} = (1 + \O ( N^{-\alpha} ) )\frac{ N^{\lambda^2/8}}{ \GB ( \lambda ; z, \beta)} \e^{ \frac{\lambda^2}{8} \kappa_4 (1 - |z|^2)^2 - \frac{\lambda \kappa_4}{4} (1- |z|^2)^2 +\1_{\beta=1} \Eone (z, \lambda, 0)}.
\eeq
The same estimate holds for the LHS replaced by $\ee[\e^{ \lambda Z_{A_1} (X, z, \eta_w)}]$ for $A_1$ sufficiently large and $\eta_w = N^{-1-\delta_w}$ with $0 < \delta_w \leq \frac{c_W}{2} \wedge \frac{1}{100}$, where $c_W>0$ is from Lemma \ref{lem:wegner}, with the $\alpha >0$ depending on $\delta_w$. 
\eet
\proof Let $Y$ be an i.i.d. matrix with Gaussian component of size $T = N^{-1/10}$ that $T$-matches $X$ in the sense of Definition \ref{def:T-match}.  Let $\eta_w = N^{-1-\delta_w}$ be as in the statement of the theorem. By Lemma \ref{lem:1-pt-apriori} we see that for $A_1$ sufficiently large
\beq \label{eqn:general-1pt-proof-a1}
 \ee \left| \e^{ \lambda Z_{A_1} (z, \eta_w ; X)}  - \e^{ \lambda \Phi_N (z, \eta_w ; X ) }  \1_{ \{ \Phi_N (z, \eta_w ; X) \leq A_1 \log N \} }  \right| \leq N^{-1}, 
\eeq
as well as the same estimate with $X$ replaced by $Y$. From the above and Lemma \ref{lem:1-pt-reg-2} we see that 
\beq \label{eqn:general-1pt-proof}
\ee \left| \frac{ | \det (X-z)|^\lambda}{\e^{ \frac{ \lambda}{2} N (|z|^2-1)}} - \e^{ \lambda Z_{A_1} ( z, \eta_w; X)} \right| \lesssim N^{-\delta_w^{3/2}/4} \cD (z, \lambda) .
\eeq
From the above,  Proposition \ref{prop:gfct-1pt}, and the estimate \eqref{eqn:general-1pt-proof-a1} applied with $X=Y$ we see that
\beq
\left| \ee\left[ \frac{ | \det (X-z)|^\lambda}{\e^{ \frac{ \lambda}{2} N (|z|^2-1)}} \right] - \ee\left[ \e^{ \lambda \Phi_N (z, \eta_w ; Y ) }  \1_{ \{ \Phi_N (z, \eta_w ; Y) \leq A_1 \log N \} } \right] \right| \leq N^{-\alpha} \cD (z, \lambda )
\eeq
for some $\alpha >0$, depending on $\delta_w >0$. By taking $\delta_w >0$ sufficiently small, we can apply Corollary \ref{cor:gde-1pt} (with the $X$ there being $Y$ here) to the second term on the LHS, which proves \eqref{eqn:general-1pt}. The claim about $Z_{A_1} (z, \eta_w ; X)$ follows from \eqref{eqn:general-1pt} and \eqref{eqn:general-1pt-proof}. \qed

\subsection{GFCT for GMC} \label{sec:GMC-GFCT}

In this section we will  establish a Green's function comparison theorem for the observables we expect to converge to the GMC (with a regularization - i.e., we use $\Phi_N (z, \eta)$ instead of $\Phi_N (z, 0)$). We then combine this result with the results of Sections \ref{sec:gmc-reg} and \ref{sec:claeys} to prove the full convergence to the GMC for a general i.i.d. matrix. 

Fix $ 0 < r <1$. Let $f : \cc \to \rr$ be a bounded measureable function supported in $\ddb_r$ and $g : \rr \to \rr$ a bounded measureable function supported in $I_r$. For an i.i.d. matrix $X$, $\eta >0$ and $\lambda \geq 0$ let us define,
\begin{align}
\Wlac (f, \eta; X) &:= \int_{\cc} f(z) \frac{\e^{ \lambda Z_A(z, \eta ; X) }  \dz}{N^{\lambda^2/8}} , \quad 
\Wlar (g, \eta ; X) :=  \int_\rr g (x) \frac{\e^{ \lambda Z_A (x, \eta ; X) } \d x}{N^{ \frac{ \lambda ( \lambda-1)}{4}}}
\end{align}
We only consider $\Wlar$ in the case $\beta=1$. 
Let $ \Qlac (\eta ; X) := \Wlac ( \1_{\ddb_r} , \eta; X)$ and $\Qlar (\eta ; X) := \Wlar ( \1_{I_r} , \eta ; X)$. By \eqref{eqn:1-pt-gfct-apriori} we have for $\eta \in [0, ( \log N)^{C_*} / N ]$,
\beq \label{eqn:GMC-gfct-a2}
\ee[\Qlaf ( \eta ; X) ] \lesssim \e^{ \lambda ( \log N)^{3/4}} , \qquad \ff \in \{ \rr, \cc \} .
\eeq
The following lets us control derivatives of $\Wlaf$ with respect to the matrix elements. 
\bel Let $L \geq 10, n \geq 1$ and $\lambda \in [0, L]$. Let $\delta \in (0, \frac{1}{100})$. 
We have with overwhelming probability, for any $\eta \in [N^{-1-\delta}, 1 ]$ and $0 \leq i, j \leq n$ that for any $\eps >0$
\beq \label{eqn:WA-deriv}
\sup_{ |\theta| \leq N^{1/10-1/2} } \left| \del_{M_{ab}}^i \del_{\bar{M}_{ab}}^j \Wlac ( f, \eta; M) \bigg\vert_{M = X(\theta_{ab} ) } \right| \leq N^{\eps+(i+j) \delta} \Qlac (\eta ; X) .
\eeq
The same estimate holds for $\Wlac$ replaced by $\Qlac$ on the LHS, as well as analogous estimates for $\Wlar, \Qlar$. 

\eel
\proof The proof follows from direct computation and \eqref{eqn:del-phi-est} and \eqref{eqn:theta-phi-est}. \qed

\

The following is a Green's function comparison theorem for the observables we expect to converge to the GMC. Note that we carry around the cut-off function $F_2$ below, as the higher moments of $\Qlaf$ are not expected to remain bounded in the large $N$ limit. This would make it difficult to apply the Lindeberg approach to $\ee[ F_1 (\Wla (f, \eta ; X))]$ directly without the cut-off function $F_2$. 
\bep \label{prop:GMC-gfct}
Let $X$ and $Y$ be two i.i.d. matrices that $T$-match in the sense of Definition \ref{def:T-match}, with $T= N^{-\omega_T}$. Let $F_1 : \rr \to \rr $ be a bounded function with bounded continuous derivatives up to order $5$. Let $\delta_c \in (0, \frac{1}{100} )$, and let $F_2: \rr \to \rr$ be a smooth function such that $F_2 (x) = 1$ for $|x| \leq N^{\delta_c}-1$ and $F_2(x) = 0$ for $|x| > N^{\delta_c}$. Let $\delta >0$ and let $\eta \in [N^{-1-\delta}, ( \log N )^{C_*} / N ]$. Then, 
\begin{align}
     &\left| \ee[ F_1 (\Wlac (f, \eta ; X) ) F_2 (\Qlac (\eta ; X) ] - \ee[ F_1 (\Wlac (f, \eta ; Y) ) F_2 (\Qlac (\eta ; Y) ]  \right| \notag\\
    \leq & N^{6(\delta+\delta_c)} (N^{-1/2} + T) .
\end{align}
as well as an analogous estimate for $\Wlar, \Qlar$. 
\eep
\proof By the support properties of $F_2$ and \eqref{eqn:WA-deriv}, for any i.i.d. matrix $W$, we have
\beq
\sup_{ |\theta| \leq N^{1/10-1/2}} \left| \del_{M_{ab}}^i \del_{ \bar{M}_{ab}}^j \left(  F_1 (\Wlac (f, \eta ; M) ) F_2 ( \Qlac (\eta ; M ) ) \right)  \big\vert_{M = W(\theta_{ab} ) } \right| \leq N^{\eps+(i+j)(\delta+\delta_c)} 
\eeq
for any $\eps >0$ 
with overwhelming probability. The claim now follows via the same proof as in Proposition \ref{prop:gfct-1pt}. The proof for $\Wlar, \Qlar$ is the same. \qed

\subsubsection{GMC for general matrices with regular functions}

We can now conclude the following convergence to the GMC for sufficiently regular functions with compact support. 
\bet \label{thm:nice-function-GMC} 
Let $X$ be an i.i.d. matrix and let $f : \cc \to \rr$ be a bounded continuous function supported in $\ddb_r$, for some $0 < r < 1$. Then for $\lambda \in (0, 2 \sqrt{2})$ we have
\beq \label{eqn:nice-GMC-C}
\lim_{N \to \infty} \int f(z) \frac{ | \det (X-z)|^{\lambda}}{ \ee[ | \det (X-z) |^\lambda]} \dz = \int f(z) \d \mu_{\GMC, \kappa_4}^{\beta,\lambda} (z)
\eeq
in distribution. If $g : \rr \to \rr$ is bounded continuous and supported in $I_r$, and $\beta=1$ then for $\lambda \in (0, \sqrt{2})$ we have
\beq \label{eqn:nice-GMC-R}
\int_{N \to \infty} \int g(x) \frac{ | \det (X-x) |^\lambda}{ \ee[ | \det (X-x ) |^\lambda] } \d x = \int g(x) \d \mu_{\GMC, \kappa_4}^{\rr, \lambda} (x) 
\eeq
in distribution. In the above, $\kappa_4$ is the fourth cumulant of $X$.
\eet
\proof We first prove \eqref{eqn:nice-GMC-C}.  We first consider the case that $f$ is bounded and Lipschitz. Let $F_1 : \rr \to \rr$ be a bounded function with bounded derivatives up to order $5$. Let $F_2$ be as in Proposition \ref{prop:GMC-gfct} for some $\delta_c >0$. Choose $\eta = N^{-1-\delta}$ for some $\delta >0$. By Theorem \ref{thm:general-1pt} the hypotheses of Proposition \ref{prop:gmc-reg-new} are satisfied for $X$. Therefore, by Proposition \ref{prop:gmc-reg-new} and \eqref{eqn:GMC-gfct-a2} we have for some $\alpha >0$,
\begin{align} \label{eqn:nice-gmc-a1}
& \left|  \ee \left[ F_1 \left( \int f(z) \frac{ | \det (X-z)|^\lambda}{ \ee[ | \det (X-z)|^\lambda]} \dz \right)  \right] - \ee[ F_1 ( \Wlac ( \tilde{f} , \eta ; X ) F_2 ( \Qlac ( \eta ; X) ] \right| \notag\\
\lesssim & N^{-\delta_c/2}  + N^{-\frac{\delta^{3/2}}{4}} + N^{-\alpha},
\end{align}
for $A \geq 1$ sufficiently large and $\delta >0$ sufficiently small. Here $\tilde{f} (z) := f(z) / \kB (z) $ where $\kB (z)$ is as in the statement of Proposition \ref{prop:gmc-reg-new}. 

Let $T = N^{-1/10}$ and let $Y$ be an i.i.d. matrix of the form $Y = \sqrt{T} G + \sqrt{1-T} W$, for $G$ an independent Ginibre matrix and $W$ an i.i.d. matrix such that $Y$ and $X$ $T$-match. Note that the fourth cumulants of $Y$ and $W$ equal that of $X$, up to an error of $\O (N^{-1/20})$, and so Assumption \ref{ass:claeys-kappa} is satisfied by $Y$.  From Proposition \ref{prop:gmc-reg-new} and \eqref{eqn:GMC-gfct-a2} we have that,
\begin{align} \label{eqn:nice-gmc-a2}
    & \left| \ee \left[ F_1 \left( \int f(z) \frac{ \e^{ \lambda Z_{A} (z, \eta ; Y)}}{ \ee[ \e^{ \lambda Z_A (z, \eta ; Y)}]} \dz \right)  \right] - \ee[ F_1 ( \Wlac ( \tilde{f} , \eta ; Y ) F_2 ( \Qlac ( \eta ; Y) ] \right| 
    \lesssim  N^{-\delta_c/2} + N^{-\alpha} .
\end{align}
By Theorem \ref{thm:GMC-gde} we have, after possibly taking $\delta >0$ smaller,
\beq \label{eqn:nice-gmc-a3}
\lim_{N \to \infty}  \ee \left[ F_1 \left( \int f(z) \frac{ \e^{ \lambda Z_{A} (z, \eta ; Y)}}{ \ee[ \e^{ \lambda Z_A (z, \eta ; Y)}]} \dz \right)  \right] = \ee\left[ F_1 \left( \int f(z) \d\mu_{\GMC, \kappa_4}^{\beta,\lambda} (z)\right) \right] .
\eeq
Here we also used \eqref{eqn:t-scaling} and the fact that $f$ is Lipschitz. 
The claim for Lipschitz $f$ now follows from \eqref{eqn:nice-gmc-a1}, \eqref{eqn:nice-gmc-a2}, \eqref{eqn:nice-gmc-a3} and Proposition \ref{prop:GMC-gfct}. The extension to general continuous functions $f$ supported in $\ddb_r$ is straightforward: one approximates $f$ by Lipschitz $\tilde{f}$ in the $L^{\infty}$ norm and uses the fact that,
\beq
\ee\left[ \int |f(z) - \tilde{f} (z) | \frac{ | \det (X-z)|^{\lambda}}{ \ee[ | \det (X-z)|^{\lambda} ]} \dz  \right] \lesssim \| f - \tilde{f} \|_\infty .
\eeq
The proof of \eqref{eqn:nice-GMC-R} is similar;  one uses Proposition \ref{prop:GMC-reg-R} instead of Proposition \ref{prop:gmc-reg-new} and Theorem \ref{thm:GMC-gde-line} instead of Theorem \ref{thm:GMC-gde}. 
\qed

\vspace{5 pt}

\noindent{\bf Proofs of Theorems \ref{theo:maingmcconv} and \ref{thm:gmc-real}}. The deduction of these results from Theorem \ref{thm:nice-function-GMC} is identical to the proof of \cite[Theorem 1.4]{claeys2021much}. I.e., see the proof of \cite[Proposition 2.12]{claeys2021much} and then the arguments in the first few paragraphs of \cite[Section 2.7]{claeys2021much}. \qed

\subsection{GFCT for K-point function} \label{sec:gfct-kpt}

We will now develop estimates leading to the computation of the $K$-point function for general i.i.d. matrices. This section will be somewhat more involved than the analogous section for the $1$-point function Section \ref{sec:gfct-1pt}. This is due to the fact that the a priori estimate Lemma \ref{lem:k-pt-apriori} and the regularization estimate Lemma \ref{lem:kpt-reg} so far hold only for Gaussian divisible matrices (compare with Lemmas~\ref{lem:1-pt-apriori} and~\ref{lem:1-pt-reg-2}) and we first need to establish them for general i.i.d. matrices via a Green's function comparison theorem. 

We now fix a  constant $b >0$ and $K$ points $z_i$ such that $z_i \in \ddhb_r$ and $|z_i -z_j| \geq N^{b-1/2}$ for all $i \neq j$. Let us fix $\blam = ( \lambda_1, \dots, \lambda_K ) \in [0, L]^K$ some $L \geq 10$. We also let the constants $\omega_1, \mfb$ and $\mfa >0$ be as in Theorem \ref{thm:local-global}. In particular,  the constant $\mfa >0$ will appear in some of the estimates below.

For an $A>0$ and a matrix $M$ we define the quantity,
\beq
\Pone_A( \eta; M ) = \prod_{i=1}^K \e^{ \lambda_i Z_A( z_i, \eta ; M) } .
\eeq
Let $X$ and $Y$ be two matrices that $T$-match. Let $\gamma(a, b ) : [\![1, N]\!]^2 \to  [\![1, N^2 ]\!] $ be some enumeration of the $N^2$ indices of an $N \times N$ matrix.  Let $U_i$, for $0 \leq i \leq N^2$, be a sequence of matrices such that $(U_i) _{ab} = Y_{ab}$ if $\gamma(a, b)  \leq i$ and $(U_i)_{ab} = X_{ab}$ if $\gamma(a, b) > i$. Then $U_0 = X$ and $U_{N^2} = Y$. 

Let $F_3$ be a bounded function with bounded derivatives up to order 5. Define for some $\eta >0$,
\beq
\Ptwo (\eta; M) = \sum_{i=1}^K F_3 \left( N \eta \Im[ \langle (H^{z_i} (M) - \i \eta)^{-1} \rangle ] \right) .
\eeq

The following is the main technical GFCT result for the $k$-point function. 
\bep \label{prop:gfct-kpt-compare} Let $\delta \in (0, \frac{1}{100} )$.  Let $\eta_1, \eta_2 \in [N^{-1-\delta}, 1]$. 
For any $D>0$ we have,
\begin{align} \label{eqn:gfct-kpt-a1}
 & \left| \ee[\Pone_A ( \eta_1; U_i) \Ptwo (\eta_2 ; U_i) ]   - \ee[\Pone_A (\eta_1; Y) \Ptwo ( \eta_2 ; Y) ]\right| \notag\\
\lesssim & N^{-D} + N^{6\delta} ( N^{-1/2} + T ) \sup_{0 \leq k \leq N^2} \ee[\Pone_A ( \eta_1; U_k) ],
\end{align} 
for all $0 \leq i \leq N^2$ and for $N$ large enough.
\eep
\proof From \eqref{eqn:del-phi-est}, \eqref{eqn:del-G-est}, and \eqref{eqn:theta-phi-est}, one sees that for any i.i.d. matrix $W$ with overwhelming probability,
\beq
\sup_{ |\theta| \leq N^{1/10-1/2}} \left| \del_{M_{ab}}^i \del_{\bar{M}_{ab}}^j [ \Pone_A (\eta_1; M) \Ptwo (\eta_2; M ) ] \bigg\vert_{M = W(\theta_{ab})} \right| \leq N^{\eps + (i+j) \delta} \Pone_A ( \eta_1 ; W) .
\eeq
Therefore, proceeding as in the proof of Proposition \ref{prop:gfct-1pt} one finds that for $0 \leq j < N^2$, we have
\begin{align} 
  & \left| \ee[\Pone_A ( \eta_1; U_j) \Ptwo (\eta_2 ; U_j) ]   - \ee[\Pone_A (\eta_1 ; U_{j+1}) \Ptwo ( \eta_2 ; U_{j+1}) ]\right| \notag\\
 \leq & N^{-D} + N^{-2} N^{\eps + 5 \delta}(T+N^{-1/2} ) \ee[\Pone_A (\eta_1 ; U_j ) ] \notag\\
 \leq & N^{-D} + N^{-2} N^{\eps + 5 \delta}(T+N^{-1/2} ) \sup_{0 \leq k \leq N^2} \ee[\Pone_A ( \eta_1 ; U_k ) ]
\end{align}
The claim now follows from writing the difference on the LHS of \eqref{eqn:gfct-kpt-a1} as a telescoping sum of the above terms for $i \leq j < N^2$. \qed

\

We now establish the following a priori bounds for general i.i.d. matrices which will allow us to deduce the $K$-point asymptotics. 
\bep Let $C_* \geq 100$.  \label{prop:gfct-kpt-apriori}
For any i.i.d. matrix $X$, uniformly for $ \eta \in [0, ( \log N)^{C_*} / N]$, we have
\beq \label{eqn:gfct-kpt-a3}
\ee[ \Pone_A (\eta ; X) ] \leq \e^{ ( \log N)^{7/8}} \cD(\bz, \blam ) ,
\eeq
for $N$ sufficiently large. Fix $\eta_w = N^{-1-\delta_w}$, with $0 < \delta_w \leq \frac{c_W}{2} \wedge \frac{1}{100} \wedge \mfa$. Then, for all $ \eta \in [0, \eta_w]$ and for all $A$ sufficiently large, we have
\beq \label{eqn:gfct-kpt-a2} 
\ee\left[ \Pone_A (\eta ; X) \sum_{i=1}^K \1_{ \{ \lambda_1^{z_i}(X) \leq \eta_w\}}  \right] \lesssim N^{-\frac{ \delta_w^{3/2}}{2}} \cD (\bz, \blam).
\eeq
\eep
\proof  Let $Y = \sqrt{T} G + (1- T)^{1/2} Y_1$ be an i.i.d. matrix that $T$-matches $X$ with $G$ an independent Ginibre matrix of the same symmetry class, with $T = N^{-1/10}$. 

We begin with the first estimate \eqref{eqn:gfct-kpt-a3}. By Lemma \ref{lem:a-priori-1} it suffices to prove it for the case $\eta = \eta_* := ( \log N)^{C_*}/N$ and $A$ sufficiently large with the $( \log N )^{7/8}$ term replaced by $\| \blam \|_1 ( \log N)^{3/4}$. Let $U_i$ be matrices interpolating between $X$ and $Y$ as defined at the beginning of Section \ref{sec:gfct-kpt}. The estimate \eqref{eqn:kpt-apriori-b1} and Lemma \ref{lem:1-pt-apriori} imply that the estimate \eqref{eqn:gfct-kpt-a3} holds for $X$ replaced by $Y$, if $A$ is sufficiently large. Therefore, applying \eqref{eqn:gfct-kpt-a1} with $F_3 = \frac{1}{K}$, we obtain
\beq
\sup_{0 \leq i \leq N^2} \ee[\Pone_{A} ( \eta_* ; U_i) ] \lesssim \e^{ \| \blam \|_1 ( \log N)^{3/4}} \cD(\bz, \blam ) + N^{-D} + N^{-1/20}  \sup_{0 \leq i \leq N^2} \ee[\Pone_{A} (\eta_* ; U_i) ] ,
\eeq
and so the estimate \eqref{eqn:gfct-kpt-a3}  follows by moving the last term on the RHS back to the LHS.

We now turn to the second estimate \eqref{eqn:gfct-kpt-a2}. Let $F_3 : \rr \to \rr$ be a smooth function such that $F_3 (x) = 0$ for $x < \frac{1}{20}$ and $F_3(x) = 1$ for $x > \frac{1}{10}$. By Lemma \ref{lem:smoothed-indicator} and Lemma \ref{lem:a-priori-1} we have 
\beq
\ee\left[ \Pone_A (\eta; X) \sum_{i=1}^K \1_{ \{ \lambda_1^{z_i} \leq \eta_w\}}  \right]  \lesssim  \ee\left[ \Pone_A (\eta_w ; X) \Ptwo ( \eta_w ; X ) \right] + N^{-D}.
\eeq
By \eqref{eqn:gfct-kpt-a1} and \eqref{eqn:gfct-kpt-a3} we then see that
\beq
\left| \ee\left[ \Pone_A ( \eta_w; X) \Ptwo (\eta_w ; X ) \right] - \ee\left[ \Pone_A (\eta_w; Y) \Ptwo (\eta_w ; Y ) \right] \right| \lesssim N^{-1/20} \cD ( \bz, \blam ) .
\eeq
By Lemma \ref{lem:smoothed-indicator} and Theorem \ref{thm:ll-2} we see that
\beq
\ee\left[ \Pone_A (\eta_w; Y) \Ptwo ( \eta_w ; Y ) \right] \leq \ee\left[ \Pone_A ( \eta_w ; Y) \sum_{i=1}^K \1_{ \{ \lambda_1^{z_i} (Y) \leq ( \log N)^3 \eta_w \}} \right] + N^{-D} .
\eeq
Then by Lemma \ref{lem:kpt-reg} (and also Lemma \ref{lem:1-pt-apriori} after taking $A$ possibly larger) we have,
\begin{align}
\ee\left[ \Pone_A (Y, \eta_w) \sum_{i=1}^K \1_{ \{ \lambda_1^{z_i} (Y) \leq ( \log N)^3 \eta_w \}} \right] \lesssim N^{- \frac{ \delta_w^{3/2}}{2}} \cD ( \bz, \blam) ,
\end{align}
and so we conclude the proof of \eqref{eqn:gfct-kpt-a2}. \qed 

\

The above results finally let us compare the $K$-point function to its regularized version for general i.i.d. matrices. 
\bec \label{cor:kpt-reg} Let $X$ be an i.i.d. matrix. 
Let $\eta_m = N^{-1-\delta_m}$ with $\delta_m \leq \frac{c_W}{2} \wedge \frac{1}{100} \wedge \mfa$. We have,
\beq \label{eqn:gfct-kpt-a4} 
\left| \ee[ \e^{ \sum_{i=1}^K \lambda_i \Phi_N (z_i,0 ; X) } ] - \ee[ \Pone (\eta_m; X) ] \right| \lesssim N^{- \frac{ \delta_m^{3/2}}{10}} \cD (\bz, \blam ).
\eeq
\eec
\proof The proof is very similar to Lemma \ref{lem:1-pt-reg-2} and so we do not give all the details. Let $\S := \bigcup_{i=1}^K \{ \lambda_1^{z_i}  \leq \eta_w \}$ with $\delta_w = \frac{\delta_m}{2}$. The estimates \eqref{eqn:gfct-kpt-a3} and \eqref{eqn:gfct-kpt-a2} implies for $\eta \in [0, \eta_w]$ that
\begin{align} \label{eqn:gfct-kpt-a5}
&\ee\left[ \e^{ \sum_{i=1}^K \lambda_i \Phi_N (z_i, \eta ) } \1_{ \S \cap \bigcap_{i=1}^K \{ \Phi_N (z_i, \eta) \leq A \log N \} }  \right] \lesssim N^{- \frac{ \delta_w^{3/2}}{2}} \cD ( \bz , \blam ), \notag\\
&\ee\left[ \e^{ \sum_{i=1}^K \lambda_i \Phi_N (z_i, \eta ) } \1_{  \bigcap_{i=1}^K \{ \Phi_N (z_i, \eta) \leq A \log N \} }  \right]  \lesssim \e^{ ( \log N)^{7/8}} \cD ( \bz, \blam).
\end{align} 
In similar fashion to the proof of Lemma \ref{lem:1-pt-reg-2}, using Lemma \ref{lem:rig-reg} (as well as Corollary \ref{cor:char-poly-a-priori}, the second estimate in \eqref{eqn:gfct-kpt-a5} and Lemma \ref{lem:1-pt-apriori}) one finds for $A$ sufficiently large,
\beq
\left| \ee[ \e^{ \sum_{i=1}^K \lambda_i \Phi_N (z_i , 0) } \1_{ \S^c}]  -\ee\left[ \e^{ \sum_{i=1}^K \lambda_i \Phi_N (z_i, \eta_m ) } \1_{ \S^c \cap \bigcap_{i=1}^K \{ \Phi_N (z_i, \eta_m) \leq A \log N \} }  \right] \right| \lesssim N^{-\delta_m/2} \cD ( \bz , \blam ) .
\eeq
The above estimates along with a few uses of Lemma \ref{lem:1-pt-apriori} to compare $\Phi_N (z, \eta)$ with $Z_{A} (z, \eta)$ are sufficient to complete the proof. \qed

\vspace{5 pt}

\noindent{\bf Proof of Theorem \ref{thm:k-pt-intro}  (case of $\gamma_i \in \rr$).}  Let $T = N^{-1/10}$ and $Y = \sqrt{T} G + (1- T)^{1/2} Y_1$ be a GDE matrix that $T$-matches $X$. Let $\eta_m$ be as in Corollary \ref{cor:kpt-reg}. From \eqref{eqn:gfct-kpt-a1} and \eqref{eqn:gfct-kpt-a3} we have
\beq \label{eqn:kpt-asymptotics-a1}
\left| \ee[ \Pone (\eta_m ; X) ] - \ee[ \Pone (\eta_m ; Y) ] \right| \lesssim N^{-\alpha} \cD ( \bz, \blam ). 
\eeq
By taking $A$ sufficiently large we see from Corollary \ref{cor:kpt-reg-gde} and Lemma \ref{lem:1-pt-apriori} that the asymptotics on the RHS of \eqref{eq:momasymp} hold for $\ee[ \Pone (\eta_m ; Y) ]$. The result now follows from  \eqref{eqn:gfct-kpt-a4} and \eqref{eqn:kpt-asymptotics-a1}. \qed

\section{Extension to complex parameters} \label{sec:complex}

In this section we show how to extend the result of Theorem \ref{thm:k-pt-intro} to complex $\lambda$. For simplicity we will only discuss the case that $z_i \in \ddb_r$ for all $i$; the case where some $z_i \in \rr$ when $\beta=1$ is only slightly more complex.

For clarity we first discuss the one-point function. Fix $L \geq 1$, $r \in (0, 1)$, $z \in \ddb_r$, and $\lambda = x + \i y \in [0, L] + \i [-1, 1]$. Our method for estimating the one-point function $\ee[ \e^{ \lambda \Phi_N (z, 0) } ]$ (and its various regularizations) will be the same. The essential major difference for the case of complex $\lambda$ is that the errors in our approach will be of the form (ignoring cut-offs)
\beq
 N^{-c} \ee  \left| \e^{ \lambda \Phi_N (z, \eta) } \right|  =N^{-c}  \ee[  \e^{ x \Phi_N (z, \eta) } ] \lesssim N^{-c} N^{x^2/8} .
\eeq
The main point is $| \e^{ \lambda \Phi_N (z, \eta) } | = \e^{ \Re[ \lambda] \Phi_N (z, \eta) } $ whose expectation is much larger than $N^{\frac{x^2-y^2}{8}}$ which is the expected size of $\ee[ \e^{ (x+ \i y ) \Phi_N (z, 0) } ]$. Therefore, in order to obtain a multiplicative error as in the statement of Theorem \ref{thm:k-pt-intro}, one must ensure that $y$ is sufficiently small, depending on the exponents $c >0$ appearing in the errors (e.g., $y < \sqrt{c}$). 

For example, following the arguments in Section \ref{sec:1pt} shows that the analog of \eqref{eqn:1-pt-prod-c1} is
\beq
\ee[ \e^{ \lambda \Phi_N (z, \eta_m, t_1 ) } \1_{ \A (z, \eta_m, t_1 , A_1)} ] = \calL \times \ee[ \e^{ \lambda \tilG } \1_{ \tilG \leq A_1 \log N } ] + \O (N^{-c} N^{x^2/8} ),
\eeq
for some $c>0$, and the local factor $\calL$ involving the complex parameter $\lambda$. In order to compute the mesoscopic term, one can proceed as in Section \ref{sec:meso} and apply Stein's method. One finds, (ignoring the cut-offs for notational simplicity) with the covariance $\cC$ and expectation $\bmE$ terms as before,
\begin{align}
\del_s \ee[ \e^{ (x + \i s) \tilG }] = &  \{ ( -s + \i x) \cC + \i \bmE  \} \ee[\e^{ (x + \i s) \tilG } ] + \O ( N^{-\omega_1/5} \ee[ |\e^{ (x + \i s) \tilG } |] + N^{-D} ) \notag\\
= & \{ ( -s + \i x) \cC + \i \bmE  \} \ee[\e^{ (x + \i s) \tilG } ] + \O ( N^{-\omega_1/5} \e^{ - \frac{x^2}{8} \log (t_1) } + N^{-D} ),
\end{align}
where in the second line we used Proposition \ref{pro:jointlapdnew} with $\lambda = x$. Hence, one can integrate the above equation and for $y$ sufficiently small, compute $\ee[ \e^{ ( x+ \i y ) \tilG } ]$ with a multiplicative error $\O (N^{-\omega_1/10} )$. With this estimate in hand, one can proceed as in Section \ref{sec:1pt} to first compute the local factor $\calL$ with complex $\lambda$ as in the proof of Theorem \ref{thm:local-factor} and then find the analog of Corollary \ref{cor:gde-1pt}. 

Similar considerations now apply to the Green's function comparison theorems of Section \ref{sec:gfct-1pt}. In particular, the error on the RHS of \eqref{eqn:gfct-1pt-c1} will be of the form $N^{6\delta + \frac{ ( \Re[ \gamma] )^2}{8} } (T + N^{-1/2} )$ which will be smaller than $| N^{\frac{ \gamma^2}{8} }|$ if $\Im[ \gamma]$ is sufficiently small. With this in hand, the proof of Theorem \ref{thm:general-1pt} for complex $\lambda$ is the same, which completes discussion of the complex parameter one-point result of Theorem \ref{thm:k-pt-intro}.

For the full $K$-point function, the main point is that now the method gives errors of the form $N^{-c} \cD(\bz, \Re[ \blam] ) $ and one notes that
\beq
\cD(\bz, \Re[ \blam] ) \lesssim N^{ K^2 \| \Im[ \blam ] \|_\infty^2 } | \cD ( \bz, \blam ) |,
\eeq
with $| \cD ( \bz, \blam ) |$ the expected size of the $K$-point function. Therefore as long as $\| \Im[ \blam ] \|_\infty$ is sufficiently small, one will find the analog of \eqref{eqn:k-pt-prod} with an error on the RHS of order $| \cD ( \bz, \blam ) | N^{-c/2}$ and then the analog of \eqref{eqn:kpt-reg-gde} with  a small multiplicative error (after, of course, extending $K$-point version of Proposition \ref{pro:jointlapdnew} in an analogous way to the one-point version discussed above). Extending the arguments of Section \ref{sec:gfct-kpt} is easy; one can directly apply the results of Proposition \ref{prop:gfct-kpt-apriori} to derive an analog of \eqref{eqn:gfct-kpt-a4} with an error of $N^{ - \frac{ \delta_m^{3/2}}{20}} | \cD (z, \blam ) |$ as long as $\| \Im[ \blam] \|_\infty$ is sufficiently small, depending on the choice of $\delta_m >0$. With this established, the proof of the complex analog of Theorem \ref{theo:momcomp} is the same, with now $\| \Im [ \blam ] \|_\infty$ taken sufficiently small depending on the $\alpha >0$ appearing on the RHS of \eqref{eqn:kpt-asymptotics-a1}. This completes the discussion of the $K$-point asymptotics with complex parameters.

\appendix

\section{Proof of various technical results}
\label{app:tech}

In this section we present various proofs of technical results using throughout the paper.

\

\noindent{\bf Proof of Proposition~\ref{prop:parameter-properties}}. 
We start with the proof of the first item. Fix a small $c>0$, we start considering the case when $|z_1-z_2|\le c$ and $\eta_1, \eta_2\le c$. By \cite[Eq. (A.10)]{cipolloni2025optimal}, we have
\beq
\label{eq:otherformV}
\mathcal{V}(z_1,\eta_1,z_2,\eta_2)=-\frac{1}{4}\log\big[u_1u_2|z_1-z_2|^2-m_i^2\frac{u_2}{u_1}(1-u_1)-m_2^2\frac{u_1}{u_2}(1-u_2)+(1-u_1)(1-u_2)].
\eeq
Additionally, by \eqref{eq:exp}, for $i=1,2$, we have
\beq
\label{eq:exp-appendix}
u_i=1-\frac{\eta_i}{\sqrt{1-|z_i|^2}}+\mathcal{O}(\eta_i^2), \qquad\quad m_i=\i \sqrt{1-|z_i|^2}+\mathcal{O}(\eta_i).
\eeq
Plugging \eqref{eq:exp-appendix} into \eqref{eq:otherformV}, we thus have
\beq 
\mathcal{V}(z_1,\eta_1,z_2,\eta_2)=-\frac{1}{4}\log\big[|z_1-z_2|^2+\sqrt{1-|z_1|^2}\eta_1+\sqrt{1-|z_2|^2}\eta_2+\mathcal{O}(\eta_1^2+\eta_2^2+|z_1-z_2|^4)\big].
\eeq
On the other hand,  if $|z_1-z_2|> c$ or one among $\eta_1, \eta_2$ is bigger than $c$, then the fact that $\mathcal{V}(z_1,\eta_1,z_2,\eta_2)\sim 1$ follows from \cite[Eq. (A.5)]{cipolloni2025optimal}, as a consequence of the fact that $V_{ij}=\log[\beta_+\beta_-]$ with the notation therein (see e.g. \cite[Eq. (A.10)]{cipolloni2025optimal}). The proof of the last two items is completely analogous using again \eqref{eq:exp-appendix}. 
\qed

\

\

We now state a (non-optimal) bound for products of two resolvents which we will use in the real case ($\beta=1$) to estimate certain error terms.
\bel
\label{lem:2grealb}
Fix a large $C>0$. Let $X$ be a real i.i.d. matrix, then, for any $\eps>0$, we have
\beq
\label{eq:desres}
\big|\langle G^z(\i\eta_1)AG^{\overline{z}}(\i\eta_2) B\rangle\big|\le  \frac{N^{3/4+\eps}}{[|z-\overline{z}|^2+N^{-1}]^{1/4}}.
\eeq
with overwhelming probability uniformly in $\eta_*:=\min_i|\eta_i|\ge (\log N)^C/N$ and $\lVert A\rVert, \lVert B\rVert\lesssim 1$. The same bound holds if one or both of the resolvents $G$ is replaced by its absolute value $|G|$.
\eel

\proof
We first present the proof of \eqref{eq:desres} and then we comment on how to obtain an analogous bound when one or both resolvents are replced by their absolute value.

By \cite[Eqs. (3.20) and (3.26)]{cipolloni2025optimal}, we have that \eqref{eq:desres} holds for $\eta_*\ge N^{-1+\eps}=:\eta_\eps$, for any $\eps>0$. Now, without loss of generality, we assume that $N^{-1}(\log N)^C\le \eta_1< N^{-1+\eps}$ and $\eta_2\ge N^{-1+\eps}$. If both are smaller than $N^{-1+\eps}$ then we can proceed similarly to the computations below, using the estimate when at least one of the $\eta_i \geq N^{-1+\eps}$ as an input (see e.g. \cite[Proof of Corollary B.6]{cipolloni2025maximum} for a similar argument). We have
\beq
\label{eq:derid}
\langle G^z(\i\eta_1)AG^{\overline{z}}(\i\eta_2) B\rangle=\langle G^z(\i\eta_\eps)AG^{\overline{z}}(\i\eta_2) B\rangle+\i\int_{\eta_1}^{\eta_\eps}\langle [G^z(\i\eta)]^2AG^z(\i\eta_2) B\rangle\, \d \eta.
\eeq
We then estimate
\beq
\begin{split}
\label{eq:bneedreal}
\int_{\eta_1}^{\eta_\eps}\langle [G^z(\i\eta)]^2AG^{\overline{z}}(\i\eta_2) B\rangle\, \d \eta &\le \int_{\eta_1}^{\eta_\eps}\frac{1}{\eta\eta_2^{1/2}}\langle \Im G^z(\i\eta)A \Im G^{\overline{z}}(\i\eta_2)A^*\rangle^{1/2}\langle \Im G^z (\i\eta) \rangle^{1/2}\, \d \eta \\
&\le N^\eps \int_{\eta_1}^{\eta_\eps}\frac{\eta_\eps^{1/2}}{\eta^{3/2}\eta_2^{1/2}}\langle \Im G^z(\i\eta_\eps)A \Im G^{\overline{z}}(\i\eta_2)A^*\rangle^{1/2}\, \d \eta \\
&\le \frac{N^{1/2+\eps}}{[|z-\overline{z}|^2+\eta_\eps]^{1/2}},
\end{split}
\eeq
where in the second inequality we used the monotonicity of $\eta\mapsto \eta\Im G^z$ and the bound 
\beq
|\langle \Im G^z (\i\eta) \rangle|\le N^\epsilon
\eeq
 by Theorem~\ref{thm:ll-1}. Additionally, in the last inequality we used \cite[Eqs. (3.20) and (3.26)]{cipolloni2025optimal}.

Finally, \eqref{eq:desres} follows by combining \eqref{eq:bneedreal} with \cite[Eqs. (3.20) and (3.26)]{cipolloni2025optimal} to estimate the first term in the LHS of \eqref{eq:derid}. Note that here we used that $\eta\ge \eta_\eps$ only to estimate the first term in the LHS of \eqref{eq:derid}. The bound \eqref{eq:desres} for $|G|$ immediately follows by an integral representation of $|G|$ in terms of $\Im G$ as in \cite[Eqs. (B.23)-(B.25)]{cipolloni2025maximum}.
\qed

In various parts of the proof we will use the following martingale bound:
\bel \label{lem:mart-bd} 
Let $M(t)$ be a continuous Martingale with $M(0) = 0$. Suppose there is a constant $B>0$ and an $\eps >0$ so that,
\beq
\pp\left[ \langle M \rangle_T > B \right] \leq \eps .
\eeq
Then, for all $ u >0$, we have
\beq
\pp\left[ \sup_{0 < t < T } |M_t| > u B^{1/2} \right] \leq \eps + C \e^{ - \frac{u^2}{2} }.
\eeq
\eel
\proof Without loss of generality we assume $B=1$. Let $\tau := \inf \{ t : \langle M \rangle_t > 1 \}$. Then, $M(t) = M(t \wedge \tau )$ for all $t \leq T$ with probability at least $1-\eps$. On the other hand $t \to M( t \wedge \tau )$ is a Martingale with quadratic variation bounded by $1$ and so by the Martingale representation theorem,
\beq
\pp\left[ \sup_{ 0 < t < T } | M ( t \wedge \tau ) | > u \right] \leq \pp\left[ \sup_{ 0 \leq t \leq 1 } |B(t) | > u \right] 
\eeq
where $B(t)$ is a Brownian motion. The claim now follows. \qed

\subsection{Proofs of moment formulas} \label{a:moments}

\noindent{\bf Proof of Proposition \ref{prop:ginibre}.} The proof in the complex case is the same as in \cite{bourgade2025fisher}, we only highlight the main steps for the reader's convenience. The first equality in \eqref{eqn:exact-complex} follows from Kostlan's theorem \cite{kostlan1992spectra} and the second from an asymptotic expansion of the Barne's $G$ function \cite[Eq. (A.6)]{voros1987spectral}. The first equality may also be derived from the fact that $X X^*$ is the Laguerre Unitary Ensemble and Selberg's integral formula.

The equality in \eqref{eqn:exact-real} is \cite[Theorem 4.8]{serebryakov2025schur} with an explicit error rate. That is, the first equality follows as in Appendix C of \cite{serebryakov2025schur} by relating the expectation to the Laguerre Orthogonal Ensemble and the asymptotics in the second equality follow as in the complex case. 
 \qed

\vspace{5 pt}

\noindent{\bf Proof of Theorem \ref{theo:momcomp}}.  For the first estimate, the case $\beta=2$ follows from \cite[Eq. (1.11)]{afanasiev2019correlation}. The case $\beta=1$  and $z\in \mathbb{D}_r^1$  follows from \cite[Eqs. (1.10), (2.12)]{afanasiev2022correlation}.\footnote{One can check that the estimate \cite[Eq. (1.10)]{afanasiev2022correlation} is uniform in the $\zeta_i$ and so one can take $\zeta_i = 0$.} Next, we consider $\beta=1$ and $z\in \rr$. Fix $z\in I_r$, and let $z_j:=z+\zeta_jN^{-1/2}$, with $|\zeta_j|\lesssim 1$. From \cite[Theorem 1.1]{afanasiev2020correlation} and the argument in the paragraph below \cite[Eq. (4.7)]{serebryakov2025schur} using \cite{tribe2014ginibre}, we have
\beq
\frac{\ee\left[\prod_{j=1}^k\big|\mathrm{det}(X-z_j)|^2\right]}{\prod_{j=1}^k\ee\left[\big|\mathrm{det}(X-z_j)|^2\right]}=C_{k,z}N^{k^2-k}\e^{\frac{k^2-k}{2}(1-z^2)^2\kappa_4}F(\zeta_1,\dots,\zeta_j)(1+o(1)),
\eeq
with $F(\zeta_1,\dots,\zeta_j)$ being an explicit function given by the ratio of a Pfaffian and a Vandermonde determinant. Then, by the paragraph below \cite[Eq. (4.7)]{serebryakov2025schur}, we see that $C_{k,z}=1$. Additionally, by the argument of \cite[Corollary 4.7]{serebryakov2025schur} we obtain
\beq
\lim_{|\zeta_1|,\dots,|\zeta_k|\to 0}F(\zeta_1,\dots,\zeta_j)=\frac{1}{\prod_{j=0}^{k-1}(2j)!}.
\eeq
Combining all the above, and using that all the estimates are uniform in $|\zeta_i|\lesssim 1$, we obtain
\beq
\label{eq:needbrealfac}
\ee\left[\big|\mathrm{det}(X-z)|^{2k}\right]=\left(\ee\left[\big|\mathrm{det}(X-z_j)|^2\right]\right)^k N^{k^2-k}\e^{\frac{k^2-k}{2}(1-z^2)^2\kappa_4}\frac{1}{\prod_{j=0}^{k-1}(2j)!}(1+o(1)).
\eeq
Finally, by \cite[Eq. (2.10)]{afanasiev2020correlation}, we have
\beq
\ee\left[\big|\mathrm{det}(X-z)|^2\right]= \e^{-N(1-z^2)}(2\pi N)^{1/2}(1+o(1)),
\eeq
which together with \eqref{eq:needbrealfac} concludes \eqref{eq:momentsreal}. \qed

\section{Local laws for products of resolvents}

In this section we will state several local laws for the quantity
$
\langle G^{z_1} ( \i \eta_1) A G^{z_2} ( \i \eta_2) B \rangle ,
$
where $A$ and $B$ are general deterministic matrices. We will show that this quantity is well-approximated by a deterministic quantity which we will denote by $\langle M(z_1, \eta_1, A, z_2, \eta_2) B \rangle$. That is, $M(z_1, \eta_1, A, z_2, \eta_2)$ is a matrix depending on $z_i, \eta_i$ and $A$. It is defined by \cite[Eq. (4.10)]{cipolloni2025optimal}. Note that it does not depend on $B$. Throughout this section we will also make use of the following notation
\beq
E_+ := (E_1 + E_2) = I, \qquad\quad E_- = (E_1 - E_2 ) .
\eeq

\subsection{2-G laws for one large parameter}

We first state a local law for two resolvents with one large parameter $1\lesssim |\eta_2| \le N^{C}$. The proof follows existing techniques but is technically not stated explicitly in the literature. In the regime where $|\eta_2|$ is large, we require an error that scales like $1 / |\eta_2|$; prior results in the literature just stated results with errors scaling like $ 1 / N$ for large $\eta_i$. We find a better estimate by carefully inspecting existing proofs. 

\bel \label{lem:2G-large} Let $X$ be an i.i.d. matrix. 
Let $r >0$ and let $z_1, z_2 \in \dd_r$. Let $\eps >0$ and $C>0$.  Let $A,B$ be two matrices with $\| A \| , \| B \| \leq 1 .$ Let $N^{\eps-1} \leq |\eta_1| \leq C$ and $C^{-1} \leq |\eta_2| \leq N^{C}$. Then,
\beq \label{eqn:2G-large}
\left| \langle G^{z_1} ( \i \eta_1 ) A G^{z_2} ( \i \eta_2 ) B \rangle - \langle M(z_1, \eta_1, A, z_2, \eta_2 ) B \rangle \right| \leq \frac{N^{\eps}}{N | \eta_1 \eta_2 | },
\eeq
with overwhelming probability
\eel
\proof
We use the notations $G_i := G_i^{z_i} ( \i \eta_i)$ and $M_i = M^{z_i} ( \i \eta_i)$. 
To keep the notation simpler throughout this proof we assume that $\eta_1,\eta_2>0$.  The other cases are completely analogous. We are going to compute large moments of the quantity on the LHS of \eqref{eqn:2G-large}. We first start with some prepatory estimates. In the following we will use the notation
\[
W:=\left(\begin{matrix} 0 & X \\
X^* & 0\end{matrix}\right).
\]
From \cite[Eq. (A.16)]{cipolloni2023mesoscopic} we have
\begin{equation}
\label{eq:staticeq}
\begin{split}
\left\langle (G_1A G_2 - M_{12}^{A})B \right\rangle
&=
\left\langle
M_1 A (G_2 - M_2)\mathfrak{B}
\right\rangle  -
\left\langle
M_1 \underline{W G_1 A G_2}\mathfrak{B}
\right\rangle \\
&\quad +
\left\langle
M_1 S[G_1-M_1]G_1 A G_2\mathfrak{B}
\right\rangle +
\left\langle
M_1 S[G_1 A G_2](G_2-M_2)\mathfrak{B}
\right\rangle.
\end{split}
\end{equation}
Here we denoted
\[
\mathfrak{B}:=\left(\mathcal{B}^{-1})^{*}[B^{*}] \right)^{*}, \qquad\quad \mathcal{B}[\cdot]:=1-M_1\mathcal{S}[\cdot]M_2, \qquad\quad \mathcal{S}[R]:=\langle R\rangle-\langle RE_-\rangle E_-,
\]
and
\[
\underline{W G_1 A G_2}:=\mathcal{S}[G_1]G_1 A G_2+\mathcal{S}[G_1AG_2]G_2.
\]
First we notice that we have
\begin{equation}
\lVert \mathfrak{B}\rVert=\lVert (\mathcal{B}^{-1})^{*}[B^{*}]\rVert \lesssim \lVert B\rVert\lesssim 1.
\end{equation}
This follows by \cite[Eq. (A.10)]{cipolloni2025optimal}, and that for $\eta_2\gg 1$ we have $m^{z_2}(\i\eta_2)\asymp \eta_2^{-1}$ and $u^{z_2}(\i\eta_2)\asymp \eta_2^{-2}$.

We now start estimating the various terms in the rhs. of \eqref{eq:staticeq} one by one. Recalling the single resolvent local law\footnote{Typically the local law is stated only for $ |\eta| \lesssim 1$. The extension to larger $\eta$ is straightforward. See, e.g. \cite[Theorem 2.2]{erdHos2019random} for a local law for $N^{-1+\eps}\le\eta\le N^C$ in the Hermitian case. The proof in the current case is analogous.}
\begin{equation}
\label{eq:sres2}
|\langle (G_i-M_i) B \rangle|\le \frac{N^{\xi/2} \lVert B\rVert}{N\eta_i},
\end{equation}
which holds with overwhelming probability for $N^{-1+\eps}\le\eta_i\le N^C$, we have 
\[
|\left\langle
M_1 A (G_2 - M_2)\mathfrak{B}
\right\rangle|\le \frac{N^\xi}{N\eta_2},
\]
with overwhelming probability. Similarly, using that $\mathcal{S}[G_1-M_1]=\langle G_1-M_1\rangle$ by the chiral symmetry, we have
\[
|\left\langle
M_1 \mathcal{S}[G_1-M_1]G_1 A G_2\mathfrak{B}
\right\rangle|\le \frac{N^\xi}{N\eta_1\eta_2},
\]
where, for $C_i$ with $|\lVert C_i\rVert\lesssim 1$, we also used
\begin{equation}
\label{eq:sschwarz}
|\langle G_1C_1G_2C_2\rangle|\le \frac{1}{\eta_2}[\langle|G_1|C_1C_1^*\rangle+\langle|G_1|C_2^*C_2\rangle]\lesssim \frac{\log N}{\eta_2},
\end{equation}
with overwhelming probability. In the second inequality of \eqref{eq:sschwarz} we used that (here $\eta_v:=\sqrt{\eta_1^2+v^2}$)
\[
\langle|G_1|C\rangle=\frac{2}{\pi}\int_0^\infty \langle \Im G^{z_1}(\i\eta_v)C\rangle \frac{\d v}{\eta_v}\lesssim \log N,
\]
with overwhelming probability.
Here we used that the regime $v\ge N^{100}$ can easily be removed by a simple norm bound.

Next, using the definition of $\mathcal{S}$ we have
\[
\begin{split}
\left\langle
M_1 \mathcal{S}[G_1 A G_2](G_2-M_2)\mathfrak{B}
\right\rangle&=\langle G_1AG_2\rangle\left\langle
M_1(G_2-M_2)\mathfrak{B}
\right\rangle \\
&\quad-\langle G_1 A G_2 E_-\rangle\left\langle
M_1 E_-(G_2-M_2)\mathfrak{B}
\right\rangle.
\end{split}
\]
Thus using the single resolvent local law \eqref{eq:sres2}, the bound \eqref{eq:sschwarz}, and $\eta_2\ge \eta_1$ we conclude
\[
|\left\langle
M_1 \mathcal{S}[G_1 A G_2](G_2-M_2)\mathfrak{B}
\right\rangle|\le \frac{N^\xi}{N\eta_1\eta_2}.
\]
Putting all this together we are left with
\begin{equation}
\label{eq:mainapprox}
\left\langle (G_1A G_2 - M_{12}^{A})B \right\rangle= -
\left\langle
M_1 \underline{W G_1 A G_2}\mathfrak{B}
\right\rangle +\mathcal{O}\left(\frac{N^\xi}{N\eta_1\eta_2}\right),\\
\end{equation}
with overwhelming probability. We can now begin estimating the moments of the LHS of \eqref{eqn:2G-large}. Fix any positive integer $p\ge 1$, then, by \eqref{eq:mainapprox}, we have
\[
\begin{split}
&\ee|\left\langle (G_1A G_2 - M_{12}^{A})B \right\rangle|^{2p} \\
&=-\ee\left\langle
M_1 \underline{W G_1 A G_2}
\left( (\mathcal{B}^{-1})^{*}[B^{*}] \right)^{*}
\right\rangle \overline{\left\langle (G_1A G_2 - M_{12}^{A})B \right\rangle} |\left\langle (G_1A G_2 - M_{12}^{A})B \right\rangle|^{2p-2} \\
&\quad+ \mathcal{O}\left(\frac{N^\xi}{N\eta_1\eta_2}\right)\ee|\left\langle (G_1A G_2 - M_{12}^{A})B \right\rangle|^{2p-1}.
\end{split}
\]
Using cumulant expansion in $W$ we obtain (here we define $C:=\mathfrak{B}M_1$ and we do only the complex case for simplicity)
\begin{equation}
\label{eq:cumexpnewapp}
\begin{split}
&\ee\left\langle\underline{W G_1 A G_2}C \right\rangle \overline{\left\langle (G_1A G_2 - M_{12}^{A})B \right\rangle} |\left\langle (G_1A G_2 - M_{12}^{A})B \right\rangle|^{2p-2}\\
&=\frac{p-1}{2N^2}\tilde{\sum}_{i\ne j} \langle G_1AG_2CE_iG_1^*B^*G_2^*A^*G_1^*E_j\rangle+\frac{p-1}{2N^2}\tilde{\sum}_{i\ne j} \langle G_1AG_2CE_iG_2^*A^*G_1^*B^*G_2^*E_j\rangle \\
&\quad+\frac{p-2}{2N^2}\tilde{\sum}_{i\ne j} \langle G_1AG_2CE_iG_1AG_2BG_1E_j\rangle+\frac{p-2}{2N^2}\tilde{\sum}_{i\ne j} \langle G_1AG_2CE_iG_2BG_1AG_2E_j\rangle \\
&\quad+\sum_{k\ge 2}\sum_{ab}\sum_{{\bm \alpha}\in\{ab,ba\}^k} \bigg\{ \frac{\kappa(ab,{\bm \alpha})}{k!} \\
& \quad \times \ee \left[\partial_{\bm \alpha}\left(\left\langle\Delta^{ab} G_1 A G_2 C \right\rangle \left\langle B^*(G_2^*A^* G_1^* - (M_{12}^A)^*) \right\rangle|\left\langle (G_1A G_2 - M_{12}^{A})B \right\rangle|^{2p-2} \right)\right] \bigg\}.
\end{split}
\end{equation}
Here $\Delta^{ab}$ denotes the matrix whose only nonzero entry is the $(a,b)$ one, $\partial_{ab}$ denotes the directional derivative in the direction $W_{ab}$, $\partial_{\bm \alpha}:=\partial_{\alpha_1}\dots \partial_{\alpha_k}$, with $\alpha_i\in \{ab,ba\}^k$, and $\kappa(ab,{\bm \alpha})$ denotes the $k+1$-th cumulant of the random variables $W_{ab}, W_{\alpha_1}$, $\dots$, $W_{\alpha_k}$. We estimate the first term in the second line by 
\begin{equation}
\frac{1}{2N^2}|\langle G_1AG_2CG_1^*B^*G_2^*A^*G_1^*\rangle|\lesssim \frac{1}{(N\eta_1\eta_2)^2},
\end{equation}
where we used Schwarz inequality and the single resolvent local law \eqref{eq:sres2}. All other terms in the second and third lines of \eqref{eq:cumexpnewapp} are analogous.

We next discuss the third order terms, i.e. the term $k=2$ in the last line of \eqref{eq:cumexpnewapp}. The terms that arise are in the form of $ ( \dots ) \times |\langle (G_1A G_2 - M_{12}^{A})B \rangle|^l $ for various combinations of sums of resolvent entries which we denote here by $( \dots )$ and some power $l$. We will show how to estimate  in detail  some exemplary terms which we list in \eqref{eq:exemplterms} below. All the remaining terms can be estimated analogously.  Note that we drop the factor $|\langle (G_1A G_2 - M_{12}^{A})B \rangle|^l$  for notational simplicity. We now estimate the following terms for $k=2$ one by one:
\begin{equation}
\label{eq:exemplterms}
\begin{split}
&\frac{1}{N^{3/2}}\langle \Delta^{ab}G_1\Delta^{ba}G_1AG_2C\rangle \langle B^*G_2^*A^*G_1^*\Delta^{ab}G_1^*\rangle \\
&\frac{1}{N^{3/2}}\langle \Delta^{ab}G_1AG_2C\rangle\langle B^*G_2^*A^* G_1^*\Delta^{ab}G_1^*\Delta^{ba}G_1^*\rangle \\
&\frac{1}{N^{3/2}}\langle \Delta^{ab}G_1\Delta^{ba}G_1\Delta^{ba}G_1AG_2C\rangle \\
&\frac{1}{N^{3/2}}\langle \Delta^{ab}G_1AG_2\Delta^{ba} G_2 C\rangle \langle B^*G_2^*A^*G_1^*\Delta^{ab}G_1^*\rangle \\
&\frac{1}{N^{3/2}}\langle \Delta^{ab}G_1AG_2C\rangle\langle B^*G_2^*\Delta ^{ba} G_2^*A^* G_1^*\Delta^{ab}G_1^*\rangle \\
&\frac{1}{N^{3/2}}\langle \Delta^{ab}G_1AG_2C\rangle\langle B^*G_2^*A^*G_1^*\Delta^{ab}G_1^*\rangle\langle B^*G_2^*A^*G_1^*\Delta^{ab}G_1^*\rangle.
\end{split}
\end{equation}\
In the remaining part of the proof all the estimates hold with overwhelming probability even if not stated explicitly. In the following we will often use (recall that $\eta_2\gtrsim 1$)
\begin{equation}
\label{eq:schwarzimplarge}
\big|(G_1AG_2C)_{ab}\big|\lesssim \frac{1}{\eta_1^{1/2}\eta_2}, \qquad\quad \big|(G_1^*B^*G_2^*A^*G_1^*)_{ab}\big|\lesssim \frac{1}{\eta_1\eta_2},
\end{equation}
which follow by Schwarz inequality, and
\begin{equation}
\label{eq:schwarzimplarge2}
\sum_{ab} |(G_1^*B^*G_2^*A^*G_1^*)_{ba}|\lesssim N^{3/2}\langle[G_1^*B^*G_2^*A^*G_1^*]^* G_1^*B^*G_2^*A^*G_1^*\rangle^{1/2}  \lesssim \frac{N^{3/2}}{\eta_1^{3/2}\eta_2}.
\end{equation}
Additionally, throughout the remaining part of the proof we will use that $|G_{ab}|\lesssim 1$. Using the first bound in \eqref{eq:schwarzimplarge} and \eqref{eq:schwarzimplarge2}, we estimate the first term in \eqref{eq:exemplterms} by 
\[
\begin{split}
\frac{1}{N^{7/2}}\left|\sum_{ab} (G_1)_{bb}(G_1AG_2C)_{aa}(G_1^*B^*G_2^*A^*G_1^*)_{ba}\right|&\lesssim \frac{1}{N^2\eta_1^2\eta_2^2}.
\end{split}
\]
For the second term in \eqref{eq:exemplterms}, we have
\[
\frac{1}{N^{7/2}}\left|\sum_{ab} (G_1AG_2C)_{ab} (G_1^*)_{bb}(G_1^*B^*G_2^*A^*G_1^*)_{aa}\right|\lesssim \frac{1}{N^2\eta_1^{3/2}\eta_2^2},
\]
where we used the second bound in \eqref{eq:schwarzimplarge}, and
\begin{equation}
\label{eq:lbhop}
\sum_{ab} |(G_1AG_2C)_{ab}|\lesssim \frac{N^{3/2}}{\eta_1^{1/2}\eta_2}.
\end{equation}
Next, for the third term in \eqref{eq:exemplterms}, using \eqref{eq:lbhop}, we estimate
\[
\frac{1}{N^{5/2}} \left|\sum_{ab}(G_1)_{bb}(G_1)_{ab}(G_1AG_2C)_{aa}\right|\lesssim  \frac{1}{N\eta_1^{1/2}\eta_2}.
\]
We then bound the fourth term in \eqref{eq:exemplterms} by 
\[
\frac{1}{N^{7/2}}\left|\sum_{ab} (G_1AG_2)_{bb}(G_2C)_{aa}(G_1^*B^*G_2^*A^*G_1^*)_{ba}\right|\lesssim \frac{1}{N^{7/2}}\cdot\frac{N^{3/2}}{\eta_1^{3/2}\eta_2}\cdot \frac{1}{\eta_1^{1/2}\eta_2}\cdot\frac{1}{\eta_2}\lesssim \frac{1}{N^2\eta_1^2\eta_2^2},
\]
where we used the first bound in \eqref{eq:schwarzimplarge}, \eqref{eq:schwarzimplarge2}, and $\lVert G_2\rVert \le \eta_2^{-1}$. Then, for the fifth term in \eqref{eq:exemplterms} we have
\[
\frac{1}{N^{7/2}}\left|\sum_{ab} (G_1AG_2C)_{ab} (G_2^*A^*G_1^*)_{aa}(G_1^*B^*G_2^*)_{bb}\right|\lesssim \frac{1}{N^{7/2}}\cdot \frac{N^{3/2}}{\eta_1^{1/2}\eta_2}\cdot \frac{1}{(\eta_1^{1/2}\eta_2)^2}\lesssim \frac{1}{(N\eta_1\eta_2)^2},
\]
where we used \eqref{eq:schwarzimplarge} and \eqref{eq:lbhop}. The last term of in \eqref{eq:exemplterms} is then bounded by
\[
\frac{1}{N^{9/2}}\left|\sum_{ab} (G_1AG_2C)_{ba}(G_1^*B^*G_2^*A^*G_1^*)_{ba} (G_1^*B^*G_2^*A^*G_1^*)_{ba}\right|\lesssim \frac{1}{N^{9/2}} \cdot \frac{N^{3/2}}{\eta_1^{1/2}\eta_2}\cdot \frac{1}{(\eta_1\eta_2)^2}\lesssim \frac{1}{(N\eta_1\eta_2)^3},
\]
where we used the second bound in \eqref{eq:schwarzimplarge} and \eqref{eq:lbhop}. This concludes the estimate of the third order terms.

Next, we consider the fourth order cumulants, i.e. the terms with $k=3$ in \eqref{eq:cumexpnewapp}. Most of these terms are easier to bound since there is an additional $N^{-1/2}$-gain from the cumulant, however, there are some terms when all the resolvent chains appearing are diagonal, potentially giving rise to larger terms. We now show how to estimate such two exemplary terms (as above we have dropped the $|\langle (G_1A G_2 - M_{12}^{A})B \rangle|^l$ factor for notational simplicity)
\begin{equation}
\label{eq:fcumul}
\begin{split}
&\frac{1}{N^2}\langle \Delta^{ab}G_1AG_2\Delta^{ba} G_2C\rangle \langle B^*G_2^*\Delta^{ab} G_2^*A^*G_1^*\Delta^{ba}G_1^*\rangle \\
&\frac{1}{N^2}\langle \Delta^{ab}G_1\Delta^{ba}G_1AG_2C\rangle \langle B^*G_2^*\Delta^{ab} G_2^*A^*G_1^*\Delta^{ba}G_1^*\rangle.
\end{split}
\end{equation}
Using the first bound in \eqref{eq:schwarzimplarge} three times, we estimate the first term in \eqref{eq:fcumul} by
\[
\frac{1}{N^4}\left|\sum_{ab} (G_1AG_2)_{bb}(G_2C)_{aa}(G_2^*A^*G_1)_{bb} (G_1^*B^*G_2^*)_{aa}\right|\lesssim \frac{1}{N^2\eta_1^{3/2}\eta_2^4}.
\]
Here we also used that $\lVert G_2\rVert\le \eta_2^{-1}$. Then, for the second term in \eqref{eq:fcumul} we have
\[
\frac{1}{N^4}\left|\sum_{ab} (G_1)_{bb} (G_1AG_2C)_{aa}(G_2^*A^*G_1)_{bb} (G_1^*B^*G_2^*)_{aa}\right|\lesssim \frac{1}{N^2\eta_1^{3/2}\eta_2^3},
\]
where we used the first bound in \eqref{eq:schwarzimplarge} three times. Finally, one can see that all the terms emerging for higher order cumulants, i.e the terms for $k\ge 4$ in \eqref{eq:cumexpnewapp}, can by estimated analogously by
\[
\sum_{l=1}^{2p} \left(\frac{1}{N\eta_1\eta_2}\right)^l \ee |\left\langle (G_1A G_2 - M_{12}^{A})B \right\rangle|^{2p-l}+N^{-D},
\]
where $D>0$ is an arbitrarily large constant. The estimate for these terms is easier than the one presented above as each additional cumulant gains a factor $N^{-1/2}$; we omit these details to keep the presentation shorter. Finally, combining all the previous estimates we have thus proven
\[
\ee|\left\langle (G_1A G_2 - M_{12}^{A})B \right\rangle|^{2p}\lesssim \sum_{l=1}^{2p} \left(\frac{1}{N\eta_1\eta_2}\right)^l \ee |\left\langle (G_1A G_2 - M_{12}^{A})B \right\rangle|^{2p-l}+N^{-D},
\]
which readily implies
\[
\ee|\left\langle (G_1A G_2 - M_{12}^{A})B \right\rangle|^{2p}\lesssim \frac{1}{(N\eta_1\eta_2)^{2p}},
\]
thus concluding the proof.
\qed

\subsection{2-G laws  on short scales}

We will prove the following.
\bep \label{prop:2G-law}
Let $X$ be an i.i.d. matrix. Let $r >0$ and let $z_1, z_2 \in \dd_r$. Let $\eps >0$ and let $A, B$ satisfy $\|A\|, \| B \| \leq 1$. Then:
\begin{enumerate}[label=(\roman*), font=\normalfont]
\item \label{it:2G-1} Let $N^{\eps -1} \leq |\eta_1|, |\eta_2| \leq 1$. Then with overwhelming probability
\beq \label{eqn:2G-1}
\left| \langle G^{z_1} ( \i \eta_1 ) A G^{z_2} ( \i \eta_2 ) B \rangle - \langle M(z_1, \eta_1, A, z_2, \eta_2 ) B \rangle \right| \leq \frac{N^{\eps}}{N | \eta_1 \eta_2 | }.
\eeq
\item \label{it:2G-2} Suppose that either $z_1 = z_2$ or $|z_1 -z_2| \geq c$ for some $c>0$. Then with overwhelming probability
\beq \label{eqn:2G-2}
\left| \langle G_{z_1} ( \i \eta_1 ) A G^{z_2} ( \i \eta_2 ) B \rangle - \langle M(z_1, \eta_1, A, z_2, \eta_2 ) B \rangle \right| \leq \frac{ ( \log N)^{10}}{N | \eta_1 \eta_2| },
\eeq
as long as $( \log N)^{10} /N \leq |\eta_1| \leq N^{-\eps}$ and $( \log N)^{10} /N \leq | \eta_2 | \leq N^{C}$. 
\end{enumerate}
\eep
  We remark that the first estimate was essentially proven in \cite{cipolloni2025optimal}; however, it was stated with an error involving the quantity $(N \min\{ |\eta_1|, |\eta_2| \} )^{-1}$, which is too weak for our purposes.

We will prove the above via a dynamical method. We therefore require the following. We use the notation $G^z ( \i \eta ; M) := \left( H^z (M) - \i \eta \right)^{-1}$ with $H^z(M)$ as in \eqref{eqn:Hermit-def}. 
\bep \label{prop:2G-compare}
Let $T = N^{-\xi}$ for some $\xi \in (0, \frac{1}{100} )$. Let $X$ and $Y$ be two i.i.d. matrices that $T$-match as in Definition \ref{def:T-match}. Let $A_1, A_2$ be matrices s.t. $\|A_1\|, \| A_2 \| \leq 1$. Let $D \geq 1$. There is an $L \geq 1$ depending only on $D$ and $\xi$ so that the following holds. Let $0 < r < 1$ and $C >0$. Let $N^{-1} \leq | \eta_1| , |\eta_2| \leq N^{C}$ and let $K \geq 1$. Then,
\begin{align}
& \pp\left[ \left| \langle G^{z_1} ( \i \eta_1 ; X  ) A_1 G^{z_2} ( \i \eta_2 ; X ) A_2 \rangle - \langle M(z_1, \eta_1, A_1, z_2, \eta_2 ) A_2 \rangle \right|  \geq K + L \right] \notag\\
\lesssim &  L \times  \pp\left[ \left| \langle G^{z_1} ( \i \eta_1 ; Y  ) A_1 G^{z_2} ( \i \eta_2 ; Y ) A_2 \rangle - \langle M(z_1, \eta_1, A_1, z_2, \eta_2 ) A_2 \rangle \right|  \geq K  \right] + N^{-D}
\end{align}
\eep
\proof   This is similar to the proof of \cite[Proposition 3.8]{cipolloni2025maximum}. Let us adopt the notation and concepts of Section \ref{sec:gfct-prelim}. For some matrices $A_1, A_2$ obeying $\| A_i\| \leq 1$ let us define the function $Q(M)$ on the space of $N \times N$ matrices by
\beq
Q(M) := (N \eta_1 \eta_2 ) \langle (H^{z_1} (M) - \i \eta_1 )^{-1} A_1 (H^{z_2} (M) - \i \eta_2 )^{-1} A_2 \rangle .
\eeq

If $X$ is an i.i.d. matrix we claim that with overwhelming probability,
\beq \label{eqn:Q-der-bd}
\sup_{ | \theta| \leq N^{1/10-1/2}} \left| \left( \del_{M_{ab}}^j \del_{\bar{M}_{ab}}^k Q\right) ( X(\theta_{ab} ) ) \right| \leq N^{\eps}
\eeq
for any $\eps$ and $j, k \geq 0$. We will prove this estimate in a moment. Let us first state how \eqref{eqn:Q-der-bd} allows us to complete the proof. 

 As in the proof of Proposition \ref{prop:gfct-kpt-compare}, we let $U_i$ be a sequence of $N^2+1$ matrices such that $U_0 = X$ and $U_{N^2} = Y$. Define
\beq
p_k := \sup_{0 \leq j\leq N^2} \pp\left[ |Q(U_j) - \langle M(z_1, \eta_1, A_1, z_2, \eta_2) A_2 \rangle | \geq K + k  \right].
\eeq
By following the proof of \cite[Proposition 3.8]{cipolloni2025maximum} we obtain for any $D>0$
\beq
p_k \lesssim N^{\eps} (T + N^{-1/2} ) p_{k-1} + N^{-D} +   \pp\left[ |Q(Y) - \langle M(z_1, \eta_1, A_1, z_2, \eta_2) A_2 \rangle | > K  \right].
\eeq
Iterating this a constant order number of times depending on $\xi$ and $D$ yields the claim.

It remains to discuss \eqref{eqn:Q-der-bd}. Denoting here $G_i(M) = G_i := (H^{z_i} (M)- \i \eta_i )^{-1}$ we have by direct computation that the derivatives are linear combinations of terms of the form
\beq \label{eqn:Q-der-bd-a1}
(\eta_1 \eta_2) (G_1A_1G_2)_{x_1y_1}(G_2BG_1)_{x_2,y_2}\prod_{i=1}^k (G_{l_i} B_i)_{x_iy_i} , \qquad 
(\eta_1\eta_2) (G_1A_1G_2A_2G_1)_{x_1y_1}\prod_{i=1}^k (G_{l_i})_{x_iy_i} ,
\eeq
for some $k$. Here $x_i, y_i \in \{ a, b+N \}$, $B, B_i\in \{A_2,I\}$ and $l_i \in \{1, 2\}$ (note that only one of the $B,B_i$ is not the identity). From \cite[Appendix A]{cipolloni2025optimal2} we have that for any unit vector ${\bf v}$ and index $n \in [\![ 2N ]\!]$ 
\beq \label{eqn:Q-der-bd-a5}
\sup_{ | \theta | \leq N^{1/10-1/2}} | [G_i (X(\theta_{ab} ) ) {\bf v}]_n | \leq N^{\eps},
\eeq
with overwhelming probability. This estimates the terms in the products in \eqref{eqn:Q-der-bd-a1}. For the remaining terms we have by Cauchy-Schwarz and the Ward identity (with $\delta_x$ the vector $(\delta_x ) (i) = \delta_{ix}$)
\beq \label{eqn:Q-der-bd-a2}
\left| (G_1 A_1 G_2 )_{xy} \right| \leq \left( \| G_1 \delta_x \|_2^2 \| G_2 \delta_y \|_2^2 \right)^{1/2} \lesssim \left( \frac{ \Im[ (G_1)_{xx}] \Im[ (G_2)_{yy}]}{\eta_1 \eta_2} \right)^{1/2}
\eeq
as well as
\begin{align} \label{eqn:Q-der-bd-a3}
| (G_1 A_1 G_2 A_2 G_1)_{xy} | \leq & \frac{1}{\eta_2}  \left( \| A_1^* (G_1)^* \delta_x \|_2^2 \| A_2 G_1 \delta_y \|_2^2 \right)^{1/2} \lesssim \frac{1}{\eta_2}  \left( \|  (G_1)^* \delta_x \|_2^2 \|  G_1 \delta_y \|_2^2 \right)^{1/2} \notag\\
\leq & \frac{1}{ \eta_1 \eta_2}\left(  \Im[ (G_1)_{xx}] \Im[ (G_1)_{yy}]\right)^{1/2},
\end{align}
where we used $\| G_2 \| \lesssim \eta_2^{-1}$. The estimate \eqref{eqn:Q-der-bd} follows from estimating all of the terms in \eqref{eqn:Q-der-bd-a1} using \eqref{eqn:Q-der-bd-a5}, \eqref{eqn:Q-der-bd-a2}, and \eqref{eqn:Q-der-bd-a3}. \qed

\subsubsection{Proof of Proposition \ref{prop:2G-law}\ref{it:2G-1}}

We first establish the result for matrices with a Gaussian component. Later we will apply Proposition \ref{prop:2G-compare} to remove it.

Consider the flow
\beq
\d X_t=-\frac{1}{2}X_t\d t+\frac{\d B_t}{\sqrt{N}}, \qquad\quad X_0=X,
\eeq
where $(B_t)_{ij}$ are i.i.d. standard real/complex Brownian motion according to the symmetry class of $X_0$. Consider the characteristic equations
\beq
\label{eq:charapp}
\partial_t \eta_t=-\Im m^{z_t}(\i\eta_t), \qquad\quad \partial_t z_t=-\frac{z_t}{2}.
\eeq
Let $\xi \in (0, \frac{1}{100})$, let $T =  N^{-\xi}$, and define $\eta_{i,t}, z_{i,t}$ as the solutions of the characteristic equations \eqref{eq:charapp} so that at time $T$ we have $\eta_{i,T}=\eta_i$, $z_{i,T}=z_i$ (see e.g. \cite[Lemma 4.2]{cipolloni2025optimal}). Due to cyclicity of the trace, it is no loss of generality to assume $|\eta_1| \leq |\eta_2|$. We will furthermore assume that $N^{10\xi-1} \leq |\eta_1| \leq N^{-\xi}$.

In a moment we will compute the evolution of observables on the LHS of \eqref{eqn:2G-1} under the characteristic flow. We require some notation.
Let $G_{i,t}:=(H_t^{z_{i,t}}-\i\eta_{i,t})^{-1}$, where $H_t^{z_{i,t}}$ is the Hermitization of  $X-z_{i,t}$. For any matrix $A$ denote $M_{ij,t}^A = M(z_{i, t} , \eta_{i, t} , A , z_{j, t} , \eta_{j, t} )$ where the deterministic approximation is as above. Similarly denote $M_{i, t} = M^{z_{i, t} } ( \i \eta_{i, t} )$. In particular it is straightforward to see that $|\eta_{1, t} | \lesssim | \eta_{2, t}|$ for all $t\in[0, T]$. We will use this frequently in order to simplify various error terms.

Then by It\^{o}'s lemma we have (see e.g. \cite[Eq. (4.20)]{cipolloni2025optimal}) for any matrices $A_1, A_2$, (in a moment it will be clear that we need to consider the evolution for choices of $A_i$ other than just $A, B$)
\beq
\label{eqn:full-2G-flow}
\begin{split}
\d \langle (G_{1,t}A_1G_{2,t}-M_{12,t}^{A_1})A_2\rangle&=\sum_{a,b=1}^{2n}\partial_{ab}\langle G_{1,t}A_1G_{2,t}A_2\rangle\frac{\d \mathfrak{B}_{ab,t}}{\sqrt{n}}+\langle (G_{1,t}A_1G_{2,t}-M_{12,t}^{A_1})A_2\rangle\d t \\
&\quad+2\tilde{\sum}_{ij}\langle (G_{1,t}A_1G_{2,t}-M_{12,t}^{A_1})E_i\rangle\langle M_{21,t}^{A_2}E_j\rangle\d t \\
&\quad+2\tilde{\sum}_{ij}\langle M_{12,t}^{A_1} E_i\rangle\langle (G_{2,t}A_2G_{1,t}-M_{21,t}^{A_2})E_j\rangle\d t \\
&\quad+2\tilde{\sum}_{ij}\langle (G_{1,t}A_1G_{2,t}-M_{12,t}^{A_1})E_i\rangle\langle (G_{2,t}A_2G_{1,t}-M_{21,t}^{A_2})E_j\rangle\d t \\
&\quad+ \langle (G_{1,t}-M_{1,t})\rangle\langle G_{1,t}^2A_1G_{2,t}A_2\rangle\d t + \langle (G_{2,t}-M_{2,t})\rangle\langle G_{2,t}^2A_2G_{1,t}A_1\rangle\d t \\
&\quad+\frac{\bm1_{\{\beta=1\}}}{N}\tilde{\sum}_{ij}\langle G_{1,t}^\mathfrak{t} E_i G_{1,t}A_1G_{2,t}A_2G_{1,t}E_j\rangle \d t\\
&\quad+\frac{\bm1_{\{\beta=1\}}}{N}\tilde{\sum}_{ij}\langle [G_{1,t}A_1G_{2,t}]^\mathfrak{t}E_1 G_{2,t}A_2G_{1,t}E_j\rangle \d t\\
&\quad+\frac{\bm1_{\{\beta=1\}}}{N}\tilde{\sum}_{ij}\langle G_{2,t}^\mathfrak{t} E_i G_{2,t}A_2G_{1,t}A_1G_{2,t}E_j\rangle \d t.
\end{split}
\eeq
Here $\mathfrak{B}_t = H^0 (B_t)$ is the Hermitization of $B_t$, $\partial_{ab}$ denotes the partial derivative in the directions of the entries of the Hermitization of $X_t$, and $\mathfrak{B}_{ab,t}$ denote the entries of $\mathfrak{B}_t$. Let us denote the Martingale term by $\d \mfe_t$, and the last four lines of \eqref{eqn:full-2G-flow} by $\mfh_t \d t$. Arguing as in \cite[(5.30)-(5.37)]{cipolloni2023universality} we have the estimate, 
\beq \label{eqn:et-qv}
\d [ \mfe_t, \mfe_t] \lesssim \left( \frac{ 1}{ N |\eta_{1, t} |^3 |\eta_{2, t}|^2} \right) \d t .
\eeq
 In a similar manner, we find 
\beq \label{eqn:ht-est}
| \mfh_t | \leq \frac{ ( \log N)^2}{N | \eta_{1, t} |^2 | \eta_{2, t}|}.
\eeq
We estimate a few representative terms contributing to $\mfh_t$ to illustrate the above two estimates. In particular, in order to estimate the terms on the fourth line of \eqref{eqn:full-2G-flow} we have,
\begin{align} \label{eqn:ht-est-1}
| \langle G_{1, t}^2 A_1 G_{2, t} A_2 \rangle |& \leq \left( \langle  G_{1, t}^* A_2^* A_2 G_{1, t} \rangle \langle G_{1, t} A_1 |G_{2, t} |^2 A_1^* G_{1, t}^* \rangle \right)^{1/2} \lesssim \frac{1}{ |\eta_{1, t} \eta_{2, t}|}, \notag\\
| \langle G_{2, t}^2 A_2 G_{1, t} A_1 \rangle |  & \leq \left( \langle G_{2, t}^* A_1^* A_1 G_{2, t} \rangle \langle G_{1, t}^* A_2^* |G_{2, t}|^2 A_2 G_{1, t}  \rangle \right)^{1/2}  \lesssim \frac{1}{ |\eta_{1, t} \eta_{2, t} |}.
\end{align}
Here we used Cauchy-Schwarz and Theorem \ref{thm:ll-2}. 
The other terms contributing to $\mfh_t$ are handled via Cauchy-Schwarz and similar arguments to \eqref{eqn:ht-est-1}. For later use (i.e., for the  proof of Proposition \ref{prop:2G-law}\ref{it:2G-2}) we make the following observations about the estimates \eqref{eqn:et-qv} and \eqref{eqn:ht-est}.  First, since they rely only on Theorem \ref{thm:ll-2} both estimates are valid as long as $|\eta_1| \geq ( \log N)^2 / N$. Furthermore, by inspecting the proof, the estimates hold even if $1 \leq |\eta_2| \leq N^{C}$. 

With these estimates in mind we now write \eqref{eqn:full-2G-flow} as 
\beq
\begin{split}
\label{eqn:eqred}
&\d \langle (G_{1,t}A_1G_{2,t}-M_{12,t}^{A_1})A_2\rangle = \langle (G_{1,t}A_1G_{2,t}-M_{12,t}^{A_1})A_2\rangle\d t \\
&\qquad\quad+ \langle G_{1,t}A_1G_{2,t}-M_{12,t}^{A_1}\rangle\langle M_{21,t}^{A_2}\rangle\d t -\langle (G_{1,t}A_1G_{2,t}-M_{12,t}^{A_1})E_-\rangle\langle M_{21,t}^{A_2}E_-\rangle\d t\\
&\qquad\quad+\langle M_{12,t}^{A_1} \rangle\langle G_{2,t}A_2G_{1,t}-M_{21,t}^{A_2}\rangle\d t-\langle M_{12,t}^{A_1} E_- \rangle\langle (G_{2,t}A_2G_{1,t}-M_{21,t}^{A_2})E_-\rangle\d t \\
&\qquad\quad+2\tilde{\sum}_{ij}\langle (G_{1,t}A_1G_{2,t}-M_{12,t}^{A_1})E_i\rangle\langle (G_{2,t}A_2G_{1,t}-M_{21,t}^{A_2})E_j\rangle\d t+ \mfh_t\d t +\d \mfe_t.
\end{split}
\eeq
We are going to apply Gr{\"o}nwall inequality to the above system of equations, and so we require estimates for the coefficients above. Let us first discuss the case that $\eta_1 \eta_2 < 0$. Then for any matrix $A$ we have
\beq
\label{eqn:gronwall-coefficients-1}
\left| \langle M_{12,t}^A \rangle \right| \lesssim \frac{ \| A \|}{|z_1  -z_2 |^2 + |\eta_{1, t}| + |\eta_{2, t}| } 
\qquad \big|\langle M_{12,t}^{A}E_-\rangle\big|\lesssim \| A \| \frac{ |z_{1,t}-z_{2,t}|}{|z_{1,t}-z_{2,t}|^2+|\eta_{1,t}|+|\eta_{2,t}|}. 
\eeq
The first follows from \cite[Eqs. (4.16-4.17)]{cipolloni2025optimal}, and the second from \cite[Eqs. (4.17),(4.27)]{cipolloni2025optimal}. We also have, if $\| A \| \leq 1$, 
\beq \label{eqn:gronwall-coefficients-2}
| \langle M_{12, t}^A \rangle| \leq | \langle M_{12, t}^{E_+} \rangle| +C_1 \left( \frac{ |z_1 - z_2|}{ |z_1-z_2|^2 + |\eta_{1, t}| + |\eta_{2, t}|} + 1 \right)
\eeq
for some $C_1 >0$. This follows from \cite[Eq. (A.35)]{cipolloni2025optimal}. 
From \eqref{eqn:eqred} we see we need to control quantities involving two resolvents tested against several different matrices, simultaneously. Therefore, introduce the quantity,
\beq \label{eqn:Yt-def}
Y_t:=\max_{A_1,A_2\in\{B_1,B_2,E_-, I\}}\big|\langle (G_{1,t}A_1G_{2,t}-M_{12,t}^{A_1})A_2\rangle\big|,
\eeq
and the stopping time,
\beq
\tau := \inf\left\{ t \geq 0 : Y_t = \frac{ N^{3 \xi}}{N | \eta_{1, t} \eta_{2, t} |} \right\} \wedge T.
\eeq
Our goal is to prove that $\tau = T$ with overwhelming probability. We now use the stopping time to integrate several terms on the RHS of \eqref{eqn:eqred}. For $t < \tau$ and any choices of the $A_i$ in the definition of $Y_t$ we have
\begin{align}
& \int_0^t | \langle (G_{1,s}A_1G_{2,s}-M_{12,s}^{A_1})A_2\rangle | +\left| \langle (G_{1,s}A_1G_{2,s}-M_{12,s}^{A_1})E_-\rangle\langle M_{21,s}^{A_2}E_-\rangle \right| \d s \notag\\
\lesssim & \int_{0^t} \frac{N^{3 \xi}}{N | \eta_{1, s} \eta_{2, s}|} + \frac{ N^{3 \xi}}{N | \eta_{1, s} \eta_{2, s} | | \eta_{2, s} |^{1/2}} \d s \lesssim  \frac{( \log N) N^{3 \xi}}{N |\eta_{2, t} |^{3/2}} \leq \frac{ N^{3 \xi } N^{-\xi/3}}{ N | \eta_{1, t} \eta_{2, t} | } 
\end{align}
where we used that $| \eta_{1, t} | \lesssim N^{-\xi}$ for all $t$. We obtain a  similar estimate for the second term on the third line of \eqref{eqn:eqred}. Here we repeatedly used the fact that $t \to |\eta_{i, t}|$ is a decreasing function as well as the fact that $\int_0^t | \eta_{i, s} |^{-1} \d s \lesssim \log N$. We will repeatedly use these two facts in the rest of the proof.

In a similar manner, the  product of the second term on the RHS of \eqref{eqn:gronwall-coefficients-2} and the first terms on the second and third line of \eqref{eqn:eqred} are bounded by
\beq
 \int_0^t \left(\frac{N^{3\xi} }{ N |\eta_{1, s} \eta_{2, s}|}\right) \times \left( 1 + \frac{ |z_1-z_2|}{ ( |z_1-z_2|^2 + |\eta_{2, s} | ) } \right) \d s \lesssim \frac{ N^{3 \xi} N^{-\xi/3}}{N | \eta_{1, t} \eta_{2, t} |}
\eeq
for $t < \tau$. Finally, for $t < \tau$, the last line of \eqref{eqn:eqred} is bounded by
\beq
\int_0^t \frac{ N^{6 \xi}}{N^2 ( \eta_{1, s} \eta_{2, s} )^2} \d  s \leq \frac{N^{ 7 \xi}}{ N^2 (\eta_{2, t} )^2 | \eta_{1, t} | } \leq \frac{ N^{-\xi}}{N | \eta_{1, t} \eta_{2, t  }|}
\eeq
using that $| \eta_{1, t} | \geq N^{10 \xi-1}$. 

Collecting these estimates as well as the estimates for $\mfe_t$ and $\mfh_t$ we find that with overwhelming probability,
\beq
Y_t \leq \frac{N^{3 \xi} N^{-\xi/4}}{N | \eta_{1, t} \eta_{2, t} |} + \int_0^t \left( 2 \big|\langle M_{12,s}^{E_+}\rangle\big|  \right) Y_s \d s .
\eeq
for all $t < \tau$. Here we also used \cite[Eq. (3.26)]{cipolloni2025optimal} as well as the fact that $|\eta_1|, |\eta_2| \gtrsim N^{-\xi}$ to bound $Y_0$. From a direct application of Gronwall's inequality as well as $\langle M_{12,t}^{E_+ } \rangle = \langle M_{21, t}^{E_+ } \rangle$ and the estimate (see \cite[Eq. (4.29)]{cipolloni2025optimal})
\beq
\exp \left( 2 \int_s^t \left| \langle M_{12, r}^{E_+} \rangle \right| \d r \right) \lesssim \left( \frac{ |z_1 -z_2 |^2 + | \eta_{2, s} | }{ |z_1-z_2|^2 + | \eta_{2, t} |} \right)^2
\eeq
we find that with overwhelming probability that for all $t < \tau$,
\begin{align}
    Y_t & \lesssim \frac{ N^{3 \xi } N^{ -\xi/4}}{N | \eta_{1, t} \eta_{2, t} | } + \frac{ N^{3 \xi} N^{-\xi/4}}{N ( |z_1 -z_2|^2 + |\eta_{2, t} | )^2} \int_0^t \frac{ |z_1 -z_2|^2 + | \eta_{2, s} |}{| \eta_{1, s} \eta_{2, s} |} \d s  \lesssim \frac{ N^{3 \xi} N^{- \xi /5 }}{N | \eta_{1, t} \eta_{2, t} | },
\end{align}
with the second inequality following in a straightforward manner. We therefore conclude that $\tau =T$ with overwhelming probability. 

The case when $\eta_1 \eta_2 >0$  differs only in notation. In particular, for \eqref{eqn:gronwall-coefficients-1} the LHS of the inequalities are switched, and \eqref{eqn:gronwall-coefficients-2} holds with $\langle M_{12, t}^A \rangle $ and $\langle M^{E_+}_{12, t} \rangle$ replaced by $\langle M_{12, t}^A E_- \rangle$ and $\langle M^{E_-}_{12, t} E_- \rangle$, respectively.

We have so far proven that for an i.i.d. matrix with a Gaussian component of size $N^{-\xi}$ we have
\beq
\left|\langle (G^{z_1} ( \i \eta_1 ) A G^{z_2} ( \i \eta_2 ) - M_{12}^A ) B \rangle \right| \leq \frac{ N^{5 \xi}}{N |\eta_{1} \eta_{2} |}
\eeq
for all $N^{10\xi-1} \leq |\eta_1|, |\eta_2| \leq N^{-\xi}$, with overwhelming probability. By  \cite[(3.26)]{cipolloni2025optimal}  this also holds for all $ N^{10\xi-1} \leq |\eta_1|, |\eta_2| \leq 1$. By Proposition \ref{prop:2G-compare}, we can conclude that the above inequality holds with $N^{5 \xi}$ replaced by $N^{6 \xi}$, for any i.i.d. matrix. The claim now follows. \qed

\subsubsection{Proof of Proposition \ref{prop:2G-law}\ref{it:2G-2}}

The argument proceeds similarly to the proof of Proposition \ref{prop:2G-law}\ref{it:2G-1}. Let $\xi \in (0, \frac{1}{100})$ and $T = N^{-\xi}$.  We will assume that $|\eta_1| \leq | \eta_2|$. We assume that $N^{-2\xi} \geq |\eta_1| \geq ( \log N)^{10} / N$. As above,  consider the characteristics $\eta_{i, t}$ with $\eta_{i, T} = \eta_i$ as defined in \eqref{eq:charapp}. Let us first consider the case that $\eta_1 \eta_2 <0$.

Our starting point is then the equation \eqref{eqn:eqred}, with the estimates \eqref{eqn:et-qv} and \eqref{eqn:ht-est} in hand for the Martingale and forcing terms, respectively.  Define $Y_t$ as in \eqref{eqn:Yt-def}, and consider the stopping time,
\beq
\tau = \inf \left\{ t \geq 0 : Y_t = \frac{ ( \log N)^5}{N |\eta_{1, t} \eta_{2, t}|} + \frac{N^{\xi/4}}{N|\eta_{1, 0} \eta_{2, 0} |} \right\} \wedge T.
\eeq
We now integrate several terms on the RHS of \eqref{eqn:eqred} up to the stopping time. First, the term on the last line of \eqref{eqn:eqred} is bounded by, 
\begin{align}
\int_0^t \frac{ ( \log N)^{10}}{N^2| \eta_{1, s} \eta_{2, s} |^2} + \frac{ N^{\xi/2}}{N^2 | \eta_{1, 0} \eta_{2, 0} |^2} \d s \lesssim \frac{ ( \log N)^{11}}{N^2 |\eta_{1, t}  ( \eta_{2, t} )^2 |} + \frac{ N^{\xi/2}}{N^2 | \eta_{1, 0}  ( \eta_{2, 0} )^2|} \lesssim \frac{ \log N}{N |\eta_{1, t} \eta_{2, t} |}  
\end{align}
where  used that $|\eta_{2, t}| \geq ( \log N)^{10} / N$ and also that $|\eta_{1, 0}|, |\eta_{2, 0} | \gtrsim T$. Similarly, the first, third and fifth terms on the RHS of \eqref{eqn:eqred} are bounded by,
\begin{align} \label{eqn:2G-proof-c1}
 & \int_0^t \left( \frac{ ( \log N)^5}{N | \eta_{1, s} \eta_{2, s} | } + \frac{ N^{\xi/4}}{N |\eta_{1, 0} \eta_{2, 0} |} \right)\left( \frac{ |z_1 -z_2|}{ | z_1-z_2|^2 + |\eta_{2, s} |} + 1 \right) \d s \notag\\
\lesssim & ( \log N)^6 \left( \frac{1}{N |\eta_{2, t}|} + \frac{1}{ N | \eta_{2, t} |^{3/2}} \right) + \frac{ N^{\xi/4} T^{1/2}}{N | \eta_{1, 0} \eta_{2, 0} |} \leq \frac{1}{N | \eta_{1, t} \eta_{2, t} | }.
\end{align}
In the last inequality we used $|\eta_{1, t} | \lesssim  N^{-\xi}$, and in the first inequality we used $\int_0^T |\eta_{2, s} |^{-1/2} \d s \lesssim T^{1/2}$. 

We now turn to the second and fourth terms on the RHS of \eqref{eqn:eqred}. If either $|z_1-z_2| \geq c$ or $|\eta_2| \geq 1$, then these terms can be bounded as in \eqref{eqn:2G-proof-c1} using the first estimate in \eqref{eqn:gronwall-coefficients-1}. If $z_1 =z_2$ and $|\eta_2| \leq 1$ then we write the fourth term as
\begin{align} \label{eqn:2G-proof-c2}
    \langle G_{2, t} A_2 G_{1, t} \rangle &= \langle G_{1, t} G_{2, t} A_2 \rangle = \frac{1}{ \i ( \eta_1 - \eta_2)} \langle (G_{1, t} - G_{2, t} ) A_2 \rangle \notag\\
    &= \frac{ \langle (M_{1, t} - M_{2, t} ) A_2 \rangle }{ \i ( \eta_1 - \eta_2)} + \O \left( \frac{ (\log N)^2}{ N ( |\eta_{1, t} \eta_{2, t} | )} \right),
\end{align}
with overwhelming probability. In the last inequality we used \cite[Theorem 3.1]{cipolloni2025optimal2}, as well as the fact that $| \eta_{2, t} - \eta_{1, t}| \geq |\eta_{2, t}|$ since they have opposite signs. Together with a similar identity for $M_{12}$ we can conclude that the second and fourth terms on the RHS of \eqref{eqn:eqred} are bounded by 
\beq
\int_0^t \frac{ ( \log N)^2}{N |\eta_{1, s} ( \eta_{2, s} )^2 |} \d s \lesssim \frac{ ( \log N)^3}{N  |\eta_{1, t} \eta_{2, t} |} \d s.
\eeq
Therefore in every case we have for all $t < \tau$ that
\beq
Y_t \leq \frac{N^{\xi/5}}{N | \eta_{1, 0} \eta_{2, 0} |} + \frac{ ( \log N)^4}{N | \eta_{1, t} \eta_{2, t}|}.
\eeq
Therefore we conclude that $\tau = T$. Since $|\eta_{1}| \leq N^{-\xi/2} |\eta_{1, 0} |$ we conclude that
\beq
|\langle (G_{1, T} A_1 G_{2, T} - M_{12, T}^{A_1} ) A_2 \rangle | \leq \frac{ ( \log N)^5}{N | \eta_{1, T} \eta_{2, T} |},
\eeq
with overwhelming probability. We wrap up the proof in the same manner as in Proposition \ref{prop:2G-law}\ref{it:2G-1}. This concludes the case that $\eta_1 \eta_2 <0$. The case when either $|z_1 -z_2| \geq c$ or $|\eta_2| \geq 1$ is analogous to the above. If $z_1=z_2$, then the roles of the second and fourth terms on the RHS of \eqref{eqn:eqred} are interchanged with the roles of the third and fifth terms. As discussed in the proof of Proposition \ref{prop:2G-law}\ref{it:2G-1}, the estimates in \eqref{eqn:gronwall-coefficients-1} are swapped. Therefore, one can integrate the second and fourth terms before the stopping time finding the same estimate as in \eqref{eqn:2G-proof-c1}. The third and fifth terms can be treated by using the identity $G^z ( \i \eta) E_- = - E_- G^z ( - \i \eta)$ to bring them into a form for which we can apply the same argument as in \eqref{eqn:2G-proof-c2}. The remainder of the argument is analogous to the $\eta_1 \eta_2 <0$ case.  \qed

\section{Proof of Theorem \ref{thm:homog}}

  The items \ref{it:constrdbm} and \ref{it:homog-ind} simply follow by construction, as we will see below. We thus now focus on \ref{it:asympind} and \ref{it:inuniv}, starting with \ref{it:asympind}.  The fact that there exists a coupling between the $b_i^{z_n}$ and the $W_i^{(n)}$, for $|i|\le N^b$, so that \eqref{eqn:homog-new-bm} is satisfied, follows analogously to the display above \cite[Eq. (4.15)]{bourgade2024fluctuations} (see also \cite[Eqs. (7.20)--(7.30)]{cipolloni2023central} for the original construction). This construction automatically ensures items \ref{it:constrdbm} and \ref{it:homog-ind}. We point out that the construction in \cite{cipolloni2023central} was performed in the complex case, but in the real case the construction is completely analogous, as it can be seen from  the discussion in \cite[Section 7.1.2]{cipolloni2021fluctuation}; we stress that, using the convention in \eqref{eq:corrbms} and in the line below it, this construction in the real case does not distinguish the case when $z_n$ is real or it is order one away from the real line. The only input used in both cases is the bound on singular vector overlaps proven in\footnote{ We point out that for the purpose of this proof we do not need the optimal bound from \cite[Corollary 3.6]{cipolloni2025optimal}, but a weaker version would be enough, as it was shown in \cite{cipolloni2021fluctuation, cipolloni2023central}.} \cite[Corollary 3.6]{cipolloni2025optimal} both in the real and complex case.  In fact, following the argument in \cite[Eqs. (7.53)--(7.54)]{cipolloni2023central}, we have
\begin{align}
\d | \langle (b_i^{z_n} - W^{(n)}_i ) , (b_j^{z_n} - W^{(n)}_j ) \rangle | & \lesssim N^{2\mathfrak{b}} \times \sup_{|i|, |j|\le N^{\mathfrak{b}}} \left\{ \sup_{\substack{  z,w\in \{z_1, \dots, z_K\} \\ z\neq w }  }|\langle {\bm u}_i^z, {\bm u} _j^w\rangle\langle {\bm v}_j^w, {\bm v}_i^z\rangle|^2  \right\} \d t \notag\\
&\lesssim \frac{N^{2\mathfrak{b}+\xi}}{(N\min_{i\ne j}|z_i-z_j|^2\nc+1)^2} \d t \le N^{-\mathfrak{b}}\d t,
\end{align}
with overwhelming probability, choosing $\mathfrak{b}$ sufficiently small in terms of $b$. Here in the first inequality we also used that in the real case $z_n\in \hat{\mathbb{D}}_r^1$, i.e. $|z_i-\overline{z}_j|\gtrsim |z_i-z_j|$. Additionally,  in the last inequality we used \cite[Corollary 3.6]{cipolloni2025optimal}, which holds uniformly for $z,w\in \hat{\dd}_r^\beta$ for both the real and complex cases, that $\xi>0$ is arbitrary small, and that $|z-w|\ge N^{-1/2+b}$.  This concludes the proof of \ref{it:asympind}. 

 We now turn to the proof of \ref{it:inuniv}.  Given \eqref{eqn:homog-new-bm} as an input, the approximation in \eqref{eqn:homog} follows as in the proof of \cite[Theorem 4.1]{bourgade2024fluctuations} in the complex case. 
 
 Let us now discuss the real i.i.d. case. We have to distinguish the analysis for $z_n$'s on the real line and those which are order one away from it. As we have already discussed, the coupling as in the theorem statement can be constructed as in \cite{cipolloni2021fluctuation}, as all the formulas in \cite[Eqs. (7.14)--(7.16)]{cipolloni2021fluctuation} holds also when (some of) the $z_n$'s are real, thanks to \eqref{eq:corrbms}. Once the coupling has been constructed, the second component of the proof is to prove \eqref{eqn:homog}. At this point, the proof proceeds for each $z_n$ individually; i.e., there is no more reference to the fact that we have multiple processes for different $z_n$.
 
 We therefore distinguish two cases: first, when $|\Im z_n|\ge c$, for some $c>0$; second, when $z_n\in I_r$. 
In the first case 
we rely on \cite[Eq. (7.44)]{cipolloni2021fluctuation}. We point out that in the proofs of \cite[Section 7]{cipolloni2021fluctuation} it was used both that $|\Im z_i|\ge N^{-\delta}$ and $|z_1-z_2|\ge N^{-\delta}$, for some small $\delta>0$. The first assumption holds in our case. As for the second restriction, inspecting the arguments in \cite[Section 7]{cipolloni2021fluctuation}, we see that the assumption $|z_1-z_2|\ge N^{-\delta}$ was only used the ensure the bound in \cite[Lemma 7.4]{cipolloni2021fluctuation} for\footnote{The bound in \cite[Lemma 7.4]{cipolloni2021fluctuation} is stated for all $|i|, |j|\le N$. However, this bound is used for all the indices only in the special case when $z_1=z$ and $z_2=\overline{z}$ within the proof of \cite[Proposition 7.7]{cipolloni2021fluctuation}. In this case the same bound holds since we are still assuming that $|\Im z_i|\sim 1$. For general $z_1, z_2$ it is used only for small indices.} $|i|, |j|\le N^{\mathfrak{b}}$, for a small $\mathfrak{b}>0$, as in the complex case.
This bound now holds for any $|z_1-z_2|\ge N^{-1/2+b}$ thanks to the more recent \cite[Corollary 3.6]{cipolloni2025optimal}. Hence, the arguments of \cite[Section 7]{cipolloni2021fluctuation} carry through with only very minor changes.

In case ii), we consider the interpolation flow 
\beq
\d x_i^{z_n}(t,\alpha)=\frac{\alpha \d b_i^{z_n}(t)+(1-\alpha)\d W_i^{(n)}(t)}{\sqrt{2N(1+N^{-\omega_r})}}+\frac{1}{2N}\sum_{j\ne \pm i}\frac{1}{x_i^{z_n}(t,\alpha)-x_j^{z_n}(t,\alpha)}\d t,
\eeq
with $x_i(0,\alpha)=\lambda_i^{z_n}(0)$ and a small $\omega_r\ge 10\mathfrak{b}$ (which ensures the well-posedness of this flow). We now proceed as in \cite[Section 4.3]{bourgade2024fluctuations}, the only difference is that \cite[Eq. (4.19)]{bourgade2024fluctuations} needs to be replaced by
\beq
\begin{split}
\label{eq:newker}
\d u_i(t,\alpha)&=\d \xi_i(t)-\sum_{j\ne \pm i}c_{ij}(t)\big(u_j(t,\alpha)-u_i(t,\alpha)\big)\d t, \\
\xi_i(t):&=\frac{b_i^{z_n}(t)-W_i^{(n)}(t)}{2\sqrt{N}}, \qquad c_{ij}(t):=\frac{1}{2N(x_i^{z_n}(t,\alpha)-x_j^{z_n}(t,\alpha))^2}.
\end{split}
\eeq
Here $u_i(t,\alpha)=\partial_\alpha x_i^{z_n}(t,\alpha)$.

It is then easy to see that the fact that \eqref{eq:newker} does not contain the term $j=-i$ does not play any role to prove analogs of \cite[Lemmas 4.10 and 4.12--4.13]{bourgade2024fluctuations}. With the analog of these lemmas at hand, we immediately conclude
\[
\max_{|i|\le N^{\mfb}} |u_i(t,\alpha)|\le N^{-1-\mathfrak{a}},
\]
for some $\mathfrak{a}, \mfb>0$, with overwhelming probability. This shows
\[
\sup_{ |i| \leq N^{\mfb}} | \lambda_i^{(z_n)} (t) - y_i^{z_n} (t) | \leq N^{-1-\mfa}.
\]
Notice that the evolution of the $y_i^{z_n} (t)$'s is governed by i.i.d. Brownian motions, hence we can use the result in \cite[Theorem 3.2]{che2019universality} to conclude \eqref{eqn:homog}, and thus the proof.\nc

\qed

\section{Positive definiteness of kernel}
\label{sec:posker}

In this section we prove that the kernel $K^\beta(z,w;\kappa_4)$ from \eqref{eqn:int-cov} is non-negative definite on $\mathbb{D}^\beta:=\mathbb{D}_1^\beta$.

\noindent\textbf{Complex case ($\beta=2$):} In this case we have
\beq
K(z,w):=K^2(z,w;\kappa_4)=-\frac{1}{2}\log|z-w|+\frac{\kappa_4}{4}\psi(z)\psi(w),\qquad\quad \psi(x):=1-|x|^2,
\eeq
on $\mathbb{D}^1=\mathbb{D}$. Consider the Sobolev space $H_0^1(\mathbb{D})$, with the norm $\lVert g\rVert_{H_0^1(\mathbb{D})}:=\big(\int_\mathbb{D}|\nabla g(z)|^2\,\d^2 z\big)^{1/2}$. Let $H^{-1}(\mathbb{D})$ be the dual space of $H_0^1(\mathbb{D})$, which is naturally equipped with the norm (see e.g. \cite[Section 5.9.1]{evans2022partial})
\beq
\label{eq:normdefsob}
\lVert f\rVert_{H^{-1}(\mathbb{D})}:=\sup_{0\ne \phi\in H_0^1(\mathbb{D})}\frac{\left|\int_\mathbb{D}f(z)\overline{\phi}(z)\,\d^2 z\right|}{\lVert \phi\rVert_{H_0^1(\mathbb{D})}}.
\eeq
From this definition for $f\in H^{-1}(\mathbb{D})$, by choosing $\phi=\psi$, it immediately follows (note that $\psi\in H_0^1(\mathbb{D})$)
\beq
\left|\int_\mathbb{D} f(z)\psi(z)\,\d z\right|^2\le \lVert f\rVert_{H^{-1}(\mathbb{D})}^2\lVert \psi\rVert_{H_0^1(\mathbb{D})}^2=2\pi \lVert f\rVert_{H^{-1}(\mathbb{D})}^2=2\pi\langle f, (-\Delta_\mathbb{D})^{-1}f\rangle,
\eeq
where $\Delta_\mathbb{D}$ is the Dirichlet Laplacian on the disc.  Here in the first equality we used that $\lVert \psi\rVert_{H_0^1(\mathbb{D})}^2=2\pi$, and in the second equality follows from the weak formulation $-\Delta_\mathbb{D}u=f$ and the identity $\lVert f\rVert_{H^{-1}(\mathbb{D})}=\lVert \nabla u\rVert_{L^2(D)}$ (see e.g. \cite[Section 3]{adams2003sobolev}).

Next, we recall that for $f\in H^{-1}$ we have
\beq
2\pi\langle f, (-\Delta)^{-1}f\rangle:=-\int\int_\cc f(z)\log|z-w|f(w)\, \d^2z\d^2w.
\eeq
Additionally, for $f$ supported in the disc we have
\beq
\langle f, (-\Delta_\mathbb{D})^{-1}f\rangle\le \langle f, (-\Delta)^{-1}f\rangle.
\eeq
We thus conclude that
\beq \label{eqn:positivity-a1}
\int\int_\cc f(z)K(z,w)f(w)\, \d^2z\d^2w\ge \frac{1}{2}\left(1+\frac{\kappa_4}{2}\right)\ge \frac{1}{4}\langle f, (-\Delta)^{-1}f\rangle,
\eeq
where in the second inequality we used that $\kappa_4\ge -1$.

\noindent \textbf{Real case ($\beta=1$):} In this case we note that is $f$ is supported on $\mathbb{D}^1$ then we have
\beq
\begin{split}
\int\int_{\mathbb{D}^1}f(z)K^1(z,w;\kappa_4)f(w)\,\d^2 z\d^2w&=\int\int_\mathbb{D}(Pf)(z)K^1(z,w;\kappa_4)(Pf)(w)\,\d^2 z\d^2w, \\
&=2\int\int_\mathbb{D}(Pf)(z)K^2(z,w;\kappa_4/2)(Pf)(w)\,\d^2 z\d^2w,
\end{split}
\eeq
where $(Pf)(z):=[f(z)+f(\overline{z})]/2$ denotes the symmetrization of $f$ with respect to the real axis. The positivity in the real case is thus analogous to the complex case above.

{\small
\bibliography{gmc_bib}
\bibliographystyle{abbrv}
}

\end{document}